\documentclass{article}

\usepackage{amsmath,amssymb,amsthm,epic,eepic}
\usepackage{epsfig} 

\title{On the role of abnormal minimizers in
sub-Riemannian geometry}

\author{B.~Bonnard, E.~Tr\'elat \\ \\
Universit\'e de Bourgogne, LAAO, \\
BP47870, 21078 Dijon Cedex, France  \\
e-mail : bbonnard@u-bourgogne.fr, trelat@topolog.u-bourgogne.fr}
\date{}

\newtheorem{thm}{Theorem}[section]
\newtheorem{cor}[thm]{Corollary}
\newtheorem{prop}[thm]{Proposition}
\newtheorem{lem}[thm]{Lemma}

\theoremstyle{remark}
\newtheorem{rem}{Remark}[section]

\theoremstyle{definition}
\newtheorem{defi}{Definition}[section]

\theoremstyle{definition}

\newtheorem{conj}[thm]{Conjecture}

\renewcommand{\geq}{\geqslant}
\renewcommand{\leq}{\leqslant}
\renewcommand{\l}{\left}
\renewcommand{\r}{\right}
\newcommand{\f}[2]{\frac{#1}{#2}}
\newcommand{\inv}[1]{\frac{1}{#1}}
\renewcommand{\it}{\textit}
\renewcommand{\rm}{\textrm}
\newcommand{\tto}{\longrightarrow}

\def\N{\rm{I\kern-0.21emN}}
\def\C{\mathbb{C}}
\def\R{\rm{I\kern-0.21emR}}
\def\H{\rm{I\kern-0.21emH}}
\def\K{\rm{I\kern-0.24emK}}
\def\P{\rm{I\kern-0.21emP}}
\def\Q{\rm{l\kern-0.5emQ}}

\newcommand{\no}{\noindent}

\newcommand{\drawcenteredtext}[3]{\put(#1,#2){\makebox(0,0){#3}}}%
\newcommand{\drawlefttext}[3]{\put(#1,#2){\makebox(0,0)[l]{#3}}}%
\newcommand{\drawrighttext}[3]{\put(#1,#2){\makebox(0,0)[r]{#3}}}%

\newcommand{\drawcircle}[4]{%
	\put(#1,#2){\circle{#3}} \put(#1,#2){\makebox(0,0)[cc]{#4}}}%

\newcommand{\drawellipse}[5]{%
	\put(#1,#2){\ellipse{#3}{#4} \makebox(0,0)[cc]{#5}}}%

\newcommand{\drawpath}[4]{\path(#1,#2)(#3,#4)}%
\newcommand{\drawdotline}[4]{\dottedline[.]{1}(#1,#2)(#3,#4)}%
\newcommand{\drawdashline}[4]{\dashline{1}(#1,#2)(#3,#4)}%
\newcommand{\drawvector}[5]{\put(#1,#2){\vector(#4,#5){#3}}}%

\newcommand{\drawdot}[2]{\put(#1,#2){\circle*{0.5}}}%
\newcommand{\drawthickdot}[2]{\put(#1,#2){\circle*{1}}}%
\newcommand{\drawhollowdot}[2]{\put(#1,#2){\circle{1}}}%

\newcommand{\drawleftbrace}[3]%
{\drawcenteredtext{#1}{#2}{$\left\{ \rule[0mm]{0mm}{#3mm} \right.$}}%

\newcommand{\drawrightbrace}[3]%
{\drawcenteredtext{#1}{#2}{$\left\} \rule[0mm]{0mm}{#3mm} \right.$}}%

\newcommand{\drawoverbrace}[3]%
{\drawcenteredtext{#1}{#2}{$\overbrace{\rule[0mm]{#3mm}{0mm}}$}}%

\newcommand{\drawunderbrace}[3]%
{\drawcenteredtext{#1}{#2}{$\underbrace{\rule[0mm]{#3mm}{0mm}}$}}%

\newcommand{\drawarc}[5]%
{\put(#1,#2){\arc{#3}{#4}{#5}}}%

\newcommand{\drawshadebox}[8]%
{%
	\makebox(0,0)[cc]{\special{sh #8}\path(#1,#4)(#3,#4)(#3,#6)(#1,#6)(#1,#4)}%
	\put(#2,#5){\makebox(0,0){\makebox(0,0)[cc]{#7}}}%
}

\newcommand{\drawshadowframebox}[5]%
{%
	\put(#1,#2){\makebox(0,0){\shadowbox{\rule{0mm}{#4mm}\rule{#3mm}{0mm}}}}%
	\put(#1,#2){\makebox(0,0){#5}}%
}%

\newcommand{\rala}{\sqrt{\lambda}}
\newcommand{\esp}{\,\,\,}

\begin{document}

\maketitle

\begin{abstract} Consider a sub-Riemannian geometry $(U,D,g)$
where
$U$ is a neighborhood at $0$ in $\R^n,$ $D$ is a rank-2 smooth
$(C^\infty $ or $C^\omega )$ distribution and $g$ is a smooth
metric on $D$. The objective of this article is to explain the
role of abnormal minimizers in SR-geometry. It is based on the
analysis of the Martinet SR-geometry.
\end{abstract}

\paragraph{Key words :} optimal control, singular
trajectories,
sub-Riemannian geometry, abnormal minimizers, sphere and
wave-front with small radii.

\paragraph{R\'esum\'e}
On consid\`ere un probl\`eme sous-Riemannien $(U,D,g)$ o\`u $U$
est un
voisinage de $0$ dans $\R^n$, $D$ une distribution lisse de rang
$2$ et $g$ une m\'etrique lisse sur $D$. L'objectif de cet article
est d'expliquer le r\^ole des g\'eod\'esiques anormales minimisantes
en g\'eom\'etrie
SR. Cette analyse est fond\'ee sur le mod\`ele SR de Martinet.

\paragraph{Titre}
Le r\^ole des g\'eod\'esiques anormales minimisantes en g\'eom\'etrie
sous-Riemannienne.

\paragraph{AMS classification :}
49J15, 53C22.

\section{Introduction} 
Consider a smooth control system on $\R ^n$ :
\begin{eqnarray}  \label{syst}
{\dot q}(t)=f(q(t),u(t)) 
\end{eqnarray}

\noindent where the set of admissible controls ${\cal U}$ is
an open set of bounded measurable mappings $u$ defined on
$[0,T(u)]$ and taking their values in $\R^m.$ We fix $q(0)=q_0$
and $T(u)=T$ and we consider the \it{end-point mapping}
$E:u\in {\cal U} \longmapsto q(T,q_0,u),$ where $q(t,q_0,u)$
is the solution of (\ref{syst}) associated to $u\in {\cal U}$ and
starting from $q_0$ at $t=0.$ We endow the set of controls defined
on $[0,T]$ with the $L^\infty $-topology. A trajectory
${\widetilde q} (t, q_0, {\widetilde u})$ denoted in short
${\widetilde q}$ is said to be \it{singular or abnormal} on
$[0,T]$ if ${\widetilde u}$ is a singular point of the end-point
mapping, i.e, the Fr\'echet derivative of $E$ is not surjective
at ${\widetilde u}.$

Consider now the \it{optimal control problem}~:
$\underset {u
(.)\in {\cal U}}{\min} \int^T_0 f^0 (q(t),u(t))dt$, where
$f^0$ is smooth and $q(t)$ is a trajectory of (\ref{syst})
subject to
boundary conditions~: $q(0)\in M_0$ and $q(T)\in M_1,$ where $M_0$
and $M_1$ are smooth submanifolds of $\R^n.$
According to the \it{weak maximum principle} \cite{P}, minimizing
trajectories are among the singular trajectories of the end-point
mapping of the \it{extended system} in $\R^{n+1}$~:
\begin{eqnarray}  \label{extended}
{\dot q} (t)& = & f(q(t),u(t))  \\ 
{\dot q}^0 (t) & = & f^0 (q(t),u(t))\nonumber
\end{eqnarray}

\noindent and they are solutions of the following equations~:
\begin{eqnarray}
{\dot q}={\displaystyle {\partial H_\nu \over \partial p}}\ 
,\ \ {\dot p}=-{\displaystyle {\partial H_\nu \over \partial q}}
\ ,\ \ {\displaystyle {\partial H_\nu \over \partial u}}=0
\end{eqnarray}

\noindent where $H_\nu =<p,$ $f(q,u)>+\nu f^0 (q,u)$ is the
\it{pseudo-Hamiltonian}, $p$ is the adjoint vector, $<,>$
the standard inner product in $\R^n$ and $\nu $ is a constant
which can be normalized to $0$ or $-1/2.$ Abnormal
trajectories correspond to $\nu =0$~; their role
in the optimal control problem has to be analyzed. Their
geometric interpretation is clear~: if ${\cal C}$ denotes
the set of curves solutions of (\ref{syst}), they correspond to
singularities of this set and the analysis of those singularities
is a preliminary step in any minimization problem. This problem
was already known in the classical calculus of variations, see
for instance the discussion in \cite{Bl} and was a major problem
for
post-second war development of this discipline whose modern name
is optimal control. The main result, concerning the analysis of
those singularities in a \it{generic context} and for
\it{affine systems} where $f (q,u)=F_0 (q)+{\displaystyle 
{\sum ^m_{i=1}}} u_i F_i (q)$ are given in
\cite{BK} and \cite{AS}. The consequence of this analysis
is to get \it{rigidity}
results about singular trajectories when $m=1,$ that is under
generic conditions a singular trajectory $\gamma $ joining
$q_0 =\gamma (0) $ to $q_1 =\gamma (T)$ is the only trajectory
contained in a $C^0$-neighborhood of $\gamma $ joining $q_0$
to $q_1$ in time $T$ (and thus is minimizing).

In optimal control the main concept is the \it{value function}
$S$ defined as follows. If $q_0, q_1$ and $T$ are fixed and
$\gamma $ is a minimizer associated to $u_\gamma$ and
joining $q_0$ to $q_1$ in time $T,$ we set~:
$$S(q_0,q_1,T)=\int ^T_0\ f^0 (\gamma (t),u_\gamma(t))dt$$

\noindent The value function is solution of
Hamilton-Jacobi-Bellman
equation and one of the main questions in optimal control is to
understand the role of abnormal trajectories on the
\it{singularities of $S$}.

The objective of this article is to make this analysis in
\it{local sub-Riemannian geometry} associated to the following
optimal control problem~:
$${\ \atop {\displaystyle \min \atop {\displaystyle {u(.)}}}}\
\int^T_0\ \sum^m_{ i=1}\ u_i^2 (t) dt\ ,\ \ 
u=(u_1 ,\ldots ,u_m)$$

\noindent subject to the constraints~:
\begin{eqnarray}   \label{quatre}
{\dot q}(t)={\displaystyle {\sum^m_{i=1}}}\ u_i (t) F_i (q(t))
\end{eqnarray}

\noindent $q\in U\subset \R^n,$ where $\{ F_1,\cdots ,F_m\}$
are $m$ linearly independant vector fields generating a
distribution $D$ and the metric $g$ is defined on $D$ by
taking the $F_i$'s as \it{orthonormal vector fields}.
The \it{length}
of a curve $q$ solution of (\ref{quatre}) on $[0,T]$ and
associated to $u\in {\cal U}$ is given by~:
$L (q)={\displaystyle {\int ^T_0}}\Bigl (\sum u_i^2 (t)\Bigr
)^{1/2}dt$ and the {\it SR-distance} between $q_0$ and $q_1$
is the minimum of the length of the curves $q$ joining $q_0$
to $q_1.$ The sphere $S(q_0,r)$ with radius $r$ is the
set of points at SR-distance $r$ from $q_0$. If any pairs
$q_0, q_1$ can be joined by a minimizer, the sphere is made of
end-points of minimizers with length $r$. It is a
\it{level set} of the value function.

It is well known (see \cite{Ag})
that in SR-geometry the sphere $S(q_0,r)$
with \it{small} radius $r$ has singularities. For instance they
are described in \cite{AEGK} in the generic contact situation
in $\R^3$ and they are {\it semi-analytic.} Our aim is to give
a \it{geometric framework} to analyze the \it{singularities} of
the sphere in the abnormal directions and to compute
\it{asymptotics} of the distance in those directions. We analyze
mainly the Martinet case extending preliminary calculations from
\cite{ABCK, BC}. The calculations are intricate because the
singularities are not in the subanalytic category even if the
distribution and the metrics are analytic. Moreover they are
related to similar computations to \it{evaluate Poincar\'e return
mappings in the Hilbert's 16th problem} (see \cite{MR,R}) using
\it{singular perturbation techniques}.

The organization and the contribution of this article is the
following.

In Sections 2 and 3 we introduce the required concepts and recall
some known results from \cite{AS2, Be, BK} to make this article
self-contained.

In Section 2, we compute the singular trajectories for
single-input affine control systems~: $\dot{q}=F_0(q)+uF_1(q)$,
using the Hamiltonian formalism. Then we evaluate under generic
conditions the accessibility set near a singular trajectory to
get rigidity results and to clarify their optimality status in SR
geometry.

In Section 3 we present some generalities concerning SR geometry.

In Section 4 we analyze the role of abnormal geodesics in SR
Martinet geometry. We study the behaviour of the geodesics
starting from $0$ in the abnormal direction by taking their
successive intersections with the Martinet surface filled by
abnormal trajectories. This defines a \it{return mapping}. To
make precise computations we use a gradated form of order $0$
where the Martinet distribution is identified to $\rm{Ker
}\omega$, $\omega=dz-\f{y^2}{2}dx$, the Martinet surface to $y=0$,
the abnormal geodesic starting from $0$ to $t\mapsto (t,0,0)$ and
the metric is truncated at~: $g=(1+\alpha y)^2dx^2+(1+\beta
x+\gamma y)^2dy^2$, where $\alpha,\beta,\gamma$ are real
parameters. The geodesics equations project onto the \it{planar
foliation}~:
\begin{equation}\label{planarfoliation}
\theta''+\sin\theta+\varepsilon\beta\cos\theta\
\theta'+\varepsilon^2\alpha\sin\theta(\alpha\cos\theta-
\beta\sin\theta) =0
\end{equation}
where $\varepsilon=\inv{\rala}$ is a \it{small parameter} near
the abnormal direction which projects onto the singular points
$\theta=k\pi$.
Here the Martinet plane $y=0$ is projected onto a \it{section}
$\Sigma$ given by~:
\begin{equation}\label{section6}
\theta'=\varepsilon(\alpha\cos\theta+\beta\sin\theta)
\end{equation}

Equation (\ref{planarfoliation}) represents a \it{perturbed
pendulum}~; to evaluate the trace of the sphere $S(0,r)$ with the
Martinet plane we compute the return mapping associated to the
section $\Sigma$. Near the abnormal direction the computations
are localized to the geodesics projecting near the separatrices
of the pendulum.

In order to estimate the asymptotics of the sphere in the
abnormal direction we use the techniques developped to compute
the asymptotic expansion of the Poincar\'e return mapping for
a one-parameter family of planar vector fields. This allows to
estimate the number of limit cycles in the Hilbert's 16th
problem, see \cite{R}. Our computations split into two parts~:
\begin{enumerate}
\item A computation where we estimate the return mapping near a
saddle point of the pendulum and which corresponds to geodesics
close to the abnormal minimizer \it{in $C^1$-topology}.
\item A global computation where we estimate the return mapping
along geodesics visiting the two saddle points and which
corresponds to geodesics close to the abnormal minimizer \it{in
$C^0$-topology, but not in $C^1$-topology} (see \cite{trelatCOCV}
for a general statement).
\end{enumerate}

Our results are the following.
\begin{itemize}
\item
If $\beta=0$, the pendulum is integrable and we prove that the
sphere belongs to the \it{log-exp category} introduced in
\cite{VDD}, and we compute the asymptotics of the sphere in the
abnormal direction.
\item
If $\beta\neq 0$, we compute the asymptotics corresponding to
geodesics $C^1$ close to the abnormal one.
\end{itemize}

We end this Section by conjecturing the cut-locus in the generic
Martinet sphere using the Liu-Sussmann example \cite{LS}.

The aim of Section 5 is to extend our previous results to the
general case and to describe a \it{Martinet sector} in the
$n$-dimensional SR sphere.

First of all in the Martinet case the exponential mapping \it{is
not proper} and the sphere is \it{tangent} to the abnormal
direction. We prove that this property is still valid if the
abnormal minimizer is strict and if the sphere is
$C^1$-stratifiable.

Then we complete the analysis of SR geometry corresponding to
stable 2-dimensional distributions in $\R^3$ by analyzing the
so-called \it{tangential case}. We compute the geodesics and make
some numerical simulations and remarks about the SR spheres.

The flat Martinet case can be lifted into the \it{Engel case}
which is a left-invariant problem of a 4-dimensional Lie group.
We give an uniform parametrization of the geodesics using the
\it{Weierstrass function}. Both Heisenberg case and Martinet flat
case can be imbedded in the Engel case.

The main contribution of Section 5 is to describe a \it{Martinet
sector} in the SR sphere in any dimension using the computations
of Section 4. We use the Hamiltonian formalism (\it{Lagrangian
manifolds}) and \it{microlocal analysis}. This leads to a
\it{stratification of the Hamilton-Jacobi equation} viewed in the
cotangent bundle.

\paragraph{Acknoledgments \\}
We thank M. Chyba for her numerical simulations concerning the
sphere.


\section{Singular or abnormal trajectories}

\subsection{Basic facts} 
Consider the smooth control system~:
\begin{eqnarray}
{\dot q}(t)=f(q(t),u(t))
\end{eqnarray}

\noindent and let $q(t)$ be a trajectory defined on $[0,T]$
and associated to a control $u\in {\cal U}.$ If we set~:
$$A(t)={\partial f\over \partial q}\ (q(t),u(t))\ ,\ B(t)
={\partial f\over \partial u}\ (q(t),u(t))$$

\noindent the linear system~:
\begin{eqnarray}   \label{lin}
\delta \dot{q}(t)=A(t)\delta q(t)+B(t)\delta u(t)
\end{eqnarray}

\noindent is called the \it{linearized} or \it{variational}
system along $(q,u)$. It is well known, see \cite{BK} that the
Fr\'echet derivative in $L^\infty$-topology
of the end-point mapping $E$ is given by~:
$$E'(v)=\Phi (T) \int^T_0 \Phi ^{-1} (s)B(s)v(s)ds$$

\noindent where $\Phi $ is the matricial solution of~:
${\dot \Phi} =A\Phi ,$ with $\Phi (0)=\rm{id}$.

Hence $(q,u)$ is singular on $[0,T]$ if and only if there
exists a non-zero vector ${\bar p}$ orthogonal to $\rm{Im }E'(u),$
that is the linear system (\ref{lin}) is \it{not controllable} on
$[0,T].$

If we introduce the row vector~: $p(t)={\bar p}\Phi (T) 
\Phi ^{-1} (t)$, a standard computation shows that the triple
$(q,p,u)$ is solution for almost all $t\in [0,T]$ of the
equations~:
$${\dot q}=f(x,u)\ ,\ \ {\dot p}=-p {\partial f\over 
\partial q}\ ,\ \ p {\partial f\over \partial u}=0$$

\noindent which takes the Hamiltonian form~:
\begin{eqnarray}   \label{ham}
{\dot q}={\displaystyle {\partial H\over \partial p}}\ ,\
 {\dot p}=-{\displaystyle {\partial H\over \partial q}}\ ,\ \ 
{\displaystyle {\partial H\over \partial u}}=0
\end{eqnarray}

\noindent where $H(q,p,u)=<p, f(q,u)>.$ The function $H$
is called the \it{Hamiltonian} and $p$ is called the
\it{adjoint vector.}

This is the parametrization of the singular trajectories
using the maximum principle. We observe that for each $t\in 
]0,T[,$ the restriction of $(q,u)$ is singular on $[0,t]$
and at $t$ the adjoint vector $p(t)$ is orthogonal to the
vector space $K(t)$ image of $L^\infty [0,t]$ by the Fr\'echet
derivative of the end-point mapping evaluated on the restriction
of $u$ to $[0,t].$ The vector space $K(t)$ corresponds to the
\it{first order Pontryagin's cone} introduced in the proof of
the maximum principle. If $t\in ]0,T]$ we shall denote by $k(t)$
the codimension of $K(t)$ or in other words the codimension of
the singularity. Using the terminology of the calculus of
variations $k(t)$ is called the \it{order of abnormality.}

The parametrization by the maximum principle allows the
computation of the singular trajectories. In this article we
are concerned by systems of the form~:
\begin{eqnarray}
{\dot q}(t)=F_0 (q(t))+u(t)F_1 (q(t))
\end{eqnarray}

\noindent and the algorithm is the following.


\subsection{Determination of the singular trajectories}

\subsubsection{The single input affine case}
It is convenient to use Hamiltonian formalism. Given any
smooth function $H$ on $T^*U,$ ${\vec H}$ will denote the
Hamiltonian vector field defined by $H.$ If $H_1,$ $H_2$ are
two smooth functions, $\{ H_1, H_2\}$ will denote their
\it{Poisson bracket}~: \mbox{$\{ H_1, H_2\}=dH_1 ({\vec H}_2)$}.
If $X$ is a smooth vector field on $U,$ we set \mbox{$H=<p,X(q)>$}
and ${\vec H}$ is the \it{Hamiltonian lift} of $X$. If $X_1, X_2$
are two vector fields with \mbox{$H_i =<p, X_i(q)>$}, $i=1,2$
we have~: \mbox{$\{ H_1, H_2\}=<p, [X_1, X_2](q)>$} where
the Lie bracket is~: $[X_1, X_2](q)={\displaystyle
{\partial X_1\over \partial q}} (q) X_2 (q)-{\displaystyle
{\partial X_2\over \partial q}} (q)X_1(q).$ We shall denote
by $H_0 =<p, F_0(q)>$ and $H_1 =<p, F_1 (q)>.$

If $f(q,u)=F_0(q)+uF_1(q),$ the equation (\ref{ham}) can be
rewritten~:
$${\dot q}={\partial H_0\over \partial p}\ +\ u {\partial H_1 
\over \partial p}\ ,\ \ {\dot p}=-\Bigl ( {\partial H_0\over
 \partial q}\ +\ u {\partial H_1\over \partial q}\Bigr
) \quad \rm{a.e.}$$
$$H_1 =0 \quad \rm{for all } t\in [0,T]$$

We denote by $z=(q,p)\in T^*U$ and let $(z,u)$ be a solution of
the above equations. Using the chain rule and the constraint~:
$H_1 =0,$ we get~:
$$0={d\over dt}\ H_1(z(t))=dH_1(z(t)){\vec H}_0 (z(t))+
u(t)dH_1 (z(t))\ {\vec H}_1 (z(t)) \quad \rm{for a.e. } t$$
This implies~: $0=\{ H_1, H_0\}(z(t))$ for all $t.$ Using
the chain rule again we get~:
$$0=\{ \{ H_1,H_0\}, H_0\} (z(t))+u(t) \{ \{ H_1,H_0\}, H_1\}
 (z(t)) \quad \rm{for a.e. } t$$

This last relation enables us to compute $u(t)$ in many cases
and justifies the following definition~:
\begin{defi}
For any singular curve $(z,u):J=[0,T]\longmapsto T^*U\times 
\R ,$ ${\cal R} (z,u)$ will denote the set $\{ t\in J, \{ \{
 H_0, H_1\}, H_1\} (z(t))\not = 0\}.$ The set ${\cal R} 
(z,u)$ possibly empty is always an open subset of $J.$
\end{defi}

\begin{defi}
A singular trajectory $(z,u):J\longrightarrow T^*U\times 
\R$ is called of order two if ${\cal R} (z,u)$ is dense in $J.$
\end{defi}

The following Proposition is straightforward~:
\begin{prop}
If $(z,u):J\longmapsto T^*M\times \R $ is a singular trajectory
and ${\cal R} (z,u)$ is not empty then~:
\begin{enumerate}
\item $z$ restricted to ${\cal R} (z,u)$ is smooth~;
\item $u(t)={\displaystyle {\{ \{ H_0, H_1 \}, H_0 \}(z(t))
\over \{ \{ H_1, H_0 \}, H_1 \} (z(t))}}$ for a.e. t
\item ${\displaystyle {dz(t))\over dt}}\ = {\vec H}_0 (z(t)) +
 {\displaystyle {\{ \{ H_0, H_1 \},H_0 \} (z(t))\over \{ \{ 
H_1, H_0 \} , H_1\}  (z(t))}}
\vec{H}_1(z(t))$ for all $t \in {\cal R} (z,u).$
\end{enumerate}
\end{prop}

Conversely, let $(F_0,F_1)$ be a pair a smooth vector fields
such that the open subset $\Omega $ of all $z\in T^*U$ such
that $\{ \{ H_0, H_1 \}, H_1\} (z)\not = 0$ is not empty.
If $H:\Omega \longmapsto \R $ is the function $H_0 +
 {\displaystyle {\{ \{ H_0, H_1 \}, H_0\}\over \{ \{ H_1,
 H_0\}, H_1 \}}}$ $H_1$ then any trajectory of $\vec H$
starting at $t=0$ from the set $H_1 =\{H_1, H_0\}=0$ is a
singular trajectory of order 2.

This algorithm allows us to compute the singular trajectories
of minimal order. More generally we can extend this computation
to the general case.

\begin{defi}
For any multi-index $\alpha \in \{ 0,1\}^n,$ $\alpha =
(\alpha _1,\ldots ,\alpha _n)$ the function $H_\alpha $ is
defined by induction by~: $H_\alpha =\{ H_{(\alpha _1,\ldots 
,\alpha _{n-1})},H_{\alpha _n}\}$. A singular trajectory $(z,u)$
is said of order $k\geq 2$ if all the brackets of order
$m\leq k$~: $H_\beta ,$ with $\beta =(\beta _1,\ldots ,\beta _m)
,$ $\beta _1=1$ are $0$ along $z$ and there exists $\alpha
=(1,\alpha_2 ,\cdots ,\alpha _k)$ such that
$H_{\alpha 1}(z)$ is not identically $0.$
\end{defi}

The generic properties of singular trajectories are described
by the following Theorems of \cite{BK2}.

\begin{thm}
There exists an open dense subset $G$ of pairs of vector
fields $(F_0,F_1)$ such that for any couple $(F_0,F_1)\in G,$
the associated control system has only singular trajectories
of minimal order 2.
\end{thm}

\begin{thm}
There exists an open dense subset $G_1$ in $G$ such that for
any couple $(F_0,F_1)$ in $G_1$,
any singular trajectory has an order of
abnormality equal to one, that is corresponds to a singularity
of the end-point mapping of codimension one.
\end{thm}


\subsubsection{The case of rank two distributions}

Consider now a distribution $D$ of rank 2. In SR-geometry we
need to compute singular trajectories $t\longmapsto q(t)$ of
the distribution and it is not restrictive to assume the
following~: $t\longmapsto q(t)$ is a smooth immersion. Then
locally there exist two vector fields $F_1, F_2$ such than
$D=$ Span $\{ F_1,$ $F_2\}$ and moreover the trajectory can be
reparametrized to satisfy the associated affine system~:
$${\dot q}(t)=u_1(t)F_1 (q(t))+u_2 (t)F_2 (q(t))$$

\noindent where $u_1(t)= 1.$ It corresponds to the choice of
a \it{projective chart} on the control domain.

Now an important remark is the following. If we introduce
the Hamiltonian lifts~: $H_i =<p, F_i (q)>$ for $i=1,2$, and
$H={\displaystyle {\sum ^2_{i=1}}} u_i H_i$ the singular
trajectories are solutions of the equations~:
$${\dot p}={\partial H\over \partial q}\ ,\ {\dot q}=-
{\partial H\over \partial p}\ ,\ {\partial H\over \partial u}\ 
=0$$

Here the constraints ${\displaystyle {\partial H\over \partial u
}}=0$ means~: $H_1 =0$ and $H_2 =0.$ This
leads to the following definition~:
\begin{defi}
Consider the affine control system~: ${\dot q}=F_1 +uF_2.$ A
singular trajectory is said exceptional if it is contained on
the level set~: $H=0$, where $H=<p, F_1>+u<p, F_2>$ is the
Hamiltonian.
\end{defi}

Hence to compute the singular trajectories associated to a
distribution we can apply locally the algorithm described in the
affine case and keeping only the exceptional trajectories. An
instant of reflexion shows that those of minimal order form a
subset of \it{codimension one} in the set of all singular
trajectories because $H$ is constant along such a trajectory
and the additional constraint $H_1=0$ has to be satisfied only
at time $t=0.$ Hence we have~:

\begin{prop}
The singular arcs of $D$ are generically singular arcs of
order 2 of the associated affine system. They are exceptional
and form a subset of codimension one in the set of all singular
trajectories.
\end{prop}


\subsection{Feedback equivalence}

\begin{defi}
Consider the class ${\cal S}$ of smooth control systems of the
form~:
$${\dot q} (t)=f(q(t), u(t)),\ q\in \R^n,\ u\in \R^m$$

\noindent Two systems $f(x,u)$ and $f'(y,v)$ are called feedback
equivalent if there exists a smooth diffeomorphism of the
form~: $\Phi : (x,u)\longmapsto (y,v),$ $y=\varphi (x),$ $v=\psi
 (x,u)$ which transforms $f$ into $f'$~:
$$d\Phi (x)f(x,u)=f'(y,v)$$

\noindent and we use the notation $f'=\Phi *f.$
\end{defi}

Here we gave a global definition but there are local associated
concepts which are~:
\begin{itemize}
\item local feedback equivalence at a point
$(x_0, u_0)\in \R^n\times \R^m.$

\item local feedback equivalence at a point $x_0$
of the state-space.

\item local feedback equivalence along a given
trajectory~: $q(t)$ or $(q(t), u(t))$ of the system.
\end{itemize}

This induces a group transformation structure called the feedback
group $G_f$ on the set of such diffeomorphisms. For affine
systems we consider a sub-group of $G_f$ which stabilizes the
class. This leads to the following definition.

\begin{defi}
Consider the class of $m$-inputs affine control systems~:
$${dq\over dt}\ (t)=F_0 (q(t)+F(q(t))u(t)$$

\noindent where $F(q)u={\displaystyle {\sum ^m_{ i=1}}} 
u_i F_i (q).$ It is identified to the set ${\cal A}_{m+1}=\{
 F_0, F_1, \ldots ,F_m\}$ of $(m+1)$-uplets of vector fields.
The vector field $F_0$ is called the drift. Let $D$ be the
distribution defined by $D=\rm{Span\ }\{ F_1 (q),\ldots ,
F_m(q)\}.$ We restrict the feedback transformations to
diffeomorphisms of the form $\Phi =(\varphi (q),$ $\psi (q,u)=
\alpha (q)+\beta (q)u),$ preserving the class ${\cal A}.$ We
denote by $G$ the set of triples $(\varphi ,\alpha ,\beta )$
endowed with the group structure induced by $G_f$.
\end{defi}

We observe
the following~: take $(F_0, F)\in {\cal A}$ and $\Phi =(\varphi ,
\alpha ,\beta )\in G,$ then the image of $(F_0, F)$ by $\Phi $ is
the affine system $(F'_0, F')$ given by~:
\begin{itemize}
\item[(i)] $F'_0 =\varphi *(F_0 +F.\beta )$
\item[(ii)] $F'=\varphi * F.\beta .$
\end{itemize}

In particular the second action corresponds to the equivalence
of the two distributions $D$ and $D'$ associated to the
respective systems.

The proof of the following result is straightforward, see
\cite{Bo}.

\begin{prop}
The singular trajectories are feedback invariants.
\end{prop}

Less trivial is the assertion that for generic systems,
singular trajectories will allow to compute a \it{complete}
set of invariants, see \cite{Bo} for such a discussion.


\subsection{Local classification of rank 2 generic
distributions ${\bf D}$ in $\R ^3$}
We recall the generic classification of rank 2 distributions
in $\R ^3,$ see \cite{Zh}, with its interpretation
using
singular trajectories. Hence we consider a system~:
$${\dot q}(t)=u_1 (t)F_1 (q(t))+u_2 (t)F_2 (q(t))\ $$

\noindent $q=(x,y,z).$ We set $D=\rm{Span\ }\{ F_1 ,F_2\}$
and we assume that $D$ is of rank 2. Our classification is
localized near a point $q_0 \in \R^3$ and we can assume
$q_0=0.$ We deal only with generic situations, that is all
the cases of codimension $\leq 3.$ We have three situations which
can be distinguished using the singular trajectories.

Introducing $H_i =<p, F_i(q)>, i=1,2,$ a singular
trajectory $z=(q,p)$ must satisfy~:
$$H_1=H_2 =\{ H_1, H_2 \}=0$$

\noindent and hence they are contained in the set $M : 
\{ q\in \R ^3 \ ;\ \rm{det}\ (F_1,$ $F_2,$ $[F_1,F_2])=0 \}$
called the \it{Martinet surface.} The singular controls of
order 2 satisfy~:
$$u_1\{ \{ H_1,H_2 \}, H_1\} +u_2 \{ \{ H_1,H_2\} , H_2 \}=0$$

\noindent We define the singular set $S=S_1\cap S_2 $
where \mbox{$S_i =\{ \rm{det}(F_1,F_2,[[F_1, F_2], F_i])=0 \}$}.
We have the following situations.

\paragraph{Case 1.} Take a point $q_0\notin M,$ then
through $q_0$ there passes no singular arc. In this case $D$
is $(C^\infty $ or $C^\omega $)- isomorphic to $\rm{Ker }
\alpha$, with $\alpha=ydx+dz.$ For this normalization
$d\alpha =dy\wedge dx $ (Darboux) and ${\displaystyle 
{\partial \over \partial z}}$ is the characteristic direction.
This case is called the \it{contact case}.

\paragraph{Case 2.} (Codimension one). We take a point
$q_0\in M\backslash S.$ Since $q_0 \notin S,$ we observe
that $M$ is near $q_0$ a smooth surface. This surface is
foliated by the singular trajectories.
A smooth $(C^\infty $ or $C^\omega $)-normal form is given by
$D=\rm{Ker }\alpha ,$ where $\alpha=dz-{\displaystyle {y^2\over 2}}
 dx$. In this normal form we have the following identification~:
\begin{itemize}
\item Martinet surface $M:y=0.$

\item The singular trajectories are the solution of
$Z={\displaystyle {\partial \over \partial x}}$ restricted
to $y=0.$
\end{itemize}

\noindent This case is called the \it{Martinet case.}

\paragraph{Case 3.} (Codimension 3). We take a point
$x_0\in M\cap S$ and we assume that the point $q_0$ is a
regular point of $M.$ The analysis of \cite{Zh} shows that in
this case we have two different $C^\infty$-\it{reductions to a}
$C^\omega $-normal form depending both upon a \it{modulus $m$.}
The two cases are~:

\begin{enumerate}
\item {\bf{Hyperbolic case.}} $D=\rm{Ker } \alpha ,$ $\alpha 
=dy+(xy+x^2z+mx^3z^2)dz.$ In this representation the Martinet
surface is given by~:
$$y+2xz+3mx^2z^2=0$$

\noindent and the singular flow in $M$ is represented in the
$(x,z)$ coordinates by~:
\begin{eqnarray*}
{\dot x} & = & 2x +(6m-1)x^2 z-2mx^3z^2\\
{\dot z}& = &-(2z+6mxz^2)
\end{eqnarray*}

\noindent We observe that $0$ is a \it{resonant saddle} and the
parameter $m$ is an \it{obstruction to the
$C^\infty $-linearization.}

\item {\bf{Elliptic case.}} $D=\rm{Ker } \alpha ,$ $\alpha
 =dy+(xy+{\displaystyle {x^3\over 3}}+xz^2+mx^3z^2)dz.$ The
Martinet surface is here identified to~:
$$y+x^2+z^2+3mx^2z^2=0$$

\noindent in which the singular flow is given by~:
\begin{eqnarray*}
{\dot x}&=& 2z-{\displaystyle {2x^3\over 3}}+6mx^2z-2mx^3z^2\\
{\dot z} & =& -(2x+6mxz^2)
\end{eqnarray*}
\end{enumerate}

\noindent Hence $0$ is a center and still $m$ is an obstruction
to $C^\infty $-linearization. The analysis of \cite{Zh2} shows that
the singularity is $C^0$-equivalent to a focus.

We call the case 3 the \it{tangential situation} because $D$
is tangent to the Martinet surface at $0.$ We must stress that
it is not a \it{simple singularity} and moreover there are
\it{numerous analytic moduli.}


\subsection{Accessibility set near a singular trajectory} 
The objective of this Section is to recall briefly the
results of \cite{BK} \it{which describe geometrically the
accessibility set near a given singular trajectory satisfying
generic assumptions} (see also \cite{trelatCOCV,trelatthese}).

\subsubsection{Basic assumptions and definitions}
We consider a smooth single input smooth affine control system~:
$${\dot q}(t)=F_0(q(t))+u(t)F_1 (q(t)),\ q\in U.$$

\noindent Let $\gamma$ be a reference singular trajectory
corresponding to a control $u\in L^\infty [0,T]$ and starting
at $t=0$ from $\gamma (0)=q_0$ and we denote by $(\gamma ,p,u)$,
where $p$ is an adjoint vector for the associated solution of the
equations (\ref{ham}) from the maximum principle. We assume the
following~:
\begin{description}
\item {$(H0)$} $(\gamma ,p)$ is contained in the
set $\Omega =\{ z=(q,p)\ ;\ \{ \{ H_0, H_1\}, H_1 \} (z)\not = 
0\},$ $\gamma $ is contained in the set $\Omega '=\{ q\ ;\ 
X(q)$ and $Y(q)$ are linearly independant$\}$ and moreover
$\gamma  : [0,T]\tto U$ is one-to-one.
\end{description}

Then according to the results of Section 2.2, the curve
$z=(\gamma ,p)$ is a singular curve of order 2 solution of the
Hamiltonian vector field ${\vec H},$ with $H=H_0+\
{\displaystyle {\{ \{ H_0,H_1\}, H_0\}\over \{ \{ H_1,H_0\},
 H_1\} }}\ H_1.$ Moreover the trajectory $\gamma  :
 [0,T]\tto \Omega '$ is a smooth curve and $\gamma $
is a one-to-one immersion.

Using the feedback invariance of the singular trajectories we
may assume the following normalizations~: $u(t)=0$ for
$t\in [0,T]$ and $\gamma $ can be taken as the trajectory~:
$t\rightarrow (t,0,0,\ldots ,0).$ Since $u$ is normalized to
$0$, by successive derivations of the constraints $H_1 =0$, i.e
$<p(t), F_1 (\gamma (t))>=0$ for $t\in [0,T]$, we get the
relations~:
$$<p(t), V^k(\gamma (t))>=0 \ ,\ \ k=0,\ldots ,+\infty $$

\noindent where $V^k$ is the vector field ad$^kF_0 (F_1)$ and
ad$^k$ is defined recursively by~:
$$\rm{ad}^0F_0(F_1)=F_1,\ \rm{ad}^kF_0(F_1)=[F_0,
 \rm{ad}^{k-1} F_0(F_1)]$$

\noindent It is well known, see \cite{H}, \cite{Kre},
that for $t>0$ the space
$E(t)=\rm{Span } \{ V^k (\gamma (t)), k\in \N\}$ is
contained in the first order Pontryagin's cone $K(t)$
evaluated along $\gamma .$ We make the following assumptions~:

\begin{description}
\item{$(H1)$} For $t\in [0,t],$ the vector space $E(t)$ is of
codimension one and generated by $\{ V^0(\gamma (t)),\ldots 
,V^{(n-2)} (\gamma (t))\}.$

\item{$(H2)$} If $n\geq 3,$ for each $t\in [0,T], 
X(\gamma (t))\notin \rm{Span } \{ V^0 (\gamma (t)),\ldots
 ,V^{(n-3)} (\gamma (t))\} $
\end{description}

\begin{defi}
Let $(\gamma (t), p(t), u(t))$ be the reference trajectory
defined on $[0,T]$ and assume that the previous assumptions
$(H0),$ $(H1),$ $(H2)$ are satisfied. According to $(H1)$
the adjoint vector $p$ is unique up to a scalar. The Hamiltonian is
$H=H_0 +uH_1$ along the reference trajectory and $H_1=0.$ If
$H=0,$ we say that $\gamma $ is $G$-exceptional. Let $D=
{\displaystyle {\partial \over \partial u}}
{\displaystyle {d^2\over dt^2}} {\displaystyle {\partial
H\over
\partial u}}=<p(t), [[F_1, F_0], F_1] (\gamma (t))>.$
The trajectory $\gamma $ is said $G$-hyperbolic if $H.D>0$
along $\gamma $ and
$G$-elliptic if $H.D<0$ along $\gamma .$
\end{defi}

\begin{rem}
According to the \it{higher-order maximum principle} the
condition $H.D_{\mid \gamma }\geq 0$ called the Legendre-Clebsch
condition is a time optimality necessary condition, see
\cite{Kre}.
\end{rem}


\subsubsection{Semi-normal forms}
The main tool to evaluate the end-point mapping is to construct
semi-normal forms along the reference trajectory $\gamma $
using the assumptions $(H0,H1,H2)$ and the action of the
feedback group localized near $\gamma .$ They are given in \cite{BK}
and we must distinguish two cases.

\begin{prop}
Assume that $\gamma $ is a $G$-hyperbolic or elliptic trajectory.
Then the system is feedback equivalent in a $C^0$-neighborhood of
$\gamma $ to a system $(N_0,N_1)$ with~:
$$N_0 ={\partial \over \partial q^1} + \sum ^{n-1}_{i=2}
 q^{i+1} {\partial \over \partial q^i} + \sum ^n_{i,j =2}
 a_{ij} (q^1)\ q^iq^j {\partial \over \partial x^1} +R\ ,\ \
 N_1= {\partial \over \partial q^n}$$

\noindent where $a_{n,n}$ is strictly positive (resp. negative)
on $[0,T]$ if $\gamma $ is elliptic (resp. hyperbolic) and
$R={\displaystyle {\sum ^{n-1}_{i=1}}}$ $R_i {\displaystyle 
{\partial \over \partial q^i}}$ is a vector field such that
the weight of $R_i$ has order greater or equal to 2 (resp. 3)
for $i=2,\ldots ,n-1$ (resp. $i=1$), the weights of the
variables $q^i$ being $0$ for $i=1$, and $1$ for $i=2,\ldots ,n.$
\end{prop}

\paragraph{Geometric interpretation :}
\begin{itemize}
\item The reference trajectory $\gamma $ is identified to
$t\longmapsto (t,0,\ldots ,0)$ and the associated control is
$u\equiv 0.$ In particular $N_{0\mid \gamma }\ =\ {\displaystyle 
{\partial \over \partial q^1_{\mid \gamma }}}.$

\item We have~:
\begin{itemize}
\item[(i)] $N_1={\displaystyle {\partial \over \partial q^n}}$
\item[(ii)]
$ad^k\ N_0.N_{1\mid \gamma }  = \l\{ \begin{array}{ll} 
 {\displaystyle {\partial \over 
\partial q^{n-k}}}\ \ & \rm{if } k=1,\ldots ,n-2  \\
 0 & \rm{if } k>n-2 
\end{array}\r. $
\item[(iii)] ad$^2 N_1.N_0={\displaystyle {\partial ^2N_1\over 
\partial {q^{n}}^2}}$
\end{itemize}

\noindent and the first order Pongryagin's cone along $\gamma $
is~: $K_{\mid \gamma }=\Bigl \{ {\displaystyle {\partial \over 
\partial q^2_{\mid \gamma }}},\ldots ,{\displaystyle {\partial 
\over \partial q^n_{\mid \gamma }}}\Bigr \}.$ The linearized
system is autonomous and in the Brunovsky canonical form~:
${\dot \varphi }^1=\varphi ^2,\ldots ,{\dot \varphi }^n =u.$

\item The adjoint vector associated to $\gamma $ is
$p=(\varepsilon ,0,\ldots ,0)$ where $\varepsilon =+1$ in the
elliptic case and $\varepsilon =-1$ in the hyperbolic case,
the Hamiltonian being $\varepsilon .$ 

\item The \it{intrinsic second-order derivative} of the
end-point mapping is identified along $\gamma $ to~:
$$\varepsilon \ \int_0^T \sum_{i,j=2}^n  a_{ij}(t)\ \varphi ^i (t)
\varphi ^j (t)dt$$

\noindent with ${\dot \varphi }^2=\varphi ^3,\ldots 
,{\dot \varphi }^{n-1}=\varphi ^n,$ ${\dot \varphi }^n=u$
and the boundary conditions at $s=0$ and $T\ :\ \varphi ^2 
(s)=\cdots =\varphi ^n(s)=0.$
\end{itemize}

\begin{prop}
Let $\gamma $ be a $G$-exceptional trajectory. Then $n\geq 3$
and there exists a $C^0$-neighborhood of $\gamma $ in which
the system is feedback equivalent to a system $(N_0,N_1)$ with~:

\begin{eqnarray*}
N_0 & = & {\displaystyle {\partial \over \partial q^1}} +
{\displaystyle {\sum^{n-2}_{i=1}}} q^{i+1} {\displaystyle 
{\partial \over \partial q^i}} + {\displaystyle 
{\sum ^{n-1}_{i,j=2}}} a_{ij} (q^1) q^i q^j {\displaystyle 
{\partial \over \partial q^n}} +R\\
N_1 &=& {\displaystyle {\partial \over \partial q^{n-1}}}
\end{eqnarray*}

\noindent where $a_{n-1,n-1}$ is strictly positive on $[0,T]$
and $R={\displaystyle {\sum ^n_{i=1}}} R_i {\displaystyle 
{\partial \over \partial q^i}}\ ,$ $R_{n-1}=0$ is a vector
field such that the weight of $R_i$ has order greater or equal
to 2 (resp. 3) for $i=1,\ldots ,n-2$ (resp. $i=n$), the weights
of the variables $q^i$ being zero for $i=1$, one for $i=2,\ldots 
,n-1$ and two for $q^n.$
\end{prop}

\paragraph{Geometric interpretation}
\begin{itemize}
\item The reference trajectory $\gamma $ is identified to
$t\longmapsto (t,0,\ldots ,0)$ and the associated control is
$u\equiv 0.$

\item We have the following normalizations~:
\begin{itemize}

\item[(i)] 
$\rm{ad}^kN_0.N_{1\mid \gamma }  = 
\l\{ \begin{array}{ll}
 {\displaystyle 
{\partial \over \partial q^{n-1-k}_{\mid \gamma }}}\ \
 & \rm{for } k=0,\ldots ,n-3 \\
  0  & \rm{for } k>n-2
\end{array} \r. $

\item[(ii)] $N_{0\mid \gamma }\ ={\rm \ ad}^{n-2}\ 
N_0.N_{1\mid \gamma }=\ {\displaystyle {\partial \over
 \partial q^1_{\mid \gamma }}}.$

\item[(iii)] ad$^2N_1.N_0=\ {\displaystyle {\partial ^2N_1\over 
\partial q^{{n-1}^2}}}$
\end{itemize}

\noindent and the first order Pontryagin's cone along $\gamma $
is $K_{\mid \gamma }={\rm \ Span\ }\Bigl \{ {\displaystyle 
{\partial \over \partial q^1_{\mid \gamma }}}\ ,\ \ldots \ ,\
 {\displaystyle {\partial \over \partial q^{n-1}_{\mid \gamma
 }}}\Bigr \} .$ \it{In the exceptional case $\dot{\gamma}(t)$ is
tangent to $K_{\mid \gamma (t)}$.}

\noindent The linearized system along $\gamma $ is the system~:
${\dot \varphi }^1=\varphi ^2,\ldots ,{\dot \varphi }^{n-2}=
\varphi ^{n-1},$ ${\dot \varphi }^{n-1}=u.$

\item The adjoint vector $p$ associated to $\gamma $ can be
normalized to \mbox{$p=(0,\ldots ,0,-1).$}

\item The intrinsic second-order derivative of the end-point
mapping along $\gamma $ is identified to~:
$$-\int ^T_0 \sum ^{n-1}_{i,j =2}\ a_{ij}\ (t)\ \varphi ^i (t)\
 \varphi ^j (t)\ dt$$

\noindent with~: ${\dot \varphi^1 }=\varphi ^2,\ldots ,{\dot 
\varphi }^{n-2}=\varphi ^{n-1},\ {\dot \varphi }^n =u$ and the
boundary conditions at $s=0$ and $T\ :\ \varphi ^1 (s)=\cdots 
=\varphi ^{n-1} (s)=0.$

\end{itemize}


\subsubsection{Evaluation of the accessibility set near
$\gamma $}
We consider all trajectories $q(t,u)$ of the system starting
at time $t=0$ from $\gamma (0)=0$~; the \it{accessibility set}
at time $t$ is the set~: $A(0,t)={\displaystyle {\bigcup_{u\in 
{\cal U}}}}$ $q(t,u).$ It is the image of the end-point mapping.

We use our semi-normal forms to evaluate the accessibility set
for all trajectories of the system contained in a
$C^0$-neighborhood of $\gamma .$ We have the following,
see \cite{BK} for the details.

\paragraph{Hyperbolic-elliptic situation \\}
By truncating the semi-normal form and replacing $q^1$ by $t$
we get a \it{linear-quadratic model}~:
\begin{eqnarray*}
{\dot q}^1 &=& 1+{\displaystyle {\sum ^n_{i,j=2}}} a_{ij}(t)
 q^i q^j\\
{\dot q}^2 & = & q^3 ,\ \ldots \ ,{\dot q}^n=u
\end{eqnarray*}

\noindent and it can be \it{integrated in cascade.}

Let $0<t\leq T$ and fix the following boundary conditions~:
$q(0)=0$ and~: $q^2 (t)=\cdots =q^{n-1} (t)=q^n(t)=0,$ we get
a projection of the accessibility set $A(0,t)$ in the line $q^1$
which describes the singularity of the end-point mapping evaluated
on $u(s)=0$ for $0\leq s\leq t.$ Fig. \ref{f1}
represents this projection when $t$ varies.

\unitlength=.25mm
\makeatletter
\def\shade{\@ifnextchar[{\shade@special}{\@killglue\special{sh}\ignorespaces}}
\def\shade@special[#1]{\@killglue\special{sh #1}\ignorespaces}
\makeatother

\begin{figure}[h]
\begin{center}

\begin{picture}(426,186)(96,-5)
\thinlines
\typeout{\space\space\space eepic-ture exported by 'qfig'.}
\font\FonttenBI=cmbxti10\relax
\font\FonttwlBI=cmbxti10 scaled \magstep1\relax
\path (120,55)(280,55)
\path (120,55)(120,163)
\path (120,55)(265,144)
\path (272.482,57.736)(280,55)(272.482,52.264)
\path (117.264,155.482)(120,163)(122.736,155.482)
\put(96,154){{\rm\rm {$q^1$}}}
\put(274,35){{\rm\rm {$t$}}}
\put(105,40){{\rm\rm {$0$}}}
\path (231,128)(231,120)
\dottedline{3}(231,118)(231,55)
\put(224,38){{\rm\rm {$t_{1c}$}}}
\path (360,55)(520,55)
\path (360,55)(360,163)
\path (360,55)(505,144)
\put(336,154){{\rm\rm {$q^1$}}}
\put(514,35){{\rm\rm {$t$}}}
\put(345,40){{\rm\rm {$0$}}}
\path (471,128)(471,120)
\dottedline{3}(471,118)(471,55)
\put(464,38){{\rm\rm {$t_{1c}$}}}
\path (512.482,57.736)(520,55)(512.482,52.264)
\path (357.264,155.482)(360,163)(362.736,155.482)
\path (158,79)(158,139)
\path (154,133)(162,139)
\path (154,125)(162,131)
\path (154,117)(162,123)
\path (154,109)(162,115)
\path (154,101)(162,107)
\path (154,93)(162,99)
\path (154,85)(162,91)
\path (196,119)(204,125)
\path (200,105)(200,165)
\path (196,159)(204,165)
\path (196,151)(204,157)
\path (196,143)(204,149)
\path (196,135)(204,141)
\path (196,127)(204,133)
\path (196,111)(204,117)
\path (247,128)(255,134)
\path (251,114)(251,144)(251,174)
\path (247,168)(255,174)
\path (247,160)(255,166)
\path (247,152)(255,158)
\path (247,144)(255,150)
\path (247,136)(255,142)
\path (247,120)(255,126)
\path (251,118)(251,87)
\path (247,112)(255,118)
\path (247,104)(255,110)
\path (247,96)(255,102)
\path (247,88)(255,94)
\dottedline{3}(158,79)(158,55)
\put(163,67){{\rm\rm {$q^1=t$}}}
\put(152,38){{\rm\rm {$t$}}}
\path (395,72)(403,78)
\path (399,18)(399,78)
\path (395,64)(403,70)
\path (395,56)(403,62)
\path (395,48)(403,54)
\path (395,40)(403,46)
\path (395,32)(403,38)
\path (395,24)(403,30)
\path (434,88)(442,94)
\path (434,96)(442,102)
\path (438,42)(438,102)
\path (434,80)(442,86)
\path (434,72)(442,78)
\path (434,64)(442,70)
\path (434,56)(442,62)
\path (434,48)(442,54)
\path (492,138)(500,144)
\path (492,98)(500,104)
\path (492,130)(500,136)
\path (492,122)(500,128)
\path (492,114)(500,120)
\path (492,106)(500,112)
\path (492,90)(500,96)
\path (496,88)(496,57)
\path (492,82)(500,88)
\path (492,74)(500,80)
\path (492,66)(500,72)
\path (492,58)(500,64)
\path (496,88)(496,176)
\path (492,146)(500,152)
\path (492,154)(500,160)
\path (492,162)(500,168)
\path (492,170)(500,176)
\put(134,-1){{\rm\rm {\underline{elliptic case}}}}
\put(380,0){{\rm\rm {\underline{hyperbolic case}}}}
\end{picture}

\end{center}
\caption{}\label{f1}
\end{figure}

\paragraph{Geometric interpretation}
The reference trajectory $\gamma $ is $C^0$-time minimal
(resp. time maximal) in the hyperbolic case (resp. elliptic
case) up to a time $t_{1c}$ which corresponds to a
\it{first conjugate time} $t_{1c}>0$ along $\gamma $ for the
time minimal (resp. time maximal) control problem.

In particular we get the following Proposition~:

\begin{prop}
Assume $T<t_{1c}.$ Then the reference singular trajectory
$\gamma $ defined on $[0,T]$ is in the hyperbolic (resp.
elliptic) case the only trajectory ${\bar \gamma }$ contained
in a $C^0$-neighborhood of $\gamma $ and satisfying the boundary
conditions~: ${\bar \gamma } (0)=\gamma (0),$ ${\bar \gamma } 
({\bar T}) =\gamma (T)$ in a time ${\bar T}\leq T$ (resp. ${\bar
 T}\geq T).$

This property is called $C^0$-\it{one-side rigidity}, compare
with \cite{AS}.
\end{prop}

\paragraph{Exceptional case \\}
We proceed as before. The model is~:
\begin{eqnarray*}
{\dot q}^1& = & 1+q^2 ,\ {\dot q}^2=q^3,\ \ldots ,{\dot q}^{n-1}=u\\
{\dot q}^n & = & {\displaystyle {\sum ^{n-1}_{i,j=2}}}
 a_{ij}(t) q^i q^j
\end{eqnarray*}

Let $0<t, t'\leq T$ and consider the following boundary
conditions~: $q(0)=0$ and $q^1(t')=t, q^2(t')=\cdots =q^{n-1}
(t')=0.$ We get a projection of the accessibility set $A(0,t')$
on the line $q^n.$ It is represented on Fig. \ref{f2}.

\unitlength=.25mm
\makeatletter
\def\shade{\@ifnextchar[{\shade@special}{\@killglue\special{sh}\ignorespaces}}
\def\shade@special[#1]{\@killglue\special{sh #1}\ignorespaces}
\makeatother

\begin{figure}[h]
\begin{center}

\begin{picture}(537,184)(100,-5)
\thinlines
\typeout{\space\space\space eepic-ture exported by 'qfig'.}
\font\FonttenBI=cmbxti10\relax
\font\FonttwlBI=cmbxti10 scaled \magstep1\relax
\path (100,76)(262,76)
\path (179,76)(179,163)
\path (254.482,78.736)(262,76)(254.482,73.264)
\path (176.264,155.482)(179,163)(181.736,155.482)
\put(160,162){{\rm\rm {$q^n$}}}
\put(169,56){{\rm\rm {$t$}}}
\put(250,58){{\rm\rm {$t'$}}}
\path (432,76)(432,163)
\path (353,76)(515,76)
\put(413,162){{\rm\rm {$q^n$}}}
\put(422,56){{\rm\rm {$t$}}}
\put(503,58){{\rm\rm {$t'$}}}
\path (429.264,155.482)(432,163)(434.736,155.482)
\path (507.483,78.736)(515,76)(507.483,73.264)
\put(179,134.346){\arc{116.692}{.109}{3.033}}
\put(432.5,9.174){\arc{133.656}{3.459}{5.966}}
\path (136,128)(136,96)
\path (152,127)(152,84)
\path (167,127)(167,78)
\path (190,126)(190,79)
\path (206,127)(206,84)
\path (221,126)(221,95)
\path (371,113)(371,38)
\path (388,114)(388,59)
\path (406,114)(406,71)
\path (423,114)(423,76)
\path (442,114)(442,75)
\path (458,114)(458,71)
\path (476,114)(476,60)
\path (494,114)(494,35)
\put(397,-1){{\rm\rm {$t_{1cc}<t<t_{1cc}+\varepsilon$}}}
\put(149,-1){{\rm\rm {$t<t_{1cc}$}}}
\end{picture}

\end{center}
\caption{}\label{f2}
\end{figure}

\paragraph{Geometric interpretation}
The reference trajectory $\gamma $ is $C^0$-time optimal up to
a time $t_{1cc}$ which corresponds to a first conjugate time
$t_{1cc}>0.$ In particular we have the following result.

\begin{prop}
Assume $T<t_{1cc}$. Then the reference singular exceptional
trajectory $\gamma $ is $C^0$-isolated (or $C^0$-rigid).
\end{prop}


\subsubsection{Conclusion : the importance of singular
trajectories in optimal control}
The previous analysis shows that singular trajectories
play generically an important role in any optimal control
problem~: Min~${\displaystyle {\int ^T_0}}\ f^0(x,u)dt$ when the
transfert time $T$ is fixed. Indeed they are locally
\it{the only trajectories satisfying the boundary conditions}
and hence are optimal. If the transfert time $T$ is not fixed only
exceptional singular trajectories play a role. In
fact as observed by \cite{AS} they correspond to the singularities
of the \it{time extended
end-point mapping~:} ${\bar E}:(T, u) 
\longmapsto q(T,x_0,u).$ It is the situation
encountered in sub-Riemannian geometry.


\section{Generalities about sub-Riemannian geometry}
From now on, we work in the $C^\omega$-category.

\begin{defi}
A SR-manifold is defined as a $n$-dimensional manifold $M$
together with a distribution $D$ of constant rank $m\leq n$ and
a Riemannian metric $g$ on $D.$ An admissible curve $t\longmapsto
q(t),$ $0\leq t\leq T$ is an absolutely continuous curve such
that ${\dot q}(t)\in D(q(t))\backslash \{ 0\}$ for almost every
$t$. The length and the energy of $q$ are respectively defined
by~:
$$L(q)=\int^T_0 ({\dot q}(t), {\dot q}(t))^{1/2} dt ,
\esp E(q)=\int^T_0 ({\dot q}(t), {\dot q}(t))dt$$

\noindent where $(\ ,\ )$ is the scalar product defined by $g$
on $D.$ The SR-distance between $q_0, q_1\in M$ denoted
$d_{SR}(q_0,q_1)$ is the infimum of the lengths of the admissible
curves joining $q_0$ to $q_1.$
\end{defi}


\subsection{Optimal control formulation}
The problem can be \it{locally} restated as follows. Let $q_0 
\in M$ and choose a coordinate system $(U,q)$ centered at $q_0$
such that there exist $m$ (smooth) vector fields $\{ F_1,\ldots 
,F_m\}$ which form an orthonormal basis of $D.$
Then each admissible curve $t\rightarrow q(t)$ on $U$ is
solution of the control system~:
\begin{eqnarray}    \label{3.1}
{\dot q}(t)={\displaystyle {\sum ^m_{i=1}}} u_i (t) F_i (q(t))
\end{eqnarray}

The length of a curve does not depend on its parametrization,
hence every admissible curve can be reparametrized into a
lipschitzian curve $s\longmapsto q(s)$ \it{parametrized by
arc-length}~: $({\dot q} (s), {\dot q} (s))=1,$ see \cite{LS}.

If an admissible curve on $U$~: $t\longmapsto q(t), 0\leq t\leq
T$
is parametrized by arc-length we have almost everywhere~:
$${\dot q} (t)=\sum ^m_{i=1} u_i F_i (q(t)) ,\ \ \sum ^m_{i=1}
 u^{2}_i (t)=1$$

\noindent and $L(q)=\int ^T_0$ $(\sum u^{2}_i )^{1/2} \ dt=T.$
Hence the length minimization problem is equivalent to a
\it{time-optimal problem} for system (\ref{3.1}). This problem
is not convex because of the constraints~: ${\displaystyle 
{\sum ^m_{i=1}}} u^{2}_i = 1$, but it is well-known that the
problem is equivalent to a time optimal control problem with
the convex constraints~: ${\displaystyle {\sum ^m_{i=1}}}
 u^{2}_i (t)\leq 1.$

It is also well-known that if every curve is parametrized on a
\it{fixed} interval $[0,T]$, the length minimization problem is
equivalent to the energy minimization problem.

Introducing the extended control system~:
\begin{eqnarray} \label{3.2}
{\dot q}(t) & = & {\displaystyle {\sum ^m_{i=1}}} u_i (t)
 F_i (q(t))\nonumber \\
{\dot q}^0 (t) & = & {\displaystyle {\sum ^m_{i=1}}} u^2_i (t)
 ,\ \ q^0 (0)=0 
\end{eqnarray}

\noindent and the end-point mapping ${\widetilde E}$ of the
extended system~: $u\in {\cal U}\longmapsto {\widetilde q}
(t,u,{\widetilde q}_0)$~, ${\widetilde q}=(q,q^0), {\widetilde 
q}(0)=(q_0,0)$, from the maximum principle the minimizers can
be selected among the solutions of the maximum principle~:
$${\dot {\widetilde q}}(t)={\partial {\widetilde H}\over
\partial {\widetilde p}}\ ,\ {\dot {\widetilde p}}(t)=-
{\displaystyle {\partial {\widetilde H}\over \partial 
{\widetilde q}}}\ ,\ {\partial {\widetilde H}\over \partial u}\ =0$$

\noindent where ${\widetilde H}=<p, {\displaystyle {\sum
^m_{i=1}}}
 u_i  F_i (q)>+p_0  {\displaystyle {\sum ^m_{i=1}}} 
u_i^2$ is the pseudo-Hamiltonian and ${\widetilde p}
=(p,p_0) \in \R^{n+1}\backslash \{ 0\}$ is the adjoint vector
of the extended system. From the previous equation
$t\longmapsto p_0 (t)$ is a constant which can be normalized
to $0$ or $-1/2.$ Introduce $H_\nu =<p, 
{\displaystyle {\sum ^m_{ i=1}}}u_i  F_i (q)>+\nu 
{\displaystyle {\sum ^m_{i=1}}} u_i^2$ where $\nu =0$ or $-1/2$~;
then the previous equations are equivalent to~:
\begin{eqnarray}
{\dot q}={\displaystyle {\partial H_\nu \over \partial p}}\
 ,\ {\dot p}=-{\displaystyle {\partial H_\nu \over \partial q}}\
 ,\ {\displaystyle {\partial H_\nu \over \partial u}}\ =0.
\end{eqnarray}

The solutions of these equations correspond to the
singularities of the end-point mapping of the extended system
and are called \it{geodesics} in the framework of SR-geometry.

They split into two categories according to the following
definition.

\begin{defi}\label{defistrict}
A geodesic is said to be \it{abnormal}
if $\nu =0$, and \it{normal} if $\nu =-1/2$. Abnormal
geodesics are precisely the singular trajectories of the
original system (\ref{3.1}).

A geodesic is called \it{strict} if
the extended adjoint vector $(p,p_0=\nu )$ is unique up to a
scalar, that is corresponds to a \it{singularity of codimension
one} of the extended end-point mapping.
\end{defi}


\subsection{Computations of the geodesics}

\subsubsection{Abnormal case}
They correspond to $\nu =0,$ and are the singular trajectories
of system (\ref{3.1}). The system is symmetric and hence
$H={\displaystyle {\sum ^m_{i=1}}} u_i P_i,$ with
$P_i=<p, F_i(q)>.$ Therefore they are \it{exceptional.}
When $m=2,$ they are computed using the algorithm of Section 2.
The case $m>2$ will be excluded in our forthcoming analysis
because from \cite{AS} in order to be optimal a singular trajectory
must satisfy the following conditions, known as Goh's conditions~:
\begin{eqnarray}
<p(t), [Fv,Fw] (q(t))>=0
\end{eqnarray}

\noindent $\forall  v, w\in \R^m, \ \forall t\in [0,T]$ and
$Fu$ denotes $\sum u_i F_i.$ If $m=2,$ this reduces to the
condition $\{ P_1, P_2\}=0$ deduced from the conditions
$P_1=P_2=0$ but if $m>2$ it is a very restrictive condition
which should not be generic (conjecture \cite{H}).


\subsubsection{Normal case}
They correspond to $\nu =-1/2.$ If the system of the
$F_i$'s is orthonormal
then ${\displaystyle {\partial H_\nu \over \partial u}}=0$ and
hence $u_i=P_i$ and $H_\nu $ reduces to $H_n ={\displaystyle
 {1\over 2}}$ ${\displaystyle {\sum ^m_{i=1}}} P^2_i.$ The
trajectories parametrized by arc length are on the level set
$H_n =1/2$ and the normal geodesics are solutions of the
following \it{Hamiltonian differential equations}~:
\begin{equation}  \label{3.4}
{\dot q}={\partial H_n\over \partial p}\ ,\ \ {\dot p}=-
{\partial H_n\over \partial q}
\end{equation}
On the domain chart $U,$ we can complete the $m$-vector
fields $\{ F_1,\ldots ,F_m\}$ to form a smooth frame $\{
F_1,\ldots
 ,F_n\}$ of $TU.$ The SR-metric $g$ can be extended into a
Riemannian metric by taking the system of the
$F_i$'s as an orthonormal frame.
We set $P_i=<p, F_i(q)>$ for $i=1,\ldots ,n$ and let
$P=(P_1,\cdots ,P_n).$ In the coordinate system $(q,P)$ the
normal geodesics are solutions of the following equations~:
\begin{eqnarray}   \label{3.5}
{\dot q}& = & {\displaystyle {\sum ^m_{i=1}}} P_iF_i (q)\nonumber \\
{\dot P}_i & = & \{ P_i,H_n\}={\displaystyle {\sum ^m_{j=1}}}
 \{ P_i,P_j\} P_j 
\end{eqnarray}

We observe that $\{ P_i ,P_j\} =<p,[F_i,F_j](q)>$ and since
the $F_i$'s form a frame we can write~:
$$[F_i ,F_j](q)={\displaystyle {\sum ^n_{k=1}}}
 c^k_{ij}(q) F_k(q)$$

\noindent where the $c^k_{ij}$'s are smooth functions.


\subsection{Exponential mapping - Conjugate and cut loci} 
Assume that the curves are parametrized by arc-length. If
$t\longmapsto q(t)$ is any geodesic, the first point where
$q(.)$ ceases to be minimizing is called a \it{cut-point}
and the set of such points when we consider all the geodesics
with $q(0)=q_0$ will form the \it{cut-locus} $L(q_0).$

The \it{sub-Riemannian sphere} with radius $r>0$ is the set
$S(q_0,r)$ of points which are at SR-distance $r$ from $q_0$.
The \it{wave front} of length $r$ is the set $W(q_0,r)$ of
end-points of geodesics with length $r$ starting from $q_0.$
If $D_{A.L.} (q_0)$ is of rank $n$ where $D_{A.L}$ is
the Lie algebra generated by $D,$ then according to
\it{Filippov's
existence Theorem} \cite{LM} if $r$ is small enough each point of
distance $r$ from $q_0$ is the end-point of a minimizing geodesic
and $S(q_0,r)$ is a subset of $W(q_0,r).$
We fix $q_0 \in U$ and let $(q(t, q_0, p_0), p(t, q_0, p_0))$ be
the normal geodesic, solution of (\ref{3.4}) and starting
from $(q_0,p_0)$ at $t=0.$ The \it{exponential mapping} is the map~:
$$\rm{exp}_{q_0}\ :\ (p_0,t)\longmapsto q(t,p_0,q_0)$$

\noindent Its domain is the set $C\times \R $ where $C$ is
$\{ p_0\ ;\ {\displaystyle {\sum ^m_{i=1}}} P^2_i (p_0,q_0)=1\}.$
If $m<n$ it is a (non compact) \it{cylinder} contrarily to the
Riemannian case~: $m=n,$ where it is a \it{sphere}.

A \it{conjugate point} along a normal geodesic is defined as
follows. Let $(p_0,t_1)$ with $t_1>0$ be a point where exp$_{q_0}$
is not an immersion. Then $t_1$ is called a conjugate time along
the normal geodesic and the image is called a conjugate point.
The \it{conjugate locus} $C(q_0)$ is the set of \it{first
conjugate points}.


\subsection{Gradated normal form}
\subsubsection{Adapted and priviliged coordinate system}
Let $(U,q)$ be a coordinate system centered at $q_0,$ with
$D=\rm{Span}\ \{ F_1,\ldots ,F_m\}.$ Assume that $D$ satisfies
the rank condition on $U.$ We define recursively~: $D_0=\{ 0\} , 
D^1=D$ and for $p\geq 2$ $D^p=\rm{Span }\{ D^{p-1}+[D^1,
D^{p-1}]\}.$ Hence $D^p$ is generated by Lie brackets of the
$F_i$'s with length $\leq p.$ At $p$ we have an increasing
sequence of vector sub-spaces~: $\{0\} =D^0(q)\subset D^1(q)\subset
\cdots \subset D^{r(q)}$ where $r(q)$ is the smallest integer
such
that
$D^{r(q)}(q)=T_qU.$

\begin{defi}
We say that $q_0$ is a \it{regular point} if
the integers $n_p(q)=\rm{dim\ }D^p(q)$ remain constant for $q$
in some neighborhood of $q_0.$ Otherwise we say that $q_0$ is a
\it{singular point.} Consider now a coordinate system $(q^1,\ldots
,q^n)$ such that $dq^j$ vanishes identically on $D^{w_j-1}(q_0)$
and doesn't vanish identically on $D^{w_j}(q_0)$ for some integer
$w_j.$ Such a coordinate system is said to be \it{adapted} to the
flag and the integer $w_j$ is the \it{weight} of $q^j .$
\end{defi}

\begin{defi}
Consider now a SR-metric $(D,g)$ defined on the chart $(U,q)$ and
represented locally by the \it{orthonormal} vector fields
$\{ F_1,\ldots ,F_m\}.$ If $f$ is a germ of smooth function at
$q_0,$ the \it{order} of $f$ at $q_0$ is~:

\begin{itemize}
\item[(i)] if $f(q_0)\not = 0, \mu (f)=0, \mu (0)=+\infty $~;

\item[(ii)] otherwise~: $\mu (f)=\rm{inf}\ \{ p\ / \  \exists 
\ V_1,\ldots ,V_p\in \{ F_1,\ldots ,F_m\} \rm{ with }
L_{V_1}\circ \cdots \circ L_{V_p} (f)(q_0)\not = 0 \}$ where $L_V$
denotes the Lie derivative. The germ $f$ is called
\it{privileged} if $\mu (f)=\min \{ p\ ;\ df(D^p(q_0))\not =
0\}.$
A coordinate system $\{ q^1,\ldots ,q^n\}$ is said to be
\it{privileged} if all the coordinates $q_i$ are privileged
at $q_0$. 
\end{itemize}
\end{defi}

We have the following very important estimation, see \cite{Be},
\cite{Ku2}~:
\begin{prop}
If $(M,D,g)$ is a SR-manifold there exists a privileged
coordinate
system $q$ at every point $q_0 =0$ of $M.$ If $w_i$ is the
order (or weight) of the coordinate $q^j$ we have the following
estimation for the SR-distance~:
$$d_{SR}(0,(q^1,\ldots ,q^n))\simeq |q^1|^{1/w_1}+\cdots
 +|q^n|^{1/w_n}.$$
\end{prop}
\begin{defi}
Let $(U,q)$ be a privileged coordinate system for the
SR-structure given locally by the $m$-orhonormal vector
fields~: $\{ F_1,\ldots ,F_m\}.$ If $w_j$ is the weight of
$q^j,$ the weight of ${\displaystyle {\partial \over \partial
 q^j }}$ is taken by convention as $-w_j.$ Every vector field
$F_i$ can be expanded into a Taylor series using the previous
gradation and we denote by ${\hat F}_i$ the homogeneous term with
lowest order $-1.$ The polysystem $\{ {\hat F}_1,\ldots ,
{\hat F}_m\}$ is called the \it{principal part} of the SR-structure.
\end{defi}

\no We have the following result, see \cite{Be}.

\begin{prop}
The vector fields ${\widehat F}_i, i=1,\ldots ,m$ generate
a nilpotent Lie algebra which satisfies the rank condition.
This Lie algebra is independant of the privileged coordinate system.
\end{prop}


\subsubsection{Gauge classification} 

Given a local SR-geometry $(U,D,g)$ represented as the optimal
control problem~:
$${\dot q}={\sum ^m_{i=1}} u_i F_i (q) $$
$$\underset{u(.)}{\min} \int^T_0\
 \Bigl (\sum ^m_{i=1} u_i^2(t)\Bigr )dt\ ,$$

\noindent there exists a pseudo-group of transformations called
the \it{gauge group} which is the subgroup of the feedback group
defined by the following transformations~:

\begin{itemize}
\item[(i)] germs of diffeomorphisms $\varphi : q\longmapsto Q$
on $U,$ preserving $q_0$~; 

\item[(ii)] feedback transformations $u=\beta (q)v$ preserving
the metric $g$ i.e, $\beta (q)\in \theta (m,\R)$ (orthogonal group).
\end{itemize}

The \it{invariants} of the associated classification problem
are the geodesics. They split into two categories~: abnormal
geodesics which are feedback invariants and normal geodesics.

If $q$ is an adapted coordinate system, a \it{gradated normal
form of order $p\geq -1$} is the polysystem $\{ F^p_1,\ldots
,F^p_m\}$ obtained by truncating the vector fields $F_i$ at
order $p$ using the weight system defined by the adapted
coordinates.


\section{The role of abnormal minimizers in SR
Martinet geometry}
In this Section we analyze the role of abnormal minimizers in
SR Martinet geometry which is the prototype of the generic rank
2 situation. Before to present this analysis it is important to
make a short visit to the contact situation in $\R^3.$

\subsection{The contact situation in $\R^3$} 

The contact situation in $\R^3$ has been analyzed in details in
several articles \cite{AEG,AEGK}. This analysis is based on
computations about the exponential mapping using a gradated
normal form. To understand the remaining of this article it is
important to make the contact situation fit into the following
framework.

First, without losing any generality we can use to understand
a generic contact SR-problem a gradated form of order 1 computed
in \cite{AEGK} where the SR-metric is defined by the two
orthonormal
vector fields~: $F_1, F_2$ where~:
$$F_1={\partial \over \partial x}\ +\ {y\over 2}(1+Q)
{\partial \over \partial z}\ ,\ \ F_2 ={\partial \over 
\partial y}\ -{x\over 2}(1+Q) {\partial \over \partial z}$$

\noindent where $Q$ is a quadratic form~: $ax^2+2bxy+cy^2$
depending on 3 parameters. The weight of $x,y$ is one and the
weight of $z$ is two. When $a=b=c=0,$ it corresponds to the
contact situation of order -1 which is the well-known
\it{Heisenberg case} but \it{also a gradated normal form of order
0.} 

To get an adapted frame we complete $F_1, F_2$ by
$F_3={\displaystyle {\partial \over \partial z}}.$ Computing
we get~: $[F_1, F_2]=(1+2Q){\displaystyle {\partial \over 
\partial z}}.$ Using $P_i=<p, F_i(q)>,$ the geodesics equations are~:
\begin{eqnarray*}
{\dot x}&=&P_1 \\
{\dot y}&=&P_2 \\ 
{\dot z}&=&{P_1y(1+Q)\over 2} \ 
 -\ {P_2 x(1+Q)\over 2} \\
{\dot P}_1 & = & \{ P_1,P_2\} P_2 \ =\ (1+2Q)\ P_2P_3\\
{\dot P}_2 & = & \{ P_2,P_1\} P_1 \ =\  -(1+2Q)P_1P_3\\
{\dot P}_3 & = & 0
\end{eqnarray*}

In the Heisenberg case we have $Q=0$, and if we set $P_3=\lambda$
we get~: $\ddot{P}_1+\lambda^2 P_1=0$, which is a \it{linear
pendulum}.

Using the \it{cylindric coordinates~:} $P_1=\sin \theta ,
P_2=\sin\theta,
P_3 =\lambda ,$ where $\theta \not = k\pi ,$ the geodesics
\it{parametrized by arc-length} are solutions of the following
equations~:
\begin{eqnarray}  \label{4.1}
{\dot x} = P_1\  \  \ \ {\dot y}&=&P_2 \  \  \ {\dot z}=
{\displaystyle {P_1y(1+Q)\over 2}}\ -{\displaystyle {P_2
 x(1+Q)\over 2}}\nonumber \\
\\
\dot{\theta} & = & (1+2Q)\lambda \nonumber
\end{eqnarray}

\noindent where $\lambda $ is a constant. The important
behavior is when $\lambda \rightarrow \infty .$ We may assume
$\lambda >0.$ By making the following \it{reparametrization}~:
\begin{eqnarray}  \label{4.2}
ds=\lambda (1+2Q)dt
\end{eqnarray}

\noindent the angle equation takes the trivial form~:
${\displaystyle {d\theta \over ds}}=1.$ Hence it is
\it{integrable}
and we obtain $\theta (s)=s+\theta _0.$

The remaining equations take the form~:
\begin{eqnarray}  \label{rema}
{\displaystyle {dx\over ds}}& = & {\displaystyle {{\sin}\ 
\theta (s)\over (1+2Q)\lambda }}\nonumber \\
{\displaystyle {dy\over ds}} & = & {\displaystyle {{\cos}\ 
\theta (s)\over (1+2Q)\lambda }}\\
{\displaystyle {dz\over ds}}& =& {\displaystyle {{\sin} 
\theta (s)(y(1+Q))-{\rm cos\ }\theta (s)(x(1+Q))\over 
(1+2Q)\lambda }}\nonumber
\end{eqnarray}

\no For large $\lambda ,$ they can be integrated as follows.
We set $\varepsilon =1/\lambda :$ small parameter,
$x=\varepsilon X, y=\varepsilon Y,$ $z=\varepsilon ^2Z, 
{\displaystyle {1\over 1+2Q}}=1+{\widetilde Q}=1+Ax^2+2Bxy+Cy^2+
\cdots $ and we get~:
\begin{eqnarray*}
{\dot X}& =&{\sin}\ (s+\theta _0)[1+\varepsilon ^2 
{\widetilde Q}(X,Y)+\rm{o} (\varepsilon ^2)]\\
{\dot Y}& =& {\cos}\ (s+\theta _0)[1+\varepsilon ^2 
{\widetilde Q}(X,Y)+\rm{o} (\varepsilon ^2)]\\
{\dot Z}& = &{\displaystyle {{\sin}\ (s+\theta  _0)Y-
{\cos}(s+\theta_0)X\over 2}}\ +\rm{o} (\varepsilon )
\end{eqnarray*}

The previous equations \it{can be integrated by quadratures}
by setting~:
\begin{eqnarray*}
X & = & X_0 +\varepsilon ^2X_1 +\rm{o} (\varepsilon ^2)\\
Y & = & Y_0 +\varepsilon ^2Y_1+\rm{o} (\varepsilon ^2 )\\
Z& =& Z_0 +\rm{o} (\varepsilon )
\end{eqnarray*}

\noindent and we get in particular
\begin{eqnarray*}
{\dot X}_0 & = &{\sin} (s+\theta _0)\\
{\dot Y}_0 & = & {\cos} (s+\theta _0)\\
{\dot Z}_0 & = & {\displaystyle {{\sin} (s+\theta_0)Y_0(s)-
{\cos\ }(s+\theta _0)\ X_0 (s)\over 2}}\\
{\dot X}_1 & = & {\sin } (s+\theta _0 )\ {\widetilde Q}(X_0,Y_0)\\
{\dot Y}_1 &= & {\cos }(s+\theta_0)\ {\widetilde Q}(X_0,Y_0).
\end{eqnarray*}

\noindent The solutions are computed in the $s$-parametrization
and the arc-length $t$ can be computed by integrating (\ref{rema})
by quadratures.

If we want to mimic this procedure in the Martinet situation,
we shall encounter \it{integrability obstructions due to the
existence of abnormal geodesics.}

\begin{figure}[h]
\epsfxsize=7cm
\epsfysize=7cm
\centerline{\epsffile{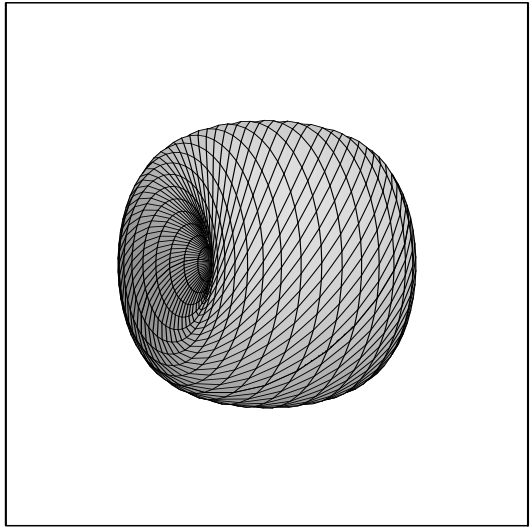}}
\caption{SR sphere in the flat contact case}\label{figcontact}
\end{figure}

The sphere in the flat contact situation is represented on Fig.
\ref{figcontact}.


\subsection{The Martinet situation}
\subsubsection{Normal forms and invariants}

The Martinet SR-geometry is rather intricate and it is difficult
to make a priori normalizations. It will appear later that a good
starting point to make the computations is to use the following
normal form computed in \cite{ABCK}~:
\begin{itemize}
\item The distribution $D$ is taken in the Martinet-Zhitomirski
normal form~: $D=\rm{Ker }\omega, \ \omega =dz-{\displaystyle 
{y^2\over 2}}dx$.

\item The metric on $D$ is taken as a \it{sum of squares}~:
$a(q)dx^2+c(q)dy^2.$
\end{itemize}
In this representation the Martinet surface containing the
abnormal geodesics is the plane~: $y=0$ and the abnormal
geodesics
are the straight-lines~: $z=z_0.$ The abnormal line passing
through $0$ is given by $\gamma : t\longmapsto (\pm t,0,0).$

The computations in \cite{ABCK} show that we can make an additional
normalization on the metric by taking either the restriction
of $a$ or $c$ to the Martinet plane $y=0$ equal to $0.$

The variables are gradated according to the following weights~:
the weight of $x,y$ is one and the weight of $z$ is three.
By identifying by convention at order $p$ two normal forms
where the Taylor series of $a$ and $c$ coincide at order $p$
we end up with the following representatives of \it{order $0$}~:

\no either
$$g=(1+\alpha y)^2 dx^2+(1+\beta x+\gamma y)^2 dy^2$$

\no or
$$g=(1+{\bar \alpha }x+{\bar \beta }y)^2 dx^2+
(1+{\bar \gamma }y)^2 dy^2.$$

\noindent In each of those representations the \it{three
parameters are,
up to sign, invariants.} They can be used to compute the
exponential mapping in the generic situation. If we truncate
$g$ to $dx^2+dy^2$ it corresponds to the principal part of
order $-1$ of the SR-structure defined previously. In the
sequel it will be called the \it{flat case.}


\subsubsection{Geodesics equations}
The distribution $D$ is generated by~:
$$G_1 ={\partial \over \partial x}\ +\ {y^2\over 2} 
{\partial \over \partial z}\ \ \ \rm{and}\ \ \ G_2 =
{\partial \over \partial y}$$

\noindent and the metric is given by $g=adx^2+cdy^2.$
We introduce the frame~:
$$F_1 ={1\over \sqrt{a}} G_1\ ,\ \ F_2={1\over \sqrt{c}}
 G_2\ ,\ \ F_3={\partial \over \partial z}$$

\noindent and $P_i=<p, F_i(q)>$ for $i=1,2,3,$ i.e
$$P_1 ={p_x+p_z y^2/2\over \sqrt{a}}\ ,\ \ P_2={p_y 
\over \sqrt{c}}\ ,\ \ P_3=p_z.$$

First, we assume that $g$ is not depending on $z$~; this is
the case for the gradated normal form of order 0. It corresponds
to an \it{isoperimetric situation}, that is the existence of
a vector field $Z$ identified here to ${\displaystyle {\partial 
\over \partial z}}$ transverse at 0 to $D(0)$ and the metric
$g$ does not depend on $z.$

The system is written~:
$${\dot x}=u_1\ ,\ \ {\dot y}=u_2\ ,\ \ {\dot z}={y^2\over 2} u_1$$

\noindent and the Hamiltonian associated to normal geodesics is~:
$$H_n (q,p)=\ {1\over 2} (u^2_1 a+u^2_2c)$$

\noindent and the geodesics controls are~:
$$u_1 ={1\over a}\ (p_x+p_z y^2/2)\ ,\ \ u_2={p_y\over c}.$$

Normal geodesics are solutions of the following equations~:
\begin{equation}  \label{4.4}
\begin{split}
{\dot x} & =  {\displaystyle {1\over a}}
 (p_x+p_z\ y^2/2),\
{\dot y}  =  {\displaystyle {p_y\over c}},\ 
{\dot z}  =  {\displaystyle {y^2\over 2a}}
 (p_x+p_z \ y^2/2) \\
{\dot p}_x & =  {\displaystyle {p_y^2 c_x\over 2c^2}}\ +\ 
{\displaystyle {(p_x +p_z\ y^2/2)^2\over 2a^2}} a_x \\
{\dot p}_y & =  {\displaystyle {p_y^2 c_y\over 2c^2}}\ +\ 
{\displaystyle {(p_x +p_z\ y^2/2)^2\over 2a^2}} a_y-
{\displaystyle {(p_x+p_z\ y^2/2)\over a}} p_zy  \\
{\dot p}_z & =  0
\end{split}
\end{equation}

In the $(q,P)$ representation the previous equations take the form~:
\begin{equation}   \label{4.5}
\begin{split}
{\dot x} & =  {\displaystyle {P_1\over \sqrt{a}}},\
{\dot y}  =  {\displaystyle {P_2\over \sqrt{c}}}, \
{\dot z}  = {\displaystyle {y^2\over 2}}\ {\displaystyle
 {P_1\over \sqrt{a}}} \\
{\dot P}_1 & =  {\displaystyle {P_2\over \sqrt{a}\sqrt{c}}}\
 \Bigl (yP_3 -{\displaystyle {a_y\over 2\sqrt{a}}} P_1 +
{\displaystyle {c_x\over 2\sqrt{c}}}\ P_2\Bigr ) \\
{\dot P}_2 & =  -{\displaystyle {P_1\over \sqrt{a}\sqrt{c}}}\ 
\Bigl ( yP_3 -{\displaystyle {a_y\over 2\sqrt{a}}} P_1 +
{\displaystyle {c_x\over 2\sqrt{c}}}\ P_2\Bigr ) \\
{\dot P}_3 & =  0 
\end{split}
\end{equation}

If we parametrize by arc-length and if we introduce the
cylindric coordinates~: $P_1=\cos \theta , P_2 =\sin \theta ,
P_3=\lambda ,$ we end up with the following equations~:
\begin{eqnarray}  \label{4.6}
{\dot x} & = & {\displaystyle {P_1\over \sqrt{a}}}\ ,\ \ 
{\dot y}={\displaystyle {P_2\over \sqrt{c}}}\ ,\ \ {\dot z} 
={\displaystyle {y^2\over 2}} {\displaystyle {P_1\over 
\sqrt{a}}}\nonumber \\
{\dot \theta } & = & -{\displaystyle {1\over \sqrt{a}\ 
\sqrt{c}}}\ \Bigl [ y\ P_3 -{\displaystyle {a_y\over 2\sqrt
{a}}}\ P_1\ +\ {\displaystyle {c_x\over 2\sqrt{c}}}\ P_2\Bigr ]\\
P_3 & = & \lambda \nonumber
\end{eqnarray}

It is proved in \cite{BH} that for a generic SR-problem, each
geodesic is strict. In our representation we have the following
result.
\begin{lem}
The abnormal geodesic $\gamma :t\longmapsto (\pm t, 0, 0)$
is strict if and only if the restriction of $a_y$ to the
Martinet plane $y=0$ is $0.$
\end{lem}
 
Using the gradated normal form of order 0 with the normalizations~:
$$a=(1+\alpha y)^2\ \ \ ,\ \ \ c=(1+\beta x+\gamma y)^2$$

\noindent the equations (\ref{4.6}) reduce to~:
\begin{equation}   \label{4.7}
\begin{split}
{\dot x} & =  {\displaystyle {{\cos}\ \theta \over
 \sqrt{a}}}\ \ \ {\dot y} = {\displaystyle {{\sin}\ \theta \over 
\sqrt{c}}}\ \ \ {\dot z}={\displaystyle {y^2\over 2}}
{\displaystyle {{\cos}\ \theta \over \sqrt{a}}} \\
{\dot \theta } & =  -{\displaystyle {1\over \sqrt{a} 
\sqrt{c}}}\ [y\lambda -\alpha {\cos\ }\theta +\beta 
 {\sin\ }\theta ] 
\end{split}
\end{equation}

The previous equation defines a \it{foliation} $({\cal F})$ of
codimension one in the plane $(y,\theta ).$ Indeed using the
parametrization~: $\sqrt{a}\ \sqrt{c}\ {\displaystyle
 {d\over dt}}={\displaystyle {d\over d\tau}}$ and denoting $'$ the
derivative with respect to $\tau ,$ the equations can be written~:
\begin{eqnarray}   \label{4.8}
x' & = & \sqrt{c}\ {\cos}\ \theta \ \ \ \ 
z' =\sqrt{c}\ {\displaystyle {y^2\over 2}}\ {\cos}\ 
\theta \nonumber \\
\\
y' & = & \sqrt{a}\ {\sin}\ \theta  \ \ \ \ 
\theta ' = - [y\lambda - \alpha {\cos}\ \theta +\beta 
{\sin\ }\theta ]\nonumber
\end{eqnarray}

\noindent and they can \it{be projected} onto the plane
$(y,\theta ).$ The last two equations are e\-qui\-va\-lent to~:
\begin{eqnarray}   \label{4.9}
\theta ''+\lambda {\sin\ }\theta +\alpha ^2{\sin\ }\theta 
 {\cos}\ \theta -\alpha \beta {\sin}^2 \theta +\beta \theta '
{\cos\ }\theta =0.
\end{eqnarray}

This equation will be used in the sequel to study the SR-Martinet
geometry in the generic case of order $0$.
Unfortunately it depends on the
choice of coordinates.
Note that in the flat case where $a=c=1$ the equation reduces to
$\theta''+\lambda\sin\theta=0$ which is a \it{nonlinear
pendulum}.


\subsubsection{Conservative case}
The analysis of Subsection 4.1 shows that in the contact case
the equation (\ref{4.1}) associated to the evolution of $\theta $
defines an integrable foliation. In the Martinet case the
foliation defined by equation (\ref{4.9}) is not in general
integrable. This leads to the following definition which is
independant of the choice of coordinates.

\begin{defi}
Let $e(t,\theta ,\lambda )$ be a normal geodesic parametrized
by arc-length starting from $q(0)=0$ and associated to $\theta
 (0)=\theta _0,$ $P_3(0)=\lambda .$ The problem is said
\it{conservative} if there exists a coordinate $y$ transverse
to the Martinet surface such that for a dense set of initial
conditions $(\theta_0,\lambda )$ the trajectory $t\longrightarrow
 y(t)$ is periodic up to reparametrization. The equation
describing the evolution of $y$ is called the \it{characteristic
equation.}
\end{defi}


\subsubsection{Analysis of the foliation ${\cal F}$}

The foliation $({\cal F})$ is described by equation (\ref{4.9})~:
$$\theta ''+\lambda {\sin\ }\theta +\alpha ^2 {\sin}\ \theta \ 
{\cos}\ \theta -\alpha \beta \ {\sin}^2 \theta +\beta \ 
\theta ' {\cos\ } \theta =0.$$

\noindent Moreover recall the relation~: $y'=(1+\alpha y)
{\cos\ }\theta , \theta '=-(y\lambda -\alpha {\cos\ }\theta 
+\beta \ {\sin\ }\theta )$

The singular line project onto $\theta =k\pi $ which correspond
to the \it{singularities} of (\ref{4.9})~: $\theta  =k\pi , 
\theta '=0.$ 

Among the solutions of (\ref{4.9}), only those satisfying the
relation~:
\begin{eqnarray}   \label{4.10}
\theta ' = \alpha \ {\cos\ }\theta +\beta \ {\sin\ }\theta 
\end{eqnarray}

\noindent at $\tau =0$ correspond to projections of geodesics
starting at $t=0$ from $q(0)=0.$

Using an energy-balance relation we can represent the solutions
of $({\cal F})$ for $\mid \lambda \mid \gg \mid \alpha \mid , 
\mid \beta \mid , \mid \gamma \mid ,$ see \cite{BC}. We may suppose
$\lambda >0.$ Introducing the small parameter~: $\varepsilon 
=1/\sqrt{\lambda }$ and the parametrization $s=\tau \sqrt
{\lambda }$ we get the equation~:
\begin{eqnarray}    \label{4.11}
{\displaystyle {d^2\theta \over ds^2}}\ +\ {\sin\ }\theta \ +\ 
\varepsilon \beta {\cos\ }\theta \ {\displaystyle {d\theta 
\over ds}}\ +\ \varepsilon ^2  \alpha \ {\sin\ }\theta (\alpha 
\ {\cos\ }\theta -\beta \ {\sin\ }\theta )=0
\end{eqnarray}

\noindent and equation (\ref{4.10}) takes the form~:
\begin{eqnarray}   \label{4.12}
{\displaystyle {d\theta \over ds}}\ =\ \varepsilon 
(\alpha \ {\cos}\ \theta  +\beta \ {\sin \ }\theta )
\end{eqnarray}

The \it{flat case} corresponds to $\alpha =\beta =0,$
i.e~: ${\displaystyle {d^2\theta \over ds^2}}+{\sin}\ 
\theta =0,$ ${\displaystyle {d\theta \over ds}}=0$ at $s=0$
and is \it{also the limit case $\varepsilon \longrightarrow 0.$}

The following result is straightforward.
\begin{lem}
The problem is conservative if and only if $\beta =0.$
\end{lem}

We represent below the trajectories of $({\cal F})$ for
$\lambda \gg \mid \alpha \mid ,$ $\mid \beta \mid ,$ $\mid 
\gamma \mid ,$ on the phase space~: $(\theta ,{\dot \theta })$
but geometrically it corresponds to a foliation on the
\it{cylinder}~: $(e^{i\theta },{\dot \theta }).$

\begin{itemize}
\item \underline{Flat case} $(\alpha =\beta =0).$
It corresponds to a \it{pendulum}, see Fig. \ref{f3}.


\setlength{\unitlength}{0.5mm}
\begin{figure}[h]
\begin{center} 

\begin{picture}(180,100)
\thinlines
\drawvector{20.0}{40.0}{142.0}{1}{0}
\drawvector{90.0}{2.0}{92.0}{0}{1}
\drawdashline{24.0}{6.0}{24.0}{90.0}
\drawdashline{156.0}{90.0}{156.0}{6.0}
\drawlefttext{92.0}{92.0}{$\frac{d\theta}{ds}$}
\drawcenteredtext{166.0}{44.0}{$\theta$}
\drawlefttext{156.6}{36.15}{$+\pi$}
\drawrighttext{22.7}{36.84}{$-\pi$}
\drawlefttext{90.0}{38.0}{$0$}
\thicklines
\drawarc{90.0}{86.0}{160.88}{0.59}{2.52}
\drawarc{90.0}{-6.0}{160.88}{3.75}{5.67}
\drawvector{98.0}{74.0}{2.0}{1}{0}
\thinlines
\drawellipse{90.0}{41.0}{92.0}{42.0}{}
\drawvector{88.0}{62.0}{2.0}{1}{0}
\path(24.0,52.0)(24.0,52.0)(24.45,52.18)(24.94,52.36)(25.38,52.56)(25.85,52.77)(26.29,52.97)(26.72,53.2)(27.17,53.4)(27.6,53.65)
\path(27.6,53.65)(28.02,53.88)(28.43,54.13)(28.85,54.38)(29.25,54.63)(29.64,54.88)(30.04,55.15)(30.45,55.4)(30.85,55.68)(31.22,55.97)
\path(31.22,55.97)(31.62,56.25)(32.0,56.54)(32.4,56.83)(32.77,57.11)(33.15,57.4)(33.54,57.72)(33.93,58.02)(34.31,58.34)(34.68,58.65)
\path(34.68,58.65)(35.06,58.97)(35.45,59.29)(35.84,59.61)(36.22,59.93)(36.63,60.25)(37.02,60.59)(37.4,60.9)(37.81,61.25)(38.22,61.59)
\path(38.22,61.59)(38.63,61.93)(39.04,62.27)(39.45,62.61)(39.88,62.95)(40.31,63.29)(40.75,63.63)(41.18,63.97)(41.63,64.3)(42.09,64.66)
\path(42.09,64.66)(42.56,65.01)(43.02,65.36)(43.5,65.69)(43.99,66.05)(44.49,66.4)(44.99,66.73)(45.5,67.08)(46.02,67.43)(46.56,67.76)
\path(46.56,67.76)(47.11,68.12)(47.65,68.44)(48.24,68.8)(48.81,69.12)(49.4,69.47)(50.02,69.8)(50.63,70.12)(51.27,70.47)(51.9,70.8)
\path(51.9,70.8)(52.56,71.12)(53.24,71.44)(53.93,71.76)(54.63,72.08)(55.34,72.4)(56.08,72.69)(56.83,73.01)(57.59,73.3)(58.38,73.62)
\path(58.38,73.62)(59.18,73.91)(59.99,74.19)(60.83,74.48)(61.68,74.76)(62.54,75.05)(63.45,75.33)(64.36,75.58)(65.29,75.86)(66.23,76.12)
\path(66.23,76.12)(67.19,76.37)(68.19,76.62)(69.19,76.87)(70.23,77.11)(71.3,77.33)(72.37,77.55)(73.48,77.76)(74.58,78.0)(75.73,78.19)
\path(75.73,78.19)(76.91,78.4)(78.11,78.58)(79.33,78.76)(80.55,78.94)(81.83,79.12)(83.12,79.29)(84.44,79.44)(85.8,79.58)(87.16,79.73)
\path(87.16,79.73)(88.55,79.87)(89.98,80.0)(90.0,80.0)
\drawvector{85.8}{79.58}{4.19}{1}{0}
\path(156.0,52.0)(156.0,52.0)(155.46,52.18)(154.94,52.36)(154.41,52.56)(153.91,52.77)(153.44,52.97)(152.96,53.2)(152.47,53.4)(152.02,53.65)
\path(152.02,53.65)(151.57,53.88)(151.11,54.13)(150.69,54.38)(150.25,54.63)(149.83,54.88)(149.41,55.15)(149.0,55.4)(148.6,55.68)(148.19,55.97)
\path(148.19,55.97)(147.8,56.25)(147.41,56.54)(147.02,56.83)(146.63,57.11)(146.25,57.4)(145.86,57.72)(145.5,58.02)(145.11,58.34)(144.75,58.65)
\path(144.75,58.65)(144.36,58.97)(144.0,59.29)(143.61,59.61)(143.25,59.93)(142.88,60.25)(142.5,60.59)(142.11,60.9)(141.74,61.25)(141.36,61.59)
\path(141.36,61.59)(140.97,61.93)(140.58,62.27)(140.19,62.61)(139.77,62.95)(139.38,63.29)(138.97,63.63)(138.55,63.97)(138.13,64.3)(137.72,64.66)
\path(137.72,64.66)(137.27,65.01)(136.83,65.36)(136.38,65.69)(135.94,66.05)(135.47,66.4)(135.0,66.73)(134.5,67.08)(134.02,67.43)(133.5,67.76)
\path(133.5,67.76)(132.99,68.12)(132.47,68.44)(131.91,68.8)(131.38,69.12)(130.8,69.47)(130.22,69.8)(129.63,70.12)(129.02,70.47)(128.41,70.8)
\path(128.41,70.8)(127.79,71.12)(127.12,71.44)(126.47,71.76)(125.79,72.08)(125.08,72.4)(124.37,72.69)(123.65,73.01)(122.9,73.3)(122.12,73.62)
\path(122.12,73.62)(121.33,73.91)(120.54,74.19)(119.69,74.48)(118.87,74.76)(118.0,75.05)(117.12,75.33)(116.19,75.58)(115.26,75.86)(114.33,76.12)
\path(114.33,76.12)(113.36,76.37)(112.36,76.62)(111.33,76.87)(110.3,77.11)(109.23,77.33)(108.12,77.55)(107.01,77.76)(105.87,78.0)(104.69,78.19)
\path(104.69,78.19)(103.51,78.4)(102.26,78.58)(101.01,78.76)(99.76,78.94)(98.44,79.12)(97.11,79.29)(95.75,79.44)(94.36,79.58)(92.93,79.73)
\path(92.93,79.73)(91.48,79.87)(90.0,80.0)(90.0,80.0)
\thicklines
\drawvector{88.0}{6.0}{2.0}{-1}{0}
\thinlines
\drawlefttext{98.0}{70.0}{$\Sigma$}
\end{picture}

\end{center}
\caption{}\label{f3}
\end{figure}

\noindent The main properties are the following.
We have two singularities~:
\begin{itemize}
\item $0$ is a center.

\item $(\pi ,0)$ is a saddle and the separatrix
$\Sigma $ is a saddle connection.
\end{itemize}

Only the \it{oscillating trajectories} correspond to
geodesics starting from 0.

\item \underline{Conservative case} $(\beta =0)$
The equation reduces to~:
$${d^2\theta \over ds^2}\ +{\sin}\ \theta +\varepsilon ^2 
\alpha ^2\ {\sin\ }\theta \ {\cos\ }\theta =0$$

\noindent Multiplying both sides by ${\displaystyle 
{d\theta \over ds}}$ and integrating on $[0,s]$ we get~:
$$\Bigl [ {1\over 2}\ \Bigl ( {d\theta \over ds}\Bigr )^2 
\Bigr ]^s_0\ =\ \Bigl [ {\cos\ }\theta  \ +\ {\varepsilon ^2
\alpha ^2 {\cos}^2 \ \theta \over 2}\Bigr ]^s_0$$

\noindent and the system has a \it{global $C^\omega $ first
integral~:}
\begin{eqnarray}   \label{4.13}
V(\theta ,{\dot \theta })={\displaystyle {1\over 2}}
 {\dot \theta }^2-\Bigl ( {\cos\ }\theta  +{\displaystyle 
{\varepsilon ^2 \alpha ^2\over 2}} {\cos}^2 \theta \Bigr ).
\end{eqnarray}

\noindent The phase portrait is similar to the one in the flat
case but the section defined by (\ref{4.12}) and corresponding
to $y=0$ is here~: ${\displaystyle {d\theta \over ds}}\ =\ 
\varepsilon \alpha {\cos\ }\theta .$ In particular if $\alpha 
\not = 0$ (strict case) there exist both oscillating and
rotating trajectories corresponding to projections of geodesics
starting from 0, see Fig. \ref{f4}.


\setlength{\unitlength}{0.5mm}
\begin{figure}[h]
\begin{center}

\begin{picture}(180,100)
\thinlines
\drawvector{20.0}{48.0}{142.0}{1}{0}
\drawvector{90.0}{4.0}{90.0}{0}{1}
\drawdashline{40.0}{10.0}{40.0}{90.0}
\drawdashline{140.0}{10.0}{140.0}{90.0}
\drawlefttext{90.9}{44.7}{$0$}
\drawcenteredtext{164.47}{51.45}{$\theta$}
\drawlefttext{91.55}{91.23}{$\frac{d\theta}{ds}$}
\drawlefttext{143.1}{44.7}{$+\pi$}
\drawdotline{40.0}{48.0}{60.0}{18.0}
\thicklines
\drawarc{90.0}{88.0}{128.05}{0.67}{2.46}
\drawarc{90.0}{8.0}{128.05}{3.75}{5.59}
\thinlines
\drawdotline{20.0}{18.0}{40.0}{48.0}
\drawdotline{40.0}{48.0}{60.0}{78.0}
\drawdotline{20.0}{78.0}{40.0}{48.0}
\drawdotline{140.0}{48.0}{160.0}{78.0}
\drawdotline{120.0}{78.0}{140.0}{48.0}
\drawdotline{140.0}{48.0}{120.0}{18.0}
\drawdotline{160.0}{18.0}{140.0}{48.0}
\drawrighttext{36.9}{44.49}{$-\pi$}
\thicklines
\drawarc{179.77}{75.16}{96.33}{1.9}{2.5}
\drawarc{178.0}{22.29}{91.69}{3.73}{4.4}
\thinlines
\path(31.04,30.78)(31.04,30.78)(31.2,31.09)(31.37,31.37)(31.53,31.69)(31.7,31.97)(31.86,32.27)(32.02,32.56)(32.18,32.83)(32.34,33.11)
\drawarc{90.0}{10.0}{87.72}{3.95}{5.46}
\drawarc{144.0}{62.0}{62.47}{0.86}{2.44}
\drawarc{34.4}{64.76}{68.58}{0.71}{2.13}
\path(32.34,33.11)(32.5,33.36)(32.68,33.63)(32.84,33.88)(33.0,34.13)(33.15,34.38)(33.34,34.61)(33.5,34.84)(33.65,35.06)(33.84,35.29)
\path(33.84,35.29)(34.0,35.5)(34.15,35.7)(34.34,35.9)(34.5,36.11)(34.65,36.29)(34.84,36.49)(35.0,36.65)(35.18,36.84)(35.34,37.0)
\path(35.34,37.0)(35.52,37.15)(35.68,37.31)(35.86,37.47)(36.02,37.61)(36.2,37.75)(36.36,37.88)(36.54,38.02)(36.7,38.13)(36.88,38.25)
\path(36.88,38.25)(37.04,38.36)(37.22,38.47)(37.4,38.56)(37.56,38.65)(37.74,38.75)(37.9,38.83)(38.09,38.9)(38.25,38.97)(38.43,39.04)
\path(38.43,39.04)(38.61,39.09)(38.77,39.15)(38.95,39.2)(39.13,39.24)(39.31,39.27)(39.47,39.31)(39.65,39.33)(39.83,39.34)(40.0,39.36)
\path(40.0,39.36)(40.18,39.36)(40.36,39.36)(40.54,39.36)(40.72,39.36)(40.9,39.34)(41.06,39.31)(41.25,39.29)(41.43,39.25)(41.61,39.22)
\path(41.61,39.22)(41.79,39.18)(41.97,39.11)(42.15,39.06)(42.33,39.0)(42.5,38.93)(42.68,38.86)(42.86,38.79)(43.04,38.7)(43.22,38.61)
\path(43.22,38.61)(43.4,38.5)(43.59,38.4)(43.77,38.29)(43.95,38.18)(44.13,38.06)(44.31,37.93)(44.5,37.81)(44.68,37.68)(44.86,37.52)
\path(44.86,37.52)(45.04,37.38)(45.24,37.22)(45.4,37.06)(45.61,36.9)(45.79,36.72)(45.97,36.56)(46.15,36.36)(46.34,36.18)(46.52,35.99)
\path(46.52,35.99)(46.72,35.79)(46.9,35.59)(47.09,35.36)(47.27,35.15)(47.47,34.93)(47.65,34.7)(47.84,34.47)(48.02,34.22)(48.22,33.97)
\path(48.22,33.97)(48.4,33.72)(48.59,33.47)(48.59,33.47)
\path(132.3,33.25)(132.3,33.25)(132.44,33.5)(132.61,33.75)(132.77,34.0)(132.91,34.24)(133.08,34.47)(133.24,34.7)(133.38,34.93)(133.55,35.15)
\path(133.55,35.15)(133.72,35.36)(133.88,35.56)(134.05,35.77)(134.21,35.97)(134.36,36.15)(134.52,36.34)(134.69,36.52)(134.86,36.7)(135.02,36.86)
\path(135.02,36.86)(135.16,37.02)(135.35,37.18)(135.5,37.34)(135.66,37.47)(135.83,37.61)(136.0,37.75)(136.16,37.88)(136.33,38.0)(136.5,38.13)
\path(136.5,38.13)(136.66,38.24)(136.83,38.34)(137.0,38.43)(137.16,38.54)(137.33,38.63)(137.5,38.7)(137.66,38.79)(137.85,38.86)(138.02,38.93)
\path(138.02,38.93)(138.19,38.99)(138.36,39.04)(138.52,39.09)(138.69,39.13)(138.88,39.18)(139.05,39.2)(139.22,39.24)(139.38,39.25)(139.57,39.27)
\path(139.57,39.27)(139.74,39.29)(139.91,39.29)(140.08,39.29)(140.25,39.29)(140.44,39.27)(140.61,39.25)(140.77,39.24)(140.96,39.2)(141.13,39.15)
\path(141.13,39.15)(141.32,39.13)(141.49,39.09)(141.66,39.04)(141.85,38.97)(142.02,38.9)(142.21,38.84)(142.38,38.77)(142.57,38.7)(142.74,38.61)
\path(142.74,38.61)(142.91,38.52)(143.11,38.43)(143.27,38.31)(143.47,38.22)(143.63,38.09)(143.83,37.99)(144.0,37.86)(144.19,37.72)(144.38,37.59)
\path(144.38,37.59)(144.55,37.45)(144.75,37.31)(144.91,37.15)(145.11,37.0)(145.3,36.84)(145.47,36.65)(145.66,36.49)(145.86,36.31)(146.02,36.13)
\path(146.02,36.13)(146.22,35.93)(146.41,35.74)(146.6,35.52)(146.77,35.31)(146.97,35.11)(147.16,34.88)(147.36,34.65)(147.55,34.43)(147.74,34.2)
\path(147.74,34.2)(147.91,33.95)(148.11,33.7)(148.3,33.45)(148.5,33.2)(148.69,32.93)(148.88,32.65)(149.07,32.38)(149.27,32.09)(149.46,31.8)
\path(149.46,31.8)(149.63,31.52)(149.85,31.21)(149.85,31.21)
\drawvector{41.15}{39.31}{0.88}{-1}{0}
\path(55.79,37.97)(55.79,37.97)(55.15,38.36)(54.56,38.77)(54.0,39.18)(53.45,39.59)(52.95,40.0)(52.47,40.4)(52.04,40.81)(51.63,41.22)
\path(51.63,41.22)(51.24,41.63)(50.88,42.04)(50.56,42.45)(50.27,42.86)(50.0,43.27)(49.75,43.68)(49.54,44.09)(49.36,44.5)(49.2,44.9)
\path(49.2,44.9)(49.06,45.31)(48.95,45.72)(48.86,46.13)(48.81,46.54)(48.77,46.93)(48.77,47.34)(48.79,47.74)(48.83,48.13)(48.88,48.54)
\path(48.88,48.54)(48.97,48.93)(49.08,49.33)(49.2,49.72)(49.36,50.11)(49.52,50.5)(49.72,50.88)(49.93,51.25)(50.15,51.65)(50.4,52.02)
\path(50.4,52.02)(50.68,52.4)(50.97,52.77)(51.27,53.13)(51.61,53.5)(51.95,53.86)(52.31,54.22)(52.68,54.59)(53.08,54.95)(53.49,55.29)
\path(53.49,55.29)(53.9,55.63)(54.36,55.99)(54.81,56.33)(55.29,56.65)(55.77,56.99)(56.27,57.31)(56.77,57.63)(57.29,57.95)(57.84,58.27)
\path(57.84,58.27)(58.38,58.59)(58.95,58.88)(59.52,59.18)(60.09,59.49)(60.7,59.77)(61.29,60.06)(61.9,60.34)(62.52,60.61)(63.15,60.88)
\path(63.15,60.88)(63.81,61.15)(64.44,61.4)(65.11,61.65)(65.76,61.9)(66.44,62.15)(67.12,62.38)(67.8,62.61)(68.5,62.83)(69.19,63.04)
\path(69.19,63.04)(69.87,63.25)(70.58,63.45)(71.3,63.65)(72.01,63.84)(72.73,64.01)(73.44,64.19)(74.18,64.37)(74.91,64.51)(75.62,64.69)
\path(75.62,64.69)(76.36,64.83)(77.08,64.97)(77.8,65.08)(78.55,65.22)(79.26,65.33)(80.01,65.44)(80.75,65.55)(81.48,65.62)(82.19,65.72)
\path(82.19,65.72)(82.93,65.8)(83.65,65.86)(84.37,65.91)(85.08,65.97)(85.8,66.01)(86.51,66.04)(87.22,66.05)(87.91,66.05)(88.62,66.08)
\path(88.62,66.08)(89.3,66.05)(89.98,66.05)(90.0,66.05)
\drawvector{87.91}{66.05}{2.04}{1}{0}
\path(124.41,37.75)(124.41,37.75)(125.08,38.15)(125.73,38.59)(126.33,39.02)(126.91,39.43)(127.44,39.86)(127.97,40.29)(128.44,40.72)(128.88,41.13)
\path(128.88,41.13)(129.3,41.56)(129.69,41.97)(130.02,42.4)(130.36,42.83)(130.66,43.25)(130.91,43.65)(131.13,44.08)(131.36,44.5)(131.52,44.9)
\path(131.52,44.9)(131.66,45.33)(131.8,45.74)(131.88,46.15)(131.96,46.56)(132.0,46.95)(132.02,47.36)(132.0,47.77)(131.97,48.15)(131.91,48.56)
\path(131.91,48.56)(131.83,48.97)(131.72,49.36)(131.61,49.75)(131.46,50.15)(131.27,50.54)(131.08,50.9)(130.86,51.29)(130.63,51.68)(130.38,52.06)
\path(130.38,52.06)(130.1,52.43)(129.8,52.79)(129.47,53.15)(129.13,53.52)(128.77,53.88)(128.41,54.25)(128.02,54.59)(127.62,54.95)(127.19,55.29)
\path(127.19,55.29)(126.76,55.63)(126.3,55.97)(125.83,56.31)(125.33,56.65)(124.83,56.97)(124.33,57.29)(123.8,57.61)(123.26,57.93)(122.69,58.24)
\path(122.69,58.24)(122.12,58.54)(121.55,58.84)(120.98,59.13)(120.37,59.43)(119.76,59.7)(119.12,59.99)(118.51,60.27)(117.86,60.54)(117.19,60.79)
\path(117.19,60.79)(116.55,61.06)(115.87,61.31)(115.19,61.56)(114.51,61.81)(113.83,62.04)(113.12,62.27)(112.44,62.5)(111.73,62.72)(111.01,62.93)
\path(111.01,62.93)(110.3,63.13)(109.58,63.34)(108.83,63.52)(108.12,63.72)(107.37,63.9)(106.66,64.08)(105.91,64.23)(105.18,64.4)(104.44,64.55)
\path(104.44,64.55)(103.69,64.69)(102.94,64.83)(102.19,64.98)(101.47,65.08)(100.73,65.22)(99.98,65.33)(99.25,65.44)(98.51,65.51)(97.76,65.62)
\path(97.76,65.62)(97.05,65.69)(96.3,65.76)(95.58,65.83)(94.87,65.87)(94.16,65.94)(93.44,65.98)(92.75,66.01)(92.05,66.04)(91.36,66.05)
\path(91.36,66.05)(90.66,66.05)(90.0,66.05)(90.0,66.05)
\drawvector{140.61}{39.31}{1.11}{-1}{0}
\thicklines
\drawvector{88.0}{72.0}{2.0}{1}{0}
\drawvector{92.0}{24.0}{2.0}{-1}{0}
\thinlines
\drawcenteredtext{166.0}{40.0}{section}
\drawcenteredtext{166.0}{26.0}{$\Sigma$}
\end{picture}

\end{center}

\caption{$\alpha >0$}\label{f4}
\end{figure}

\item \underline{General case} $(\beta \not = 0)$
The two main differences are the following~:

\begin{itemize}
\item the center $0$ becomes a focus~;

\item the saddle connection is broken.
\end{itemize}

The trajectories are represented on Fig. \ref{f5}.


\setlength{\unitlength}{0.4mm}
\begin{figure}[h]
\begin{center} 

\begin{picture}(180,100)
\thinlines
\drawvector{20.0}{40.0}{142.0}{1}{0}
\drawvector{90.0}{4.0}{90.0}{0}{1}
\drawdashline{24.0}{10.0}{24.0}{90.0}
\drawdashline{156.0}{90.0}{156.0}{10.0}
\drawlefttext{92.0}{92.0}{$\frac{d\theta}{ds}$}
\drawcenteredtext{166.0}{44.0}{$\theta$}
\drawlefttext{156.0}{36.0}{$+\pi$}
\drawrighttext{22.0}{36.0}{$-\pi$}
\drawlefttext{92.0}{36.0}{$0$}
\path(24.0,44.0)(24.0,44.0)(24.79,44.31)(25.59,44.63)(26.37,44.93)(27.17,45.25)(27.95,45.54)(28.73,45.86)(29.53,46.15)(30.3,46.45)
\path(30.3,46.45)(31.07,46.75)(31.85,47.02)(32.63,47.31)(33.38,47.59)(34.15,47.88)(34.91,48.15)(35.68,48.43)(36.43,48.7)(37.18,48.97)
\path(37.18,48.97)(37.93,49.24)(38.68,49.5)(39.43,49.75)(40.18,50.0)(40.91,50.25)(41.65,50.5)(42.38,50.75)(43.11,51.0)(43.84,51.22)
\path(43.84,51.22)(44.56,51.47)(45.29,51.7)(46.02,51.93)(46.74,52.15)(47.45,52.38)(48.15,52.59)(48.86,52.81)(49.58,53.02)(50.27,53.22)
\path(50.27,53.22)(50.97,53.43)(51.68,53.63)(52.36,53.84)(53.06,54.04)(53.75,54.22)(54.43,54.43)(55.13,54.61)(55.81,54.79)(56.47,54.97)
\path(56.47,54.97)(57.15,55.15)(57.83,55.33)(58.5,55.5)(59.16,55.66)(59.83,55.83)(60.49,55.99)(61.15,56.15)(61.81,56.31)(62.45,56.45)
\path(62.45,56.45)(63.11,56.61)(63.75,56.75)(64.4,56.9)(65.05,57.04)(65.69,57.16)(66.31,57.31)(66.94,57.43)(67.58,57.56)(68.2,57.68)
\path(68.2,57.68)(68.83,57.79)(69.45,57.91)(70.08,58.02)(70.69,58.15)(71.3,58.25)(71.91,58.36)(72.52,58.45)(73.12,58.54)(73.73,58.65)
\path(73.73,58.65)(74.33,58.74)(74.93,58.83)(75.52,58.9)(76.12,58.99)(76.7,59.06)(77.29,59.15)(77.87,59.22)(78.45,59.29)(79.02,59.34)
\path(79.02,59.34)(79.61,59.41)(80.18,59.47)(80.75,59.52)(81.31,59.59)(81.87,59.63)(82.44,59.68)(83.0,59.72)(83.55,59.75)(84.11,59.79)
\path(84.11,59.79)(84.65,59.83)(85.19,59.86)(85.75,59.88)(86.29,59.91)(86.81,59.93)(87.36,59.95)(87.88,59.97)(88.41,59.97)(88.94,59.99)
\path(88.94,59.99)(89.47,59.99)(89.98,60.0)(90.0,60.0)
\drawvector{88.41}{59.97}{1.58}{1}{0}
\path(24.0,60.0)(24.0,60.0)(24.87,60.18)(25.75,60.38)(26.62,60.59)(27.47,60.77)(28.34,60.97)(29.2,61.15)(30.04,61.34)(30.88,61.52)
\path(30.88,61.52)(31.73,61.7)(32.56,61.9)(33.4,62.06)(34.24,62.25)(35.06,62.43)(35.88,62.59)(36.7,62.77)(37.5,62.93)(38.31,63.11)
\path(38.31,63.11)(39.11,63.27)(39.91,63.43)(40.72,63.59)(41.5,63.75)(42.29,63.9)(43.06,64.06)(43.84,64.22)(44.61,64.37)(45.38,64.51)
\path(45.38,64.51)(46.15,64.66)(46.9,64.8)(47.65,64.94)(48.4,65.08)(49.15,65.23)(49.9,65.37)(50.63,65.51)(51.36,65.63)(52.09,65.76)
\path(52.09,65.76)(52.81,65.9)(53.54,66.02)(54.25,66.15)(54.97,66.26)(55.66,66.4)(56.38,66.51)(57.06,66.62)(57.77,66.75)(58.45,66.86)
\path(58.45,66.86)(59.13,66.97)(59.81,67.08)(60.5,67.19)(61.16,67.29)(61.83,67.38)(62.49,67.5)(63.15,67.58)(63.81,67.69)(64.45,67.79)
\path(64.45,67.79)(65.09,67.87)(65.73,67.97)(66.37,68.05)(67.01,68.15)(67.62,68.23)(68.26,68.3)(68.87,68.4)(69.48,68.47)(70.09,68.55)
\path(70.09,68.55)(70.69,68.62)(71.3,68.69)(71.9,68.76)(72.48,68.83)(73.08,68.91)(73.66,68.97)(74.23,69.02)(74.8,69.08)(75.37,69.15)
\path(75.37,69.15)(75.94,69.2)(76.51,69.26)(77.06,69.31)(77.62,69.37)(78.16,69.41)(78.7,69.47)(79.25,69.51)(79.77,69.55)(80.3,69.58)
\path(80.3,69.58)(80.83,69.62)(81.36,69.66)(81.87,69.7)(82.38,69.73)(82.9,69.76)(83.4,69.8)(83.9,69.83)(84.4,69.84)(84.88,69.87)
\path(84.88,69.87)(85.37,69.9)(85.86,69.91)(86.33,69.93)(86.8,69.94)(87.27,69.95)(87.73,69.97)(88.19,69.98)(88.66,69.98)(89.11,69.98)
\path(89.11,69.98)(89.55,69.98)(89.98,70.0)(90.0,70.0)
\drawvector{88.66}{69.98}{1.33}{1}{0}
\path(24.0,78.0)(24.0,78.0)(24.67,78.11)(25.35,78.23)(26.03,78.34)(26.7,78.45)(27.38,78.58)(28.06,78.69)(28.75,78.8)(29.42,78.9)
\path(29.42,78.9)(30.1,79.01)(30.77,79.11)(31.45,79.22)(32.13,79.31)(32.79,79.41)(33.47,79.51)(34.15,79.62)(34.81,79.7)(35.5,79.8)
\path(35.5,79.8)(36.16,79.9)(36.84,79.98)(37.52,80.06)(38.18,80.16)(38.86,80.25)(39.52,80.33)(40.2,80.41)(40.86,80.5)(41.54,80.56)
\path(41.54,80.56)(42.2,80.65)(42.88,80.73)(43.54,80.8)(44.22,80.87)(44.88,80.94)(45.54,81.01)(46.22,81.08)(46.88,81.15)(47.54,81.22)
\path(47.54,81.22)(48.22,81.27)(48.88,81.33)(49.54,81.4)(50.2,81.45)(50.86,81.51)(51.54,81.56)(52.2,81.62)(52.86,81.68)(53.52,81.73)
\path(53.52,81.73)(54.18,81.76)(54.84,81.81)(55.5,81.87)(56.16,81.91)(56.83,81.94)(57.49,82.0)(58.15,82.02)(58.81,82.06)(59.47,82.11)
\path(59.47,82.11)(60.13,82.13)(60.79,82.18)(61.45,82.2)(62.11,82.23)(62.75,82.26)(63.41,82.29)(64.06,82.3)(64.73,82.33)(65.38,82.36)
\path(65.38,82.36)(66.04,82.37)(66.69,82.4)(67.34,82.41)(68.0,82.43)(68.66,82.44)(69.3,82.45)(69.95,82.47)(70.61,82.47)(71.26,82.48)
\path(71.26,82.48)(71.91,82.48)(72.56,82.48)(73.22,82.48)(73.87,82.5)(74.51,82.48)(75.16,82.48)(75.81,82.48)(76.47,82.48)(77.11,82.48)
\path(77.11,82.48)(77.76,82.47)(78.41,82.45)(79.05,82.44)(79.69,82.43)(80.34,82.41)(81.0,82.4)(81.63,82.37)(82.29,82.36)(82.93,82.33)
\path(82.93,82.33)(83.56,82.31)(84.22,82.29)(84.86,82.26)(85.51,82.23)(86.15,82.2)(86.79,82.18)(87.43,82.13)(88.06,82.11)(88.7,82.06)
\path(88.7,82.06)(89.34,82.02)(89.98,82.0)(90.0,82.0)
\drawvector{88.06}{82.11}{1.93}{1}{0}
\path(90.0,70.0)(90.0,70.0)(90.75,69.94)(91.51,69.9)(92.26,69.84)(93.01,69.79)(93.76,69.73)(94.51,69.66)(95.26,69.58)(96.01,69.51)
\path(96.01,69.51)(96.75,69.41)(97.5,69.33)(98.23,69.23)(98.97,69.13)(99.7,69.04)(100.44,68.93)(101.16,68.8)(101.9,68.69)(102.62,68.55)
\path(102.62,68.55)(103.34,68.43)(104.06,68.3)(104.8,68.15)(105.51,68.01)(106.23,67.86)(106.94,67.69)(107.66,67.54)(108.37,67.37)(109.08,67.19)
\path(109.08,67.19)(109.79,67.01)(110.48,66.83)(111.19,66.65)(111.9,66.44)(112.58,66.26)(113.29,66.05)(113.98,65.83)(114.68,65.62)(115.37,65.41)
\path(115.37,65.41)(116.05,65.19)(116.75,64.95)(117.43,64.72)(118.11,64.48)(118.79,64.23)(119.47,63.97)(120.15,63.72)(120.83,63.47)(121.5,63.2)
\path(121.5,63.2)(122.16,62.93)(122.83,62.65)(123.51,62.36)(124.16,62.08)(124.83,61.79)(125.48,61.5)(126.15,61.18)(126.8,60.88)(127.47,60.56)
\path(127.47,60.56)(128.11,60.25)(128.77,59.93)(129.41,59.59)(130.07,59.27)(130.71,58.93)(131.35,58.58)(131.99,58.24)(132.63,57.88)(133.27,57.52)
\path(133.27,57.52)(133.91,57.15)(134.53,56.79)(135.16,56.4)(135.8,56.02)(136.42,55.63)(137.05,55.25)(137.66,54.86)(138.28,54.45)(138.91,54.04)
\path(138.91,54.04)(139.52,53.63)(140.14,53.22)(140.75,52.79)(141.36,52.36)(141.97,51.93)(142.58,51.5)(143.19,51.06)(143.78,50.61)(144.38,50.15)
\path(144.38,50.15)(144.99,49.7)(145.58,49.22)(146.19,48.75)(146.77,48.29)(147.36,47.81)(147.96,47.33)(148.55,46.84)(149.13,46.34)(149.71,45.84)
\path(149.71,45.84)(150.28,45.34)(150.86,44.81)(151.44,44.31)(152.02,43.79)(152.6,43.25)(153.16,42.72)(153.74,42.18)(154.3,41.65)(154.86,41.09)
\path(154.86,41.09)(155.42,40.54)(155.99,40.0)(156.0,40.0)
\path(90.0,82.0)(90.0,82.0)(90.83,81.87)(91.66,81.75)(92.5,81.62)(93.33,81.48)(94.15,81.36)(94.97,81.22)(95.79,81.08)(96.59,80.94)
\path(96.59,80.94)(97.41,80.8)(98.2,80.66)(99.01,80.51)(99.81,80.34)(100.61,80.19)(101.4,80.04)(102.19,79.87)(102.97,79.72)(103.75,79.55)
\path(103.75,79.55)(104.52,79.37)(105.3,79.2)(106.08,79.04)(106.83,78.86)(107.59,78.68)(108.36,78.48)(109.12,78.3)(109.87,78.12)(110.62,77.93)
\path(110.62,77.93)(111.36,77.73)(112.09,77.54)(112.83,77.33)(113.56,77.12)(114.3,76.93)(115.02,76.72)(115.75,76.51)(116.47,76.3)(117.19,76.08)
\path(117.19,76.08)(117.9,75.86)(118.61,75.63)(119.31,75.41)(120.01,75.19)(120.7,74.94)(121.41,74.72)(122.09,74.48)(122.79,74.25)(123.47,74.0)
\path(123.47,74.0)(124.15,73.76)(124.83,73.51)(125.5,73.26)(126.16,73.01)(126.83,72.75)(127.48,72.5)(128.14,72.23)(128.8,71.97)(129.46,71.69)
\path(129.46,71.69)(130.11,71.43)(130.75,71.16)(131.38,70.87)(132.02,70.61)(132.66,70.33)(133.28,70.04)(133.91,69.76)(134.53,69.47)(135.16,69.16)
\path(135.16,69.16)(135.77,68.87)(136.38,68.58)(136.99,68.27)(137.58,67.98)(138.19,67.66)(138.78,67.36)(139.38,67.05)(139.97,66.73)(140.55,66.41)
\path(140.55,66.41)(141.13,66.09)(141.72,65.76)(142.3,65.44)(142.86,65.12)(143.44,64.79)(144.0,64.44)(144.55,64.12)(145.11,63.77)(145.66,63.43)
\path(145.66,63.43)(146.22,63.09)(146.77,62.74)(147.3,62.38)(147.85,62.04)(148.38,61.68)(148.91,61.31)(149.44,60.95)(149.97,60.59)(150.5,60.22)
\path(150.5,60.22)(151.0,59.86)(151.52,59.47)(152.03,59.11)(152.55,58.72)(153.05,58.34)(153.55,57.95)(154.05,57.56)(154.53,57.18)(155.02,56.79)
\path(155.02,56.79)(155.5,56.38)(155.99,56.0)(156.0,56.0)
\path(90.0,60.0)(90.0,60.0)(90.3,59.95)(90.62,59.9)(90.94,59.86)(91.25,59.81)(91.55,59.75)(91.86,59.7)(92.16,59.63)(92.47,59.56)
\path(92.47,59.56)(92.76,59.5)(93.05,59.43)(93.34,59.36)(93.62,59.27)(93.91,59.2)(94.19,59.11)(94.48,59.04)(94.76,58.95)(95.02,58.84)
\path(95.02,58.84)(95.3,58.75)(95.56,58.65)(95.83,58.54)(96.09,58.45)(96.36,58.34)(96.61,58.22)(96.87,58.11)(97.12,58.0)(97.37,57.86)
\path(97.37,57.86)(97.61,57.75)(97.86,57.61)(98.09,57.49)(98.33,57.36)(98.56,57.22)(98.8,57.08)(99.02,56.93)(99.26,56.79)(99.48,56.63)
\path(99.48,56.63)(99.69,56.47)(99.91,56.31)(100.12,56.15)(100.34,56.0)(100.55,55.84)(100.76,55.66)(100.97,55.49)(101.16,55.31)(101.36,55.13)
\path(101.36,55.13)(101.55,54.95)(101.75,54.77)(101.94,54.58)(102.12,54.38)(102.3,54.18)(102.5,54.0)(102.66,53.79)(102.84,53.59)(103.01,53.38)
\path(103.01,53.38)(103.19,53.16)(103.36,52.95)(103.51,52.74)(103.69,52.52)(103.84,52.29)(104.0,52.06)(104.15,51.84)(104.3,51.59)(104.44,51.36)
\path(104.44,51.36)(104.59,51.11)(104.73,50.88)(104.87,50.63)(105.01,50.38)(105.15,50.13)(105.27,49.88)(105.41,49.61)(105.54,49.36)(105.66,49.09)
\path(105.66,49.09)(105.77,48.81)(105.88,48.54)(106.01,48.27)(106.12,48.0)(106.23,47.7)(106.33,47.43)(106.44,47.13)(106.54,46.84)(106.62,46.56)
\path(106.62,46.56)(106.73,46.25)(106.81,45.95)(106.91,45.65)(107.0,45.34)(107.08,45.04)(107.16,44.72)(107.23,44.4)(107.3,44.08)(107.38,43.75)
\path(107.38,43.75)(107.44,43.43)(107.51,43.11)(107.58,42.77)(107.65,42.43)(107.69,42.09)(107.76,41.75)(107.8,41.4)(107.86,41.06)(107.91,40.7)
\path(107.91,40.7)(107.94,40.34)(108.0,40.0)(108.0,40.0)
\path(108.0,40.0)(108.0,40.0)(107.94,39.79)(107.91,39.59)(107.86,39.4)(107.8,39.2)(107.76,39.02)(107.69,38.84)(107.65,38.65)(107.58,38.47)
\path(107.58,38.47)(107.51,38.29)(107.45,38.11)(107.38,37.93)(107.3,37.77)(107.23,37.59)(107.16,37.43)(107.08,37.27)(107.0,37.09)(106.91,36.93)
\path(106.91,36.93)(106.81,36.77)(106.73,36.63)(106.62,36.47)(106.54,36.31)(106.44,36.18)(106.33,36.02)(106.23,35.88)(106.12,35.75)(106.01,35.61)
\path(106.01,35.61)(105.88,35.47)(105.77,35.34)(105.66,35.2)(105.54,35.08)(105.41,34.95)(105.27,34.81)(105.15,34.7)(105.01,34.58)(104.87,34.47)
\path(104.87,34.47)(104.73,34.34)(104.59,34.24)(104.44,34.13)(104.3,34.02)(104.15,33.91)(104.0,33.81)(103.84,33.7)(103.69,33.61)(103.51,33.52)
\path(103.51,33.52)(103.36,33.43)(103.19,33.33)(103.01,33.25)(102.84,33.15)(102.66,33.08)(102.5,33.0)(102.3,32.91)(102.12,32.84)(101.94,32.77)
\path(101.94,32.77)(101.75,32.68)(101.55,32.63)(101.36,32.56)(101.16,32.49)(100.97,32.43)(100.76,32.36)(100.55,32.31)(100.34,32.25)(100.12,32.2)
\path(100.12,32.2)(99.91,32.15)(99.69,32.11)(99.48,32.06)(99.26,32.02)(99.02,31.97)(98.8,31.94)(98.56,31.9)(98.33,31.87)(98.09,31.84)
\path(98.09,31.84)(97.86,31.81)(97.61,31.79)(97.37,31.77)(97.12,31.75)(96.87,31.72)(96.61,31.7)(96.36,31.7)(96.09,31.68)(95.83,31.68)
\path(95.83,31.68)(95.56,31.67)(95.3,31.65)(95.02,31.65)(94.76,31.65)(94.48,31.67)(94.19,31.67)(93.91,31.68)(93.62,31.69)(93.34,31.7)
\path(93.34,31.7)(93.05,31.7)(92.76,31.72)(92.47,31.75)(92.16,31.77)(91.86,31.79)(91.55,31.81)(91.25,31.85)(90.94,31.88)(90.62,31.92)
\path(90.62,31.92)(90.3,31.95)(90.0,31.98)(90.0,32.0)
\path(90.0,32.0)(90.0,32.0)(89.87,32.04)(89.76,32.08)(89.63,32.11)(89.51,32.15)(89.41,32.2)(89.3,32.25)(89.18,32.29)(89.06,32.34)
\path(89.06,32.34)(88.95,32.38)(88.86,32.43)(88.75,32.47)(88.63,32.52)(88.54,32.58)(88.43,32.63)(88.33,32.68)(88.23,32.74)(88.12,32.79)
\path(88.12,32.79)(88.02,32.84)(87.93,32.9)(87.83,32.95)(87.73,33.0)(87.65,33.06)(87.55,33.13)(87.45,33.18)(87.37,33.25)(87.27,33.31)
\path(87.27,33.31)(87.19,33.36)(87.11,33.43)(87.01,33.49)(86.93,33.54)(86.84,33.61)(86.76,33.68)(86.69,33.75)(86.61,33.81)(86.52,33.88)
\path(86.52,33.88)(86.44,33.95)(86.37,34.02)(86.3,34.09)(86.23,34.15)(86.15,34.22)(86.08,34.31)(86.01,34.38)(85.94,34.45)(85.87,34.52)
\path(85.87,34.52)(85.8,34.61)(85.73,34.68)(85.68,34.75)(85.62,34.84)(85.55,34.91)(85.5,35.0)(85.44,35.08)(85.37,35.15)(85.31,35.24)
\path(85.31,35.24)(85.26,35.31)(85.2,35.4)(85.16,35.49)(85.09,35.56)(85.05,35.65)(85.0,35.75)(84.94,35.83)(84.91,35.91)(84.86,36.0)
\path(84.86,36.0)(84.81,36.09)(84.76,36.18)(84.73,36.27)(84.69,36.38)(84.65,36.47)(84.61,36.56)(84.56,36.65)(84.54,36.75)(84.5,36.84)
\path(84.5,36.84)(84.47,36.95)(84.43,37.04)(84.4,37.15)(84.37,37.24)(84.33,37.34)(84.3,37.45)(84.29,37.54)(84.26,37.65)(84.23,37.75)
\path(84.23,37.75)(84.2,37.86)(84.19,37.95)(84.16,38.06)(84.15,38.18)(84.12,38.27)(84.11,38.38)(84.09,38.5)(84.08,38.61)(84.06,38.72)
\path(84.06,38.72)(84.05,38.83)(84.04,38.95)(84.02,39.06)(84.01,39.16)(84.01,39.29)(84.01,39.4)(84.0,39.52)(84.0,39.63)(84.0,39.75)
\path(84.0,39.75)(84.0,39.88)(84.0,39.99)(84.0,40.0)
\path(84.0,40.0)(84.0,40.0)(84.04,40.11)(84.08,40.22)(84.12,40.34)(84.16,40.45)(84.19,40.56)(84.23,40.68)(84.27,40.79)(84.33,40.88)
\path(84.33,40.88)(84.37,40.99)(84.41,41.09)(84.45,41.18)(84.5,41.29)(84.55,41.38)(84.58,41.47)(84.63,41.56)(84.69,41.65)(84.73,41.75)
\path(84.73,41.75)(84.77,41.83)(84.83,41.9)(84.87,42.0)(84.91,42.06)(84.97,42.15)(85.01,42.22)(85.06,42.29)(85.12,42.36)(85.16,42.43)
\path(85.16,42.43)(85.22,42.5)(85.26,42.56)(85.31,42.63)(85.37,42.7)(85.43,42.75)(85.48,42.81)(85.52,42.86)(85.58,42.91)(85.63,42.97)
\path(85.63,42.97)(85.69,43.02)(85.75,43.06)(85.8,43.11)(85.86,43.15)(85.91,43.2)(85.97,43.22)(86.02,43.27)(86.08,43.31)(86.13,43.34)
\path(86.13,43.34)(86.19,43.36)(86.26,43.4)(86.31,43.43)(86.37,43.45)(86.44,43.47)(86.5,43.5)(86.55,43.5)(86.62,43.52)(86.68,43.54)
\path(86.68,43.54)(86.73,43.56)(86.8,43.56)(86.86,43.58)(86.91,43.59)(86.98,43.59)(87.05,43.59)(87.11,43.59)(87.18,43.59)(87.23,43.59)
\path(87.23,43.59)(87.3,43.59)(87.37,43.58)(87.44,43.56)(87.51,43.56)(87.56,43.54)(87.63,43.52)(87.7,43.5)(87.76,43.5)(87.83,43.47)
\path(87.83,43.47)(87.91,43.45)(87.98,43.43)(88.05,43.4)(88.12,43.36)(88.19,43.34)(88.26,43.31)(88.33,43.27)(88.4,43.22)(88.47,43.2)
\path(88.47,43.2)(88.55,43.15)(88.62,43.11)(88.69,43.06)(88.76,43.02)(88.83,42.97)(88.91,42.91)(88.98,42.86)(89.05,42.81)(89.13,42.75)
\path(89.13,42.75)(89.2,42.7)(89.29,42.63)(89.37,42.56)(89.44,42.5)(89.51,42.43)(89.59,42.36)(89.68,42.29)(89.76,42.22)(89.83,42.15)
\path(89.83,42.15)(89.91,42.06)(90.0,42.0)(90.0,42.0)
\path(90.0,42.0)(90.0,42.0)(90.02,41.95)(90.06,41.91)(90.11,41.88)(90.15,41.84)(90.18,41.79)(90.22,41.75)(90.26,41.72)(90.29,41.68)
\path(90.29,41.68)(90.31,41.65)(90.36,41.61)(90.38,41.58)(90.41,41.54)(90.44,41.5)(90.48,41.47)(90.51,41.43)(90.52,41.4)(90.55,41.36)
\path(90.55,41.36)(90.58,41.34)(90.61,41.31)(90.62,41.27)(90.66,41.24)(90.68,41.2)(90.69,41.18)(90.72,41.15)(90.75,41.11)(90.76,41.09)
\path(90.76,41.09)(90.77,41.06)(90.8,41.02)(90.81,41.0)(90.83,40.97)(90.84,40.95)(90.87,40.91)(90.87,40.88)(90.88,40.86)(90.9,40.84)
\path(90.9,40.84)(90.91,40.81)(90.93,40.79)(90.94,40.75)(90.94,40.74)(90.94,40.72)(90.95,40.68)(90.97,40.66)(90.98,40.63)(90.98,40.61)
\path(90.98,40.61)(90.98,40.59)(90.98,40.58)(90.98,40.56)(90.98,40.54)(90.98,40.52)(91.0,40.5)(90.98,40.47)(90.98,40.45)(90.98,40.43)
\path(90.98,40.43)(90.98,40.41)(90.98,40.4)(90.98,40.38)(90.98,40.36)(90.97,40.34)(90.95,40.33)(90.94,40.31)(90.94,40.29)(90.94,40.27)
\path(90.94,40.27)(90.93,40.27)(90.91,40.25)(90.9,40.24)(90.88,40.22)(90.87,40.2)(90.87,40.2)(90.84,40.18)(90.83,40.18)(90.81,40.15)
\path(90.81,40.15)(90.8,40.15)(90.77,40.13)(90.76,40.13)(90.75,40.11)(90.72,40.11)(90.69,40.09)(90.68,40.09)(90.66,40.08)(90.62,40.08)
\path(90.62,40.08)(90.61,40.06)(90.58,40.06)(90.55,40.04)(90.52,40.04)(90.51,40.04)(90.48,40.02)(90.44,40.02)(90.41,40.02)(90.38,40.02)
\path(90.38,40.02)(90.36,40.02)(90.31,40.0)(90.29,40.0)(90.26,40.0)(90.22,40.0)(90.19,40.0)(90.15,40.0)(90.11,40.0)(90.06,40.0)
\path(90.06,40.0)(90.02,40.0)(90.0,40.0)(90.0,40.0)
\end{picture}

\end{center}
\caption{$\beta >0$}\label{f5}
\end{figure}

The respective generic behaviors of $t\longmapsto y(t)$ are
represented on Fig. \ref{f6}.


\setlength{\unitlength}{0.5mm}
\begin{figure}[h]
\begin{center} 

\begin{picture}(180,100)
\thinlines
\drawvector{20.0}{48.0}{60.0}{1}{0}
\drawvector{100.0}{48.0}{60.0}{1}{0}
\drawvector{22.0}{24.0}{54.0}{0}{1}
\drawvector{102.0}{24.0}{54.0}{0}{1}
\drawcenteredtext{160.0}{52.0}{$t$}
\drawcenteredtext{80.0}{52.0}{$t$}
\drawcenteredtext{106.0}{78.0}{$y$}
\drawcenteredtext{26.0}{78.0}{$y$}
\drawcenteredtext{48.0}{20.0}{$\beta =0$}
\drawcenteredtext{130.0}{20.0}{$\beta \neq 0$}
\path(22.0,48.0)(22.0,48.0)(22.22,48.83)(22.46,49.63)(22.7,50.43)(22.95,51.22)(23.19,51.97)(23.42,52.72)(23.65,53.45)(23.88,54.18)
\path(23.88,54.18)(24.12,54.86)(24.36,55.56)(24.59,56.22)(24.81,56.86)(25.04,57.5)(25.28,58.11)(25.51,58.7)(25.72,59.27)(25.95,59.84)
\path(25.95,59.84)(26.19,60.38)(26.4,60.91)(26.62,61.43)(26.86,61.93)(27.07,62.4)(27.29,62.86)(27.52,63.31)(27.75,63.75)(27.95,64.16)
\path(27.95,64.16)(28.18,64.55)(28.39,64.93)(28.62,65.29)(28.84,65.62)(29.04,65.95)(29.27,66.26)(29.47,66.56)(29.69,66.83)(29.89,67.11)
\path(29.89,67.11)(30.12,67.34)(30.32,67.58)(30.54,67.79)(30.75,67.98)(30.95,68.15)(31.15,68.3)(31.37,68.45)(31.57,68.58)(31.78,68.69)
\path(31.78,68.69)(31.97,68.79)(32.18,68.86)(32.38,68.91)(32.59,68.95)(32.79,68.98)(32.99,69.0)(33.18,68.98)(33.38,68.95)(33.59,68.91)
\path(33.59,68.91)(33.79,68.86)(33.97,68.79)(34.18,68.69)(34.38,68.58)(34.56,68.45)(34.75,68.3)(34.95,68.16)(35.15,67.98)(35.34,67.79)
\path(35.34,67.79)(35.52,67.58)(35.72,67.34)(35.9,67.11)(36.09,66.83)(36.27,66.56)(36.47,66.26)(36.65,65.95)(36.83,65.63)(37.02,65.29)
\path(37.02,65.29)(37.2,64.93)(37.38,64.55)(37.56,64.16)(37.74,63.75)(37.91,63.31)(38.09,62.86)(38.27,62.4)(38.45,61.93)(38.63,61.43)
\path(38.63,61.43)(38.81,60.91)(38.99,60.38)(39.15,59.84)(39.33,59.27)(39.5,58.7)(39.68,58.11)(39.84,57.5)(40.02,56.86)(40.18,56.22)
\path(40.18,56.22)(40.34,55.56)(40.52,54.86)(40.68,54.18)(40.86,53.45)(41.02,52.72)(41.18,51.99)(41.34,51.22)(41.5,50.43)(41.66,49.63)
\path(41.66,49.63)(41.83,48.83)(41.99,48.0)(42.0,48.0)
\path(42.0,48.0)(42.0,48.0)(42.18,47.15)(42.4,46.34)(42.59,45.54)(42.79,44.77)(43.0,44.0)(43.18,43.25)(43.4,42.52)(43.59,41.81)
\path(43.59,41.81)(43.79,41.11)(44.0,40.43)(44.18,39.77)(44.4,39.11)(44.59,38.49)(44.79,37.88)(45.0,37.29)(45.18,36.7)(45.4,36.13)
\path(45.4,36.13)(45.59,35.59)(45.79,35.06)(46.0,34.54)(46.18,34.06)(46.4,33.58)(46.59,33.11)(46.79,32.66)(47.0,32.25)(47.18,31.82)
\path(47.18,31.82)(47.4,31.44)(47.59,31.05)(47.79,30.7)(48.0,30.36)(48.18,30.03)(48.4,29.71)(48.59,29.42)(48.79,29.14)(49.0,28.88)
\path(49.0,28.88)(49.18,28.63)(49.4,28.4)(49.59,28.2)(49.79,28.01)(50.0,27.84)(50.18,27.68)(50.38,27.53)(50.59,27.4)(50.79,27.29)
\path(50.79,27.29)(51.0,27.2)(51.18,27.12)(51.38,27.06)(51.59,27.03)(51.79,27.0)(51.99,27.0)(52.18,27.0)(52.38,27.03)(52.59,27.06)
\path(52.59,27.06)(52.79,27.12)(52.99,27.2)(53.18,27.29)(53.38,27.4)(53.59,27.53)(53.79,27.68)(53.99,27.82)(54.18,28.01)(54.38,28.2)
\path(54.38,28.2)(54.59,28.4)(54.79,28.63)(54.99,28.87)(55.18,29.14)(55.38,29.42)(55.59,29.71)(55.79,30.03)(55.99,30.35)(56.18,30.7)
\path(56.18,30.7)(56.38,31.05)(56.59,31.44)(56.79,31.82)(56.99,32.24)(57.18,32.66)(57.38,33.11)(57.59,33.58)(57.79,34.06)(57.99,34.54)
\path(57.99,34.54)(58.18,35.06)(58.38,35.59)(58.59,36.13)(58.79,36.7)(58.99,37.27)(59.18,37.88)(59.38,38.49)(59.59,39.11)(59.79,39.77)
\path(59.79,39.77)(59.99,40.43)(60.18,41.11)(60.38,41.81)(60.59,42.52)(60.79,43.25)(60.99,44.0)(61.18,44.77)(61.38,45.54)(61.59,46.34)
\path(61.59,46.34)(61.79,47.15)(61.99,47.99)(62.0,48.0)
\path(62.0,48.0)(62.0,48.0)(62.15,48.7)(62.31,49.4)(62.47,50.09)(62.63,50.77)(62.79,51.43)(62.95,52.09)(63.11,52.72)(63.27,53.36)
\path(63.27,53.36)(63.43,53.97)(63.59,54.56)(63.75,55.15)(63.9,55.74)(64.08,56.31)(64.23,56.86)(64.4,57.4)(64.55,57.93)(64.72,58.43)
\path(64.72,58.43)(64.87,58.95)(65.04,59.43)(65.19,59.91)(65.36,60.38)(65.51,60.83)(65.68,61.27)(65.83,61.7)(66.0,62.11)(66.16,62.52)
\path(66.16,62.52)(66.31,62.91)(66.47,63.29)(66.62,63.65)(66.8,64.01)(66.95,64.36)(67.11,64.69)(67.26,65.0)(67.44,65.3)(67.59,65.59)
\path(67.59,65.59)(67.75,65.87)(67.91,66.15)(68.08,66.4)(68.23,66.63)(68.38,66.87)(68.55,67.08)(68.72,67.3)(68.87,67.48)(69.02,67.66)
\path(69.02,67.66)(69.19,67.83)(69.36,68.0)(69.51,68.13)(69.66,68.26)(69.83,68.38)(70.0,68.5)(70.16,68.58)(70.3,68.66)(70.47,68.73)
\path(70.47,68.73)(70.62,68.8)(70.8,68.83)(70.94,68.87)(71.11,68.88)(71.26,68.9)(71.44,68.88)(71.58,68.87)(71.75,68.83)(71.91,68.8)
\path(71.91,68.8)(72.08,68.75)(72.23,68.68)(72.38,68.59)(72.55,68.51)(72.72,68.4)(72.87,68.29)(73.02,68.16)(73.19,68.01)(73.36,67.86)
\path(73.36,67.86)(73.51,67.69)(73.66,67.51)(73.83,67.31)(74.0,67.12)(74.15,66.9)(74.3,66.68)(74.47,66.43)(74.62,66.18)(74.79,65.91)
\path(74.79,65.91)(74.94,65.63)(75.11,65.34)(75.26,65.04)(75.43,64.73)(75.58,64.4)(75.75,64.05)(75.91,63.7)(76.06,63.34)(76.23,62.95)
\path(76.23,62.95)(76.38,62.58)(76.55,62.16)(76.7,61.75)(76.87,61.33)(77.02,60.88)(77.19,60.43)(77.34,59.97)(77.51,59.5)(77.66,59.0)
\path(77.66,59.0)(77.83,58.5)(77.98,58.0)(78.0,58.0)
\path(102.0,48.0)(102.0,48.0)(102.19,48.83)(102.38,49.63)(102.59,50.43)(102.8,51.22)(103.0,51.97)(103.19,52.72)(103.38,53.45)(103.59,54.18)
\path(103.59,54.18)(103.8,54.86)(104.0,55.56)(104.19,56.22)(104.38,56.86)(104.59,57.5)(104.8,58.11)(105.0,58.7)(105.19,59.27)(105.38,59.84)
\path(105.38,59.84)(105.59,60.38)(105.8,60.91)(106.0,61.43)(106.19,61.93)(106.38,62.4)(106.59,62.86)(106.8,63.31)(107.0,63.75)(107.19,64.16)
\path(107.19,64.16)(107.38,64.55)(107.59,64.93)(107.8,65.29)(108.0,65.62)(108.19,65.95)(108.38,66.26)(108.59,66.56)(108.8,66.83)(109.0,67.11)
\path(109.0,67.11)(109.19,67.34)(109.38,67.58)(109.58,67.79)(109.8,67.98)(110.0,68.15)(110.19,68.3)(110.38,68.45)(110.58,68.58)(110.8,68.69)
\path(110.8,68.69)(111.0,68.79)(111.19,68.86)(111.38,68.91)(111.58,68.95)(111.8,68.98)(112.0,69.0)(112.19,68.98)(112.38,68.95)(112.58,68.91)
\path(112.58,68.91)(112.8,68.86)(113.0,68.79)(113.19,68.69)(113.38,68.58)(113.58,68.45)(113.79,68.3)(114.0,68.16)(114.19,67.98)(114.38,67.79)
\path(114.38,67.79)(114.58,67.58)(114.79,67.34)(115.0,67.11)(115.19,66.83)(115.38,66.56)(115.58,66.26)(115.79,65.95)(116.0,65.63)(116.19,65.29)
\path(116.19,65.29)(116.38,64.93)(116.58,64.55)(116.79,64.16)(116.98,63.75)(117.19,63.31)(117.38,62.86)(117.58,62.4)(117.79,61.93)(117.98,61.43)
\path(117.98,61.43)(118.19,60.91)(118.38,60.38)(118.58,59.84)(118.79,59.27)(118.98,58.7)(119.19,58.11)(119.38,57.5)(119.58,56.86)(119.79,56.22)
\path(119.79,56.22)(119.98,55.56)(120.19,54.86)(120.38,54.18)(120.58,53.45)(120.79,52.72)(120.98,51.99)(121.19,51.22)(121.38,50.43)(121.58,49.63)
\path(121.58,49.63)(121.79,48.83)(121.98,48.0)(122.0,48.0)
\path(122.0,48.0)(122.0,48.0)(122.16,47.36)(122.31,46.74)(122.48,46.13)(122.63,45.54)(122.8,44.95)(122.97,44.38)(123.12,43.83)(123.3,43.27)
\path(123.3,43.27)(123.47,42.75)(123.62,42.22)(123.8,41.72)(123.97,41.24)(124.13,40.75)(124.3,40.29)(124.48,39.84)(124.66,39.38)(124.83,38.95)
\path(124.83,38.95)(125.0,38.54)(125.18,38.15)(125.36,37.75)(125.52,37.38)(125.7,37.0)(125.88,36.65)(126.06,36.31)(126.25,36.0)(126.43,35.68)
\path(126.43,35.68)(126.61,35.38)(126.79,35.09)(126.97,34.81)(127.16,34.56)(127.33,34.31)(127.51,34.06)(127.7,33.84)(127.9,33.63)(128.08,33.43)
\path(128.08,33.43)(128.27,33.25)(128.46,33.08)(128.64,32.91)(128.83,32.77)(129.02,32.63)(129.22,32.5)(129.41,32.4)(129.61,32.31)(129.8,32.22)
\path(129.8,32.22)(130.0,32.15)(130.19,32.09)(130.39,32.04)(130.6,32.02)(130.8,32.0)(131.0,32.0)(131.19,32.0)(131.39,32.02)(131.6,32.04)
\path(131.6,32.04)(131.8,32.09)(132.0,32.15)(132.21,32.22)(132.41,32.31)(132.61,32.4)(132.83,32.5)(133.02,32.63)(133.24,32.77)(133.44,32.91)
\path(133.44,32.91)(133.66,33.08)(133.86,33.25)(134.08,33.43)(134.3,33.63)(134.5,33.84)(134.72,34.06)(134.94,34.31)(135.14,34.54)(135.36,34.81)
\path(135.36,34.81)(135.58,35.09)(135.8,35.38)(136.02,35.68)(136.25,35.99)(136.47,36.31)(136.69,36.65)(136.91,37.0)(137.13,37.38)(137.35,37.75)
\path(137.35,37.75)(137.58,38.15)(137.8,38.54)(138.02,38.95)(138.25,39.38)(138.47,39.83)(138.71,40.29)(138.94,40.75)(139.16,41.24)(139.39,41.72)
\path(139.39,41.72)(139.63,42.22)(139.86,42.75)(140.1,43.27)(140.33,43.83)(140.57,44.38)(140.8,44.95)(141.03,45.54)(141.27,46.13)(141.52,46.74)
\path(141.52,46.74)(141.75,47.36)(141.99,47.99)(142.0,48.0)
\path(142.0,48.0)(142.0,48.0)(142.08,48.18)(142.16,48.38)(142.24,48.58)(142.32,48.77)(142.41,48.95)(142.5,49.13)(142.58,49.31)(142.66,49.49)
\path(142.66,49.49)(142.75,49.66)(142.86,49.84)(142.94,50.0)(143.03,50.15)(143.13,50.31)(143.22,50.47)(143.33,50.63)(143.42,50.79)(143.52,50.93)
\path(143.52,50.93)(143.63,51.08)(143.72,51.22)(143.83,51.36)(143.94,51.49)(144.05,51.61)(144.14,51.75)(144.25,51.86)(144.36,52.0)(144.47,52.11)
\path(144.47,52.11)(144.58,52.22)(144.71,52.34)(144.82,52.45)(144.94,52.56)(145.05,52.65)(145.16,52.75)(145.28,52.84)(145.41,52.95)(145.52,53.04)
\path(145.52,53.04)(145.64,53.11)(145.77,53.2)(145.89,53.27)(146.02,53.36)(146.16,53.43)(146.27,53.5)(146.41,53.56)(146.53,53.63)(146.67,53.7)
\path(146.67,53.7)(146.8,53.75)(146.94,53.81)(147.08,53.86)(147.22,53.9)(147.36,53.95)(147.5,53.99)(147.63,54.02)(147.77,54.06)(147.91,54.09)
\path(147.91,54.09)(148.05,54.13)(148.21,54.15)(148.36,54.18)(148.5,54.2)(148.64,54.2)(148.8,54.22)(148.94,54.22)(149.11,54.24)(149.25,54.24)
\path(149.25,54.24)(149.41,54.24)(149.57,54.24)(149.72,54.24)(149.88,54.22)(150.05,54.2)(150.21,54.2)(150.36,54.18)(150.52,54.15)(150.69,54.13)
\path(150.69,54.13)(150.86,54.09)(151.02,54.06)(151.19,54.02)(151.36,54.0)(151.53,53.95)(151.71,53.9)(151.88,53.86)(152.05,53.81)(152.22,53.75)
\path(152.22,53.75)(152.41,53.7)(152.58,53.63)(152.77,53.56)(152.94,53.5)(153.13,53.43)(153.3,53.36)(153.5,53.27)(153.67,53.2)(153.86,53.11)
\path(153.86,53.11)(154.05,53.04)(154.24,52.95)(154.42,52.84)(154.61,52.75)(154.82,52.65)(155.0,52.56)(155.19,52.45)(155.39,52.34)(155.6,52.22)
\path(155.6,52.22)(155.8,52.11)(156.0,52.0)(156.0,52.0)
\end{picture}

\end{center}
\caption{$\beta \geq 0$}\label{f6}
\end{figure}
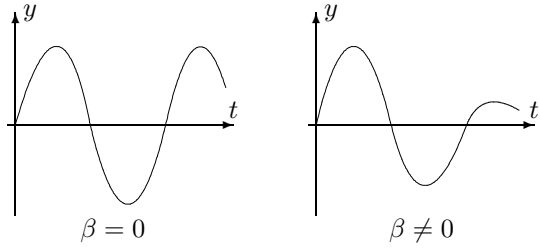

\noindent \it{It is important to observe that our
description of the behavior of
$t\longmapsto y(t)$ is true for the gradated
form of order 0, but also of any order when $\lambda 
\longrightarrow \infty .$}

\end{itemize}


\subsubsection{Characteristic equation}
If $\beta =0,$ the Hamiltonian $H_n={\displaystyle {1\over 2}}
 (P^2_1+P^2_2)$ has two cyclic coordinates~: $x$ and $z$ and
therefore $p_x={\cos\ }\theta (0)$ and $p_z=\lambda $ are first
integrals. The equation $H_n=1/2,$ with $P_1={\displaystyle
 {p_x+p_z y^2/2\over \sqrt{a}}}$ and $P_2={\displaystyle 
{p_y\over \sqrt{c}}}$ takes the form~:
$$(\sqrt{c}\ {\dot y})^2+\Bigl ( {p_x+p_z  y^2/2\over 
\sqrt{a}}\Bigr )^2=1$$

\no Introducing~: $d\tau ={\displaystyle {dt\over \sqrt{a}
 \sqrt{c}}}$ it becomes~:
$$\Bigl ( {dy \over d\tau }\Bigr )^2\ +\ (p_x+p_z y^2/2)^2=a$$

\noindent where $a=(1+\alpha y)^2.$ Hence we get~:
\begin{eqnarray}  \label{4.111}
\Bigl ( {\displaystyle {dy\over d\tau }}\Bigr ) ^2 = F(y)
\end{eqnarray}

\noindent where $F(y)=(1+\alpha y)^2-(p_x+p_z y^2/2)^2.$
The analysis is based on the \it{roots of the quartic} $F(y).$
We assume $\lambda >0.$

We observe that $F$ can be factorized as $F_1F_2$ with~:
$$F_1 =(1+\alpha y)-(p_x+p_z y^2/2)\ \ ,\ \ \ F_2 =(1+\alpha y)-
(p_z+p_z y^2/2)$$

\noindent and we can write~:
$$F(y)=\Bigl (2m^2 -{\lambda \over 2} (y-{\alpha \over \lambda 
})^2\Bigr ) \Bigl ( 2m'' +{\lambda \over 2} (y+{\alpha \over
 \lambda })^2\Bigr )$$

\noindent where~:
$2m^2 =1-p_x +{\alpha ^2\over 2\lambda }\ ,\ \ \ 2m''=1+p_x-
 {\alpha ^2\over 2\lambda }$

\noindent and~:
$m^2+m''=1\ \ ,\ \ m\geq 0.$

\noindent If we set~:
$\eta =\ {\sqrt{\lambda } y\over 2m} - {\alpha \over 
2m\sqrt{\lambda }}\ \ ,\ \ {\bar \eta}\ =\ {\sqrt{\lambda }
y\over 2m} + {\alpha \over 2m\sqrt{\lambda }}$

\noindent we can write~:
\begin{eqnarray}  \label{4.12'}
F(y)=4m^2(1-\eta ^2)(m''+m^2 {\bar \eta }^2)
\end{eqnarray}

\noindent $F$ is a quartic whose roots on $\C $ are $\eta =\pm 
1\ ,$ \ ${\bar \eta }=\pm \ {\displaystyle {\sqrt{m''}\over m}}.$

The case $m''=0$ is called critical and it corresponds to a
double root for $F.$ We have~:
\begin{lem}
In the strict case $\alpha \not = 0,$ there exist geodesics
starting from $0$ which are critical.
\end{lem}
\paragraph{Geometric interpretation \\}
The critical geodesics project in the $(\theta ,$ ${\dot 
\theta })$ phase space onto a separatrix, see Fig. \ref{f4}.

The characteristic equation can be put into a normal form
using an \it{homographic transformation} to normalize the
roots of $F.$ The procedure is standard, see \cite{L}. We
proceed as follows~; $F$ is factorized into $F_1 F_2$ and we
consider the \it{pencil} $F_1+\nu F_2$ of two quadratic forms.
If $\alpha \not = 0,$ there exist two distinct real numbers
$\nu_1, \nu_2$ such that $F_1+\nu F_2$ is a \it{perfect square}~:
$K_1(y-p)^2, K_2 (y-q)^2.$ Using the homographic transformation~:
\begin{eqnarray}   \label{4.13'}
u={\displaystyle {y-p\over y-q}}\ \ ,
\end{eqnarray}

\noindent the characteristic equation can be written in the
normal form~:
\begin{eqnarray}    \label{4.14}
{\displaystyle {dy\over \sqrt{F(y)}}}\ =\ {\displaystyle 
{(p-q)^{-1}\ du\over \sqrt{(A_1 u^2+B_1)(A_2u^2+B_2)}}}.
\end{eqnarray}

The right hand side corresponds to an integrand of an
\it{elliptic integral of the first kind.} More precisely,
excepted the critical case $m''=0,$ the solution $y$ in the
$u$-coordinate can be computed as follows~:
\begin{itemize}
\item if the quartic $F$ admits two real roots,
$u$ can be parametrized using the cn \it{Jacobi function}~;

\item if the quartic $F$ admits four real roots, $u$ can be
parametrized using the dn \it{Jacobi function.}
\end{itemize}

If $\alpha =0,$ the analysis is simpler, indeed $F(y)$ can be
written~:
$$F(y)=4k^2 (1 -\eta ^2)\ ({k'}^2 +k^2 \ \eta ^2)$$

\noindent where $\eta ={\displaystyle {\sqrt{\lambda }y
\over 2k}}$ and $\eta $ can be computed using only the cn
function.

\begin{prop}
We have two cases~:
\begin{itemize}
\item[(i)] If $\alpha =0,$ \ $y={\displaystyle {2k\over 
\sqrt{\lambda }}} \eta $ where $\eta $ is the \rm{cn}
Jacobi function.

\item[(ii)] If $\alpha \not =0,$ $y$ is generically the
image by an homography of the \rm{cn} or \rm{dn} Jacobi function.
\end{itemize}
\end{prop}

\paragraph{Geometric interpretation}
If $\alpha =0,$ the motion of $y$ is a cn whose amplitude is
${\displaystyle {2k\over \sqrt{\lambda }}}.$ The motion is
symmetric with respect to $y=0$ and the amplitude tends to $0$
when $\lambda $ tends to the infinity, see Fig. \ref{f7}.

\noindent If $\alpha \not = 0,$ we can expand~:
$y={\displaystyle {ua-p\over u-1}}$ near $u=0.$ The motion of
$y$ is no more symmetric with respect to $y=0$ and there is a
\it{shift.} Hence $y$ can be approximated by a constant plus
a cn or dn motion.


\setlength{\unitlength}{0.5mm}
\begin{figure}[h]
\begin{center} 

\begin{picture}(180,100)
\thinlines
\drawvector{20.0}{48.0}{60.0}{1}{0}
\drawvector{100.0}{48.0}{60.0}{1}{0}
\drawvector{22.0}{24.0}{54.0}{0}{1}
\drawvector{102.0}{24.0}{54.0}{0}{1}
\drawcenteredtext{160.0}{52.0}{$\tau$}
\drawcenteredtext{80.0}{52.0}{$\tau$}
\drawcenteredtext{106.0}{78.0}{$y$}
\drawcenteredtext{26.0}{78.0}{$y$}
\drawcenteredtext{48.0}{20.0}{$\alpha =0$}
\drawcenteredtext{130.0}{20.0}{$\alpha \neq 0$}
\path(22.0,48.0)(22.0,48.0)(22.2,48.83)(22.45,49.63)(22.7,50.43)(22.95,51.22)(23.19,51.97)(23.42,52.72)(23.63,53.45)(23.87,54.18)
\path(23.87,54.18)(24.12,54.86)(24.36,55.56)(24.59,56.22)(24.79,56.86)(25.04,57.5)(25.28,58.11)(25.51,58.7)(25.7,59.27)(25.95,59.84)
\path(25.95,59.84)(26.19,60.38)(26.38,60.9)(26.62,61.43)(26.86,61.93)(27.05,62.4)(27.29,62.86)(27.52,63.31)(27.75,63.75)(27.95,64.16)
\path(27.95,64.16)(28.18,64.55)(28.37,64.93)(28.62,65.29)(28.84,65.62)(29.04,65.94)(29.27,66.26)(29.45,66.55)(29.69,66.83)(29.87,67.11)
\path(29.87,67.11)(30.12,67.33)(30.3,67.58)(30.54,67.79)(30.75,67.98)(30.95,68.15)(31.13,68.3)(31.37,68.44)(31.55,68.58)(31.78,68.69)
\path(31.78,68.69)(31.95,68.79)(32.18,68.86)(32.38,68.91)(32.59,68.94)(32.79,68.98)(32.99,69.0)(33.18,68.98)(33.38,68.94)(33.59,68.91)
\path(33.59,68.91)(33.79,68.86)(33.97,68.79)(34.18,68.69)(34.38,68.58)(34.56,68.44)(34.75,68.3)(34.95,68.16)(35.15,67.98)(35.34,67.79)
\path(35.34,67.79)(35.52,67.58)(35.72,67.33)(35.9,67.11)(36.09,66.83)(36.27,66.55)(36.47,66.26)(36.65,65.94)(36.83,65.62)(37.02,65.29)
\path(37.02,65.29)(37.2,64.93)(37.38,64.55)(37.56,64.16)(37.74,63.75)(37.9,63.31)(38.09,62.86)(38.27,62.4)(38.45,61.93)(38.63,61.43)
\path(38.63,61.43)(38.81,60.9)(38.99,60.38)(39.15,59.84)(39.33,59.27)(39.5,58.7)(39.68,58.11)(39.84,57.5)(40.02,56.86)(40.18,56.22)
\path(40.18,56.22)(40.34,55.56)(40.52,54.86)(40.68,54.18)(40.86,53.45)(41.02,52.72)(41.18,51.99)(41.34,51.22)(41.5,50.43)(41.65,49.63)
\path(41.65,49.63)(41.83,48.83)(41.99,48.0)(42.0,48.0)
\path(42.0,48.0)(42.0,48.0)(42.18,47.15)(42.4,46.34)(42.59,45.54)(42.79,44.77)(43.0,44.0)(43.18,43.25)(43.4,42.52)(43.59,41.81)
\path(43.59,41.81)(43.79,41.11)(44.0,40.43)(44.18,39.77)(44.4,39.11)(44.59,38.49)(44.79,37.88)(45.0,37.29)(45.18,36.7)(45.4,36.13)
\path(45.4,36.13)(45.59,35.59)(45.79,35.06)(46.0,34.54)(46.18,34.06)(46.4,33.58)(46.59,33.11)(46.79,32.65)(47.0,32.25)(47.18,31.8)
\path(47.18,31.8)(47.4,31.44)(47.59,31.04)(47.79,30.7)(48.0,30.36)(48.18,30.03)(48.4,29.7)(48.59,29.42)(48.79,29.12)(49.0,28.87)
\path(49.0,28.87)(49.18,28.62)(49.4,28.38)(49.59,28.2)(49.79,28.01)(50.0,27.84)(50.18,27.68)(50.38,27.53)(50.59,27.38)(50.79,27.29)
\path(50.79,27.29)(51.0,27.2)(51.18,27.12)(51.38,27.04)(51.59,27.03)(51.79,27.0)(51.99,27.0)(52.18,27.0)(52.38,27.03)(52.59,27.04)
\path(52.59,27.04)(52.79,27.12)(52.99,27.2)(53.18,27.29)(53.38,27.38)(53.59,27.53)(53.79,27.68)(53.99,27.8)(54.18,28.01)(54.38,28.2)
\path(54.38,28.2)(54.59,28.38)(54.79,28.62)(54.99,28.87)(55.18,29.12)(55.38,29.42)(55.59,29.7)(55.79,30.03)(55.99,30.35)(56.18,30.7)
\path(56.18,30.7)(56.38,31.04)(56.59,31.44)(56.79,31.8)(56.99,32.24)(57.18,32.65)(57.38,33.11)(57.59,33.58)(57.79,34.06)(57.99,34.54)
\path(57.99,34.54)(58.18,35.06)(58.38,35.59)(58.59,36.13)(58.79,36.7)(58.99,37.27)(59.18,37.88)(59.38,38.49)(59.59,39.11)(59.79,39.77)
\path(59.79,39.77)(59.99,40.43)(60.18,41.11)(60.38,41.81)(60.59,42.52)(60.79,43.25)(60.99,44.0)(61.18,44.77)(61.38,45.54)(61.59,46.34)
\path(61.59,46.34)(61.79,47.15)(61.99,47.99)(62.0,48.0)
\path(62.0,48.0)(62.0,48.0)(62.15,48.7)(62.31,49.4)(62.47,50.09)(62.63,50.77)(62.79,51.43)(62.95,52.09)(63.11,52.72)(63.27,53.36)
\path(63.27,53.36)(63.43,53.97)(63.59,54.56)(63.75,55.15)(63.9,55.74)(64.08,56.31)(64.23,56.86)(64.4,57.4)(64.55,57.93)(64.72,58.43)
\path(64.72,58.43)(64.87,58.95)(65.04,59.43)(65.19,59.9)(65.36,60.38)(65.51,60.83)(65.68,61.27)(65.83,61.7)(66.0,62.11)(66.16,62.52)
\path(66.16,62.52)(66.3,62.9)(66.47,63.29)(66.62,63.65)(66.8,64.01)(66.94,64.36)(67.11,64.69)(67.26,65.0)(67.44,65.3)(67.58,65.58)
\path(67.58,65.58)(67.75,65.87)(67.91,66.15)(68.08,66.4)(68.23,66.62)(68.37,66.87)(68.55,67.08)(68.72,67.3)(68.87,67.48)(69.01,67.66)
\path(69.01,67.66)(69.19,67.83)(69.36,68.0)(69.51,68.12)(69.66,68.26)(69.83,68.37)(70.0,68.5)(70.16,68.58)(70.3,68.66)(70.47,68.73)
\path(70.47,68.73)(70.62,68.8)(70.8,68.83)(70.94,68.87)(71.11,68.87)(71.26,68.9)(71.44,68.87)(71.58,68.87)(71.75,68.83)(71.91,68.8)
\path(71.91,68.8)(72.08,68.75)(72.23,68.68)(72.37,68.58)(72.55,68.51)(72.72,68.4)(72.87,68.29)(73.01,68.16)(73.19,68.01)(73.36,67.86)
\path(73.36,67.86)(73.51,67.69)(73.66,67.51)(73.83,67.3)(74.0,67.12)(74.15,66.9)(74.3,66.68)(74.47,66.43)(74.62,66.18)(74.79,65.91)
\path(74.79,65.91)(74.94,65.62)(75.11,65.33)(75.26,65.04)(75.43,64.73)(75.58,64.4)(75.75,64.05)(75.91,63.7)(76.05,63.34)(76.23,62.95)
\path(76.23,62.95)(76.37,62.58)(76.55,62.15)(76.69,61.75)(76.87,61.33)(77.01,60.88)(77.19,60.43)(77.33,59.97)(77.51,59.5)(77.66,59.0)
\path(77.66,59.0)(77.83,58.5)(77.98,58.0)(78.0,58.0)
\path(102.0,54.0)(102.0,54.0)(102.19,54.79)(102.37,55.56)(102.58,56.31)(102.8,57.06)(103.0,57.79)(103.19,58.5)(103.37,59.2)(103.58,59.88)
\path(103.58,59.88)(103.8,60.54)(104.0,61.2)(104.19,61.83)(104.37,62.43)(104.58,63.04)(104.8,63.63)(105.0,64.19)(105.19,64.75)(105.37,65.26)
\path(105.37,65.26)(105.58,65.8)(105.8,66.3)(106.0,66.8)(106.19,67.26)(106.37,67.72)(106.58,68.16)(106.8,68.58)(107.0,69.0)(107.19,69.37)
\path(107.19,69.37)(107.37,69.76)(107.58,70.12)(107.8,70.47)(108.0,70.8)(108.19,71.11)(108.37,71.4)(108.58,71.68)(108.8,71.94)(109.0,72.19)
\path(109.0,72.19)(109.19,72.43)(109.37,72.62)(109.58,72.83)(109.8,73.01)(110.0,73.19)(110.19,73.33)(110.37,73.48)(110.58,73.58)(110.8,73.69)
\path(110.8,73.69)(111.0,73.79)(111.19,73.87)(111.37,73.91)(111.58,73.94)(111.8,73.98)(112.0,74.0)(112.19,73.98)(112.37,73.94)(112.58,73.91)
\path(112.58,73.91)(112.8,73.87)(113.0,73.8)(113.19,73.69)(113.37,73.58)(113.58,73.48)(113.79,73.33)(114.0,73.19)(114.19,73.01)(114.37,72.83)
\path(114.37,72.83)(114.58,72.62)(114.79,72.43)(115.0,72.19)(115.19,71.94)(115.37,71.68)(115.58,71.4)(115.79,71.11)(116.0,70.8)(116.19,70.47)
\path(116.19,70.47)(116.37,70.12)(116.58,69.76)(116.79,69.37)(116.98,69.0)(117.19,68.58)(117.37,68.16)(117.58,67.72)(117.79,67.26)(117.98,66.8)
\path(117.98,66.8)(118.19,66.3)(118.37,65.8)(118.58,65.26)(118.79,64.75)(118.98,64.19)(119.19,63.63)(119.37,63.04)(119.58,62.43)(119.79,61.83)
\path(119.79,61.83)(119.98,61.2)(120.19,60.54)(120.37,59.88)(120.58,59.2)(120.79,58.5)(120.98,57.79)(121.19,57.06)(121.37,56.31)(121.58,55.56)
\path(121.58,55.56)(121.79,54.79)(121.98,54.0)(122.0,54.0)
\path(122.0,54.0)(122.0,54.0)(122.19,53.11)(122.37,52.27)(122.58,51.43)(122.8,50.61)(123.0,49.81)(123.19,49.02)(123.37,48.27)(123.58,47.52)
\path(123.58,47.52)(123.8,46.79)(124.0,46.08)(124.19,45.38)(124.37,44.7)(124.58,44.04)(124.8,43.4)(125.0,42.77)(125.19,42.15)(125.37,41.58)
\path(125.37,41.58)(125.58,41.0)(125.8,40.45)(126.0,39.9)(126.19,39.4)(126.37,38.88)(126.58,38.4)(126.8,37.93)(127.0,37.5)(127.19,37.06)
\path(127.19,37.06)(127.37,36.65)(127.58,36.25)(127.79,35.88)(128.0,35.52)(128.19,35.15)(128.38,34.84)(128.6,34.54)(128.8,34.25)(129.0,33.97)
\path(129.0,33.97)(129.19,33.72)(129.38,33.47)(129.58,33.25)(129.8,33.06)(130.0,32.88)(130.19,32.7)(130.38,32.56)(130.58,32.43)(130.8,32.31)
\path(130.8,32.31)(131.0,32.22)(131.19,32.13)(131.38,32.06)(131.58,32.02)(131.8,32.0)(132.0,32.0)(132.19,32.0)(132.38,32.02)(132.58,32.06)
\path(132.58,32.06)(132.77,32.13)(133.0,32.2)(133.19,32.31)(133.38,32.43)(133.58,32.56)(133.77,32.7)(134.0,32.86)(134.19,33.06)(134.38,33.25)
\path(134.38,33.25)(134.6,33.47)(134.77,33.72)(135.0,33.97)(135.19,34.25)(135.38,34.54)(135.6,34.84)(135.77,35.15)(136.0,35.5)(136.19,35.88)
\path(136.19,35.88)(136.38,36.25)(136.6,36.65)(136.77,37.06)(137.0,37.49)(137.19,37.93)(137.38,38.4)(137.6,38.88)(137.77,39.4)(138.0,39.9)
\path(138.0,39.9)(138.19,40.45)(138.38,41.0)(138.58,41.58)(138.77,42.15)(139.0,42.77)(139.19,43.4)(139.38,44.04)(139.58,44.7)(139.77,45.38)
\path(139.77,45.38)(140.0,46.06)(140.19,46.79)(140.38,47.52)(140.58,48.27)(140.77,49.02)(141.0,49.81)(141.19,50.61)(141.38,51.43)(141.58,52.27)
\path(141.58,52.27)(141.77,53.11)(142.0,53.99)(142.0,54.0)
\path(142.0,54.0)(142.0,54.0)(142.16,54.65)(142.32,55.33)(142.47,55.97)(142.63,56.61)(142.77,57.25)(142.94,57.86)(143.11,58.47)(143.27,59.06)
\path(143.27,59.06)(143.44,59.65)(143.6,60.2)(143.75,60.77)(143.91,61.31)(144.07,61.84)(144.22,62.38)(144.38,62.88)(144.55,63.38)(144.72,63.88)
\path(144.72,63.88)(144.88,64.36)(145.02,64.8)(145.19,65.26)(145.35,65.72)(145.5,66.15)(145.66,66.55)(145.83,66.97)(146.0,67.37)(146.16,67.75)
\path(146.16,67.75)(146.32,68.12)(146.47,68.48)(146.63,68.83)(146.77,69.18)(146.94,69.5)(147.11,69.8)(147.27,70.12)(147.44,70.41)(147.6,70.69)
\path(147.6,70.69)(147.75,70.94)(147.91,71.19)(148.07,71.44)(148.22,71.69)(148.38,71.91)(148.55,72.12)(148.72,72.3)(148.88,72.51)(149.02,72.69)
\path(149.02,72.69)(149.19,72.83)(149.35,73.0)(149.5,73.12)(149.66,73.26)(149.83,73.37)(150.0,73.5)(150.16,73.58)(150.32,73.66)(150.47,73.73)
\path(150.47,73.73)(150.63,73.8)(150.77,73.83)(150.94,73.87)(151.11,73.91)(151.27,73.91)(151.44,73.93)(151.6,73.91)(151.75,73.87)(151.91,73.86)
\path(151.91,73.86)(152.07,73.8)(152.22,73.76)(152.38,73.69)(152.55,73.61)(152.72,73.51)(152.88,73.41)(153.02,73.3)(153.19,73.18)(153.35,73.04)
\path(153.35,73.04)(153.5,72.87)(153.66,72.73)(153.83,72.55)(154.0,72.37)(154.16,72.16)(154.32,71.97)(154.47,71.75)(154.63,71.51)(154.77,71.26)
\path(154.77,71.26)(154.94,71.01)(155.11,70.76)(155.27,70.48)(155.44,70.19)(155.6,69.87)(155.75,69.58)(155.91,69.25)(156.07,68.91)(156.22,68.55)
\path(156.22,68.55)(156.38,68.22)(156.55,67.83)(156.72,67.44)(156.88,67.05)(157.02,66.66)(157.19,66.25)(157.35,65.8)(157.5,65.37)(157.66,64.93)
\path(157.66,64.93)(157.83,64.47)(158.0,64.0)(158.0,64.0)
\end{picture}

\end{center}
\caption{$\alpha \geq 0$}\label{f7}
\end{figure}
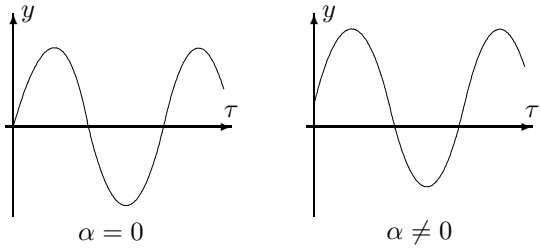


\subsubsection{Integral formulas in the general conservative case} 
If the metric $g$ does not depend on $x,$ it is convenient to
use the following integral formulas from \cite{Ku} to compute $x$
and $z$ in terms of $y.$

We denote by $e(t), t\in [0,T]$ a normal geodesic starting from
0 and we assume that the component~: $t\longmapsto y(t)$
oscillates periodically with period ${\cal P}.$ We denote by
$0<t_1\ <\cdots <t_N \leq T$ the successive times such that
$y(t_i)=0.$ We introduce~:
$$\sigma = \left \{ 
\begin{array}{lll}
{\rm{sign}}\ {\dot y}(0) \ \ \ & \rm{if}\ \ \ &  
{\dot y}(0)\not = 0\\
{\rm{sign}}\ {\ddot y}(0)\ \ & \rm{if} \ \ \ & {\dot y}(0)=0
\end{array}
\right .$$

\noindent and we set~:
$$y_+=\underset{t\in [0,P]}{\max} y(t) \ \ ,\ \ \ y_-=
\underset{t\in [0,P]}{\min} y(t) $$

Parametrizing the geodesics by $y$ we must integrate the equations~:
$${dx\over dy}\ =\ \f{\sqrt{c}}{\sqrt{a}}\f{P_1}{P_2}
\ \ ,\ \ {dz\over dy}\ =\ {y^2\over 2} {\sqrt{c}\over \sqrt{a}}
 {P_1\over P_2}\ \ ,\ \ \ dt={\sqrt{c}\over P_2} dy$$

\noindent where $P_2 (y)=\sigma \sqrt{1-P_1^2 (y)}$ for
$t\in [0,t_1].$

This allows to get explicit integral formulas. In particular
if $y(T)=0$ for $T=t_N$ we get~:
\begin{itemize}
\item \underline{$N$ odd}
\begin{equation}  \label{param}
\begin{split}
x(T) & =  2{\displaystyle {\int }^{y\sigma }_0}\ \sigma 
{\displaystyle {\sqrt{c}\over \sqrt{a}}}\ {\displaystyle 
{P_1 (y)\over \sqrt{1-P_1^2(y)}}}\ dy+(N-1)\ {\displaystyle 
{\int ^{y_+}_{y_-}}}\ {\displaystyle {\sqrt{c}\over \sqrt{a}}}\ 
{\displaystyle {P_1(y)\over \sqrt{1-P^2_1(y)}}}\ dy \\
z(T) & =  {\displaystyle {\int ^{y\sigma }_0}}\ \sigma
{\displaystyle {\sqrt{c}\over \sqrt{a}}}\ {\displaystyle 
{y^2\ P_1(y)\over \sqrt{1-P_1^2(y)}}}\ dy+(N-1) {\displaystyle \ 
{\int ^{y_+}_{y_-}}}\ {\displaystyle {\sqrt{c}\over 2\sqrt{a}}}\ 
{\displaystyle {y^2\ P_1 (y)\over \sqrt{1-P^2_1(y)}}}\ dy
\end{split}
\end{equation}

\item \underline{$N$ even}
\begin{equation} \label{even}
\begin{split}
x(T) & =  N\ {\displaystyle {\int ^{y_+}_{y_-}}}\ 
{\displaystyle {\sqrt{c}\over \sqrt{a}}}\ {\displaystyle 
{P_1 (y)\over \sqrt{1-P_1^2 (y)}}}\ dy \\
z(T) & =  N\ {\displaystyle {\int ^{y_+}_{y_-}}}\ 
{\displaystyle {\sqrt{c}\over 2\sqrt{a}}}\ {\displaystyle 
{y^2P_1 (y)\over \sqrt{1-P^2_1(y)}}}\ dy
\end{split}
\end{equation}
\end{itemize}

\noindent and the period is given by~:
\begin{eqnarray}
{\cal P} = 2\ {\displaystyle {\int ^{y_+}_{y_-}}}\ 
{\displaystyle {\sqrt{c}\over \sqrt{1-P^2_1 (y)}}}\ dy.
\end{eqnarray}

The integrands have simple poles when $P_1 (y)=\pm 1$ so the
integrals are well-defined.


\subsubsection{The return mapping}
The main geometric object to understand the role of abnormal
trajectories in the problem is the \it{return mapping.} Indeed
if we consider the trace of the sphere and the wave front
in the plane $y=0$~:
$${\widetilde S}(0,r)=S(0,r)\cap (y=0)\ ,\ \ \ 
{\widetilde W}(0,r)=W(0,r)\cap (y=0) \ ,$$

\noindent they are in the image of the following mappings.

\begin{defi}
Let $e\ :\ (t\in [0,T],$ $\theta (0),\lambda )\longmapsto
 (x(t),$ $y(t),$ $z(t))$ be a normal geodesic, parametrized
by arc-length. If $y(t) \neq 0,$ we can define
$0<t_1 <\cdots <t_N\leq T$ as the times corresponding to
$y(t_i )=0.$ The first return mapping is~:
$$R_1 \ :\ (\lambda ,\ \theta (0))\in D_1 \longmapsto 
 (x(t_1),z(t_1))$$

\noindent and more generally the n-th return mapping is the map~:
$$R_n \ :\ (\lambda ,\theta (0))\in D_n \longmapsto 
(x(t_n),z(t_n))$$ 

\noindent where $D_i$ are the domains.
\end{defi}

If the length is fixed to $r$, we observe that ${\widetilde W}
(0,r)$ is the union of the image of the
return mappings and $(\pm r,0)$ which are the end-points of
the abnormal geodesics.

The remaining of this Section is devoted to the analysis of
the return mapping. We proceed by perturbations of the flat case.
We shall estimate the asymptotic expansions of $\widetilde{S}$
and $\widetilde{W}$ in the abnormal direction. They are an union
of curves in the plane. Such a curve is subanalytic if and only
if it admits a Puiseux expansion. \it{It is a practical criterion to
measure the transcendence of the sphere and wave front in the
abnormal direction.}


\subsubsection{The pendulum and the elastica in the flat case}
In the flat case the equation (\ref{4.11}) is a simple pendulum~:
$${d^2\theta \over ds^2}+{\sin\ }\theta =0$$

\noindent where $s=t\sqrt{\lambda }$, $t$ is the arc-length parameter
and $y=-{\displaystyle {d\theta \over \sqrt{\lambda }ds}}.$
In particular if $y(0)=0,$ we have ${\displaystyle 
{d\theta \over ds}}=0.$ We get~:
$$\Bigl ( {d\theta \over ds}\Bigr ) ^2 \ =\ 2({\cos\ }\theta 
-{\cos\ }\theta  (0))$$

The integration is standard using elliptic integrals \cite{L}.
The characteristic equation takes the form~:
$${\dot y}^2=\Bigl ( 1-p_x-{y^2\over 2} p_z\Bigr )\ \Bigl ( 1
+p_x+{y^2\over 2} p_z\Bigr )$$

\noindent and we introduce $k, k'\in [0,1]$ by setting~:
$$2k^2=1-p_x\ ,\ \ 2{k'}^2=1+p_x$$

\noindent where $p_x={\cos\ }\theta (0).$ We set~:
$\eta ={\displaystyle {y\sqrt{\lambda }\over 2k}}$ and we get
the equation~:
$${\displaystyle {{\dot \eta }^2\over \lambda }}\ =
\ (1-\eta ^2)\ ({k'}^2+k^2\eta ^2)$$

\noindent We integrate with $\eta (0)=y(0)=0$ and we choose
the branch ${\dot \eta } (0)>0$ 
co\-rres\-ponding to ${\dot y}(0)
={\sin\ }\theta (0)>0.$ We get using the cn \it{Jacobi function~:}
$$\eta (t) = -\rm{cn }(K(k)+t\sqrt{\lambda },k)$$

\noindent where $4K(k)$ is the period, $K$ being the complete
elliptic integral of the first kind~:
$$K(k)=\int ^1_0\ {d\eta \over \sqrt{(1-\eta ^2)({k'}^2+
k^2\eta ^2)}}\ =\ \int ^{\pi /2}_0\ (1-k^2 {\sin}^2\theta
)^{-\inv{2}} d\theta $$

\noindent Hence
\begin{eqnarray}  \label{4.18}
y(t)=-{\displaystyle {2k\over \sqrt{\lambda }}}\ \rm{cn }(u,k)\ ,
\ \ u=K+t\sqrt{\lambda }
\end{eqnarray}

\noindent which coincides with the formula obtained by integrating
the pendulum.

The components $y$ and $z$ can be computed by quadratures and we get~:
\begin{equation}  \label{4.20}
\begin{split}
x(t)& = -t+{\displaystyle {2\over \sqrt{\lambda }}}\ (E(u)-E)\\
z(t) & =  {\displaystyle {2\over 3\lambda ^{3/2}}}\ [(2k^2-1)\ 
(E(u)-E(K))+{k'}^2 t\sqrt{\lambda }+2k^2\rm{\rm{sn }u\ \rm{cn }u\
\rm{dn }u}]
\end{split}
\end{equation}

\noindent where $E$ is the complete elliptic integral of the
second kind~:
$$E(k)={\displaystyle {\int ^{K(k)}_0}}\ \rm{dn }^2u
\ du={\displaystyle
 {\int ^{\pi /2}_0}}\ (1-k^2\sin^2\theta )^{1/2}d\theta $$

The previous parametrization corresponds to geodesics with
$\lambda>0, \theta (0)\in ]0,\pi [.$ The solutions corresponding
to $\lambda >0, \theta (0)\in ]-\pi ,0[$ are deduced using the
symmetry~: $S_1 : (x,y,z)\longmapsto (x,-y,z).$ The solutions
corresponding to $\lambda <0$ are deduced using the symmetry~:
$S_2 \ :\ (x,y,z)\longmapsto (-x,y,-z).$ The solutions with
$\lambda =0$ play no role in our analysis.

\paragraph{Elastica \\ }
The projections of the geodesics on the plane $(x,y)$ are
parametrized by~:
$$y(t)=-{2k\over \sqrt{\lambda }} \rm{cn }(u,k)\ ,\ \ x(t)=-t+
{2\over \sqrt{\lambda }} (E(u)-E)$$

\noindent They are precisely the \it{inflexional elastica}
described in \cite{Love}.

They take various shapes whose typical ones are represented
on Fig. \ref{f8}.


\setlength{\unitlength}{0.4mm}
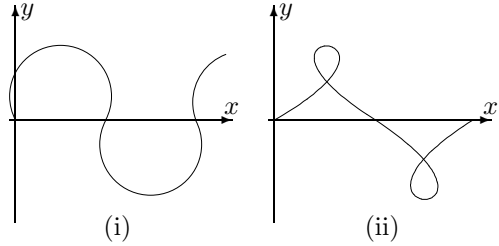
\begin{figure}[h]
\begin{center} 

\begin{picture}(180,100)
\thinlines
\drawvector{10.0}{48.0}{74.0}{1}{0}
\drawvector{96.0}{48.0}{74.0}{1}{0}
\drawvector{12.0}{14.0}{72.0}{0}{1}
\drawvector{98.0}{14.0}{72.0}{0}{1}
\drawcenteredtext{170.0}{52.0}{$x$}
\drawcenteredtext{84.0}{52.0}{$x$}
\drawcenteredtext{102.0}{84.0}{$y$}
\drawcenteredtext{16.0}{84.0}{$y$}
\drawcenteredtext{46.0}{12.0}{(i)}
\drawcenteredtext{134.0}{12.0}{(ii)}
\drawarc{27.0}{56.0}{34.0}{2.65}{0.48}
\drawarc{57.0}{40.0}{34.0}{5.79}{3.63}
\drawarc{88.0}{54.0}{34.16}{2.78}{4.35}
\path(98.0,48.0)(98.0,48.0)(99.63,48.95)(101.19,49.88)(102.69,50.79)(104.09,51.68)(105.44,52.56)(106.7,53.43)(107.9,54.27)(109.02,55.09)
\path(109.02,55.09)(110.09,55.9)(111.09,56.68)(112.02,57.45)(112.91,58.2)(113.72,58.93)(114.47,59.65)(115.16,60.34)(115.8,61.02)(116.4,61.68)
\path(116.4,61.68)(116.93,62.31)(117.41,62.93)(117.84,63.54)(118.23,64.12)(118.58,64.69)(118.87,65.23)(119.12,65.76)(119.33,66.27)(119.51,66.76)
\path(119.51,66.76)(119.63,67.23)(119.73,67.69)(119.8,68.12)(119.83,68.52)(119.81,68.93)(119.79,69.3)(119.73,69.66)(119.63,70.0)(119.51,70.3)
\path(119.51,70.3)(119.37,70.61)(119.2,70.88)(119.02,71.15)(118.81,71.38)(118.59,71.61)(118.36,71.8)(118.11,71.98)(117.84,72.16)(117.56,72.3)
\path(117.56,72.3)(117.27,72.41)(116.98,72.52)(116.68,72.61)(116.37,72.68)(116.05,72.72)(115.75,72.75)(115.43,72.75)(115.12,72.73)(114.8,72.69)
\path(114.8,72.69)(114.5,72.65)(114.19,72.56)(113.9,72.48)(113.62,72.36)(113.33,72.23)(113.08,72.06)(112.83,71.9)(112.58,71.69)(112.37,71.48)
\path(112.37,71.48)(112.16,71.25)(112.0,71.0)(111.83,70.72)(111.69,70.43)(111.58,70.11)(111.51,69.76)(111.44,69.41)(111.41,69.04)(111.41,68.63)
\path(111.41,68.63)(111.45,68.22)(111.52,67.77)(111.62,67.31)(111.77,66.83)(111.95,66.33)(112.18,65.8)(112.44,65.26)(112.75,64.7)(113.09,64.12)
\path(113.09,64.12)(113.5,63.52)(113.94,62.88)(114.44,62.24)(114.98,61.56)(115.58,60.88)(116.23,60.16)(116.94,59.43)(117.72,58.68)(118.55,57.91)
\path(118.55,57.91)(119.43,57.11)(120.37,56.29)(121.4,55.47)(122.48,54.61)(123.62,53.72)(124.83,52.83)(126.12,51.9)(127.48,50.95)(128.91,49.99)
\path(128.91,49.99)(130.41,49.0)(131.99,48.0)(132.0,48.0)
\path(132.0,48.0)(132.0,48.0)(133.58,46.93)(135.08,45.88)(136.5,44.86)(137.86,43.86)(139.14,42.88)(140.36,41.91)(141.5,40.99)(142.58,40.08)
\path(142.58,40.08)(143.61,39.18)(144.55,38.33)(145.44,37.49)(146.27,36.66)(147.03,35.86)(147.74,35.09)(148.39,34.34)(149.0,33.61)(149.53,32.9)
\path(149.53,32.9)(150.03,32.2)(150.47,31.54)(150.88,30.9)(151.22,30.29)(151.52,29.69)(151.78,29.11)(152.0,28.55)(152.17,28.03)(152.32,27.52)
\path(152.32,27.52)(152.41,27.03)(152.49,26.55)(152.52,26.12)(152.5,25.69)(152.47,25.29)(152.41,24.9)(152.33,24.54)(152.21,24.21)(152.07,23.89)
\path(152.07,23.89)(151.89,23.61)(151.71,23.34)(151.5,23.09)(151.27,22.86)(151.02,22.64)(150.77,22.45)(150.5,22.29)(150.22,22.15)(149.91,22.03)
\path(149.91,22.03)(149.61,21.93)(149.3,21.85)(148.97,21.79)(148.64,21.76)(148.32,21.73)(148.0,21.75)(147.66,21.77)(147.33,21.81)(147.0,21.87)
\path(147.0,21.87)(146.69,21.96)(146.38,22.07)(146.07,22.2)(145.77,22.36)(145.49,22.53)(145.21,22.72)(144.96,22.94)(144.71,23.17)(144.47,23.43)
\path(144.47,23.43)(144.27,23.7)(144.08,24.0)(143.91,24.3)(143.77,24.64)(143.66,25.01)(143.57,25.37)(143.5,25.78)(143.47,26.2)(143.47,26.62)
\path(143.47,26.62)(143.5,27.09)(143.57,27.56)(143.66,28.06)(143.8,28.59)(143.97,29.12)(144.19,29.69)(144.46,30.27)(144.77,30.86)(145.11,31.47)
\path(145.11,31.47)(145.5,32.11)(145.94,32.77)(146.44,33.45)(146.99,34.15)(147.58,34.88)(148.25,35.61)(148.94,36.36)(149.72,37.15)(150.55,37.93)
\path(150.55,37.93)(151.42,38.75)(152.38,39.59)(153.39,40.45)(154.47,41.33)(155.61,42.22)(156.83,43.13)(158.11,44.06)(159.47,45.02)(160.91,45.99)
\path(160.91,45.99)(162.41,46.97)(163.99,47.99)(164.0,48.0)
\end{picture}

\end{center}
\caption{elastica}\label{f8}
\end{figure}

\noindent When $k'\rightarrow 0$ the limit behavior is
represented on Fig. \ref{f8} (ii), see also Fig. \ref{f9}
(behaviour on the separatrix).
 

\setlength{\unitlength}{0.3mm}
\begin{figure}[h]
\begin{center} 

\begin{picture}(180,100)
\thinlines
\drawvector{6.0}{2.0}{92.0}{0}{1}
\drawvector{2.0}{8.0}{172.0}{1}{0}
\drawcenteredtext{10.0}{90.0}{$y$}
\drawcenteredtext{170.0}{4.0}{$x$}
\path(6.0,8.0)(6.0,8.0)(7.51,8.59)(9.01,9.18)(10.5,9.77)(11.97,10.36)(13.43,10.94)(14.88,11.52)(16.3,12.1)(17.72,12.68)
\path(17.72,12.68)(19.13,13.25)(20.53,13.81)(21.92,14.38)(23.28,14.93)(24.63,15.48)(25.97,16.04)(27.3,16.59)(28.62,17.12)(29.93,17.67)
\path(29.93,17.67)(31.21,18.2)(32.49,18.75)(33.75,19.28)(35.0,19.79)(36.24,20.31)(37.45,20.84)(38.66,21.36)(39.86,21.87)(41.04,22.37)
\path(41.04,22.37)(42.22,22.87)(43.38,23.37)(44.52,23.87)(45.65,24.37)(46.77,24.87)(47.88,25.35)(48.97,25.82)(50.04,26.3)(51.11,26.79)
\path(51.11,26.79)(52.15,27.26)(53.2,27.72)(54.22,28.2)(55.24,28.65)(56.22,29.11)(57.22,29.56)(58.18,30.02)(59.15,30.46)(60.09,30.9)
\path(60.09,30.9)(61.02,31.35)(61.95,31.79)(62.86,32.22)(63.75,32.65)(64.62,33.06)(65.48,33.49)(66.34,33.9)(67.19,34.33)(68.01,34.74)
\path(68.01,34.74)(68.83,35.15)(69.62,35.54)(70.41,35.95)(71.19,36.34)(71.94,36.74)(72.69,37.13)(73.43,37.5)(74.16,37.9)(74.86,38.27)
\path(74.86,38.27)(75.55,38.65)(76.23,39.02)(76.91,39.38)(77.56,39.75)(78.2,40.11)(78.83,40.47)(79.44,40.83)(80.05,41.16)(80.63,41.52)
\path(80.63,41.52)(81.22,41.86)(81.77,42.2)(82.33,42.54)(82.87,42.86)(83.38,43.2)(83.9,43.52)(84.4,43.84)(84.87,44.15)(85.34,44.47)
\path(85.34,44.47)(85.8,44.79)(86.26,45.09)(86.69,45.38)(87.11,45.68)(87.51,45.99)(87.9,46.27)(88.27,46.56)(88.63,46.86)(89.0,47.13)
\path(89.0,47.13)(89.33,47.4)(89.66,47.68)(89.97,47.95)(90.26,48.22)(90.55,48.49)(90.83,48.75)(91.08,49.0)(91.33,49.25)(91.56,49.5)
\path(91.56,49.5)(91.79,49.75)(91.98,49.99)(92.0,50.0)
\path(92.0,50.0)(92.0,50.0)(92.26,50.43)(92.54,50.86)(92.8,51.29)(93.05,51.7)(93.3,52.13)(93.55,52.54)(93.77,52.95)(94.0,53.34)
\path(94.0,53.34)(94.22,53.74)(94.44,54.13)(94.63,54.52)(94.83,54.9)(95.02,55.27)(95.2,55.65)(95.37,56.0)(95.55,56.36)(95.7,56.72)
\path(95.7,56.72)(95.87,57.06)(96.01,57.41)(96.16,57.75)(96.29,58.09)(96.41,58.41)(96.52,58.74)(96.63,59.06)(96.75,59.36)(96.83,59.68)
\path(96.83,59.68)(96.93,59.97)(97.01,60.27)(97.08,60.56)(97.15,60.86)(97.22,61.13)(97.26,61.4)(97.3,61.68)(97.34,61.95)(97.37,62.2)
\path(97.37,62.2)(97.41,62.47)(97.43,62.72)(97.44,62.95)(97.44,63.2)(97.44,63.43)(97.41,63.65)(97.4,63.88)(97.37,64.11)(97.34,64.31)
\path(97.34,64.31)(97.3,64.52)(97.26,64.73)(97.19,64.93)(97.13,65.12)(97.06,65.3)(97.0,65.5)(96.91,65.66)(96.81,65.83)(96.72,66.01)
\path(96.72,66.01)(96.62,66.16)(96.51,66.33)(96.38,66.48)(96.26,66.62)(96.12,66.76)(95.98,66.9)(95.83,67.02)(95.68,67.16)(95.51,67.27)
\path(95.51,67.27)(95.34,67.4)(95.16,67.51)(94.98,67.61)(94.79,67.7)(94.58,67.8)(94.38,67.88)(94.18,67.98)(93.95,68.05)(93.73,68.12)
\path(93.73,68.12)(93.48,68.19)(93.25,68.26)(93.0,68.31)(92.75,68.37)(92.48,68.41)(92.2,68.45)(91.93,68.5)(91.65,68.52)(91.36,68.55)
\path(91.36,68.55)(91.05,68.58)(90.75,68.58)(90.43,68.59)(90.11,68.61)(89.79,68.61)(89.44,68.61)(89.11,68.59)(88.76,68.58)(88.4,68.55)
\path(88.4,68.55)(88.04,68.54)(87.66,68.5)(87.27,68.47)(86.9,68.43)(86.51,68.37)(86.11,68.33)(85.69,68.26)(85.27,68.2)(84.86,68.13)
\path(84.86,68.13)(84.43,68.06)(84.0,68.0)(84.0,68.0)
\path(84.0,68.0)(84.0,68.0)(83.68,67.87)(83.37,67.75)(83.06,67.62)(82.79,67.48)(82.5,67.36)(82.23,67.22)(81.97,67.08)(81.72,66.93)
\path(81.72,66.93)(81.47,66.79)(81.23,66.62)(81.01,66.48)(80.79,66.31)(80.58,66.16)(80.37,66.0)(80.18,65.83)(80.0,65.66)(79.83,65.48)
\path(79.83,65.48)(79.66,65.31)(79.5,65.13)(79.36,64.94)(79.22,64.76)(79.08,64.58)(78.95,64.38)(78.84,64.19)(78.75,64.0)(78.65,63.79)
\path(78.65,63.79)(78.55,63.59)(78.48,63.38)(78.41,63.16)(78.36,62.95)(78.3,62.74)(78.26,62.52)(78.23,62.29)(78.19,62.06)(78.19,61.84)
\path(78.19,61.84)(78.18,61.59)(78.18,61.36)(78.19,61.11)(78.2,60.88)(78.23,60.63)(78.26,60.38)(78.31,60.13)(78.37,59.88)(78.43,59.61)
\path(78.43,59.61)(78.51,59.36)(78.58,59.09)(78.66,58.81)(78.76,58.54)(78.87,58.27)(78.98,58.0)(79.12,57.7)(79.25,57.43)(79.38,57.13)
\path(79.38,57.13)(79.55,56.84)(79.69,56.56)(79.87,56.25)(80.05,55.95)(80.23,55.65)(80.43,55.34)(80.62,55.04)(80.84,54.72)(81.06,54.4)
\path(81.06,54.4)(81.3,54.08)(81.54,53.75)(81.77,53.43)(82.04,53.11)(82.3,52.77)(82.58,52.43)(82.86,52.09)(83.15,51.75)(83.45,51.4)
\path(83.45,51.4)(83.76,51.06)(84.08,50.7)(84.41,50.34)(84.73,50.0)(85.08,49.63)(85.44,49.27)(85.8,48.9)(86.18,48.52)(86.55,48.15)
\path(86.55,48.15)(86.94,47.77)(87.33,47.4)(87.75,47.0)(88.16,46.63)(88.58,46.24)(89.01,45.84)(89.45,45.43)(89.91,45.04)(90.37,44.63)
\path(90.37,44.63)(90.83,44.24)(91.3,43.83)(91.8,43.4)(92.29,43.0)(92.79,42.58)(93.3,42.15)(93.83,41.72)(94.34,41.29)(94.88,40.86)
\path(94.88,40.86)(95.44,40.43)(95.98,40.0)(96.0,40.0)
\path(96.0,40.0)(96.0,40.0)(96.27,39.75)(96.56,39.5)(96.87,39.27)(97.19,39.04)(97.51,38.79)(97.83,38.54)(98.18,38.31)(98.52,38.08)
\path(98.52,38.08)(98.88,37.84)(99.26,37.59)(99.62,37.36)(100.01,37.11)(100.41,36.86)(100.81,36.63)(101.23,36.4)(101.65,36.15)(102.08,35.9)
\path(102.08,35.9)(102.52,35.68)(102.98,35.43)(103.44,35.18)(103.9,34.95)(104.37,34.72)(104.87,34.47)(105.36,34.22)(105.87,34.0)(106.37,33.75)
\path(106.37,33.75)(106.91,33.52)(107.44,33.27)(107.98,33.04)(108.54,32.79)(109.09,32.56)(109.66,32.31)(110.23,32.08)(110.83,31.84)(111.43,31.6)
\path(111.43,31.6)(112.04,31.36)(112.65,31.12)(113.27,30.87)(113.91,30.63)(114.55,30.39)(115.2,30.15)(115.87,29.92)(116.54,29.68)(117.22,29.44)
\path(117.22,29.44)(117.91,29.2)(118.61,28.95)(119.31,28.71)(120.02,28.47)(120.76,28.23)(121.48,28.0)(122.23,27.76)(122.98,27.52)(123.76,27.28)
\path(123.76,27.28)(124.52,27.04)(125.3,26.79)(126.09,26.55)(126.9,26.31)(127.7,26.07)(128.52,25.84)(129.35,25.6)(130.19,25.36)(131.03,25.12)
\path(131.03,25.12)(131.88,24.87)(132.75,24.63)(133.63,24.39)(134.5,24.15)(135.39,23.92)(136.3,23.68)(137.22,23.44)(138.13,23.2)(139.05,22.95)
\path(139.05,22.95)(140.0,22.71)(140.94,22.47)(141.89,22.23)(142.86,22.0)(143.83,21.76)(144.83,21.52)(145.82,21.28)(146.82,21.04)(147.83,20.79)
\path(147.83,20.79)(148.86,20.55)(149.88,20.31)(150.91,20.07)(151.97,19.84)(153.02,19.6)(154.1,19.36)(155.16,19.12)(156.25,18.87)(157.35,18.63)
\path(157.35,18.63)(158.44,18.39)(159.57,18.15)(160.69,17.92)(161.82,17.68)(162.96,17.44)(164.11,17.2)(165.27,16.95)(166.44,16.71)(167.61,16.47)
\path(167.61,16.47)(168.8,16.23)(169.99,16.0)(170.0,16.0)
\end{picture}

\end{center}
\caption{behaviour on the separatrix}\label{f9}
\end{figure}
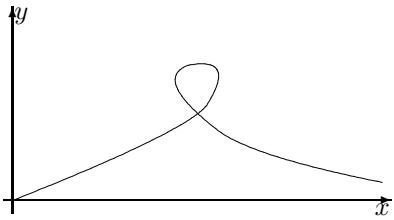

In this representation $\theta $ is up to a constant the angle
of the normal with respect to a fixed direction. The rotating
trajectories of the pendulum corres\-pond to geodesics not starting
from 0. They project on the space $(x,y)$ onto \it{non
inflexional elastica}, see Fig. \ref{f10} (ii).


\setlength{\unitlength}{0.4mm}
\begin{figure}[h]
\begin{center} 

\begin{picture}(180,100)
\thinlines
\drawvector{20.0}{48.0}{60.0}{1}{0}
\drawvector{100.0}{48.0}{60.0}{1}{0}
\drawvector{22.0}{24.0}{54.0}{0}{1}
\drawvector{102.0}{24.0}{54.0}{0}{1}
\drawcenteredtext{160.0}{52.0}{$t$}
\drawcenteredtext{80.0}{52.0}{$t$}
\drawcenteredtext{106.0}{78.0}{$y$}
\drawcenteredtext{26.0}{78.0}{$y$}
\path(22.0,56.0)(22.0,56.0)(22.2,56.75)(22.39,57.47)(22.59,58.2)(22.79,58.9)(23.0,59.61)(23.2,60.27)(23.39,60.93)(23.59,61.59)
\path(23.59,61.59)(23.79,62.22)(24.0,62.84)(24.2,63.43)(24.39,64.01)(24.59,64.58)(24.79,65.15)(25.0,65.69)(25.2,66.2)(25.39,66.72)
\path(25.39,66.72)(25.59,67.2)(25.79,67.69)(26.0,68.16)(26.2,68.59)(26.39,69.04)(26.6,69.44)(26.79,69.86)(27.0,70.25)(27.2,70.62)
\path(27.2,70.62)(27.39,70.97)(27.59,71.31)(27.79,71.63)(28.0,71.94)(28.19,72.25)(28.38,72.52)(28.59,72.8)(28.79,73.05)(29.0,73.29)
\path(29.0,73.29)(29.19,73.51)(29.38,73.7)(29.59,73.9)(29.79,74.08)(29.98,74.23)(30.19,74.37)(30.38,74.51)(30.59,74.62)(30.79,74.72)
\path(30.79,74.72)(30.98,74.8)(31.19,74.87)(31.38,74.93)(31.59,74.95)(31.79,74.98)(31.98,75.0)(32.18,74.98)(32.38,74.95)(32.59,74.93)
\path(32.59,74.93)(32.79,74.87)(32.99,74.8)(33.18,74.72)(33.38,74.62)(33.59,74.51)(33.79,74.37)(33.99,74.23)(34.18,74.08)(34.38,73.9)
\path(34.38,73.9)(34.59,73.7)(34.79,73.51)(34.99,73.29)(35.18,73.05)(35.38,72.8)(35.59,72.52)(35.79,72.25)(35.99,71.95)(36.18,71.63)
\path(36.18,71.63)(36.38,71.31)(36.59,70.97)(36.79,70.62)(36.99,70.25)(37.18,69.86)(37.38,69.44)(37.59,69.04)(37.79,68.59)(37.99,68.16)
\path(37.99,68.16)(38.18,67.69)(38.38,67.2)(38.59,66.72)(38.79,66.2)(38.99,65.69)(39.18,65.15)(39.38,64.58)(39.59,64.01)(39.79,63.43)
\path(39.79,63.43)(39.99,62.84)(40.18,62.22)(40.38,61.59)(40.59,60.93)(40.79,60.27)(40.99,59.61)(41.18,58.9)(41.38,58.2)(41.59,57.47)
\path(41.59,57.47)(41.79,56.75)(41.99,56.0)(42.0,56.0)
\path(42.0,56.0)(42.0,56.0)(42.18,55.2)(42.4,54.43)(42.59,53.66)(42.79,52.91)(43.0,52.2)(43.18,51.47)(43.4,50.79)(43.59,50.11)
\path(43.59,50.11)(43.79,49.43)(44.0,48.79)(44.18,48.15)(44.4,47.54)(44.59,46.95)(44.79,46.36)(45.0,45.79)(45.18,45.24)(45.4,44.7)
\path(45.4,44.7)(45.59,44.18)(45.79,43.68)(46.0,43.2)(46.18,42.72)(46.4,42.27)(46.59,41.83)(46.79,41.4)(47.0,41.0)(47.18,40.59)
\path(47.18,40.59)(47.4,40.22)(47.59,39.86)(47.79,39.52)(48.0,39.2)(48.18,38.88)(48.4,38.59)(48.59,38.31)(48.79,38.04)(49.0,37.79)
\path(49.0,37.79)(49.18,37.56)(49.4,37.34)(49.59,37.15)(49.79,36.95)(50.0,36.79)(50.18,36.63)(50.38,36.5)(50.59,36.38)(50.79,36.27)
\path(50.79,36.27)(51.0,36.2)(51.18,36.11)(51.38,36.06)(51.59,36.02)(51.79,36.0)(51.99,36.0)(52.18,36.0)(52.38,36.02)(52.59,36.06)
\path(52.59,36.06)(52.79,36.11)(52.99,36.18)(53.18,36.27)(53.38,36.38)(53.59,36.5)(53.79,36.63)(53.99,36.79)(54.18,36.95)(54.38,37.15)
\path(54.38,37.15)(54.59,37.34)(54.79,37.56)(54.99,37.79)(55.18,38.04)(55.38,38.31)(55.59,38.59)(55.79,38.88)(55.99,39.18)(56.18,39.52)
\path(56.18,39.52)(56.38,39.86)(56.59,40.22)(56.79,40.59)(56.99,40.99)(57.18,41.4)(57.38,41.83)(57.59,42.27)(57.79,42.72)(57.99,43.18)
\path(57.99,43.18)(58.18,43.68)(58.38,44.18)(58.59,44.7)(58.79,45.24)(58.99,45.79)(59.18,46.36)(59.38,46.95)(59.59,47.54)(59.79,48.15)
\path(59.79,48.15)(59.99,48.79)(60.18,49.43)(60.38,50.11)(60.59,50.79)(60.79,51.47)(60.99,52.18)(61.18,52.91)(61.38,53.66)(61.59,54.43)
\path(61.59,54.43)(61.79,55.2)(61.99,55.99)(62.0,56.0)
\path(62.0,56.0)(62.0,56.0)(62.15,56.66)(62.31,57.33)(62.47,57.97)(62.63,58.61)(62.79,59.25)(62.95,59.86)(63.11,60.45)(63.27,61.04)
\path(63.27,61.04)(63.43,61.63)(63.59,62.18)(63.75,62.75)(63.9,63.29)(64.08,63.81)(64.23,64.33)(64.4,64.83)(64.55,65.33)(64.72,65.81)
\path(64.72,65.81)(64.87,66.29)(65.04,66.75)(65.19,67.19)(65.36,67.62)(65.51,68.05)(65.68,68.45)(65.83,68.86)(66.0,69.25)(66.16,69.62)
\path(66.16,69.62)(66.31,69.98)(66.47,70.33)(66.62,70.66)(66.8,71.0)(66.95,71.3)(67.11,71.61)(67.26,71.9)(67.44,72.18)(67.59,72.44)
\path(67.59,72.44)(67.75,72.69)(67.91,72.94)(68.08,73.16)(68.23,73.38)(68.38,73.58)(68.55,73.79)(68.72,73.97)(68.87,74.13)(69.02,74.3)
\path(69.02,74.3)(69.19,74.44)(69.36,74.58)(69.51,74.69)(69.66,74.8)(69.83,74.91)(70.0,75.0)(70.16,75.06)(70.3,75.12)(70.47,75.18)
\path(70.47,75.18)(70.62,75.22)(70.8,75.25)(70.94,75.26)(71.11,75.26)(71.26,75.25)(71.44,75.23)(71.58,75.19)(71.75,75.15)(71.91,75.08)
\path(71.91,75.08)(72.08,75.01)(72.23,74.94)(72.38,74.84)(72.55,74.73)(72.72,74.62)(72.87,74.48)(73.02,74.34)(73.19,74.19)(73.36,74.02)
\path(73.36,74.02)(73.51,73.84)(73.66,73.66)(73.83,73.45)(74.0,73.25)(74.15,73.01)(74.3,72.77)(74.47,72.52)(74.62,72.26)(74.79,72.0)
\path(74.79,72.0)(74.94,71.7)(75.11,71.41)(75.26,71.09)(75.43,70.77)(75.58,70.44)(75.75,70.09)(75.91,69.73)(76.06,69.37)(76.23,68.98)
\path(76.23,68.98)(76.38,68.59)(76.55,68.19)(76.7,67.76)(76.87,67.33)(77.02,66.9)(77.19,66.44)(77.34,65.98)(77.51,65.5)(77.66,65.01)
\path(77.66,65.01)(77.83,64.51)(77.98,64.0)(78.0,64.0)
\path(102.0,54.0)(102.0,54.0)(102.87,54.83)(103.7,55.63)(104.5,56.43)(105.26,57.22)(105.97,58.0)(106.65,58.75)(107.29,59.49)(107.88,60.2)
\path(107.88,60.2)(108.44,60.9)(108.98,61.59)(109.48,62.27)(109.94,62.93)(110.37,63.56)(110.76,64.19)(111.13,64.8)(111.48,65.38)(111.79,65.95)
\path(111.79,65.95)(112.06,66.51)(112.33,67.05)(112.55,67.58)(112.76,68.08)(112.94,68.56)(113.08,69.04)(113.23,69.5)(113.33,69.93)(113.43,70.34)
\path(113.43,70.34)(113.5,70.75)(113.55,71.13)(113.58,71.5)(113.61,71.84)(113.61,72.18)(113.58,72.48)(113.55,72.79)(113.51,73.06)(113.47,73.33)
\path(113.47,73.33)(113.4,73.56)(113.31,73.79)(113.23,73.98)(113.13,74.18)(113.04,74.34)(112.91,74.5)(112.8,74.62)(112.68,74.73)(112.55,74.83)
\path(112.55,74.83)(112.41,74.91)(112.29,74.95)(112.15,75.0)(112.01,75.01)(111.87,75.01)(111.75,75.0)(111.61,74.95)(111.48,74.9)(111.34,74.81)
\path(111.34,74.81)(111.23,74.72)(111.11,74.59)(111.0,74.47)(110.9,74.3)(110.8,74.12)(110.7,73.94)(110.63,73.72)(110.56,73.48)(110.51,73.23)
\path(110.51,73.23)(110.47,72.95)(110.44,72.66)(110.43,72.33)(110.43,72.01)(110.44,71.65)(110.48,71.26)(110.52,70.87)(110.61,70.45)(110.69,70.01)
\path(110.69,70.01)(110.8,69.55)(110.94,69.08)(111.09,68.58)(111.27,68.05)(111.48,67.51)(111.7,66.94)(111.95,66.37)(112.25,65.76)(112.55,65.12)
\path(112.55,65.12)(112.88,64.48)(113.26,63.81)(113.66,63.11)(114.08,62.4)(114.55,61.66)(115.05,60.9)(115.58,60.13)(116.15,59.33)(116.75,58.5)
\path(116.75,58.5)(117.37,57.66)(118.05,56.79)(118.77,55.9)(119.52,55.0)(120.31,54.06)(121.16,53.11)(122.04,52.13)(122.95,51.13)(123.93,50.11)
\path(123.93,50.11)(124.94,49.06)(125.98,48.0)(126.0,48.0)
\path(138.0,48.0)(138.0,48.0)(139.05,49.0)(140.05,50.0)(141.02,50.97)(141.96,51.91)(142.83,52.86)(143.67,53.77)(144.49,54.66)(145.25,55.54)
\path(145.25,55.54)(145.97,56.4)(146.66,57.24)(147.3,58.06)(147.91,58.86)(148.49,59.63)(149.02,60.4)(149.53,61.13)(150.0,61.86)(150.44,62.56)
\path(150.44,62.56)(150.86,63.25)(151.24,63.91)(151.6,64.56)(151.91,65.19)(152.21,65.8)(152.47,66.4)(152.72,66.97)(152.92,67.52)(153.11,68.05)
\path(153.11,68.05)(153.28,68.58)(153.42,69.06)(153.55,69.55)(153.64,70.01)(153.72,70.44)(153.8,70.87)(153.83,71.26)(153.86,71.66)(153.86,72.01)
\path(153.86,72.01)(153.86,72.37)(153.83,72.69)(153.78,73.0)(153.74,73.29)(153.67,73.55)(153.6,73.81)(153.5,74.05)(153.41,74.26)(153.3,74.45)
\path(153.3,74.45)(153.17,74.63)(153.05,74.8)(152.92,74.94)(152.78,75.05)(152.63,75.16)(152.5,75.25)(152.33,75.3)(152.19,75.36)(152.02,75.38)
\path(152.02,75.38)(151.88,75.4)(151.72,75.38)(151.55,75.36)(151.39,75.3)(151.25,75.25)(151.1,75.16)(150.96,75.06)(150.8,74.94)(150.67,74.8)
\path(150.67,74.8)(150.55,74.65)(150.42,74.48)(150.3,74.27)(150.21,74.06)(150.11,73.83)(150.02,73.58)(149.96,73.33)(149.89,73.04)(149.85,72.73)
\path(149.85,72.73)(149.8,72.41)(149.8,72.06)(149.78,71.72)(149.8,71.33)(149.83,70.94)(149.88,70.52)(149.96,70.09)(150.05,69.65)(150.14,69.18)
\path(150.14,69.18)(150.28,68.69)(150.44,68.19)(150.61,67.66)(150.82,67.12)(151.03,66.55)(151.28,65.98)(151.57,65.38)(151.88,64.77)(152.21,64.13)
\path(152.21,64.13)(152.57,63.49)(152.97,62.81)(153.38,62.13)(153.85,61.43)(154.33,60.7)(154.86,59.97)(155.41,59.2)(156.0,58.43)(156.63,57.63)
\path(156.63,57.63)(157.28,56.83)(157.99,56.0)(158.0,56.0)
\path(126.0,48.0)(126.0,48.0)(126.11,47.88)(126.23,47.75)(126.36,47.65)(126.47,47.52)(126.59,47.43)(126.72,47.31)(126.83,47.2)(126.94,47.11)
\path(126.94,47.11)(127.08,47.0)(127.19,46.91)(127.31,46.81)(127.44,46.72)(127.55,46.63)(127.66,46.54)(127.8,46.45)(127.91,46.38)(128.02,46.29)
\path(128.02,46.29)(128.16,46.22)(128.27,46.15)(128.38,46.06)(128.52,46.0)(128.63,45.93)(128.75,45.86)(128.88,45.81)(129.0,45.75)(129.11,45.68)
\path(129.11,45.68)(129.24,45.63)(129.36,45.58)(129.47,45.52)(129.6,45.47)(129.72,45.43)(129.83,45.38)(129.96,45.34)(130.07,45.29)(130.19,45.27)
\path(130.19,45.27)(130.32,45.22)(130.44,45.2)(130.55,45.16)(130.66,45.13)(130.78,45.11)(130.91,45.09)(131.03,45.06)(131.16,45.04)(131.27,45.04)
\path(131.27,45.04)(131.38,45.02)(131.5,45.0)(131.63,45.0)(131.75,45.0)(131.88,45.0)(132.0,45.0)(132.11,45.0)(132.22,45.0)(132.35,45.0)
\path(132.35,45.0)(132.47,45.0)(132.6,45.02)(132.72,45.04)(132.83,45.04)(132.94,45.06)(133.07,45.09)(133.19,45.11)(133.32,45.13)(133.44,45.16)
\path(133.44,45.16)(133.55,45.2)(133.66,45.22)(133.78,45.27)(133.91,45.29)(134.03,45.34)(134.16,45.38)(134.27,45.43)(134.38,45.47)(134.5,45.52)
\path(134.5,45.52)(134.63,45.58)(134.75,45.63)(134.88,45.68)(135.0,45.74)(135.11,45.81)(135.22,45.86)(135.35,45.93)(135.47,46.0)(135.6,46.06)
\path(135.6,46.06)(135.72,46.15)(135.83,46.22)(135.94,46.29)(136.07,46.38)(136.19,46.45)(136.32,46.54)(136.44,46.63)(136.55,46.72)(136.66,46.81)
\path(136.66,46.81)(136.78,46.9)(136.91,47.0)(137.02,47.11)(137.16,47.2)(137.27,47.31)(137.38,47.41)(137.5,47.52)(137.63,47.65)(137.75,47.75)
\path(137.75,47.75)(137.88,47.88)(138.0,47.99)(138.0,48.0)
\drawcenteredtext{48.0}{18.0}{(i)}
\drawcenteredtext{130.0}{18.0}{(ii)}
\end{picture}

\end{center}
\caption{}\label{f10}
\end{figure}


\subsubsection{Trace of $S(0,r)$ and $W(0,r)$ in $y=0$ in
the flat case}
The successive intersection times with $y=0$ are given by~:
$t_i ={\displaystyle {2K\over \sqrt{\lambda }}}$, \mbox{$
i=1,\ldots 
,N.$} If we fix the length to $t_i=r,$ we get the following curves~:
\begin{equation*}
\begin{split}
x & =  -r+{\displaystyle {2\over \sqrt{\lambda }}}\ 
(E(K+i2K)-E ) \\
z & =  {\displaystyle {2\over 3\lambda ^{3/2}}}\ \Bigl [ 
(2k^2-1)\ (E (K+i2K)-E )+2Ki{k'}^2\Bigr ]\nonumber
\end{split}
\end{equation*}

It represents a parametric curve, where the parameter is
$k\in ]0,1[.$ Using the relation~: $E(K+i2K)=(2i+1)E$ we obtain
for each $i$ the following curves~: $k\longmapsto C_i(k)=
(x_i(k), z_i(k)),$
\begin{equation*}
\begin{split}
x_i (k) & =  -r +2 r \ {\displaystyle {E\over K}} \\
z_i (k) & =  {\displaystyle {r^3\over 6i ^2K^3}}\ \Bigl [
 (2k ^2 -1)E+{k'}^2 K\Bigr ]
\end{split}
\end{equation*}

\noindent where $k\in ]0,1[.$ We can easily draw those curves
using the standard package about elliptic functions in
Mathematica, see Fig. \ref{f11}.


\unitlength=.25mm
\makeatletter
\def\shade{\@ifnextchar[{\shade@special}{\@killglue\special{sh}\ignorespaces}}
\def\shade@special[#1]{\@killglue\special{sh #1}\ignorespaces}
\makeatother

\begin{figure}[h]
\begin{center}
\begin{picture}(215,115)(232,-5)
\thinlines
\typeout{\space\space\space eepic-ture exported by 'qfig'.}
\font\FonttenBI=cmbxti10\relax
\font\FonttwlBI=cmbxti10 scaled \magstep1\relax
\path (232,18)(400,18)
\path (315,6)(315,99)
\path (392.483,20.736)(400,18)(392.483,15.264)
\path (312.264,91.482)(315,99)(317.736,91.482)
\put(405,19){{\rm\rm x}}
\put(295,93){{\rm\rm z}}
\path (245,18)(246.5,18.155)(248,18.32)(249.5,18.495)(251,18.68)
(252.5,18.875)(254,19.08)(255.5,19.295)(257,19.52)(258.5,19.755)
(260,20)(260.87,20.075)(261.88,20.2)(263.03,20.375)(264.32,20.6)
(265.75,20.875)(267.32,21.2)(269.03,21.575)(270.88,22)(272.87,22.475)
(275,23)(278.395,23.71)(281.68,24.44)(284.855,25.19)(287.92,25.96)
(290.875,26.75)(293.72,27.56)(296.455,28.39)(299.08,29.24)(301.595,30.11)
(304,31)(305.89,31.775)(307.76,32.6)(309.61,33.475)(311.44,34.4)
(313.25,35.375)(315.04,36.4)(316.81,37.475)(318.56,38.6)(320.29,39.775)
(322,41)(323.51,42.905)(325.04,44.72)(326.59,46.445)(328.16,48.08)
(329.75,49.625)(331.36,51.08)(332.99,52.445)(334.64,53.72)(336.31,54.905)
(338,56)(340.025,57.095)(342,58.08)(343.925,58.955)(345.8,59.72)
(347.625,60.375)(349.4,60.92)(351.125,61.355)(352.8,61.68)(354.425,61.895)
(356,62)(357.525,61.77)(359,61.48)(360.425,61.13)(361.8,60.72)
(363.125,60.25)(364.4,59.72)(365.625,59.13)(366.8,58.48)(367.925,57.77)
(369,57)(369.8,56.035)(370.6,55.04)(371.4,54.015)(372.2,52.96)
(373,51.875)(373.8,50.76)(374.6,49.615)(375.4,48.44)(376.2,47.235)
(377,46)(378.025,44.6)(379,43.2)(379.925,41.8)(380.8,40.4)
(381.625,39)(382.4,37.6)(383.125,36.2)(383.8,34.8)(384.425,33.4)
(385,32)(385.525,30.6)(386,29.2)(386.425,27.8)(386.8,26.4)
(387.125,25)(387.4,23.6)(387.625,22.2)(387.8,20.8)(387.925,19.4)
(388,18)
\path (246,18)(247.595,18.43)(249.28,18.92)(251.055,19.47)(252.92,20.08)
(254.875,20.75)(256.92,21.48)(259.055,22.27)(261.28,23.12)(263.595,24.03)
(266,25)(268.99,25.94)(271.96,26.96)(274.91,28.06)(277.84,29.24)
(280.75,30.5)(283.64,31.84)(286.51,33.26)(289.36,34.76)(292.19,36.34)
(295,38)(297.88,40.1)(300.72,42.2)(303.52,44.3)(306.28,46.4)
(309,48.5)(311.68,50.6)(314.32,52.7)(316.92,54.8)(319.48,56.9)
(322,59)(324.435,61.73)(326.84,64.32)(329.215,66.77)(331.56,69.08)
(333.875,71.25)(336.16,73.28)(338.415,75.17)(340.64,76.92)(342.835,78.53)
(345,80)(347.18,81.645)(349.32,83.08)(351.42,84.305)(353.48,85.32)
(355.5,86.125)(357.48,86.72)(359.42,87.105)(361.32,87.28)(363.18,87.245)
(365,87)(367.05,86.14)(369,85.16)(370.85,84.06)(372.6,82.84)
(374.25,81.5)(375.8,80.04)(377.25,78.46)(378.6,76.76)(379.85,74.94)
(381,73)(381.78,70.535)(382.52,68.04)(383.22,65.515)(383.88,62.96)
(384.5,60.375)(385.08,57.76)(385.62,55.115)(386.12,52.44)(386.58,49.735)
(387,47)(387.38,44.235)(387.72,41.44)(388.02,38.615)(388.28,35.76)
(388.5,32.875)(388.68,29.96)(388.82,27.015)(388.92,24.04)(388.98,21.035)
(389,18)
\put(370,90){{\rm\rm {i=1}}}
\put(349,63){{\rm\rm {i=2}}}
\put(375,-1){{\rm\rm {$(r,0)$}}}
\put(233,0){{\rm\rm {$(-r,0)$}}}
\end{picture}
\end{center}
\caption{}\label{f11}
\end{figure}

The \it{exterior} curve obtained for $i=1$ represents the
intersection of the sphere $S(0,r)$ with the Martinet plane
in the domain $z>0.$ Each point of this curve is the end-point
of \it{two distinct minimizers} and by obvious geometric
reasoning we have~:

\begin{prop}
The cut locus $L(0,r)$ is $C_1 \cup -C_1.$
\end{prop}

Moreover by inspecting Fig. \ref{f11} we deduce the following~:

\begin{prop}
The abnormal geodesics are minimizers.
\end{prop}

This result is not new but here the proof is based on the
analysis of the geodesic flow. \it{The main property is that
at each intersection with $y=0,$ the variable $z$ has non zero
drift which can be easily evaluated using (\ref{4.20}). This
will lead to optimality results for the general metric, by
stability.}

This is an alternative proof to the optimality results
presented in Section 3 or in \cite{AS}, \cite{LS},
where we consider all
the trajectories of the system.

\begin{rem}
We observe that $(-r,0)$ is a ramified point of the trace of the
wave front on the Martinet plane with an \it{infinite number of
branches}. This gives us a precise geometric interpretation on
the structure of the geodesics of fixed length with respect to
the abnormal line. \it{Indeed for every neighborhood $U$ of
$(-r,0,0)$ and every $n\in \N$, there exists a geodesic of length
$r$ with end-point in $U$, with $n$ oscillations}.
\end{rem}

We represent on Fig. \ref{f12} the first and second return mapping,
the length being fixed to $r,$ and by restricting the domain
to $\lambda >0, \theta (0)\in [0,\pi ]$.


\unitlength=.25mm
\makeatletter
\def\shade{\@ifnextchar[{\shade@special}{\@killglue\special{sh}\ignorespaces}}
\def\shade@special[#1]{\@killglue\special{sh #1}\ignorespaces}
\makeatother

\begin{figure}[h]
\begin{center}
\begin{picture}(409,170)(112,-5)
\thinlines
\typeout{\space\space\space eepic-ture exported by 'qfig'.}
\font\FonttenBI=cmbxti10\relax
\font\FonttwlBI=cmbxti10 scaled \magstep1\relax
\path (340,41)(508,41)
\path (423,29)(423,122)
\path (500.482,43.736)(508,41)(500.482,38.264)
\path (420.264,114.482)(423,122)(425.736,114.482)
\put(513,42){{\rm\rm x}}
\put(403,116){{\rm\rm z}}
\path (353,41)(354.5,41.155)(356,41.32)(357.5,41.495)(359,41.68)
(360.5,41.875)(362,42.08)(363.5,42.295)(365,42.52)(366.5,42.755)
(368,43)(368.87,43.075)(369.88,43.2)(371.03,43.375)(372.32,43.6)
(373.75,43.875)(375.32,44.2)(377.03,44.575)(378.88,45)(380.87,45.475)
(383,46)(386.395,46.71)(389.68,47.44)(392.855,48.19)(395.92,48.96)
(398.875,49.75)(401.72,50.56)(404.455,51.39)(407.08,52.24)(409.595,53.11)
(412,54)(413.89,54.775)(415.76,55.6)(417.61,56.475)(419.44,57.4)
(421.25,58.375)(423.04,59.4)(424.81,60.475)(426.56,61.6)(428.29,62.775)
(430,64)(431.51,65.905)(433.04,67.72)(434.59,69.445)(436.16,71.08)
(437.75,72.625)(439.36,74.08)(440.99,75.445)(442.64,76.72)(444.31,77.905)
(446,79)(448.025,80.095)(450,81.08)(451.925,81.955)(453.8,82.72)
(455.625,83.375)(457.4,83.92)(459.125,84.355)(460.8,84.68)(462.425,84.895)
(464,85)(465.525,84.77)(467,84.48)(468.425,84.13)(469.8,83.72)
(471.125,83.25)(472.4,82.72)(473.625,82.13)(474.8,81.48)(475.925,80.77)
(477,80)(477.8,79.035)(478.6,78.04)(479.4,77.015)(480.2,75.96)
(481,74.875)(481.8,73.76)(482.6,72.615)(483.4,71.44)(484.2,70.235)
(485,69)(486.025,67.6)(487,66.2)(487.925,64.8)(488.8,63.4)
(489.625,62)(490.4,60.6)(491.125,59.2)(491.8,57.8)(492.425,56.4)
(493,55)(493.525,53.6)(494,52.2)(494.425,50.8)(494.8,49.4)
(495.125,48)(495.4,46.6)(495.625,45.2)(495.8,43.8)(495.925,42.4)
(496,41)
\path (354,41)(355.595,41.43)(357.28,41.92)(359.055,42.47)(360.92,43.08)
(362.875,43.75)(364.92,44.48)(367.055,45.27)(369.28,46.12)(371.595,47.03)
(374,48)(376.99,48.94)(379.96,49.96)(382.91,51.06)(385.84,52.24)
(388.75,53.5)(391.64,54.84)(394.51,56.26)(397.36,57.76)(400.19,59.34)
(403,61)(405.88,63.1)(408.72,65.2)(411.52,67.3)(414.28,69.4)
(417,71.5)(419.68,73.6)(422.32,75.7)(424.92,77.8)(427.48,79.9)
(430,82)(432.435,84.73)(434.84,87.32)(437.215,89.77)(439.56,92.08)
(441.875,94.25)(444.16,96.28)(446.415,98.17)(448.64,99.92)(450.835,101.53)
(453,103)(455.18,104.645)(457.32,106.08)(459.42,107.305)(461.48,108.32)
(463.5,109.125)(465.48,109.72)(467.42,110.105)(469.32,110.28)(471.18,110.245)
(473,110)(475.05,109.14)(477,108.16)(478.85,107.06)(480.6,105.84)
(482.25,104.5)(483.8,103.04)(485.25,101.46)(486.6,99.76)(487.85,97.94)
(489,96)(489.78,93.535)(490.52,91.04)(491.22,88.515)(491.88,85.96)
(492.5,83.375)(493.08,80.76)(493.62,78.115)(494.12,75.44)(494.58,72.735)
(495,70)(495.38,67.235)(495.72,64.44)(496.02,61.615)(496.28,58.76)
(496.5,55.875)(496.68,52.96)(496.82,50.015)(496.92,47.04)(496.98,44.035)
(497,41)
\put(483,22){{\rm\rm {$(r,0)$}}}
\put(341,23){{\rm\rm {$(-r,0)$}}}
\path (112,41)(265,41)
\path (250,132)(250,41)
\path (257.482,43.736)(265,41)(257.482,38.264)
\path (121,41)(121,132)
\path (187,34)(187,149)
\path (184.264,141.482)(187,149)(189.736,141.482)
\put(193,148){{\rm\rm {$\lambda$}}}
\put(263,25){{\rm\rm {$\cos \ \theta (0)$}}}
\put(239,20){{\rm\rm {$+1$}}}
\put(179,20){{\rm\rm {$0$}}}
\put(114,21){{\rm\rm {$-1$}}}
\path (250,56)(246.75,55.905)(243.4,55.92)(239.95,56.045)(236.4,56.28)
(232.75,56.625)(229,57.08)(225.15,57.645)(221.2,58.32)(217.15,59.105)
(213,60)(207.355,61.005)(201.92,62.12)(196.695,63.345)(191.68,64.68)
(186.875,66.125)(182.28,67.68)(177.895,69.345)(173.72,71.12)(169.755,73.005)
(166,75)(162.815,77.105)(159.76,79.32)(156.835,81.645)(154.04,84.08)
(151.375,86.625)(148.84,89.28)(146.435,92.045)(144.16,94.92)(142.015,97.905)
(140,101)(138.115,104.205)(136.36,107.52)(134.735,110.945)(133.24,114.48)
(131.875,118.125)(130.64,121.88)(129.535,125.745)(128.56,129.72)(127.715,133.805)
(127,138)
\path (250,77)(245.965,77.14)(241.96,77.36)(237.985,77.66)(234.04,78.04)
(230.125,78.5)(226.24,79.04)(222.385,79.66)(218.56,80.36)(214.765,81.14)
(211,82)(206.725,82.895)(202.6,83.88)(198.625,84.955)(194.8,86.12)
(191.125,87.375)(187.6,88.72)(184.225,90.155)(181,91.68)(177.925,93.295)
(175,95)(172.54,97.065)(170.16,99.16)(167.86,101.285)(165.64,103.44)
(163.5,105.625)(161.44,107.84)(159.46,110.085)(157.56,112.36)(155.74,114.665)
(154,117)(152.34,119.365)(150.76,121.76)(149.26,124.185)(147.84,126.64)
(146.5,129.125)(145.24,131.64)(144.06,134.185)(142.96,136.76)(141.94,139.365)
(141,142)
\path (282,98)(283.62,99.605)(285.28,101.12)(286.98,102.545)(288.72,103.88)
(290.5,105.125)(292.32,106.28)(294.18,107.345)(296.08,108.32)(298.02,109.205)
(300,110)(302.155,110.615)(304.32,111.16)(306.495,111.635)(308.68,112.04)
(310.875,112.375)(313.08,112.64)(315.295,112.835)(317.52,112.96)(319.755,113.015)
(322,113)(324.795,112.825)(327.48,112.6)(330.055,112.325)(332.52,112)
(334.875,111.625)(337.12,111.2)(339.255,110.725)(341.28,110.2)(343.195,109.625)
(345,109)(346.695,108.325)(348.28,107.6)(349.755,106.825)(351.12,106)
(352.375,105.125)(353.52,104.2)(354.555,103.225)(355.48,102.2)(356.295,101.125)
(357,100)
\path (313,116)(320,113)(313,110)
\path (164,100)(162,106)(169,106)
\path (151,82)(149,89)(156,87)
\put(191,90){{\rm\rm {Dom $R_2$}}}
\put(136,58){{\rm\rm {Dom $R_1$}}}
\put(291,121){{\rm\rm {$R_i, i=1,2$}}}
\path (407,70)(404,62)(413,62)
\path (408,57)(402,52)(408,49)
\put(448,115){{\rm\rm {Im $R_1=\widetilde{S}$}}}
\put(442,62){{\rm\rm {Im $R_2$}}}
\put(149,-1){{\rm\rm {Dom $R_i$}}}
\end{picture}
\end{center}
\caption{}\label{f12}
\end{figure}

\noindent In the phase space $(\theta ,{\dot \theta }),$ $R_1$
corresponds to the symmetry~: $(\theta ,0)\longmapsto 
(-\theta ,0)$ and $R_2$ corresponds to the identity~:
$(\theta ,0)\longmapsto (\theta ,0).$

We represent on Fig. \ref{f13} the two branches $C_1$ and
${\bar C}_1$ in ${\widetilde S}(0,r)$ ending at $(-r,0)$
and $(r,0)$ and corresponding respectively to the behaviors
of the geodesics near the center $0$ and the separatrix $\Sigma $.


\unitlength=.25mm
\makeatletter
\def\shade{\@ifnextchar[{\shade@special}{\@killglue\special{sh}\ignorespaces}}
\def\shade@special[#1]{\@killglue\special{sh #1}\ignorespaces}
\makeatother

\begin{figure}[h]
\begin{center}
\begin{picture}(479,148)(80,-5)
\thinlines
\typeout{\space\space\space eepic-ture exported by 'qfig'.}
\font\FonttenBI=cmbxti10\relax
\font\FonttwlBI=cmbxti10 scaled \magstep1\relax
\path (80,69)(241,69)
\path (159,0)(159,129)
\path (233.482,71.736)(241,69)(233.482,66.264)
\path (156.264,121.482)(159,129)(161.736,121.482)
\put(168,118){{\rm\rm {$\dot{\theta}$}}}
\put(159,72.25){\arc{38.552}{.169}{2.972}}
\put(159.5,78.206){\arc{120.416}{.154}{2.988}}
\thicklines
\put(158,78.526){\arc{137.328}{.139}{3.002}}
\thinlines
\put(220,73){{\rm\rm {$\pi$}}}
\path (162,56)(156,53)(162,49)
\path (155,23)(147,19)(154,15)
\put(189,5){{\rm\rm {$\Sigma$}}}
\path (355,69)(519,69)
\path (437,30)(437,125)
\path (434.264,117.482)(437,125)(439.736,117.482)
\path (511.483,71.736)(519,69)(511.483,66.264)
\put(443,115){{\rm\rm z}}
\path (365,69)(366.29,69.155)(367.56,69.32)(368.81,69.495)(370.04,69.68)
(371.25,69.875)(372.44,70.08)(373.61,70.295)(374.76,70.52)(375.89,70.755)
(377,71)(378.09,71.3)(379.16,71.6)(380.21,71.9)(381.24,72.2)
(382.25,72.5)(383.24,72.8)(384.21,73.1)(385.16,73.4)(386.09,73.7)
(387,74)(387.71,74.12)(388.44,74.28)(389.19,74.48)(389.96,74.72)
(390.75,75)(391.56,75.32)(392.39,75.68)(393.24,76.08)(394.11,76.52)
(395,77)(395.91,77.52)(396.84,78.08)(397.79,78.68)(398.76,79.32)
(399.75,80)(400.76,80.72)(401.79,81.48)(402.84,82.28)(403.91,83.12)
(405,84)
\path (502,69)(501.945,69.645)(501.88,70.28)(501.805,70.905)(501.72,71.52)
(501.625,72.125)(501.52,72.72)(501.405,73.305)(501.28,73.88)(501.145,74.445)
(501,75)(501.025,75.365)(501,75.76)(500.925,76.185)(500.8,76.64)
(500.625,77.125)(500.4,77.64)(500.125,78.185)(499.8,78.76)(499.425,79.365)
(499,80)(498.255,80.98)(497.52,81.92)(496.795,82.82)(496.08,83.68)
(495.375,84.5)(494.68,85.28)(493.995,86.02)(493.32,86.72)(492.655,87.38)
(492,88)(491.355,88.58)(490.72,89.12)(490.095,89.62)(489.48,90.08)
(488.875,90.5)(488.28,90.88)(487.695,91.22)(487.12,91.52)(486.555,91.78)
(486,92)
\put(359,52){{\rm\rm {$(-r,0)$}}}
\put(517,73){{\rm\rm x}}
\put(487,52){{\rm\rm {$(r,0)$}}}
\put(383,87){{\rm\rm {$C_1$}}}
\put(468,98){{\rm\rm {$\bar{C_1}$}}}
\path (285,101)(286.075,102.125)(287.2,103.2)(288.375,104.225)(289.6,105.2)
(290.875,106.125)(292.2,107)(293.575,107.825)(295,108.6)(296.475,109.325)
(298,110)(299.71,110.67)(301.44,111.28)(303.19,111.83)(304.96,112.32)
(306.75,112.75)(308.56,113.12)(310.39,113.43)(312.24,113.68)(314.11,113.87)
(316,114)(318.27,114.115)(320.48,114.16)(322.63,114.135)(324.72,114.04)
(326.75,113.875)(328.72,113.64)(330.63,113.335)(332.48,112.96)(334.27,112.515)
(336,112)(337.67,111.415)(339.28,110.76)(340.83,110.035)(342.32,109.24)
(343.75,108.375)(345.12,107.44)(346.43,106.435)(347.68,105.36)(348.87,104.215)
(350,103)
\path (346.619,110.251)(350,103)(342.75,106.381)
\put(303,124){{\rm\rm {$R_1$}}}
\end{picture}
\end{center}
\caption{}\label{f13}
\end{figure}

Inspection of Fig. \ref{f12} leads to the following.

\begin{prop}
For each $n\geq 1,$ the return mapping $R_n$ is not proper.
\end{prop}

\begin{proof}
The inverse image of a compact ball centered at $(-r,0)$
corresponds to an asymptotic branch in the parameter space
$(\theta (0),\lambda ).$ The transcendence of this branch can
be easily computed. Indeed when $k'\longrightarrow 0,$ $K(k')
\simeq \rm{ln }1/k'$ and the branch is \it{logarithmic.}
\end{proof}


\subsubsection{Asymptotics of the sphere and wave front near
(r,0) and (-r,0)}
We can estimate the branches ${\bar C}_1$ and $C_1.$
The computations are geometrically different. Indeed the
computation of ${\bar C}_1$ requires the estimation of the
leaves of the foliation ${\cal F}$, localized near the center
but the computation of $C_1$ requires the estimation of the
leaves near the separatrix $\Sigma $ connecting the saddle
points $(-\pi ,0)$ and $(\pi ,0).$ To make the estimation we
use the parametric representation~:
\begin{equation*}
\begin{split}
x(k) & =  -r +2r \ {\displaystyle {E\over K}} \\
z(k) & =  {\displaystyle {r^3\over 6K^3}} [(2k^2-1)E+{k'}^2 K]
\end{split}
\end{equation*}

\noindent where $k\in ]0,1[$ and ${\bar C}_1$ (resp. $C_1$)
is obtained by making $k\rightarrow 0$ (resp. $k\rightarrow 1$).

The transcendence of the branches is related to the properties
of the complete integrals~:
$$K=\int^1_0   {d\eta \over \sqrt{(1-\eta ^2)({k'}^2 +{k'}^2 
+k^2 \eta ^2)}}\ =\ \int ^{\pi /2}_0 (1-k^2\ {\sin\ }^2
\theta )^{-1/2}\ d\theta $$

\noindent and
$$E=\int ^K_0 dn^2u\ du=\int ^{\pi /2}_0 (1-k^2\ {\sin\ }^2
\theta )^{1/2}d\theta $$

Both $E$ and $K$ are solutions of hypergeometric equations
whose singular points are located at $k=0$ and $1.$ Using
this properties we deduce the following \cite{D}.

\begin{lem}   \label{lem4.9}
When $k\rightarrow 0,$ $E$ and $K$ are given by the following
converging asymptotic expansions~:
\begin{equation*}
\begin{split}
K(k) & =  {\displaystyle {\pi \over 2}}\ \Bigl [ 1+\Bigl ( 
{\displaystyle {1\over 2}}\Bigr ) ^2k^2 \ +\ \Bigl ( 
{\displaystyle {1\over 2}}\ .\ {\displaystyle {3\over 4}}
\Bigr ) ^2\ k^4+\cdots \Bigr ]\\
E(k) & =  {\displaystyle {\pi \over 2}}\ \Bigl [1-\Bigl ( 
{\displaystyle {1\over 2}}\Bigr ) ^2 k^2 -{\displaystyle 
{1\over 3}}\ \Bigl ( {\displaystyle {1\over 2}}\ 
{\displaystyle {3\over 4}}\Bigr ) ^2\ k^4+\cdots \Bigr ].
\end{split}
\end{equation*}
\end{lem}

\begin{lem}
When $k'=\sqrt{1-k^2}\rightarrow 0$ we have~:
\begin{equation*}
\begin{split}
E(k) & =  u_1 ({k'}^2)\rm{ln } {\displaystyle {4\over k'}}\ 
+\ u_2({k'}^2)\\
K(k) & =  u_3 ({k'}^2)\rm{ln } {\displaystyle {4\over k'}}\ 
+\ u_4 ({k'}^2)
\end{split}
\end{equation*}

\noindent where the $u_i$'s are analytic near $0$ and can be
written as~:

$ \begin{array}{rclrcl}
u_1({k'}^2) & = & {\displaystyle {{k'}^2\over 2}}\ +\ 
\rm{O}({k'}^4)\ , & u_2 ({k'}^2) &= & 1-{\displaystyle 
{{k'}^2\over 4}}\ +\ \rm{O}({k'}^4)\\
u_3 ({k'}^2) & = & 1+{\displaystyle {{k'}^2\over 4}}\ +\ 
\rm{O}({k'}^4)\ , & u_4 ({k'}^2) &= & -{\displaystyle 
{{k'}^2\over 4}}+\ \rm{O}({k'}^4).
\end{array} $
\end{lem}

\begin{rem}
The complete expansions are given in \cite{D}. The general theory
about Fuchsian differential equations guarantees the convergence
of the previous expansions and the coefficients can be
recursively computed using the ODE. Another method which can
be applied in the general conservative case is to use the
integral formulas.
\end{rem}

\paragraph{Estimation of ${\bar C}_1$ \\}
When $k\rightarrow 0,$ $E$ and $1/K$ are analytic and we have
the following estimations using Lemma \ref{lem4.9}~:
\begin{eqnarray*}
{\displaystyle {E\over K}} & = & 1-{\displaystyle {k^2\over 2}}+
\rm{o} (k^2)\\
x(k)-r & = & -rk^2 +\rm{o} (k^2)\\
z(k) & = & {\displaystyle {2r^3\over 3\pi ^2}}\ k^2 +\rm{o} (k^2)
\end{eqnarray*}

\noindent In particular we deduce the following~:
\begin{prop}   \label{p414}
When $k\rightarrow 0,$ the branch ${\bar C}_1$ is semi-analytic
and is given by a graph of the form~:
$$z = -{2r^2\over 3\pi ^2}\ (x-r)+\rm{o} (x-r)\ ,\ \ x\leq r$$
\end{prop}

\paragraph{Estimation of $C_1$ \\}
When $k'\rightarrow 0,$
we cannot work in the analytic category but in the
\it{log-exp category} introduced in \cite{VDD}. Using
\cite{LR}, the elimination of the parameter $k'$ is allowed in
this category and will lead to a log-exp graph. The precise
algorithm to evaluate $C_1$ has been established in \cite{ABCK}
and we proceed as follows.

We set $X={\displaystyle {x+r\over 2r}} , \ Z={\displaystyle
 {z\over r^3}}$,
and we get~:
$$X = {E\over K}\ ,\ \ Z = {1\over 6K^3}\ [(2k^2-1)E+{k'}^2K]$$
If we introduce~:
$X_1 =k'\ ,\ \ X_2 ={1\over \rm{ln}\ 4/k'}$,
we have $X_1,$ $X_2\rightarrow 0$ when
$k'\rightarrow 0^+$ and both $X$ and $Z$ are analytic
functions of $X_1$ and $X_2.$

An easy computation shows that~:
$$X_1\simeq 4 \rm{e}^{-\f{1}{X}} \ ,\ \ X_2 \simeq X\ \ \ \rm{when\ \ }
 X\rightarrow 0^+$$

\noindent and we can write~:
$$X_1 =4 \rm{e}^{-\f{1}{X}} \ (1+Y_1 (X))\ ,\ \ X_2 =X (1+Y_2(X))$$

\noindent where $Y_1, Y_2\rightarrow 0$ when $X\rightarrow 0^+.$

Both $Y_1$ and $Y_2$ can be compared and a computation gives us~:
$$Y_2 =XA_1 (X,Y_1)\ ,\ \ Y_1 \simeq {Y_2\over X}\ 
\rm{when\ }X\rightarrow 0^+$$

\noindent where $A_1$ is a germ of an analytic function at $0$.

Now the equation $X=E/K$ can be solved in the variables
$Y_1, X_1, X_2$ using the Implicit Function Theorem in the
analytic category and the computations show that~:
$$Y_1=A_2 ( X, \f{\rm{e}^{-\f{1}{X}}}{X}  ) $$

\noindent where $A_2$ is a germ of an analytic function at 0.
Using this relation we end with~:
$$Z=F ( X,\f{\rm{e}^{-\f{1}{X}}}{X}  ) $$

\noindent where $F$ is a germ of an analytic function at 0.

This is the constructive algorithm to compute the branch $C_1$
as a graph in the log-exp category. Hence $Z$ can be expanded as~:
$$Z=\sum ^{+\infty}_{p=0}\ u_p(X)\ \l( \f{\rm{e}^{-\f{1}{X}}}{X}
 \r) ^p $$

\noindent To ensure that $C_1$ is not semi-analytic we must
check that there exists a non zero term of the form
$a_{k,p}X^k\ ( {\displaystyle \rm{e}^{-\f{1}{X}} }
  ) ^p\ ,$ $p>0$ in the expansion. For this we compute
the \it{first} non zero coefficient according to the
lexicographic order on the pair $(p,k).$ The simplest
computation made in \cite{ABCK} is to observe that~:
$$X=\f{E}{K}\ ,\ \ 6Z={1-2{k'}^2\over E^2}\ \l(\f{E}{K}\r)^3+
\ {{k'}^2 
\over K^2}$$

\noindent but the algorithm which can be generalized is the
following. We use the approxi\-mations~:
\begin{eqnarray*}
E & = & 1+{\displaystyle {{k'}^2\over 2}}\ \rm{ln}\ 4/k'\ 
-\ {\displaystyle {{k'}^2\over 4}}\ +\ \rm{o} ({k'}^2)\\
K & = & \rm{ln}\ 4/k'+\ {\displaystyle {{k'}^2\over 4}}\ 
\rm{ln}\ 4/k'+\rm{o} ({k'}^2 \rm{ln}\ 4/k').
\end{eqnarray*}

\noindent Easy computations lead to the formula~:
$$6Z =X^3 -{5\over 4}\ {{k'}^2\over ({\ln}\ 4/k')^4}\ +\rm{o}
 \Bigl ( {{k'}^2\over {\ln} \ 4/k'}\Bigr )$$

\noindent Using $k'\simeq 4\rm{e}^{-\f{1}{X}} ,$ $
 {1\over {\ln}\ 4/k'
}\simeq X$ we obtain~:
$$Z =\inv{6}X^3 -4\rm{e}^{-\f{1}{X}} \ X^3 +\rm{o} (X^3
\rm{e}^{-\f{1}{X}} )$$

\begin{rem}
We observe the following~:
\begin{itemize}
\item $u_0 (X)=X^3/6$ is algebraic.

\item There is a phenomenon of compensation and the first
non zero flat term is of the form $X^3\rm{e}^{-\f{1}{X}}$ and
not $X^2\rm{e}^{-\f{1}{X}}$~; that's why we need three terms in $E$
and two terms in $K.$

\item In general the computation  of the first non zero $a_{p,k}$
can be done in a finite number of steps, for instance using a finite
number of coefficients of $u_0 (X).$
\end{itemize}
\end{rem}


\subsubsection{Numerical aspects}
\begin{figure}[h]
\epsfxsize=7cm
\epsfysize=7cm
\centerline{\epsffile{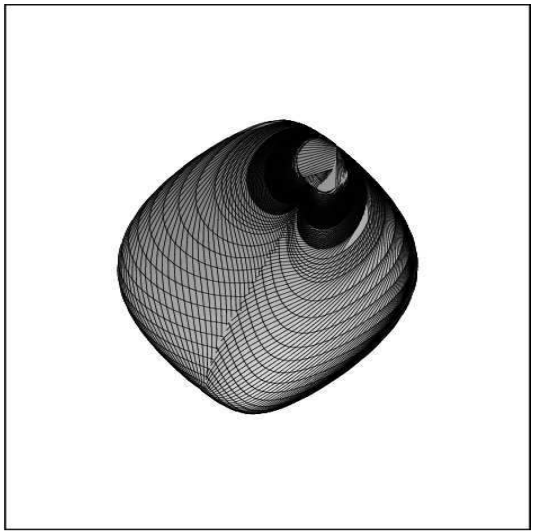}}
\caption{Flat Martinet sphere}\label{f14}
\end{figure}

Fig. \ref{f14} represents the numerical simulation of the flat
Martinet sphere. We observe a
numerical problem when computing near the abnormal direction.


\subsubsection{Asymptotics of the sphere and wave front in the
abnormal direction in the conservative case}
\label{section4.2.13}
\paragraph{Geometric preliminaries \\ }
We can estimate the sphere and the wave front in the abnormal
direction when $g=(1+\alpha y)^2 dx^2+(1+\gamma y)^2dy^2$ (or
in the general case) using the integral formulas (\ref{param}).
We observe that the geometry remains invariant for the following
symmetry~: $S_1:(x,y,z)\longmapsto (-x,y,-z)$ and in our study
we can assume $\lambda \geq 0.$ Another symmetry is the
following. Adding to the geodesics the equations~: ${\dot \alpha
 }=0,$ ${\dot \gamma }=0$ we can observe that the geodesics
equations are left invariant by the transformation~:
$(x,y,z,p_x,p_y,p_z,\alpha ,\gamma )\longmapsto (x,-y,z,p_x,-p_y, 
p_z, -\alpha, -\gamma ).$ Hence we can fix the sign of one of the
parameters and we shall make the following choice~:
$\alpha \geq 0.$

Let $e(t)=(x(t),$ $y(t),$ $z(t))$ be a normal geodesic starting
from $0$ and associated to $p_y (0)={\sin}\ \theta (0), 
p_x ={\cos\ }\theta (0)$ and $p_z =\lambda .$ We observe the
following.
If $\lambda $ is non zero the $y$ component of a geodesic
oscillates periodically unless it corresponds to a separatrix
$\Sigma $ between two values $y_-$ and $y_+$ and we have $y_-
<0<y_+$ if ${\dot y}(0)\not = 0.$ If ${\dot y}(0)=0,$ then
sign ${\ddot y}(0)=$ sign~$\alpha >0$ when $\alpha >0.$

Moreover using Fig. \ref{f4} or the integral formulas (\ref{even}),
we deduce the following Proposition.

\begin{prop}
Let $e(t)=(x(t),$ $y(t),$ $z(t))$ be a geodesic starting from
$0$ such that $y$ oscillates periodically, ${\dot y} (0)\not =
0$ and corresponding to the initial conditions ${\dot y}(0), 
p_x$ and $p_z.$ Let ${\widetilde e} (t)=({\widetilde x}(t), 
{\widetilde y}(t), {\widetilde z}(t))$ be the geodesic
associated to $-{\dot y} (0),$ $p_x$ and $p_z.$ Then $e$
and ${\widetilde e}$ are distinct but their even intersections
with the plane $y=0$ are identical and have the same length.
In particular $e(.)$ is not a minimizer beyond its second
intersection with the plane $y=0.$
\end{prop}

This is illustrated on Fig. \ref{f15} where we project a geodesic
in the plane $(x,y)$.


\setlength{\unitlength}{0.6mm}
\begin{figure}[h]
\begin{center} 

\begin{picture}(180,100)
\thinlines
\drawvector{19.5}{35.95}{60.0}{1}{0}
\drawvector{97.56}{35.95}{60.0}{1}{0}
\drawvector{22.0}{24.0}{54.0}{0}{1}
\drawvector{102.0}{24.0}{54.0}{0}{1}
\drawcenteredtext{158.0}{40.0}{$x$}
\drawcenteredtext{79.19}{40.43}{$t$}
\drawcenteredtext{26.0}{76.0}{$y$}
\drawcenteredtext{106.0}{76.0}{$y$}
\path(102.72,36.02)(102.94,36.06)(103.19,36.09)(103.44,36.13)(103.66,36.18)(103.91,36.25)(104.16,36.31)(104.38,36.38)(104.63,36.47)(104.88,36.56)
\path(104.88,36.56)(105.11,36.66)(105.35,36.77)(105.6,36.88)(105.83,37.02)(106.07,37.15)(106.3,37.29)(106.55,37.43)(106.78,37.59)(107.02,37.75)
\path(107.02,37.75)(107.27,37.93)(107.5,38.11)(107.75,38.29)(108.0,38.49)(108.22,38.7)(108.47,38.9)(108.72,39.13)(108.94,39.36)(109.19,39.59)
\path(109.19,39.59)(109.44,39.84)(109.66,40.09)(109.91,40.34)(110.14,40.61)(110.38,40.88)(110.63,41.18)(110.86,41.47)(111.11,41.77)(111.35,42.08)
\path(111.35,42.08)(111.58,42.38)(111.83,42.72)(112.07,43.04)(112.3,43.38)(112.55,43.74)(112.78,44.09)(113.02,44.45)(113.27,44.83)(113.5,45.2)
\path(113.5,45.2)(113.75,45.59)(113.99,45.99)(114.0,46.0)
\path(114.0,46.0)(114.0,46.0)(114.22,46.43)(114.47,46.86)(114.69,47.29)(114.91,47.72)(115.13,48.13)(115.35,48.56)(115.55,48.97)(115.75,49.36)
\path(115.75,49.36)(115.96,49.77)(116.16,50.18)(116.33,50.56)(116.52,50.95)(116.71,51.34)(116.88,51.72)(117.05,52.09)(117.22,52.47)(117.38,52.84)
\path(117.38,52.84)(117.53,53.2)(117.69,53.56)(117.83,53.91)(117.97,54.25)(118.11,54.61)(118.25,54.95)(118.36,55.29)(118.5,55.61)(118.61,55.95)
\path(118.61,55.95)(118.72,56.27)(118.83,56.59)(118.94,56.9)(119.03,57.2)(119.13,57.52)(119.22,57.81)(119.3,58.11)(119.38,58.4)(119.44,58.7)
\path(119.44,58.7)(119.52,58.97)(119.58,59.25)(119.64,59.54)(119.69,59.81)(119.75,60.06)(119.8,60.34)(119.83,60.59)(119.88,60.84)(119.91,61.09)
\path(119.91,61.09)(119.94,61.34)(119.96,61.58)(119.97,61.81)(119.99,62.04)(119.99,62.27)(120.0,62.49)(119.99,62.7)(119.99,62.93)(119.97,63.13)
\path(119.97,63.13)(119.96,63.34)(119.94,63.54)(119.91,63.74)(119.88,63.93)(119.83,64.11)(119.8,64.3)(119.75,64.47)(119.69,64.65)(119.64,64.81)
\path(119.64,64.81)(119.58,64.98)(119.52,65.13)(119.46,65.3)(119.38,65.44)(119.3,65.59)(119.22,65.73)(119.13,65.87)(119.03,66.01)(118.94,66.13)
\path(118.94,66.13)(118.83,66.26)(118.72,66.38)(118.61,66.51)(118.5,66.62)(118.36,66.73)(118.25,66.83)(118.11,66.93)(117.97,67.01)(117.83,67.11)
\path(117.83,67.11)(117.69,67.19)(117.53,67.27)(117.38,67.36)(117.22,67.43)(117.05,67.5)(116.88,67.55)(116.71,67.62)(116.52,67.68)(116.33,67.73)
\path(116.33,67.73)(116.16,67.76)(115.96,67.81)(115.75,67.84)(115.55,67.88)(115.35,67.91)(115.13,67.94)(114.91,67.95)(114.69,67.98)(114.47,67.98)
\path(114.47,67.98)(114.22,67.98)(114.0,68.0)(114.0,68.0)
\path(114.0,68.0)(114.0,68.0)(113.75,67.98)(113.52,67.98)(113.3,67.98)(113.07,67.95)(112.86,67.94)(112.63,67.91)(112.42,67.88)(112.22,67.84)
\path(112.22,67.84)(112.02,67.81)(111.83,67.76)(111.64,67.73)(111.46,67.68)(111.27,67.62)(111.11,67.55)(110.94,67.5)(110.77,67.43)(110.61,67.36)
\path(110.61,67.36)(110.44,67.27)(110.3,67.19)(110.16,67.12)(110.0,67.01)(109.88,66.93)(109.74,66.83)(109.61,66.73)(109.5,66.62)(109.38,66.51)
\path(109.38,66.51)(109.25,66.38)(109.16,66.26)(109.05,66.13)(108.94,66.01)(108.86,65.87)(108.77,65.73)(108.69,65.59)(108.61,65.44)(108.53,65.3)
\path(108.53,65.3)(108.47,65.13)(108.39,64.98)(108.33,64.81)(108.28,64.65)(108.24,64.48)(108.19,64.3)(108.14,64.11)(108.11,63.93)(108.08,63.74)
\path(108.08,63.74)(108.05,63.54)(108.02,63.34)(108.02,63.13)(108.0,62.93)(108.0,62.7)(108.0,62.5)(108.0,62.27)(108.0,62.04)(108.02,61.81)
\path(108.02,61.81)(108.02,61.58)(108.05,61.34)(108.08,61.09)(108.11,60.84)(108.14,60.59)(108.19,60.34)(108.24,60.08)(108.28,59.81)(108.33,59.54)
\path(108.33,59.54)(108.39,59.25)(108.47,58.97)(108.52,58.7)(108.61,58.4)(108.69,58.11)(108.77,57.81)(108.86,57.52)(108.94,57.22)(109.05,56.9)
\path(109.05,56.9)(109.16,56.59)(109.25,56.27)(109.38,55.95)(109.5,55.61)(109.61,55.29)(109.74,54.95)(109.88,54.61)(110.0,54.25)(110.14,53.91)
\path(110.14,53.91)(110.3,53.56)(110.44,53.2)(110.61,52.84)(110.77,52.47)(110.92,52.09)(111.11,51.72)(111.27,51.34)(111.46,50.95)(111.64,50.56)
\path(111.64,50.56)(111.83,50.18)(112.02,49.77)(112.22,49.36)(112.42,48.97)(112.63,48.56)(112.85,48.13)(113.07,47.72)(113.3,47.29)(113.52,46.86)
\path(113.52,46.86)(113.75,46.43)(113.99,46.0)(114.0,46.0)
\path(114.06,45.77)(114.05,45.77)(114.33,45.38)(114.58,45.0)(114.84,44.61)(115.11,44.25)(115.37,43.9)(115.62,43.54)(115.88,43.2)(116.15,42.86)
\path(116.15,42.86)(116.41,42.52)(116.66,42.2)(116.93,41.9)(117.19,41.59)(117.44,41.29)(117.7,41.0)(117.98,40.72)(118.23,40.45)(118.5,40.18)
\path(118.5,40.18)(118.76,39.93)(119.01,39.68)(119.27,39.43)(119.54,39.2)(119.8,38.97)(120.05,38.75)(120.31,38.54)(120.58,38.34)(120.83,38.15)
\path(120.83,38.15)(121.09,37.97)(121.36,37.79)(121.62,37.61)(121.87,37.45)(122.12,37.29)(122.38,37.15)(122.65,37.0)(122.91,36.88)(123.16,36.75)
\path(123.16,36.75)(123.43,36.63)(123.69,36.52)(123.94,36.43)(124.2,36.33)(124.47,36.25)(124.73,36.16)(124.98,36.09)(125.25,36.04)(125.51,35.99)
\path(125.51,35.99)(125.76,35.93)(126.01,35.9)(126.27,35.88)(126.54,35.86)(126.8,35.84)(127.05,35.84)(127.31,35.84)(127.58,35.84)(127.83,35.86)
\drawdashline{20.0}{68.0}{78.0}{68.0}
\drawdashline{20.0}{28.0}{78.0}{28.0}
\path(22.0,36.0)(22.0,36.0)(23.28,37.18)(24.51,38.34)(25.67,39.5)(26.78,40.63)(27.82,41.74)(28.81,42.81)(29.75,43.88)(30.62,44.93)
\path(30.62,44.93)(31.45,45.95)(32.22,46.95)(32.93,47.93)(33.61,48.9)(34.24,49.83)(34.81,50.75)(35.34,51.63)(35.83,52.5)(36.27,53.36)
\path(36.27,53.36)(36.65,54.18)(37.02,54.99)(37.34,55.77)(37.61,56.52)(37.86,57.27)(38.06,57.97)(38.22,58.66)(38.36,59.34)(38.47,59.97)
\path(38.47,59.97)(38.54,60.59)(38.59,61.2)(38.59,61.77)(38.59,62.33)(38.54,62.86)(38.49,63.36)(38.4,63.84)(38.29,64.3)(38.15,64.73)
\path(38.15,64.73)(38.0,65.13)(37.84,65.51)(37.65,65.87)(37.45,66.2)(37.22,66.51)(36.99,66.8)(36.75,67.06)(36.49,67.3)(36.22,67.51)
\path(36.22,67.51)(35.95,67.69)(35.66,67.86)(35.38,67.98)(35.09,68.09)(34.79,68.19)(34.5,68.25)(34.18,68.27)(33.88,68.27)(33.59,68.26)
\path(33.59,68.26)(33.29,68.22)(33.0,68.15)(32.72,68.05)(32.43,67.91)(32.15,67.76)(31.89,67.58)(31.63,67.38)(31.39,67.16)(31.17,66.9)
\path(31.17,66.9)(30.95,66.61)(30.75,66.3)(30.55,65.95)(30.38,65.58)(30.23,65.19)(30.11,64.76)(30.01,64.31)(29.92,63.84)(29.86,63.33)
\path(29.86,63.33)(29.81,62.79)(29.8,62.24)(29.81,61.65)(29.87,61.02)(29.95,60.38)(30.04,59.7)(30.19,59.0)(30.36,58.25)(30.56,57.5)
\path(30.56,57.5)(30.8,56.7)(31.1,55.88)(31.42,55.04)(31.78,54.15)(32.18,53.25)(32.61,52.31)(33.11,51.34)(33.63,50.34)(34.22,49.31)
\path(34.22,49.31)(34.84,48.25)(35.52,47.15)(36.25,46.04)(37.02,44.9)(37.84,43.72)(38.72,42.5)(39.66,41.27)(40.65,40.0)(41.72,38.68)
\path(41.72,38.68)(42.83,37.36)(43.99,36.0)(44.0,36.0)
\path(52.0,36.0)(52.0,36.0)(53.16,37.29)(54.27,38.58)(55.34,39.83)(56.34,41.04)(57.27,42.25)(58.18,43.41)(59.02,44.56)(59.81,45.66)
\path(59.81,45.66)(60.56,46.75)(61.25,47.81)(61.9,48.84)(62.5,49.86)(63.06,50.84)(63.56,51.79)(64.04,52.7)(64.47,53.61)(64.86,54.47)
\path(64.86,54.47)(65.2,55.31)(65.52,56.13)(65.8,56.91)(66.04,57.68)(66.25,58.41)(66.43,59.13)(66.56,59.81)(66.68,60.45)(66.76,61.09)
\path(66.76,61.09)(66.81,61.68)(66.84,62.27)(66.84,62.81)(66.83,63.34)(66.77,63.83)(66.7,64.3)(66.62,64.73)(66.51,65.16)(66.37,65.55)
\path(66.37,65.55)(66.23,65.91)(66.06,66.25)(65.9,66.55)(65.69,66.83)(65.5,67.09)(65.27,67.33)(65.05,67.54)(64.81,67.72)(64.56,67.87)
\path(64.56,67.87)(64.31,68.0)(64.05,68.09)(63.79,68.16)(63.52,68.22)(63.25,68.25)(63.0,68.25)(62.72,68.22)(62.45,68.16)(62.18,68.08)
\path(62.18,68.08)(61.93,67.98)(61.66,67.84)(61.41,67.69)(61.16,67.51)(60.93,67.3)(60.7,67.06)(60.49,66.8)(60.29,66.52)(60.09,66.22)
\path(60.09,66.22)(59.91,65.87)(59.75,65.51)(59.61,65.13)(59.47,64.73)(59.36,64.29)(59.27,63.83)(59.2,63.34)(59.15,62.83)(59.13,62.29)
\path(59.13,62.29)(59.13,61.74)(59.16,61.15)(59.22,60.54)(59.31,59.9)(59.41,59.24)(59.56,58.54)(59.74,57.84)(59.95,57.11)(60.18,56.34)
\path(60.18,56.34)(60.45,55.56)(60.77,54.75)(61.13,53.91)(61.52,53.06)(61.95,52.16)(62.41,51.27)(62.93,50.33)(63.49,49.38)(64.08,48.4)
\path(64.08,48.4)(64.73,47.38)(65.43,46.36)(66.18,45.29)(66.97,44.22)(67.8,43.11)(68.7,41.99)(69.65,40.84)(70.65,39.66)(71.7,38.47)
\path(71.7,38.47)(72.81,37.24)(73.98,36.0)(74.0,36.0)
\path(22.0,36.0)(22.0,36.0)(22.28,35.68)(22.55,35.36)(22.84,35.06)(23.12,34.77)(23.39,34.47)(23.68,34.18)(23.95,33.9)(24.25,33.63)
\path(24.25,33.63)(24.53,33.36)(24.8,33.11)(25.1,32.86)(25.37,32.61)(25.67,32.38)(25.95,32.13)(26.23,31.9)(26.53,31.69)(26.8,31.47)
\path(26.8,31.47)(27.1,31.27)(27.38,31.06)(27.68,30.87)(27.95,30.69)(28.25,30.5)(28.54,30.32)(28.82,30.15)(29.12,30.0)(29.4,29.84)
\path(29.4,29.84)(29.7,29.69)(29.98,29.54)(30.28,29.4)(30.56,29.28)(30.87,29.14)(31.15,29.03)(31.45,28.92)(31.75,28.8)(32.04,28.71)
\path(32.04,28.71)(32.33,28.62)(32.63,28.54)(32.91,28.45)(33.22,28.37)(33.5,28.31)(33.81,28.25)(34.11,28.2)(34.4,28.14)(34.7,28.11)
\path(34.7,28.11)(35.0,28.06)(35.29,28.04)(35.59,28.02)(35.9,28.01)(36.2,28.0)(36.49,28.0)(36.79,28.0)(37.09,28.01)(37.4,28.02)
\path(37.4,28.02)(37.7,28.04)(38.0,28.06)(38.29,28.11)(38.59,28.14)(38.9,28.2)(39.2,28.25)(39.5,28.3)(39.81,28.37)(40.11,28.45)
\path(40.11,28.45)(40.43,28.54)(40.72,28.62)(41.04,28.7)(41.34,28.8)(41.65,28.92)(41.95,29.03)(42.27,29.14)(42.56,29.27)(42.88,29.4)
\path(42.88,29.4)(43.18,29.54)(43.5,29.69)(43.81,29.84)(44.11,29.98)(44.43,30.15)(44.74,30.32)(45.04,30.5)(45.36,30.69)(45.66,30.87)
\path(45.66,30.87)(45.99,31.06)(46.29,31.27)(46.61,31.47)(46.93,31.69)(47.24,31.9)(47.54,32.13)(47.86,32.38)(48.18,32.61)(48.5,32.86)
\path(48.5,32.86)(48.81,33.11)(49.13,33.36)(49.45,33.63)(49.75,33.9)(50.08,34.18)(50.4,34.47)(50.72,34.77)(51.04,35.06)(51.36,35.36)
\path(51.36,35.36)(51.68,35.68)(51.99,35.99)(52.0,36.0)
\path(44.0,36.0)(44.0,36.0)(44.31,35.68)(44.63,35.36)(44.95,35.06)(45.27,34.77)(45.59,34.47)(45.9,34.18)(46.22,33.9)(46.54,33.63)
\path(46.54,33.63)(46.86,33.36)(47.18,33.11)(47.49,32.86)(47.81,32.61)(48.11,32.38)(48.43,32.13)(48.75,31.9)(49.06,31.69)(49.38,31.47)
\path(49.38,31.47)(49.68,31.27)(50.0,31.06)(50.31,30.87)(50.63,30.69)(50.93,30.5)(51.25,30.32)(51.56,30.15)(51.86,30.0)(52.18,29.84)
\path(52.18,29.84)(52.49,29.69)(52.79,29.54)(53.11,29.4)(53.4,29.28)(53.72,29.14)(54.02,29.03)(54.34,28.92)(54.63,28.8)(54.95,28.71)
\path(54.95,28.71)(55.25,28.62)(55.56,28.54)(55.86,28.45)(56.16,28.37)(56.47,28.31)(56.77,28.25)(57.08,28.2)(57.38,28.14)(57.68,28.11)
\path(57.68,28.11)(57.99,28.06)(58.29,28.04)(58.59,28.02)(58.88,28.01)(59.18,28.0)(59.49,28.0)(59.79,28.0)(60.09,28.01)(60.38,28.02)
\path(60.38,28.02)(60.68,28.04)(60.99,28.06)(61.29,28.11)(61.59,28.14)(61.88,28.2)(62.18,28.25)(62.47,28.3)(62.77,28.37)(63.06,28.45)
\path(63.06,28.45)(63.36,28.54)(63.65,28.62)(63.95,28.7)(64.23,28.8)(64.54,28.92)(64.83,29.03)(65.12,29.14)(65.41,29.27)(65.7,29.4)
\path(65.7,29.4)(66.0,29.54)(66.29,29.69)(66.58,29.84)(66.87,29.98)(67.16,30.15)(67.44,30.32)(67.73,30.5)(68.02,30.69)(68.3,30.87)
\path(68.3,30.87)(68.59,31.06)(68.88,31.27)(69.18,31.47)(69.45,31.69)(69.75,31.9)(70.04,32.13)(70.31,32.38)(70.61,32.61)(70.88,32.86)
\path(70.88,32.86)(71.16,33.11)(71.45,33.36)(71.73,33.63)(72.02,33.9)(72.3,34.18)(72.58,34.47)(72.87,34.77)(73.15,35.06)(73.43,35.36)
\path(73.43,35.36)(73.7,35.68)(73.98,35.99)(74.0,36.0)
\drawcenteredtext{48.0}{14.0}{$\dot{y}(0) \neq 0$}
\drawcenteredtext{130.0}{14.0}{$\dot{y}(0)=0$}
\path(127.83,35.86)(128.08,35.88)(128.35,35.91)(128.61,35.95)(128.86,36.0)(129.13,36.06)(129.38,36.13)(129.64,36.2)(129.89,36.29)(130.16,36.36)
\path(130.16,36.36)(130.41,36.47)(130.67,36.56)(130.94,36.68)(131.19,36.81)(131.44,36.93)(131.71,37.06)(131.97,37.2)(132.22,37.36)(132.49,37.52)
\path(132.49,37.52)(132.74,37.68)(133.0,37.86)(133.25,38.04)(133.52,38.24)(133.77,38.43)(134.02,38.63)(134.28,38.86)(134.55,39.08)(134.8,39.31)
\path(134.8,39.31)(135.05,39.54)(135.32,39.79)(135.58,40.04)(135.83,40.31)(136.08,40.56)(136.35,40.84)(136.61,41.13)(136.86,41.43)(137.11,41.72)
\path(137.11,41.72)(137.38,42.04)(137.63,42.34)(137.88,42.68)(138.14,43.0)(138.39,43.34)(138.66,43.7)(138.91,44.04)(139.16,44.41)(139.42,44.79)
\path(139.42,44.79)(139.69,45.16)(139.94,45.54)(139.94,45.56)
\path(140.0,46.0)(140.0,46.0)(140.22,46.43)(140.47,46.86)(140.69,47.29)(140.91,47.72)(141.13,48.13)(141.35,48.56)(141.55,48.97)(141.75,49.36)
\path(141.75,49.36)(141.96,49.77)(142.16,50.18)(142.33,50.56)(142.52,50.95)(142.71,51.34)(142.88,51.72)(143.05,52.09)(143.22,52.47)(143.38,52.84)
\path(143.38,52.84)(143.53,53.2)(143.69,53.56)(143.83,53.91)(143.97,54.25)(144.11,54.61)(144.25,54.95)(144.36,55.29)(144.5,55.61)(144.61,55.95)
\path(144.61,55.95)(144.72,56.27)(144.83,56.59)(144.94,56.9)(145.03,57.2)(145.13,57.52)(145.22,57.81)(145.3,58.11)(145.38,58.4)(145.44,58.7)
\path(145.44,58.7)(145.52,58.97)(145.58,59.25)(145.64,59.54)(145.69,59.81)(145.75,60.06)(145.8,60.34)(145.83,60.59)(145.88,60.84)(145.91,61.09)
\path(145.91,61.09)(145.94,61.34)(145.96,61.58)(145.97,61.81)(145.99,62.04)(145.99,62.27)(146.0,62.49)(145.99,62.7)(145.99,62.93)(145.97,63.13)
\path(145.97,63.13)(145.96,63.34)(145.94,63.54)(145.91,63.74)(145.88,63.93)(145.83,64.11)(145.8,64.3)(145.75,64.47)(145.69,64.65)(145.64,64.81)
\path(145.64,64.81)(145.58,64.98)(145.52,65.13)(145.46,65.3)(145.38,65.44)(145.3,65.59)(145.22,65.73)(145.13,65.87)(145.03,66.01)(144.94,66.13)
\path(144.94,66.13)(144.83,66.26)(144.72,66.38)(144.61,66.51)(144.5,66.62)(144.36,66.73)(144.25,66.83)(144.11,66.93)(143.97,67.01)(143.83,67.11)
\path(143.83,67.11)(143.69,67.19)(143.53,67.27)(143.38,67.36)(143.22,67.43)(143.05,67.5)(142.88,67.55)(142.71,67.62)(142.52,67.68)(142.33,67.73)
\path(142.33,67.73)(142.16,67.76)(141.96,67.81)(141.75,67.84)(141.55,67.88)(141.35,67.91)(141.13,67.94)(140.91,67.95)(140.69,67.98)(140.47,67.98)
\path(140.47,67.98)(140.22,67.98)(140.0,68.0)(140.0,68.0)
\path(140.0,68.0)(140.0,68.0)(139.75,67.98)(139.52,67.98)(139.3,67.98)(139.07,67.95)(138.86,67.94)(138.63,67.91)(138.42,67.88)(138.22,67.84)
\path(138.22,67.84)(138.02,67.81)(137.83,67.76)(137.64,67.73)(137.46,67.68)(137.27,67.62)(137.11,67.55)(136.94,67.5)(136.77,67.43)(136.61,67.36)
\path(136.61,67.36)(136.44,67.27)(136.3,67.19)(136.16,67.12)(136.0,67.01)(135.88,66.93)(135.74,66.83)(135.61,66.73)(135.5,66.62)(135.38,66.51)
\path(135.38,66.51)(135.25,66.38)(135.16,66.26)(135.05,66.13)(134.94,66.01)(134.86,65.87)(134.77,65.73)(134.69,65.59)(134.61,65.44)(134.53,65.3)
\path(134.53,65.3)(134.47,65.13)(134.39,64.98)(134.33,64.81)(134.28,64.65)(134.24,64.48)(134.19,64.3)(134.14,64.11)(134.11,63.93)(134.08,63.74)
\path(134.08,63.74)(134.05,63.54)(134.02,63.34)(134.02,63.13)(134.0,62.93)(134.0,62.7)(134.0,62.5)(134.0,62.27)(134.0,62.04)(134.02,61.81)
\path(134.02,61.81)(134.02,61.58)(134.05,61.34)(134.08,61.09)(134.11,60.84)(134.14,60.59)(134.19,60.34)(134.24,60.08)(134.28,59.81)(134.33,59.54)
\path(134.33,59.54)(134.39,59.25)(134.47,58.97)(134.52,58.7)(134.61,58.4)(134.69,58.11)(134.77,57.81)(134.86,57.52)(134.94,57.22)(135.05,56.9)
\path(135.05,56.9)(135.16,56.59)(135.25,56.27)(135.38,55.95)(135.5,55.61)(135.61,55.29)(135.74,54.95)(135.88,54.61)(136.0,54.25)(136.14,53.91)
\path(136.14,53.91)(136.3,53.56)(136.44,53.2)(136.61,52.84)(136.77,52.47)(136.92,52.09)(137.11,51.72)(137.27,51.34)(137.46,50.95)(137.64,50.56)
\path(137.64,50.56)(137.83,50.18)(138.02,49.77)(138.22,49.36)(138.42,48.97)(138.63,48.56)(138.85,48.13)(139.07,47.72)(139.3,47.29)(139.52,46.86)
\path(139.52,46.86)(139.75,46.43)(139.99,46.0)(140.0,46.0)
\path(140.0,46.0)(140.0,46.0)(140.07,45.79)(140.16,45.59)(140.22,45.4)(140.32,45.2)(140.38,45.02)(140.47,44.83)(140.55,44.63)(140.63,44.45)
\path(140.63,44.45)(140.72,44.27)(140.78,44.09)(140.88,43.91)(140.94,43.74)(141.03,43.56)(141.11,43.38)(141.19,43.22)(141.27,43.04)(141.35,42.88)
\path(141.35,42.88)(141.44,42.72)(141.5,42.56)(141.6,42.4)(141.66,42.24)(141.75,42.08)(141.83,41.91)(141.91,41.77)(142.0,41.61)(142.07,41.47)
\path(142.07,41.47)(142.16,41.31)(142.22,41.18)(142.32,41.04)(142.38,40.9)(142.47,40.75)(142.55,40.61)(142.63,40.47)(142.72,40.34)(142.78,40.22)
\path(142.78,40.22)(142.88,40.09)(142.94,39.95)(143.03,39.84)(143.11,39.72)(143.19,39.59)(143.27,39.47)(143.35,39.36)(143.44,39.24)(143.5,39.13)
\path(143.5,39.13)(143.6,39.02)(143.66,38.9)(143.75,38.79)(143.83,38.7)(143.91,38.59)(144.0,38.5)(144.07,38.4)(144.16,38.29)(144.22,38.2)
\path(144.22,38.2)(144.32,38.11)(144.38,38.02)(144.47,37.93)(144.55,37.84)(144.63,37.75)(144.72,37.68)(144.78,37.59)(144.88,37.52)(144.94,37.43)
\path(144.94,37.43)(145.03,37.36)(145.11,37.29)(145.19,37.22)(145.27,37.15)(145.35,37.08)(145.44,37.02)(145.5,36.95)(145.6,36.9)(145.66,36.84)
\path(145.66,36.84)(145.75,36.77)(145.83,36.72)(145.91,36.66)(146.0,36.61)(146.07,36.56)(146.16,36.52)(146.22,36.47)(146.32,36.43)(146.38,36.4)
\path(146.38,36.4)(146.47,36.36)(146.55,36.31)(146.63,36.27)(146.72,36.25)(146.78,36.22)(146.88,36.18)(146.94,36.15)(147.03,36.13)(147.11,36.11)
\path(147.11,36.11)(147.19,36.09)(147.27,36.08)(147.35,36.06)(147.44,36.04)(147.5,36.02)(147.6,36.02)(147.66,36.0)(147.75,36.0)(147.83,36.0)
\path(147.83,36.0)(147.91,36.0)(148.0,36.0)(148.0,36.0)
\end{picture}

\end{center}
\caption{}\label{f15}
\end{figure}

\paragraph{Conclusion \\ }
The previous Proposition tells us that except when
$p_y (0)={\sin\ }\theta (0)=0,$ the sphere is contained in
the image of $R_1$ and $R_2.$ The others cases can by studied
by continuity or using a numerical algorithm developped in \cite{Ch}
to compute the conjugate points.

We shall now estimate the image of $R_1$ and $R_2$ near the
two singularities of the foliation ${\cal F}.$

\paragraph{Estimation of $R_1$ \\ }
The constraint $y=0$ takes the form $S:\ {\displaystyle 
{d\theta \over ds}}=\varepsilon \alpha {\cos\ }\theta $ where
${\cos\ }\theta $ can be approximated by $\pm 1$ near
$\theta =0,\pi .$ Contrarily to the flat case we must
distinguish the case $\theta (0)\in ]-\pi ,0[$ where
$\sigma =\rm{sign\ }{\dot y}(0)=+1$ from the case $\theta (0)
\in ]0,\pi [$ where $\sigma =-1.$ We use following notations~:

\begin{itemize}
\item $C(D)$ branches corresponding to an oscillating (resp.
rotating) pendulum or $CD$~: mixed behaviors.

\item Symbols without bars~: behavior near the separatrix,
symbols with bar~: behaviours near the focus.

\item When $\sigma =+1$, we use the symbol '.
\end{itemize}

\noindent They are images by $R_1$ of curves in the parameters
$\lambda ,\theta (0)$ denoted by the same but minuscule symbol.
We obtain the Fig. \ref{f16}.

\unitlength=.25mm
\makeatletter
\def\shade{\@ifnextchar[{\shade@special}{\@killglue\special{sh}\ignorespaces}}
\def\shade@special[#1]{\@killglue\special{sh #1}\ignorespaces}
\makeatother

\begin{figure}[h]
\begin{center}

\begin{picture}(629,234)(20,-5)
\thinlines
\typeout{\space\space\space eepic-ture exported by 'qfig'.}
\font\FonttenBI=cmbxti10\relax
\font\FonttwlBI=cmbxti10 scaled \magstep1\relax
\path (160,223)(160,0)
\path (76,215)(76,15)
\path (244,216)(244,12)
\path (157.264,215.482)(160,223)(162.736,215.482)
\thicklines
\put(160,59){\arc{186.075}{3.586}{5.839}}
\thinlines
\thicklines
\put(160,143){\arc{189.652}{.483}{2.659}}
\thinlines
\path (163,155)(171,151)(163,147)
\path (169,45)(161,49)(169,53)
\path (159,140)(167,136)(159,132)
\put(45,83){{\rm\rm {$-\pi$}}}
\put(251,84){{\rm\rm {$+\pi$}}}
\put(274,106){{\rm\rm {$\theta$}}}
\put(165,212){{\rm\rm {$\f{d\theta}{ds}$}}}
\put(162,84){{\rm\rm {$0$}}}
\dottedline{3}(288,164)(283.715,157.505)(279.36,150.92)(274.935,144.245)(270.44,137.48)
(265.875,130.625)(261.24,123.68)(256.535,116.645)(251.76,109.52)(246.915,102.305)
(242,95)(237.015,87.605)(231.96,80.12)(226.835,72.545)(221.64,64.88)
(216.375,57.125)(211.04,49.28)(205.635,41.345)(200.16,33.32)(194.615,25.205)
(189,17)
\dottedline{3}(299,18)(294.715,24.495)(290.36,31.08)(285.935,37.755)(281.44,44.52)
(276.875,51.375)(272.24,58.32)(267.535,65.355)(262.76,72.48)(257.915,79.695)
(253,87)(248.015,94.395)(242.96,101.88)(237.835,109.455)(232.64,117.12)
(227.375,124.875)(222.04,132.72)(216.635,140.655)(211.16,148.68)(205.615,156.795)
(200,165)
\dottedline{3}(130,19)(125.715,25.495)(121.36,32.08)(116.935,38.755)(112.44,45.52)
(107.875,52.375)(103.24,59.32)(98.535,66.355)(93.76,73.48)(88.915,80.695)
(84,88)(79.015,95.395)(73.96,102.88)(68.835,110.455)(63.64,118.12)
(58.375,125.875)(53.04,133.72)(47.635,141.655)(42.16,149.68)(36.615,157.795)
(31,166)
\dottedline{3}(119,165)(114.715,158.505)(110.36,151.92)(105.935,145.245)(101.44,138.48)
(96.875,131.625)(92.24,124.68)(87.535,117.645)(82.76,110.52)(77.915,103.305)
(73,96)(68.015,88.605)(62.96,81.12)(57.835,73.545)(52.64,65.88)
(47.375,58.125)(42.04,50.28)(36.635,42.345)(31.16,34.32)(25.615,26.205)
(20,18)
\path (57,100)(306,100)(306,100)
\path (298.483,102.736)(306,100)(298.483,97.264)
\thicklines
\put(352.206,28.941){\arc{255.57}{3.731}{4.195}}
\thinlines
\thicklines
\put(380.674,180.152){\arc{315.162}{2.123}{2.608}}
\thinlines
\path (199,69)(302,69)
\path (121,68)(28,68)
\put(245.5,68.25){\arc{23.049}{3.207}{6.218}}
\put(75.5,65.9){\arc{25.35}{3.308}{6.117}}
\path (249,75)(241,79)(249,83)
\path (81,74)(73,78)(81,82)
\put(132,95.5){\arc{65.734}{2.15}{4.133}}
\put(161.5,45){\arc{182.003}{4.15}{5.274}}
\put(192,96){\arc{62.769}{5.247}{7.319}}
\put(300,41){{\rm\rm {$\Sigma$}}}
\put(305,63){{\rm\rm {section}}}
\put(161.5,104.056){\arc{28.754}{3.493}{5.932}}
\path (174.148,116.954)(175,109)(169.234,114.546)
\path (143,109)(179,109)
\put(108,77){{\rm\rm {$c_1$}}}
\put(78,53){{\rm\rm {$d_1$}}}
\put(178,110){{\rm\rm {$\bar{c}_1$}}}
\path (439,100)(596,100)
\path (465,100)(466.97,100.065)(468.88,100.16)(470.73,100.285)(472.52,100.44)
(474.25,100.625)(475.92,100.84)(477.53,101.085)(479.08,101.36)(480.57,101.665)
(482,102)(483.235,102.545)(484.44,103.08)(485.615,103.605)(486.76,104.12)
(487.875,104.625)(488.96,105.12)(490.015,105.605)(491.04,106.08)(492.035,106.545)
(493,107)(493.845,107.22)(494.68,107.48)(495.505,107.78)(496.32,108.12)
(497.125,108.5)(497.92,108.92)(498.705,109.38)(499.48,109.88)(500.245,110.42)
(501,111)(501.745,111.62)(502.48,112.28)(503.205,112.98)(503.92,113.72)
(504.625,114.5)(505.32,115.32)(506.005,116.18)(506.68,117.08)(507.345,118.02)
(508,119)
\path (469,100)(467.575,99.7)(466.2,99.4)(464.875,99.1)(463.6,98.8)
(462.375,98.5)(461.2,98.2)(460.075,97.9)(459,97.6)(457.975,97.3)
(457,97)(456.21,96.745)(455.44,96.48)(454.69,96.205)(453.96,95.92)
(453.25,95.625)(452.56,95.32)(451.89,95.005)(451.24,94.68)(450.61,94.345)
(450,94)(449.455,93.78)(448.92,93.52)(448.395,93.22)(447.88,92.88)
(447.375,92.5)(446.88,92.08)(446.395,91.62)(445.92,91.12)(445.455,90.58)
(445,90)(444.555,89.38)(444.12,88.72)(443.695,88.02)(443.28,87.28)
(442.875,86.5)(442.48,85.68)(442.095,84.82)(441.72,83.92)(441.355,82.98)
(441,82)
\path (568,137)(568.89,137.035)(569.76,137.04)(570.61,137.015)(571.44,136.96)
(572.25,136.875)(573.04,136.76)(573.81,136.615)(574.56,136.44)(575.29,136.235)
(576,136)(576.735,135.6)(577.44,135.2)(578.115,134.8)(578.76,134.4)
(579.375,134)(579.96,133.6)(580.515,133.2)(581.04,132.8)(581.535,132.4)
(582,132)(582.39,131.78)(582.76,131.52)(583.11,131.22)(583.44,130.88)
(583.75,130.5)(584.04,130.08)(584.31,129.62)(584.56,129.12)(584.79,128.58)
(585,128)(585.19,127.38)(585.36,126.72)(585.51,126.02)(585.64,125.28)
(585.75,124.5)(585.84,123.68)(585.91,122.82)(585.96,121.92)(585.99,120.98)
(586,120)
\path (344,146)(346.79,148.04)(349.56,149.96)(352.31,151.76)(355.04,153.44)
(357.75,155)(360.44,156.44)(363.11,157.76)(365.76,158.96)(368.39,160.04)
(371,161)(373.545,161.705)(376.08,162.32)(378.605,162.845)(381.12,163.28)
(383.625,163.625)(386.12,163.88)(388.605,164.045)(391.08,164.12)(393.545,164.105)
(396,164)(399.075,163.445)(402,162.88)(404.775,162.305)(407.4,161.72)
(409.875,161.125)(412.2,160.52)(414.375,159.905)(416.4,159.28)(418.275,158.645)
(420,158)(421.575,157.345)(423,156.68)(424.275,156.005)(425.4,155.32)
(426.375,154.625)(427.2,153.92)(427.875,153.205)(428.4,152.48)(428.775,151.745)
(429,151)
\path (425.619,158.251)(429,151)(421.75,154.381)
\put(372,176){{\rm\rm {$R_1$}}}
\put(485,122){{\rm\rm {$C_1$}}}
\put(561,144){{\rm\rm {$\bar{C_1}$}}}
\put(438,66){{\rm\rm {$D_1$}}}
\put(418,10){{\rm\rm {\underline{case $\sigma=-1$}}}}
\end{picture}

\begin{picture}(634,252)(20,-5)
\thinlines
\typeout{\space\space\space eepic-ture exported by 'qfig'.}
\font\FonttenBI=cmbxti10\relax
\font\FonttwlBI=cmbxti10 scaled \magstep1\relax
\path (160,241)(160,0)
\path (76,233)(76,16)
\path (244,234)(244,17)
\path (157.264,233.482)(160,241)(162.736,233.482)
\thicklines
\put(160,77){\arc{186.075}{3.586}{5.839}}
\thinlines
\thicklines
\put(160,161){\arc{189.652}{.483}{2.659}}
\thinlines
\path (163,173)(171,169)(163,165)
\path (169,62)(161,66)(169,70)
\put(45,101){{\rm\rm {$-\pi$}}}
\put(251,102){{\rm\rm {$+\pi$}}}
\put(274,124){{\rm\rm {$\theta$}}}
\put(165,230){{\rm\rm {$\f{d\theta}{ds}$}}}
\dottedline{3}(288,182)(283.715,175.505)(279.36,168.92)(274.935,162.245)(270.44,155.48)
(265.875,148.625)(261.24,141.68)(256.535,134.645)(251.76,127.52)(246.915,120.305)
(242,113)(237.015,105.605)(231.96,98.12)(226.835,90.545)(221.64,82.88)
(216.375,75.125)(211.04,67.28)(205.635,59.345)(200.16,51.32)(194.615,43.205)
(189,35)
\dottedline{3}(299,36)(294.715,42.495)(290.36,49.08)(285.935,55.755)(281.44,62.52)
(276.875,69.375)(272.24,76.32)(267.535,83.355)(262.76,90.48)(257.915,97.695)
(253,105)(248.015,112.395)(242.96,119.88)(237.835,127.455)(232.64,135.12)
(227.375,142.875)(222.04,150.72)(216.635,158.655)(211.16,166.68)(205.615,174.795)
(200,183)
\dottedline{3}(130,37)(125.715,43.495)(121.36,50.08)(116.935,56.755)(112.44,63.52)
(107.875,70.375)(103.24,77.32)(98.535,84.355)(93.76,91.48)(88.915,98.695)
(84,106)(79.015,113.395)(73.96,120.88)(68.835,128.455)(63.64,136.12)
(58.375,143.875)(53.04,151.72)(47.635,159.655)(42.16,167.68)(36.615,175.795)
(31,184)
\dottedline{3}(119,183)(114.715,176.505)(110.36,169.92)(105.935,163.245)(101.44,156.48)
(96.875,149.625)(92.24,142.68)(87.535,135.645)(82.76,128.52)(77.915,121.305)
(73,114)(68.015,106.605)(62.96,99.12)(57.835,91.545)(52.64,83.88)
(47.375,76.125)(42.04,68.28)(36.635,60.345)(31.16,52.32)(25.615,44.205)
(20,36)
\path (57,118)(306,118)(306,118)
\path (298.483,120.736)(306,118)(298.483,115.264)
\thicklines
\put(352.206,46.941){\arc{255.57}{3.731}{4.195}}
\thinlines
\thicklines
\put(380.674,198.152){\arc{315.162}{2.123}{2.608}}
\thinlines
\path (200,87)(302,87)
\path (121,86)(77,86)
\put(300,59){{\rm\rm {$\Sigma$}}}
\put(305,81){{\rm\rm {section}}}
\path (427,118)(599,118)
\put(160.005,159.302){\arc{206.909}{.774}{2.354}}
\put(211,53){{\rm\rm {$c'_1d'_1$}}}
\path (142,126)(179,126)
\put(140,98){{\rm\rm {$\bar{c}_1'$}}}
\path (319,175)(320.73,176.76)(322.52,178.44)(324.37,180.04)(326.28,181.56)
(328.25,183)(330.28,184.36)(332.37,185.64)(334.52,186.84)(336.73,187.96)
(339,189)(341.555,190.005)(344.12,190.92)(346.695,191.745)(349.28,192.48)
(351.875,193.125)(354.48,193.68)(357.095,194.145)(359.72,194.52)(362.355,194.805)
(365,195)(368.06,195.15)(371.04,195.2)(373.94,195.15)(376.76,195)
(379.5,194.75)(382.16,194.4)(384.74,193.95)(387.24,193.4)(389.66,192.75)
(392,192)(394.26,191.15)(396.44,190.2)(398.54,189.15)(400.56,188)
(402.5,186.75)(404.36,185.4)(406.14,183.95)(407.84,182.4)(409.46,180.75)
(411,179)
\path (407.619,186.251)(411,179)(403.75,182.381)
\path (456,130)(463,130)(474,133)(484,140)(490,148)
\path (565,160)(578,158)(583,151)(587,142)
\put(346,208){{\rm\rm {$R_1$}}}
\put(454,149){{\rm\rm {$C'_1D'_1$}}}
\put(558,168){{\rm\rm {$\bar{C}_1'$}}}
\put(160.5,119.773){\arc{31.561}{5.878}{9.83}}
\path (163,108)(158,105)(163,101)
\put(160.168,162.061){\arc{178.644}{1.019}{2.101}}
\path (155,77)(148,75)(155,71)
\put(142,78){{\rm\rm {$c'_1$}}}
\put(527,154.556){\arc{44.933}{3.783}{5.642}}
\put(516,180){{\rm\rm {$C'_1$}}}
\put(419,51){{\rm\rm {\underline{case $\sigma=+1$}}}}
\end{picture}

\end{center}

\caption{}    \label{f16}
\end{figure}

\paragraph{Estimation of $R_2$ \\ }
The analysis is simpler because the branches corresponding
to $\sigma =+1$ and $\sigma =-1$ are similar.

We get the Fig. \ref{f18}.


\setlength{\unitlength}{0.7mm}
\begin{figure}[h]
\begin{center} 

\begin{picture}(180,100)
\thinlines
\drawpath{10.0}{48.0}{98.0}{48.0}
\thicklines
\drawarc{33.0}{54.0}{43.68}{0.26}{2.85}
\drawarc{75.0}{54.0}{43.68}{0.26}{2.85}
\drawarc{75.0}{42.0}{43.68}{3.41}{6.0}
\thinlines
\path(94.0,64.0)(94.0,64.0)(94.08,64.06)(94.16,64.15)(94.23,64.23)(94.31,64.3)(94.41,64.37)(94.5,64.45)(94.58,64.54)(94.66,64.61)
\path(94.66,64.61)(94.76,64.68)(94.86,64.76)(94.94,64.83)(95.04,64.9)(95.13,64.97)(95.23,65.04)(95.33,65.11)(95.43,65.16)(95.52,65.23)
\path(95.52,65.23)(95.62,65.3)(95.73,65.37)(95.83,65.44)(95.94,65.5)(96.05,65.55)(96.15,65.62)(96.26,65.68)(96.37,65.75)(96.48,65.8)
\path(96.48,65.8)(96.58,65.86)(96.7,65.91)(96.81,65.98)(96.93,66.04)(97.05,66.08)(97.16,66.15)(97.29,66.19)(97.41,66.25)(97.52,66.3)
\path(97.52,66.3)(97.65,66.36)(97.77,66.41)(97.9,66.45)(98.02,66.51)(98.15,66.55)(98.27,66.59)(98.41,66.65)(98.54,66.69)(98.68,66.73)
\path(98.68,66.73)(98.8,66.79)(98.94,66.83)(99.08,66.87)(99.22,66.91)(99.36,66.94)(99.5,67.0)(99.63,67.02)(99.77,67.06)(99.91,67.11)
\path(99.91,67.11)(100.05,67.15)(100.2,67.18)(100.36,67.22)(100.5,67.26)(100.65,67.29)(100.8,67.31)(100.94,67.36)(101.11,67.38)(101.26,67.41)
\path(101.26,67.41)(101.41,67.44)(101.56,67.48)(101.73,67.51)(101.88,67.52)(102.05,67.55)(102.2,67.58)(102.37,67.61)(102.52,67.62)(102.69,67.66)
\path(102.69,67.66)(102.87,67.68)(103.02,67.69)(103.19,67.72)(103.37,67.75)(103.54,67.76)(103.7,67.77)(103.88,67.8)(104.05,67.81)(104.23,67.83)
\path(104.23,67.83)(104.41,67.84)(104.58,67.87)(104.76,67.87)(104.94,67.88)(105.12,67.91)(105.3,67.91)(105.5,67.93)(105.68,67.94)(105.87,67.94)
\path(105.87,67.94)(106.05,67.94)(106.23,67.95)(106.43,67.97)(106.62,67.98)(106.81,67.98)(107.01,67.98)(107.19,67.98)(107.4,67.98)(107.59,67.98)
\path(107.59,67.98)(107.8,67.98)(107.98,68.0)(108.0,68.0)
\drawvector{106.81}{67.98}{1.18}{1}{0}
\drawcenteredtext{108.0}{72.0}{$R_2$}
\drawellipse{75.0}{48.0}{38.0}{24.0}{}
\drawellipse{75.0}{48.0}{14.0}{8.0}{}
\drawdashline{8.0}{40.0}{20.0}{40.0}
\drawdashline{46.0}{40.0}{62.0}{40.0}
\drawdashline{68.0}{50.0}{82.0}{50.0}
\drawpath{76.05}{36.84}{74.47}{35.95}
\drawpath{74.47}{35.95}{76.26}{35.04}
\drawpath{75.8}{44.93}{74.91}{44.25}
\drawpath{74.91}{44.25}{76.05}{43.59}
\path(56.7,39.77)(56.7,39.77)(56.63,39.9)(56.58,40.04)(56.52,40.15)(56.45,40.29)(56.4,40.43)(56.34,40.54)(56.29,40.65)(56.22,40.79)
\path(56.22,40.79)(56.15,40.9)(56.11,41.02)(56.06,41.13)(56.0,41.24)(55.93,41.34)(55.88,41.43)(55.83,41.54)(55.77,41.63)(55.72,41.72)
\path(55.72,41.72)(55.65,41.81)(55.61,41.9)(55.54,42.0)(55.49,42.08)(55.43,42.15)(55.38,42.22)(55.33,42.31)(55.27,42.38)(55.22,42.45)
\path(55.22,42.45)(55.15,42.5)(55.11,42.58)(55.04,42.63)(55.0,42.68)(54.93,42.75)(54.88,42.79)(54.84,42.84)(54.77,42.88)(54.72,42.93)
\path(54.72,42.93)(54.65,42.97)(54.61,43.0)(54.56,43.04)(54.5,43.08)(54.45,43.11)(54.4,43.13)(54.34,43.15)(54.29,43.18)(54.25,43.2)
\path(54.25,43.2)(54.2,43.22)(54.13,43.22)(54.09,43.24)(54.04,43.24)(53.99,43.25)(53.93,43.25)(53.88,43.25)(53.83,43.24)(53.77,43.24)
\path(53.77,43.24)(53.72,43.22)(53.68,43.22)(53.63,43.2)(53.58,43.18)(53.52,43.15)(53.47,43.13)(53.4,43.11)(53.36,43.08)(53.31,43.04)
\path(53.31,43.04)(53.27,43.0)(53.22,42.97)(53.15,42.93)(53.11,42.88)(53.06,42.84)(53.02,42.79)(52.97,42.75)(52.93,42.68)(52.88,42.63)
\path(52.88,42.63)(52.83,42.58)(52.77,42.5)(52.72,42.45)(52.68,42.38)(52.63,42.31)(52.59,42.22)(52.54,42.15)(52.49,42.08)(52.43,42.0)
\path(52.43,42.0)(52.38,41.9)(52.34,41.81)(52.29,41.72)(52.25,41.63)(52.2,41.54)(52.15,41.43)(52.11,41.34)(52.06,41.24)(52.02,41.13)
\path(52.02,41.13)(51.97,41.02)(51.93,40.9)(51.88,40.79)(51.83,40.65)(51.79,40.54)(51.74,40.43)(51.7,40.29)(51.65,40.15)(51.61,40.04)
\path(51.75,39.77)(51.75,39.77)(51.63,39.59)(51.52,39.4)(51.4,39.25)(51.27,39.08)(51.15,38.9)(51.04,38.74)(50.9,38.58)(50.79,38.4)
\path(50.79,38.4)(50.65,38.25)(50.54,38.09)(50.4,37.93)(50.27,37.77)(50.15,37.61)(50.02,37.47)(49.88,37.31)(49.75,37.15)(49.61,37.02)
\path(49.61,37.02)(49.47,36.86)(49.33,36.72)(49.18,36.58)(49.04,36.43)(48.9,36.29)(48.75,36.15)(48.61,36.02)(48.45,35.88)(48.31,35.75)
\path(48.31,35.75)(48.15,35.61)(48.0,35.49)(47.84,35.36)(47.7,35.24)(47.54,35.11)(47.38,34.99)(47.22,34.86)(47.04,34.74)(46.88,34.61)
\path(46.88,34.61)(46.72,34.5)(46.56,34.38)(46.38,34.27)(46.22,34.15)(46.04,34.04)(45.88,33.93)(45.7,33.83)(45.52,33.72)(45.34,33.61)
\path(45.34,33.61)(45.15,33.52)(44.99,33.4)(44.81,33.31)(44.63,33.22)(44.43,33.11)(44.25,33.02)(44.06,32.93)(43.88,32.84)(43.68,32.75)
\path(43.68,32.75)(43.5,32.65)(43.29,32.56)(43.11,32.49)(42.9,32.4)(42.7,32.31)(42.5,32.24)(42.31,32.15)(42.11,32.08)(41.9,32.0)
\path(41.9,32.0)(41.7,31.93)(41.49,31.86)(41.29,31.79)(41.08,31.7)(40.86,31.63)(40.65,31.59)(40.43,31.52)(40.22,31.45)(40.0,31.38)
\path(40.0,31.38)(39.79,31.34)(39.56,31.28)(39.34,31.2)(39.11,31.17)(38.9,31.12)(38.65,31.05)(38.43,31.02)(38.22,30.95)(37.99,30.92)
\path(37.99,30.92)(37.75,30.87)(37.52,30.81)(37.27,30.79)(37.04,30.75)(36.81,30.7)(36.56,30.67)(36.33,30.62)(36.09,30.6)(35.84,30.55)
\path(35.84,30.55)(35.59,30.54)(35.34,30.51)(35.09,30.47)(34.84,30.45)(34.59,30.44)(34.34,30.39)(34.09,30.37)(33.84,30.37)(33.58,30.36)
\path(33.58,30.36)(33.33,30.34)(33.06,30.3)(33.06,30.31)
\drawvector{34.09}{30.37}{1.02}{-1}{0}
\path(33.52,30.31)(33.52,30.3)(33.25,30.31)(32.97,30.34)(32.7,30.35)(32.45,30.36)(32.18,30.37)(31.92,30.37)(31.64,30.38)(31.38,30.42)
\path(31.38,30.42)(31.13,30.44)(30.87,30.45)(30.62,30.46)(30.37,30.51)(30.12,30.53)(29.87,30.54)(29.62,30.59)(29.37,30.62)(29.12,30.64)
\path(29.12,30.64)(28.88,30.7)(28.63,30.71)(28.39,30.77)(28.17,30.79)(27.93,30.86)(27.7,30.88)(27.45,30.95)(27.21,31.0)(27.0,31.04)
\path(27.0,31.04)(26.77,31.1)(26.54,31.13)(26.29,31.2)(26.06,31.27)(25.86,31.31)(25.62,31.37)(25.39,31.45)(25.19,31.52)(24.96,31.56)
\path(24.96,31.56)(24.76,31.63)(24.54,31.7)(24.31,31.79)(24.12,31.87)(23.88,31.94)(23.69,32.02)(23.47,32.09)(23.28,32.18)(23.05,32.25)
\path(23.05,32.25)(22.87,32.34)(22.67,32.43)(22.45,32.52)(22.27,32.61)(22.05,32.7)(21.87,32.79)(21.68,32.88)(21.46,32.99)(21.29,33.09)
\path(21.29,33.09)(21.1,33.18)(20.89,33.29)(20.71,33.38)(20.54,33.5)(20.35,33.61)(20.17,33.72)(19.97,33.83)(19.79,33.93)(19.62,34.06)
\path(19.62,34.06)(19.45,34.15)(19.28,34.29)(19.1,34.4)(18.93,34.54)(18.76,34.65)(18.59,34.79)(18.42,34.9)(18.26,35.04)(18.09,35.15)
\path(18.09,35.15)(17.93,35.31)(17.76,35.43)(17.6,35.58)(17.44,35.72)(17.29,35.86)(17.12,36.0)(16.96,36.13)(16.8,36.29)(16.67,36.43)
\path(16.67,36.43)(16.52,36.59)(16.37,36.74)(16.2,36.88)(16.06,37.04)(15.93,37.2)(15.77,37.36)(15.64,37.52)(15.51,37.68)(15.35,37.84)
\path(15.35,37.84)(15.22,38.0)(15.1,38.18)(14.96,38.34)(14.81,38.52)(14.68,38.68)(14.56,38.86)(14.43,39.04)(14.31,39.22)(14.18,39.4)
\path(14.18,39.4)(14.06,39.58)(13.93,39.75)(13.93,39.77)
\drawcenteredtext{34.4}{26.72}{$d_2$}
\drawcenteredtext{74.47}{57.75}{$c_2$}
\drawcenteredtext{75.37}{41.11}{$\bar{c}_2$}
\path(128.0,52.0)(128.0,52.0)(128.19,52.0)(128.38,52.0)(128.58,52.0)(128.77,52.0)(128.99,52.0)(129.16,52.02)(129.38,52.02)(129.57,52.02)
\path(129.57,52.02)(129.75,52.04)(129.96,52.06)(130.13,52.06)(130.33,52.08)(130.52,52.09)(130.72,52.11)(130.91,52.13)(131.08,52.15)(131.27,52.15)
\path(131.27,52.15)(131.47,52.18)(131.63,52.2)(131.83,52.24)(132.02,52.25)(132.19,52.29)(132.38,52.31)(132.55,52.34)(132.75,52.36)(132.91,52.4)
\path(132.91,52.4)(133.1,52.43)(133.27,52.47)(133.46,52.5)(133.63,52.54)(133.8,52.56)(133.99,52.61)(134.16,52.65)(134.33,52.68)(134.5,52.72)
\path(134.5,52.72)(134.66,52.77)(134.85,52.81)(135.02,52.86)(135.19,52.9)(135.36,52.95)(135.52,53.0)(135.69,53.04)(135.86,53.09)(136.02,53.15)
\path(136.02,53.15)(136.19,53.2)(136.35,53.25)(136.5,53.31)(136.66,53.38)(136.83,53.43)(137.0,53.5)(137.13,53.56)(137.3,53.61)(137.47,53.68)
\path(137.47,53.68)(137.63,53.74)(137.77,53.81)(137.94,53.88)(138.1,53.93)(138.25,54.0)(138.38,54.08)(138.55,54.15)(138.71,54.22)(138.86,54.29)
\path(138.86,54.29)(139.0,54.38)(139.16,54.45)(139.3,54.52)(139.44,54.61)(139.6,54.68)(139.75,54.77)(139.88,54.84)(140.02,54.93)(140.16,55.02)
\path(140.16,55.02)(140.32,55.11)(140.46,55.18)(140.6,55.27)(140.75,55.36)(140.88,55.45)(141.02,55.54)(141.16,55.65)(141.3,55.74)(141.44,55.83)
\path(141.44,55.83)(141.57,55.93)(141.71,56.02)(141.83,56.13)(141.97,56.22)(142.1,56.33)(142.24,56.43)(142.36,56.54)(142.5,56.63)(142.63,56.75)
\path(142.63,56.75)(142.75,56.84)(142.88,56.95)(143.0,57.06)(143.13,57.18)(143.25,57.29)(143.38,57.4)(143.5,57.52)(143.63,57.63)(143.75,57.75)
\path(143.75,57.75)(143.86,57.88)(144.0,57.99)(144.0,58.0)
\path(128.0,52.0)(128.0,52.0)(128.22,52.04)(128.47,52.08)(128.71,52.11)(128.94,52.15)(129.16,52.2)(129.38,52.25)(129.63,52.31)(129.85,52.36)
\path(129.85,52.36)(130.07,52.4)(130.3,52.47)(130.5,52.52)(130.72,52.59)(130.94,52.65)(131.16,52.7)(131.36,52.77)(131.58,52.84)(131.77,52.9)
\path(131.77,52.9)(131.99,52.97)(132.19,53.04)(132.38,53.11)(132.58,53.18)(132.77,53.25)(132.99,53.34)(133.16,53.4)(133.36,53.5)(133.55,53.58)
\path(133.55,53.58)(133.75,53.65)(133.91,53.74)(134.11,53.83)(134.3,53.9)(134.47,54.0)(134.63,54.09)(134.83,54.18)(135.0,54.27)(135.16,54.38)
\path(135.16,54.38)(135.33,54.47)(135.5,54.56)(135.66,54.65)(135.83,54.77)(136.0,54.86)(136.13,54.97)(136.3,55.09)(136.47,55.18)(136.61,55.29)
\path(136.61,55.29)(136.77,55.4)(136.91,55.52)(137.07,55.63)(137.21,55.75)(137.35,55.88)(137.5,55.99)(137.63,56.11)(137.77,56.24)(137.91,56.36)
\path(137.91,56.36)(138.02,56.49)(138.16,56.61)(138.3,56.74)(138.41,56.86)(138.55,57.0)(138.66,57.13)(138.77,57.27)(138.91,57.4)(139.02,57.54)
\path(139.02,57.54)(139.13,57.68)(139.25,57.83)(139.36,57.97)(139.47,58.11)(139.58,58.27)(139.69,58.4)(139.77,58.56)(139.88,58.7)(139.99,58.86)
\path(139.99,58.86)(140.08,59.02)(140.19,59.18)(140.27,59.34)(140.36,59.49)(140.46,59.65)(140.55,59.81)(140.63,59.97)(140.71,60.15)(140.8,60.31)
\path(140.8,60.31)(140.86,60.47)(140.94,60.65)(141.02,60.83)(141.1,61.0)(141.16,61.15)(141.24,61.34)(141.3,61.52)(141.36,61.7)(141.41,61.88)
\path(141.41,61.88)(141.5,62.06)(141.55,62.25)(141.61,62.45)(141.66,62.63)(141.72,62.81)(141.77,63.0)(141.82,63.2)(141.86,63.4)(141.91,63.59)
\path(141.91,63.59)(141.94,63.79)(142.0,63.99)(142.0,64.0)
\path(162.0,66.0)(162.0,66.0)(162.11,65.98)(162.22,65.98)(162.35,65.98)(162.47,65.98)(162.58,65.97)(162.71,65.94)(162.83,65.94)(162.94,65.93)
\path(162.94,65.93)(163.05,65.91)(163.16,65.9)(163.27,65.87)(163.41,65.83)(163.52,65.83)(163.63,65.8)(163.75,65.76)(163.86,65.73)(163.97,65.69)
\path(163.97,65.69)(164.08,65.66)(164.19,65.62)(164.3,65.58)(164.41,65.55)(164.52,65.51)(164.63,65.47)(164.75,65.41)(164.86,65.37)(164.97,65.3)
\path(164.97,65.3)(165.08,65.26)(165.19,65.19)(165.3,65.15)(165.41,65.08)(165.52,65.01)(165.63,64.97)(165.74,64.91)(165.83,64.83)(165.94,64.76)
\path(165.94,64.76)(166.05,64.69)(166.16,64.62)(166.27,64.55)(166.36,64.47)(166.47,64.4)(166.58,64.3)(166.66,64.23)(166.77,64.15)(166.88,64.05)
\path(166.88,64.05)(166.99,63.97)(167.08,63.88)(167.19,63.79)(167.27,63.68)(167.38,63.59)(167.5,63.5)(167.58,63.38)(167.69,63.29)(167.77,63.18)
\path(167.77,63.18)(167.88,63.08)(167.99,62.97)(168.08,62.86)(168.19,62.75)(168.27,62.63)(168.38,62.5)(168.47,62.4)(168.57,62.27)(168.66,62.15)
\path(168.66,62.15)(168.75,62.02)(168.86,61.9)(168.94,61.77)(169.02,61.63)(169.13,61.5)(169.22,61.36)(169.32,61.22)(169.41,61.09)(169.5,60.95)
\path(169.5,60.95)(169.6,60.81)(169.69,60.65)(169.77,60.52)(169.86,60.36)(169.96,60.22)(170.05,60.06)(170.13,59.9)(170.22,59.75)(170.3,59.59)
\path(170.3,59.59)(170.38,59.43)(170.49,59.27)(170.58,59.11)(170.66,58.93)(170.75,58.77)(170.83,58.59)(170.91,58.43)(171.0,58.25)(171.08,58.06)
\path(171.08,58.06)(171.16,57.9)(171.25,57.7)(171.33,57.52)(171.41,57.34)(171.5,57.15)(171.58,56.97)(171.66,56.77)(171.75,56.59)(171.83,56.38)
\path(171.83,56.38)(171.91,56.18)(172.0,56.0)(172.0,56.0)
\drawcenteredtext{140.0}{66.0}{$D_2$}
\drawcenteredtext{148.0}{58.0}{$C_2$}
\drawcenteredtext{164.0}{68.0}{$\bar{C}_2$}
\path(108.0,68.0)(108.0,68.0)(108.15,67.98)(108.3,67.98)(108.47,67.98)(108.62,67.98)(108.77,67.98)(108.93,67.98)(109.08,67.98)(109.23,67.97)
\path(109.23,67.97)(109.38,67.95)(109.54,67.94)(109.68,67.94)(109.83,67.94)(109.97,67.93)(110.12,67.91)(110.26,67.91)(110.4,67.88)(110.54,67.87)
\path(110.54,67.87)(110.68,67.87)(110.81,67.84)(110.94,67.83)(111.08,67.81)(111.22,67.8)(111.36,67.77)(111.48,67.76)(111.62,67.75)(111.75,67.72)
\path(111.75,67.72)(111.87,67.69)(112.0,67.68)(112.12,67.66)(112.26,67.62)(112.37,67.61)(112.5,67.58)(112.62,67.55)(112.73,67.52)(112.86,67.51)
\path(112.86,67.51)(112.98,67.48)(113.08,67.44)(113.2,67.41)(113.31,67.38)(113.44,67.36)(113.55,67.31)(113.66,67.29)(113.76,67.26)(113.87,67.22)
\path(113.87,67.22)(113.98,67.19)(114.08,67.15)(114.19,67.11)(114.29,67.06)(114.38,67.02)(114.5,67.0)(114.58,66.94)(114.69,66.91)(114.79,66.87)
\path(114.79,66.87)(114.88,66.83)(114.98,66.79)(115.06,66.73)(115.16,66.69)(115.26,66.65)(115.34,66.59)(115.43,66.55)(115.51,66.51)(115.61,66.45)
\path(115.61,66.45)(115.69,66.41)(115.77,66.36)(115.86,66.3)(115.94,66.25)(116.01,66.19)(116.09,66.15)(116.18,66.08)(116.26,66.04)(116.33,65.98)
\path(116.33,65.98)(116.4,65.91)(116.48,65.86)(116.55,65.8)(116.62,65.75)(116.69,65.68)(116.76,65.62)(116.81,65.55)(116.88,65.5)(116.94,65.44)
\path(116.94,65.44)(117.01,65.37)(117.08,65.3)(117.13,65.23)(117.19,65.16)(117.26,65.11)(117.31,65.04)(117.37,64.97)(117.43,64.9)(117.48,64.83)
\path(117.48,64.83)(117.52,64.76)(117.58,64.68)(117.63,64.61)(117.69,64.54)(117.73,64.45)(117.77,64.38)(117.83,64.3)(117.87,64.23)(117.91,64.15)
\path(117.91,64.15)(117.94,64.06)(118.0,64.0)(118.0,64.0)
\end{picture}

\end{center}
\caption{$\sigma=\pm 1$}\label{f18}
\end{figure}
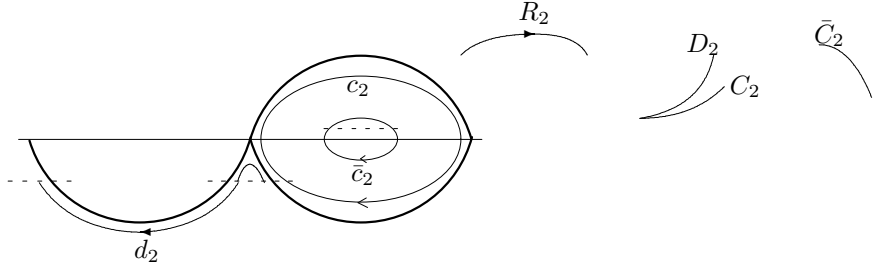

\paragraph{Estimation problems \\ }
We must estimate the branches $C_1, D_1, C'_1D'_1, {\bar C}_1,
\bar{C}'_1, C_2, D_2$ and ${\bar C}_2.$ We know a priori the
following~:

\begin{itemize}
\item The branches ${\bar C}_1, {\bar C}'_1, C'_1 D'_1$ and
${\bar C}_2$ are semi-analytic. We must check if they end on the
abnormal direction.

\item The branches $C_1, D_1, C_2$ and $D_2$ are in the
exp-log category and are ending on the abnormal direction.

\item We must compare the positions of the
branches
$C_1, D_1, C_2$ and $D_2$ to determine which ones
are in the sphere.
\end{itemize}

All the computations are made in the general integrable case,
i.e. the coefficients of the metrics $a$ and $c$ are \it{analytic
functions of $y$} so that~:
\begin{equation*}
\begin{split}
a&=1+2\alpha y+\cdots \\
c&=1+2\gamma y+\cdots
\end{split}
\end{equation*}

\no Our computations are based on the integral formulas
(\ref{even}) and lead to the following~:

\begin{prop}[Comparison of branches $C_1, C_2, D_2$]
Let $X=\f{x+r}{2r}$ and $Z=\f{z}{r^3}$. We have the estimates~:
\begin{itemize}
\item branch $C_1$ :
$Z=\inv{6}X^3+\l(\f{r^2\alpha ^2}{64}+\f{\pi
r}{32}(\alpha+\gamma) \r)X^4+\rm{o}(X^4) $
\item branch $C_2$ :
$Z=\inv{24}X^3+\rm{o}(X^3)$
\item branch $D_2$ :
$Z=\inv{6}X^3+\l(\f{r^2\alpha ^2}{64}-\f{\pi
r}{32}(\alpha+\gamma) \r)X^4+\rm{o}(X^4) $
\end{itemize}
and we can conclude :
\begin{itemize}
\item if $\gamma>-\alpha $, the branch $C_1$ is in the sphere.
\item if $\gamma<-\alpha $, the branch $D_2$ is in the sphere.
\end{itemize}
\end{prop}

\begin{rem}
At $0$ the Gauss curvature of the Riemannian metric $g_R=adx^2 +
cdy^2$ is $K=\f{\alpha (\alpha +\gamma )+\beta^2}{4}$. If $\beta
=0$, it reduces to $\alpha (\alpha +\gamma )\over 4$. Hence the
critical value $\alpha +\gamma =0$ is connected to $K=0$.
\end{rem}

If $\alpha =0$ in the gradated form of order $0$,
the section reduces to : $y=0$. Then the branch
$D_1$ does not exist (see Fig. \ref{f16}) and the branch
$\bar{C}_1=\bar{C}'_1$ ends on the abnormal direction (and is in
the sphere). Also the branches $C'_1D'_1$ and $\bar{C}_2$ end on
the abnormal direction, but are not in the sphere, as can easily
checked.

If $\alpha \neq 0 $, the branches $\bar{C}_1, \bar{C}'_1,
C'_1D'_1$ and $\bar{C}_2$ do not end on the abnormal direction. A
new branch appears~: $D_1$, which is the only branch in $z<0$
that ends on $(-r,0)$ (the same is available in $z>0$ on
$(r,0)$). Therefore \it{$D_1$ is in the sphere}.
\\

Hence we know the asymptotics of the trace of the sphere with
$y=0$ near the singularity $(-r,0)$ (resp.
$(r,0)$) in the general integrable case. Now an important
question is to check \it{in which class it is}. In \cite{ABCK}
it was proved that the sphere in the flat case is \it{not
subanalytic}. Very precise evaluations of flat terms of branch
$C_1$ lead to the following~:

\begin{thm}
In the general integrable case the sphere is not subanalytic.
\end{thm}

\begin{rem}
This result cannot be obtained by perturbation
of the flat case. The explanation is the
following.

\no We proved that in the flat case the sphere is not
subanalytic~:
$$Z=\f{1}{6}X^3-4X^3\rm{e}^{-\inv{X}}+\rm{o}(X^3\rm{e}^{-\inv{X}})$$
In the general case (not only integrable) a natural idea would be
to invoke some perturbation argument in order to check non
subanalyticity.
We may think that the previous graph is
continuous with respect to the coefficients of the metrics, or
with respect to the radius of the sphere.
But this is wrong, as shown in the following
example~:
$$F_1=(1+\varepsilon y)\f{\partial}{\partial x}+\f{y^2}{2}
\f{\partial}{\partial z}\ ,\ \ F_2=\f{\partial}{\partial y}\ , \
\ g=\inv{(1+\varepsilon y)^2}dx^2+dy^2\ \ (\varepsilon<0)$$
We obtain~:
$$Z=\inv{6}X^3+\f{r^2\varepsilon ^2}{32}X^4+\cdots+r\varepsilon
(\f{3}{4}-\f{7}{12}r^2\varepsilon ^2)X^4\rm{e}^{-\inv{X}}
+\rm{o}(X^4\rm{e}^{-\inv{X}})$$
This is actually not surprising, since in the step of elimination
of the parameter $k'$ (see \cite{BLT}), we replaced $k'$ with its
expression in function of $X$. But this step needs an
exponentiation, and we know that equivalents do not pass through
exponentiation.

However we could expect that the expansions of $X$ and $Z$ in
function of $\inv{\sqrt{\lambda }}, k'$ (see \cite{BLT}) are
continuous
with respect to the coefficients. It is still wrong~:
\begin{itemize}
\item flat case : \esp \! $Z-\inv{6}X^3=-2\f{{k'}^2}{\lambda ^{3/2}}
+\rm{o}(\f{{k'}^2}{\lambda ^{3/2}})$
\item case $\varepsilon <0$ : $Z-\inv{6}X^3=\inv{\lambda ^2}
\rm{An}(\inv{\sqrt{\lambda}})-(3+\f{r^2\varepsilon ^2}{4})
\f{{k'}^2}{\lambda ^{3/2}}+\rm{o}(\f{{k'}^2}{\lambda ^{3/2}})$
\end{itemize}
Nevertheless we can observe that \it{the analytic part of the
graph is always continuous with respect to the coefficients}.
Instability only
appears in flat terms. This can be easily explained
in the case $\varepsilon <0$~: to compute $X$ and $Z$, we need to
evaluate some integrals. To do that, the change of variable $\eta
=\f{k'}{k}\rm{sh } t$ is relevant (see \cite{BLT}) and leads to
expand $X$ and $Z$ as a sum of terms containing $\rm{Argsh } \f{
\varepsilon p_x}{2k'\sqrt{\lambda}}$.
Now if one wants to expand this last expression (using the
formula $\rm{Argsh }x=\ln (x+\sqrt{1+x^2})$), with $k'
\sqrt{\lambda} \rightarrow 0$, it is necessary to assume \it{
$\varepsilon $ fixed} (so as $r$) to get~:
\begin{equation*}
\begin{split}
 \rm{Argsh } \f{
 \varepsilon p_x}{2k'\sqrt{\lambda}}
 &=\ln \f{\varepsilon p_x}{2k'\sqrt{\lambda}}\ + \ 
  \ln \l( 1+\sqrt{1+\f{4{k'}^2\lambda}{\varepsilon ^2p_x^2}} \r) \\
 &=\ln \f{\varepsilon p_x}{2k'\sqrt{\lambda}}\ + \ 
 \ln 2 +\rm{An}\l( \f{{k'}^2\lambda}{\varepsilon ^2p_x^2} \r)
\end{split}
\end{equation*}
in order to obtain analytic expansions of $X$ and $Z$, which
prove that the sphere belongs to the log-exp category.
Unfortunately in 
this last expression, there is no sense to make $\varepsilon
\rightarrow 0$ because we needed to assume $\varepsilon $ fixed.
Moreover note that ${k'}^2\rm{sh }2\rm{Argsh }
\f{\varepsilon p_x}{2k'\sqrt{\lambda}}
=\f{\varepsilon ^2p_x^2}{2\lambda}+{k'}^2+\rm{o}({k'}^2)$, so
that this term brings new flat terms with coefficients having the
same order as unity.

We could now expect to have continuity with respect to
parameters if we do not expand the Argsh's, and try to make the
following reasoning~:
\begin{enumerate}
\item $x\longmapsto f(0,x)$ is not subanalytic. 
\item $\varepsilon \longmapsto f(\varepsilon,x)$ is continuous.
\end{enumerate}
Then for $\varepsilon \neq 0 \quad 
x\longmapsto f(\varepsilon,x)$ is not subanalytic.

\no But this is wrong, see the following example~:
\begin{equation*}
\begin{split}
f(\varepsilon,t)&=\ln t+\rm{Argsh }\f{\varepsilon}{t}
=\ln 2\varepsilon + \f{t^2}{2\varepsilon ^2}+\cdots \ \ :
\ \rm{analytic in $t$.} \\
f(0,t)&=\ln t \ \ : \ \rm{not subanalytic.}
\end{split}
\end{equation*}
\end{rem}

So the sphere is not subanalytic. Now the main question is : in
which category is the sphere ?
In \cite{BLT}, we proved that the branch $C_1$ belongs to the
log-exp category. A precise answer is the following~:

\begin{prop}
We set near the
singularity $(-r,0)$~: $X=\f{x+r}{2r}$, $Z=\f{z}{r^3}$, and we
have~:
\begin{itemize}
\item branch $C_1$ :
$Z=\rm{An}(X,X\!\ln X,X\!\ln ^2X,X\!\ln ^3X,\f{\rm{e}^{-\inv{X}}
}{X^3})=\inv{6}X^3+\cdots$ \\
where $\rm{An}(.)$ is a germ at $0$ of an analytic function. 
Moreover the analytic part of $Z(X)$ is continuous with respect
to $r$ and the coefficients of the metrics. \\
A similar result holds for $D_2$.
\item branch $D_1$ :
$Z=\rm{An}(\sqrt{-X},\sqrt{-X}\ln (-X),\rm{e}^{-\f{r\alpha}{2
\sqrt{2}\sqrt{-X}}})
=\f{-8}{r^2\alpha ^2}X^2+\cdots$ \\
Moreover the analytic part of $X(\sqrt{Z})$ is continuous
with respect
to $r$ and the coefficients of the metrics.
\end{itemize}
\end{prop}

\begin{cor}\label{corLE}
In the general integrable Martinet case the sphere belongs to the
log-exp category.
\end{cor}

\begin{proof}
Our estimations show that near the abnormal direction the sphere
is log-exp. In the other directions the sphere is subanalytic,
see \cite{Ag}.
\end{proof}


\subsubsection{Asymptotics of the sphere and wave front in the
abnormal direction in the general gradated case of order 0}
We set : $g=(1+\alpha y)^2dx^2+(1+\beta x+\gamma y)^2dy^2$ with
$\alpha, \beta \neq 0$. In this case the equation in $(\theta,
\dot{\theta})$ obtained by projection is not integrable. In order
to compute the asymptotics of the sphere we can use \it{formal
first integrals} near the saddles.
Moreover toric blowing-up allow us to evaluate the solution if
$\lambda$ is fixed, see \cite{Bru,Cha}.
The technics are similar to
the ones used by \cite{R} and others to evaluate the
Poincar\'e-Dulac return mapping near a polycycle for a
one-parameter family $(X_\varepsilon )$ of vector fields. This
computation can be reduced to the evaluation of the
Poincar\'e-Dulac
mapping \it{near a resonant saddle}~:
$$
X_\varepsilon =\lambda_1(\varepsilon )\f{\partial}{\partial x'}
+\lambda_2(\varepsilon )\f{\partial}{\partial y'}$$
$$r=\f{\lambda_1(\varepsilon )}{\lambda_2(\varepsilon )}
\qquad \lambda_1\lambda_2<0$$
$$r(0)=-1 \qquad r=-1+\alpha_1(\varepsilon )$$
In this method we use a normal form in which~:
\begin{itemize}
\item the separatrices are normalized to~: $x'=0, y'=0$.
\item $X_\varepsilon \sim x'\f{\partial}{\partial x'}
-y'\f{\partial}{\partial y'}-\l( \displaystyle{\overset{N}
{\underset{i=0}{\sum}}}
   \alpha_{i+1}(
\varepsilon )(x'y')^i\r) y'\f{\partial}{\partial y'}$
\end{itemize}
up to a change of time parametrization, and in some neighborhood
$U\times {\cal{E}}(\varepsilon ) \subset {\R}^{2} \times \R$, with
${\cal{E}}(\varepsilon )\underset{\varepsilon \rightarrow 0}
{\longrightarrow} 0$, see \cite{R} for details.

\unitlength=.25mm
\makeatletter
\def\shade{\@ifnextchar[{\shade@special}{\@killglue\special{sh}\ignorespaces}}
\def\shade@special[#1]{\@killglue\special{sh #1}\ignorespaces}
\makeatother

\begin{figure}[h]
\begin{center}

\begin{picture}(184,109)(247,-5)
\thinlines
\typeout{\space\space\space eepic-ture exported by 'qfig'.}
\font\FonttenBI=cmbxti10\relax
\font\FonttwlBI=cmbxti10 scaled \magstep1\relax
\path (260,22)(370,22)
\path (270,6)(270,93)
\path (258,73)(299,73)
\path (346,42)(346,0)
\path (362.483,24.736)(370,22)(362.483,19.264)
\path (267.264,85.482)(270,93)(272.736,85.482)
\put(351,43){{\rm\rm {$\sigma_1$}}}
\put(283,78){{\rm\rm {$\sigma_2$}}}
\put(361,2){{\rm\rm {$y'$}}}
\put(247,87){{\rm\rm {$x'$}}}
\put(342.367,89.614){\arc{119.449}{1.51}{2.86}}
\path (306,27)(298,22)(306,17)
\path (265,41)(270,48)(274,41)
\path (300,43)(298,50)(306,49)
\end{picture}

\end{center}
\end{figure}

In our case the method has to be \it{adapted}. Indeed in our case
the section depends on the parameter and is imposed by the
geometry. In particular the distance of the saddle point to the
section $\sigma_\varepsilon $ tends to $0$ when $\varepsilon 
\rightarrow 0$. The method is the following. Let $d$ be the
distance to the separatrix. Then we want to compute~:
$d\longmapsto (x(d)+r,z(d))$ when $d\approx 0$ (using the normal form
for $(X_\varepsilon , \sigma_\varepsilon )$).

This computation generalizes the computation in the conservative
case where $d$ is the distance to the root of multiplicity two of
the potential.

The algorithm to evaluate step by step this application is to
consider the $k$-jet of $(X_\varepsilon , \sigma_\varepsilon )$.
It is not clear a priori that the $k$-jet is sufficient to
compute the first $k$ terms in the expansion. However we shall
prove that the $1$-jet is sufficient to compute the first term in
the expansion. It gives us the contact of the branch $D_1$ with
the abnormal direction.

\begin{prop} \label{propcontact}
Let us suppose $a=(1+\alpha y)^2, c=(1+\beta x+\gamma y)^2$ with
$\alpha > 0$. Let $X=\f{x+r}{2r}$, $Z=\f{z}{r^3}$. Then near
$X=0$ the graph of the branch $D_1$ is the following~:
$$Z=\f{-2}{r^2\alpha ^2}X^2+\rm{o}(X^2)$$
\end{prop}

\begin{rem}
Observe that in the flat case, the abnormal geodesic is not
strict and the contact is of order $1$ (see prop \ref{p414}).
\end{rem}

\begin{proof}
The differential system is~:
\begin{equation*}
\begin{split}
\f{dx}{dt}&=\f{\cos \theta}{1+\alpha y} \ , \ \ 
\f{dy}{dt}=\f{\sin \theta}{1+\beta x+\gamma y} \ , \ \ 
\f{dz}{dt}=\f{y^2}{2}\f{\cos \theta}{1+\alpha y} \\
\f{d\theta}{dt}&=-\inv{(1+\alpha y)(1+\beta x+\gamma y)}
(\lambda y -\alpha \cos \theta+\beta \sin \theta)
\end{split}
\end{equation*}
Reparametrizing with : $ds=\sqrt{\lambda}
\inv{(1+\alpha y)(1+\beta x+\gamma y)}dt$,
we obtain~:
\begin{equation}
\begin{split}
\f{dx}{ds}&=\inv{\rala}(1+\beta x+\gamma y)\cos \theta \label
{eqx} \\
\f{dy}{ds}&=\inv{\rala}(1+\alpha y)\sin \theta \\
\f{dz}{ds}&=\inv{\rala}\f{y^2}{2}(1+\beta x
+\gamma y){\cos \theta} \\
\f{d\theta}{ds}&=
-\sqrt{\lambda} y +\f{\alpha}{\rala}
 \cos \theta-\f{\beta}{\rala} \sin \theta
\end{split}
\end{equation}
Hence the equation governing $\theta$ is~:
$$\f{d^2\theta}{ds^2}+\sin \theta +\f{\alpha^2}{\lambda}\sin
\theta \cos \theta -\f{\alpha \beta}{\lambda}\sin^2\theta
+\f{\beta}{\rala}\cos \theta \f{d\theta}{ds}=0$$
Set $u=\theta +\pi, v=\f{d\theta}{ds}$. Then~:
\begin{equation*}
\begin{split}
\f{du}{ds}&=v \\
\f{dv}{ds}&=\sin u-\f{\alpha^2}{\lambda}\sin u\cos u+\f{\alpha
\beta}{\lambda}\sin ^2u+\f{\beta}{\rala}v\cos u
\end{split}
\end{equation*}
The eigenvalues of the linearized system are solutions of :
$\mu^2-\f{\beta}{\rala}\mu-(1-\f{\alpha^2}{\lambda})=0$, hence~:
$\mu_1=1+\f{\beta}{2\rala}+\rm{O}\l(\inv{\lambda}\r) \ , \ \ 
\mu_2=-1+\f{\beta}{2\rala}+\rm{O}\l(\inv{\lambda}\r)$.

\no Let $u=u_1+v_1, v=\mu_1u_1+\mu_2u_2$. We get~:
\begin{equation*}
\begin{split}
\f{du_1}{ds}&=\mu_1u_1+\rm{O}\l(\f{u_1^2}{\lambda},
\f{v_1^2}{\lambda}, \f{u_1v_1}{\lambda}, \f{u_1^3}{\rala},
\f{u_1^2v_1}{\rala}, \f{u_1v_1^2}{\rala}, \f{v_1^3}{\rala}\r) \\
\f{dv_1}{ds}&=\mu_2v_1+\rm{O}\l(\f{u_1^2}{\lambda},
\f{v_1^2}{\lambda}, \f{u_1v_1}{\lambda}, \f{u_1^3}{\rala},
\f{u_1^2v_1}{\rala}, \f{u_1v_1^2}{\rala}, \f{v_1^3}{\rala}\r) \\
\end{split}
\end{equation*}
and after integration~:
\begin{equation*}
\begin{split}
u(s)&=A\rm{e}^{\mu_1s}+B\rm{e}^{\mu_2s}
+\rm{O}\Big(\f{A^2}{\lambda}\rm{e}^{2\mu
_1s}, \f{B^2}{\lambda}\rm{e}^{2\mu
_2s}, \f{AB}{\lambda}\rm{e}^{(\mu
_1+\mu_2)s}, \f{A^3}{\sqrt{\lambda}}\rm{e}^{3\mu
_1s},  \\
& \qquad \qquad \qquad \qquad \qquad \qquad \qquad
A^2B\rm{e}^{(2\mu
_1+\mu_2)s}, AB^2\rm{e}^{(\mu
_1+2\mu_2)s}, \f{B^3}{\sqrt{\lambda}}\rm{e}^{3\mu
_2s} \Big) \\
v(s)&=\mu_1A\rm{e}^{\mu_1s}+\mu_2B\rm{e}^{\mu_2s}
+\rm{O}(\cdots)
\end{split}
\end{equation*}
where $A$ and $B$ are constants to determine.\\
The section is $y=0$, hence~:
$v=\f{-\alpha}{\rala}\cos u+\f{\beta}{\rala}\sin u
=\f{-\alpha}{\rala}+\f{\beta}{\rala}u+\rm{O}\l(\f{u^2}
{\rala}\r)$
Let $s_f$ be the parameter corresponding to the final time $t=r$,
i.e. : $y(0)=y(s_f)=0$. Putting
these conditions in the previous equations we obtain~:
$$B=\f{\alpha}{\rala}+\rm{O}\l(\inv{\lambda}\r)\ , \ \ 
A=\f{-\alpha}{\rala}\rm{e}^{-\mu_1s_f}+\rm{O}\l(
\inv{\lambda}\rm{e}^{-\mu_1s_f}\r)$$
Hence~:
\begin{equation}  \label{theta}
\theta(s)+\pi=u(s)=\f{-\alpha}{\rala}\rm{e}^{\mu_1(s-s_f)} +
\f{\alpha}{\rala}\rm{e}^{\mu_2s} + \rm{O}\l(\inv{\lambda}
\rm{e}^{\mu_1(s-s_f)}, \inv{\lambda}\rm{e}^{\mu_2s}\r)
\end{equation}

\no To get $y$, just note that :
$y=-\inv{\rala}\f{d\theta}{ds}+\f{\alpha}{\lambda}\cos
\theta-\f{\beta}{\rala}\sin \theta$, hence~:
\begin{equation}  \label{y}
y(s)=-\f{\alpha}{\lambda}+\f{\alpha}{\lambda}
\rm{e}^{\mu_1(s-s_f)}+\f{\alpha}{\lambda}
\rm{e}^{\mu_2s}+\rm{O}\l(\inv{\lambda^{3/2}}\r)
\end{equation}

Then we have to compute $x$, which amounts to integrating
equation (\ref{eqx}). We get~:
\begin{equation}  \label{x}
1+\beta x(s)=\rm{e}^{-\f{\beta}{\rala}s}+\f{\alpha \gamma}
{\lambda}-\f{\alpha \gamma}{\lambda}\rm{e}^{-\f{\beta}{\rala}s}
+\rm{O}\l(\inv{\lambda^{3/2}}\r)
\end{equation}

The computation of $z$ is then similar and we obtain~:
\begin{equation}  \label{z}
z(s)=-\f{\alpha^2}{2\lambda^2}\f{1-\rm{e}^{-\f{\beta}{\rala}s}}
{\beta}+\rm{o}\l(\inv{\lambda^2}\r)
\end{equation}

It remains to estimate $s_f$. From the equation~:
$\f{dt}{ds}=\inv{\rala}(1+\alpha y)(1+\beta x+\gamma y)$ we get~:
\begin{equation}  \label{r}
r=\f{1-\rm{e}^{-\f{\beta}{\rala}s_f}}{\beta}\l(1-\alpha \f{\alpha
+\gamma}{\lambda} \r)+\rm{O}\l(\inv{\lambda^{3/2}}\r)
\end{equation}

This leads to the conclusion~:
$Z=\f{-2}{r^2\alpha^2}X^2+\rm{o}(X^2)$

\end{proof}

\begin{rem}
Another way to compute this expansion is to use the theory
developed in \cite{trelatCOCV}, which states that the so-called
$L^\infty$-sector has a contact of order 2 with the abnormal
direction, and moreover gives an explicit formula to estimate the
contact.
\end{rem}

The previous method cannot be applied to study the contact of
branches $C_1$ and $D_2$ with the abnormal direction, because in
the phase plane of the pendulum these branches correspond
to a \it{global computation of return mapping}, and thus the
calculations cannot be localized near a saddle as previously.
Anyway inspecting carefully the system leads to the following~:

\begin{lem}
In the general gradated case of order $0$ the contact of branches
$C_1$ and $D_2$ with the abnormal direction is~:
$$Z=(\inv{6}+\rm{O}(r))X^3+\rm{o}(X^3)$$
\end{lem}

Note that contacts are still in the polynomial category.

\begin{proof}
We have~: $\dot{y}=v/\sqrt{c}$, hence $||y||_\infty=\rm{O}(r)$.
On the other part~: $\dot{x}=u/(1+\alpha y)$, and thus~:
$\dot{x}=u(1+\rm{O}(r))$. In the same way~:
$\dot{z}=u\f{y^2}{2}(1+\rm{O}(r))$. Then the result in the flat
case leads easily to the conclusion.
\end{proof}

\begin{rem}
From our previous study we can assert that minimizing controls
steering $0$ to points of $C_1$ (resp. $D_2$) are close to the
abnormal reference control in $L^2$-topology, but not in
$L^\infty$-topology. It is a crucial difference with the branch
$D_1$.
\end{rem}

Concerning the \it{transcendance} of this branch $D_1$, the
following fact was proved in \cite{trelatCRAS}~:

\begin{prop}
In the general gradated case of order $0$, if $\alpha\neq 0$ then
the branch $D_1$ is $C^\infty$ and is not subanalytic at
$x=-r,z=0$.
\end{prop}

\begin{cor}
In the general gradated case of order $0$, if the abnormal
minimizer is strict then the SR spheres with small radii are not
subanalytic.
\end{cor}

\begin{proof}
Let $A=(-r,0,0)$ denote the end-point of the abnormal trajectory.
We shall prove that $D_1$ is not subanalytic at $A$.
The method is the following. First of all
the Maximum Principle gives a parametrization of
minimizing trajectories steering $0$ to points of $D_1$.
Then we prove that the set of Lagrange multipliers associated to
these points (i.e. end-points of the corresponding adjoint
vectors) is not subanalytic. Finally we conclude using the fact
that, roughly speaking, these vectors coincide with the gradient
of the sub-Riemannian distance (where it is well-defined). These
facts are summarized in the following~:

\begin{lem}
To each point $q$ of $D_1$ is
associated a control $u$, and we denote by $(\psi(q),\psi^0(q))$
an associated Lagrange multiplier. Then we set~:
$${\cal L}=\left\{ \left({\psi_x(q)\over \sqrt{\psi_x(q)^2
+\psi_z(q)^2}}, {\psi_z(q)\over \sqrt{\psi_x(q)^2
+\psi_z(q)^2}}\right)\ /\ q\in{D_1}\right\}$$
where $\psi_x$ (resp. $\psi_z$) is the projection on the axis $x$
(resp. on the axis $z$) of the vector $\psi$.
If the set ${\cal L}$ is not subanalytic then the curve
${D_1}$ is not subanalytic.
\end{lem}

\begin{proof}[Proof of the Lemma]
Let $(q(\tau))_{0\leq \tau< 1}$ be a parametrization of
the curve
${D_1}$ such that $q(0)=A$. For each $\tau$ let
$u_\tau$ be a control such that $E(u_\tau)=q(\tau)$, and let
$(\psi_\tau, \psi^0_\tau)$
be an associated Lagrange multiplier, i.e.~:
$$\psi_\tau .dE(u_\tau)=-\psi^0_\tau dC(u_\tau)$$
Then~:
$\psi_\tau.{d\over d\tau}q(\tau)
=\psi_\tau.dE(u_\tau).{d\over d\tau}u_\tau
=-\psi^0_\tau dC(u_\tau).{d\over d\tau}u_\tau
=-\psi^0_\tau {d\over d\tau}C(u_\tau)$.
Moreover for each $\tau$ the point $q(\tau)$ belongs to the sphere
$S(0,r)$, hence $C(u_\tau)=r$, and thus~:
$\psi_\tau.{d\over d\tau}q(\tau)=0$.
Therefore in the plane $(y=0)$ the vectors of the set
${\cal L}$ are unitary normal vectors to the curve
$D_1$. Then the conclusion is immediate.
\end{proof}

With notations of Proof of Proposition \ref{propcontact},
we are now lead to study a family of
vector fields $(X_\varepsilon)$ depending on the parameter
$\varepsilon={1\over\sqrt{\lambda}}$, in the neighborhood of a
{\it saddle point} $u=v=0$.
For the section $\Sigma$ corresponding to $y=0$
we estimate the \it{return time},
i.e. the time needed to a trajectory starting from
$\Sigma$ to reach again $\Sigma$~; then we claim that this time
is $t=r$. This gives us
a relation between $\theta(r)$ and $\lambda$, thus between $p_x(r)$
and $p_z(r)$. Then one has to show that this
relation is not subanalytic. We proceed in the following way.
First of all recall that
$$1+\beta x(s)=\rm{e}^{-{\beta\over\sqrt{\lambda}}s }+O({1
\over\lambda }),\ \ y(s)=O({1\over\lambda})$$
We need a result which is independant of the parameter
$\varepsilon={1\over\sqrt{\lambda}}$. So it is no use trying to
write an analytic normal form, since the saddle may be resonant.
On the other hand $C^k$ normal forms (see \cite{R}) are not enough
because flat terms that we aim to exhibit disappear up to a
$O(\varepsilon^k)$. However near the saddle
separatrices of $X_\varepsilon$ are analytic in
$u,v,\varepsilon$, and actually there exists an
analytic change of coordinates $(u_1,v_1)=An(u,v)$ (in the sequel
$An(.)$ denotes an analytic germ at $0$) such that in these new
coordinates separatrices are $u_1=0,v_1=0$, and the system is~:
$$\dot{u}_1=\mu_1({1\over\sqrt{\lambda}})u_1
(1+o({1\over\sqrt{\lambda}})),
\
\dot{v}_1=\mu_2({1\over\sqrt{\lambda}})v_1
(1+o({1\over\sqrt{\lambda}}))$$
where
$\mu_1({1\over\sqrt{\lambda}}),\mu_2({1\over\sqrt{\lambda}})$
are the eigenvalues of the saddle~; in particular~:
$\mu_2({1\over\sqrt{\lambda}})=-1+{\beta\over 2\sqrt{\lambda}}+
O({1\over\lambda})$.
Moreover we have~:
$u=u_1+v_1+o({1\over\sqrt{\lambda}}),
v=u_1-v_1+o({1\over\sqrt{\lambda}})$, therefore the section is
$\Sigma~:
v_1=u_1+{\alpha\over\sqrt{\lambda}}+o({1\over\sqrt{\lambda}})$.
Let $s_f$ denote the parameter corresponding
to the return time, i.e. $(u_1(s_f),v_1(s_f))\in\Sigma$.
We have~: $v_1(0)={1\over\sqrt{\lambda}}
+o({1\over\sqrt{\lambda}})$. Then~:
\begin{equation}\label{eq(1)}
s_f=\int_0^{s_f}ds=\int_{v_1(0)}^{v_1(s_f)} {dv_1\over
\mu_2({1\over\sqrt{\lambda}})v_1(1+o({1\over\sqrt{\lambda}})) }
={1\over\mu_2}(1+o({1\over\sqrt{\lambda}}))\rm{ln }{v_1(s_f)\over
v_1(0)}
\end{equation}
On the other part~:
${dt\over ds}=\inv{\sqrt{\lambda}}(1+\alpha y)(1+\beta x+\gamma
y) ={1\over \sqrt{\lambda}}\rm{e}^{-{\beta\over\sqrt{\lambda}}s }
+O({1\over\lambda^{3/2}})$. Hence we get~:
$$r=\int_0^{s_f}{dt\over ds}ds
={1-\rm{e}^{-{\beta\over\sqrt{\lambda}}s_f} \over \beta r}
+O({1\over\lambda})$$
And thus~:
$s_f=-\sqrt{\lambda}{\rm{ln }(1-\beta r)\over\beta}+O({1\over
\sqrt{\lambda}}) $. Putting into (\ref{eq(1)}) we obtain finally~:
$$v_1(s_f)\sim {1\over\sqrt{1-\beta r}} {1\over\sqrt{\lambda}}
\rm{e}^{\sqrt{\lambda} {\rm{ln }(1-\beta r)\over\beta} }$$
In particular $v_1(s_f)$ is not an analytic function in
${1\over\sqrt{\lambda}}$.

We know that $v_1(s_f)=An(u(s_f),v(s_f)) \sim {u(s_f)-v(s_f)\over
2}$.
Moreover, on the section $\Sigma$, we have~:
$v(s_f)=-{\alpha\over\sqrt{\lambda}}\cos u(s_f)+{\beta \over
\sqrt{\lambda}}\sin u(s_f)$. Hence~:
$v_1(s_f)=An(u(s_f),{1\over\sqrt{\lambda}})={u(s_f)\over 2}+\cdots$.
From the Implicit Function Theorem in the analytic class we get~:
$u(s_f)=An(v_1(s_f),{1\over\sqrt{\lambda}})$.
Therefore $u(s_f)$ is not an analytic function in
${1\over\sqrt{\lambda}}$. So the set ${\cal L}$ is not
subanalytic, which ends the proof.
\end{proof}

Unfortunately we have no general result similar to Corollary
\ref{corLE}. We think that the log-exp category is not wide
enough in the non integrable case. Indeed
due to the dissipation phenomenon observed in the
pendulum representation if $\beta\neq 0$ we cannot expect to keep
the analytic properties required in the definition of log-exp
functions. Moreover in the
phase plane of the pendulum, the foliation is not a
priori integrable in the analytic category for any value of the
parameter. We can observe that if we fix $\lambda$
to $1$ and evaluate the Poincar\'e-Dulac mapping, it is pfaffian
if and only if $X_\varepsilon $ ($\lambda=1$) is
$C^\omega$-integrable (see \cite{MM,MR}).

Hence we conjecture~:

\begin{conj}
If $\beta\neq 0$ then SR Martinet spheres are not log-exp, even
not pfaffian.
\end{conj}

Hence we should try to extend the log-exp category to a wider
category in which analyticity would be replaced with some
asymptotic properties. That is why we should be interested in
\it{Il'Yashenko's class of functions} (see \cite{Il}). Indeed
several years ago the Dulac's problem of finiteness of limit
cycles was solved independantly by Ecalle and Il'Yashenko~; in
his proof, Il'Yashenko introduces a very wide class of
\it{non-oscillating functions} to describe Poincar\'e return
mappings. Actually he needs to expand real germs of
functions into terms having not only the order of $x^n$ but also
into flat terms. His category is contructed by recurrence as
follows. Let $M_0$ be a class of functions that can be
expanded in an unique way into an ordinary Dulac's series, i.e.
series of the form~:
$$\Sigma_0 = cx^{\nu_0}+\sum_{i=1}^{\infty}
P_i(\rm{ln }x)x^{\nu_i}$$
where $c>0$, the $P_i$'s are polynomial and $(\nu_i)$ is an
increasing sequence of positive numbers going to infinity.
We only sketch the first step of the recurrence.
By definition germs of the class $M_1$ can be expanded in series
containing flat exponential terms whose coefficients belong to
$M_0$. For instance a super-accurate series of some germ $f$ may
be of the form~:
$$\Sigma_1=a_0(x)+\sum_{i=1}^\infty
a_i(x)\rm{e}^{-\f{\nu_i}{x}}$$
where $a_i\in M_0$ and $(\nu_i)$ is an increasing sequence of
positive numbers tending to infinity. It is a generalization of
ordinary series in so much as the usual Dulac's series of $f$ is
$a_0(x)$.

The complete definition of super-accurate series then goes by
recurrence (see \cite{Il}).
Their interest is all in the fact that the application $f\mapsto
\hat{f}$ is \it{one-to-one}, where $\hat{f}$ is the
super-accurate series associated to $f$.

The relation with our problem is the following.
We deal actually with Poincar\'e return mappings
in the phase plane of the pendulum, and their study is crucial to
estimate the spheres. The difference is that our pendulum depends
on parameters (namely $\lambda$)~; hence our problem is related
to Dulac's problem with parameters, i.e. the Hilbert's 16th
problem.

Hence we should try to construct a category of functions similar
to the one introduced by Il'Yashenko, but with parameters. In any
case it is a possible way to try to solve the problem of
transcendence in SR geometry.

\begin{conj}
In the general non integrable case SR spheres belong to some
extended Il'Yashenko's category.
\end{conj}


\subsubsection{Conjecture about the cut-locus~:
the Martinet sphere in the Liu-Sussmann example}
We shall construct the cut-locus in the Liu-Sussmann example
\cite{LS}, the reasoning being generalizable to compute the
generic SR-Martinet sphere. The model is the following~:
$$D=\rm{Ker }\omega\ , \ \ \omega=(1+\varepsilon y)dz
-\f{y^2}{2}dx \ , \ \ g=\f{dx^2}{(1+\varepsilon y)^2}+dy^2$$

The model is non generic because it is conservative~; moreover
the Lie algebra generated by the orthonormal frame is nilpotent.
In the cylindric coordinates the geodesics equations are~:
\begin{equation*}
\begin{split}
\dot{x}&=(1+\varepsilon y)\cos \theta\ , \ \ \dot{y}=\sin \theta\
,\ \ \dot{z}=\f{y^2}{2}\cos \theta \\
\dot{\theta}&=-(p_x\varepsilon +p_zy)\ ,\ \ P_3=\lambda
\end{split}
\end{equation*}

\no and the angle evolution is the pendulum~: $\ddot{\theta}+\lambda
\sin \theta =0$ if $\lambda \neq 0$. Using the symmetry~:
$(x,y,z) \longmapsto (-x,y,-z)$ we can assume $\lambda >0$. The
abnormal geodesic starting from $0$~: $t\longmapsto (\pm t,0,0)$,
is strict if and only if $\varepsilon \neq 0$. We may assume
$\varepsilon \leq 0$. Introducing $s=t\rala$, and denoting by $'$
the derivative with respect to $s$, the pendulum is normalized
to~:
$\ddot{\theta}+\sin \theta =0$. The constraint $y=0$ defines the
section $\Sigma : \dot{\theta}=-p_x\varepsilon $, which can be
written~: $\theta '=-\f{\varepsilon \cos \theta}{\rala}$.

The geodesics corresponding to $\lambda =0$ are globally optimal
if the length $r$ is small enough. They divide the sphere
$S(0,r)$ into two hemispheres and we compute the cut-locus in the
northern hemisphere $(\lambda>0)$.


\setlength{\unitlength}{0.5mm}
\begin{figure}[h]
\begin{center} 

\begin{picture}(180,100)
\thinlines
\drawvector{10.0}{48.0}{78.0}{1}{0}
\drawvector{48.0}{10.0}{76.0}{0}{1}
\drawdashline{12.0}{12.0}{12.0}{84.0}
\drawdashline{82.0}{12.0}{82.0}{84.0}
\drawarc{47.0}{62.0}{75.37}{0.38}{2.75}
\drawarc{47.0}{34.0}{75.37}{3.51}{5.9}
\drawlefttext{50.0}{84.0}{$\theta '$}
\drawcenteredtext{90.0}{52.0}{$\theta$}
\drawlefttext{82.8}{44.7}{$+\pi$}
\thicklines
\drawrighttext{11.25}{44.93}{$-\pi$}
\drawlefttext{48.81}{44.7}{$0$}
\thinlines
\drawvector{46.34}{71.68}{1.35}{1}{0}
\drawvector{47.9}{24.26}{1.35}{-1}{0}
\drawdashline{30.0}{76.0}{18.0}{48.0}
\drawvector{18.67}{49.43}{0.67}{-1}{-2}
\drawellipse{47.0}{48.0}{46.0}{24.0}{}
\drawvector{46.79}{60.0}{1.11}{1}{0}
\drawvector{47.9}{36.15}{1.12}{-1}{0}
\drawcircle{138.0}{48.0}{64.0}{}
\drawpath{106.0}{48.0}{170.0}{48.0}
\drawarc{104.0}{48.0}{4.0}{4.71}{1.57}
\drawarc{172.0}{48.0}{4.0}{1.57}{4.71}
\drawcenteredtext{138.0}{64.0}{$y>0$}
\drawcenteredtext{138.0}{32.0}{$y<0$}
\drawcenteredtext{157.72}{51.45}{$y=0$}
\drawrighttext{115.19}{73.93}{$\lambda =0$}
\drawcenteredtext{30.0}{78.0}{$\sigma$}
\end{picture}

\end{center}
\caption{}\label{cutf1}
\end{figure}

If $\varepsilon =0$, the SR-sphere is the Martinet flat sphere.
The abnormal line is not strict and cuts the equator $\lambda =0$
in two points. The cut-locus is the plane $y=0$
minus the abnormal line, in which, due to
the symmetry $(x,y,z)\longmapsto (x,-y,z)$, two normal geodesics
are intersecting with the same length. It is represented on
Fig. \ref{cutf1}.

When $\varepsilon \neq 0$, the section and the pendulum are
represented on Fig. \ref{cutf2}.


\unitlength=.25mm
\makeatletter
\def\shade{\@ifnextchar[{\shade@special}{\@killglue\special{sh}\ignorespaces}}
\def\shade@special[#1]{\@killglue\special{sh #1}\ignorespaces}
\makeatother

\begin{figure}[h]
\begin{center}

\begin{picture}(292,187)(200,-5)
\thinlines
\typeout{\space\space\space eepic-ture exported by 'qfig'.}
\font\FonttenBI=cmbxti10\relax
\font\FonttwlBI=cmbxti10 scaled \magstep1\relax
\path (200,82)(434,82)
\path (206,9)(206,164)
\path (394,168)(394,6)
\path (302,0)(302,177)
\path (426.482,84.736)(434,82)(426.482,79.264)
\path (299.264,169.482)(302,177)(304.736,169.482)
\put(300,30.082){\arc{214.77}{3.646}{5.779}}
\put(300,136.95){\arc{217.766}{.529}{2.613}}
\path (328,140)(334,133)(325,130)
\path (322,36)(316,30)(324,26)
\path (206,64)(206.92,63.955)(207.88,63.92)(208.88,63.895)(209.92,63.88)
(211,63.875)(212.12,63.88)(213.28,63.895)(214.48,63.92)(215.72,63.955)
(217,64)(218.41,63.83)(219.84,63.72)(221.29,63.67)(222.76,63.68)
(224.25,63.75)(225.76,63.88)(227.29,64.07)(228.84,64.32)(230.41,64.63)
(232,65)(233.835,65.52)(235.64,66.08)(237.415,66.68)(239.16,67.32)
(240.875,68)(242.56,68.72)(244.215,69.48)(245.84,70.28)(247.435,71.12)
(249,72)(250.625,73.055)(252.2,74.12)(253.725,75.195)(255.2,76.28)
(256.625,77.375)(258,78.48)(259.325,79.595)(260.6,80.72)(261.825,81.855)
(263,83)(263.675,84.38)(264.4,85.72)(265.175,87.02)(266,88.28)
(266.875,89.5)(267.8,90.68)(268.775,91.82)(269.8,92.92)(270.875,93.98)
(272,95)(273.535,96.07)(275.04,97.08)(276.515,98.03)(277.96,98.92)
(279.375,99.75)(280.76,100.52)(282.115,101.23)(283.44,101.88)(284.735,102.47)
(286,103)(287.37,103.29)(288.68,103.56)(289.93,103.81)(291.12,104.04)
(292.25,104.25)(293.32,104.44)(294.33,104.61)(295.28,104.76)(296.17,104.89)
(297,105)(297.365,105)(297.76,105)(298.185,105)(298.64,105)
(299.125,105)(299.64,105)(300.185,105)(300.76,105)(301.365,105)
(302,105)(302.845,105.18)(303.68,105.32)(304.505,105.42)(305.32,105.48)
(306.125,105.5)(306.92,105.48)(307.705,105.42)(308.48,105.32)(309.245,105.18)
(310,105)(310.52,104.78)(311.08,104.52)(311.68,104.22)(312.32,103.88)
(313,103.5)(313.72,103.08)(314.48,102.62)(315.28,102.12)(316.12,101.58)
(317,101)(318.19,100.29)(319.36,99.56)(320.51,98.81)(321.64,98.04)
(322.75,97.25)(323.84,96.44)(324.91,95.61)(325.96,94.76)(326.99,93.89)
(328,93)(328.855,92)(329.72,91)(330.595,90)(331.48,89)
(332.375,88)(333.28,87)(334.195,86)(335.12,85)(336.055,84)
(337,83)(337.82,81.91)(338.68,80.84)(339.58,79.79)(340.52,78.76)
(341.5,77.75)(342.52,76.76)(343.58,75.79)(344.68,74.84)(345.82,73.91)
(347,73)(348.355,71.93)(349.72,70.92)(351.095,69.97)(352.48,69.08)
(353.875,68.25)(355.28,67.48)(356.695,66.77)(358.12,66.12)(359.555,65.53)
(361,65)(362.365,64.665)(363.76,64.36)(365.185,64.085)(366.64,63.84)
(368.125,63.625)(369.64,63.44)(371.185,63.285)(372.76,63.16)(374.365,63.065)
(376,63)(377.665,62.965)(379.36,62.96)(381.085,62.985)(382.84,63.04)
(384.625,63.125)(386.44,63.24)(388.285,63.385)(390.16,63.56)(392.065,63.765)
(394,64)
\put(329,16){{\rm\rm {$\Sigma$}}}
\put(322,99){{\rm\rm {$\sigma$}}}
\put(420,63){{\rm\rm {$\theta$}}}
\put(310,164){{\rm\rm {$\theta '$}}}
\path (299,107)(305,102)
\path (299,102)(305,107)
\path (391,66)(397,61)
\path (391,61)(397,66)
\put(399,55){{\rm\rm {$m$}}}
\put(288,107){{\rm\rm {$M$}}}
\end{picture}

\end{center}
\caption{}\label{cutf2}
\end{figure}

\no $\sigma$ is given by $\theta '=-\f{\varepsilon \cos
\theta}{\rala}$. When $\lambda \tto +\infty$, the section tends
to $\theta '=0$, and when $\lambda \tto 0$ the points $M$ and $m$
tend to $\infty$.

We can easily construct the cut-locus on the small sphere using
the following conjectures~:
\begin{itemize}
\item Only the geodesics where the section $\sigma$ is in the
configuration of Fig. \ref{cutf2} have cut-points (this is
justified by the fact that when $\lambda =0$ the geodesics are
globally optimal if the length is small enough).
\item A separatrix has no cut-point.
\end{itemize}

The trace of the sphere with the Martinet plane has been computed
in \cite{BC}. An important property is that the rotating
trajectories of the pendulum near the separatrix have a cut-point
located in the plane $y=0$, corresponding to its second
intersection with $y=0$. We represent on Fig. \ref{cutf3} the
construction of the cut-locus.


\unitlength=.25mm
\makeatletter
\def\shade{\@ifnextchar[{\shade@special}{\@killglue\special{sh}\ignorespaces}}
\def\shade@special[#1]{\@killglue\special{sh #1}\ignorespaces}
\makeatother

\begin{figure}[h]
\begin{center}

\begin{picture}(572,195)(57,-5)
\thinlines
\typeout{\space\space\space eepic-ture exported by 'qfig'.}
\font\FonttenBI=cmbxti10\relax
\font\FonttwlBI=cmbxti10 scaled \magstep1\relax
\path (72,99)(306,99)
\path (78,26)(78,181)
\path (266,185)(266,23)
\path (298.483,101.736)(306,99)(298.483,96.264)
\thicklines
\put(172,47){\arc{214.849}{3.647}{5.778}}
\thinlines
\thicklines
\put(172,154){\arc{217.816}{.529}{2.612}}
\thinlines
\path (200,157)(206,150)(197,147)
\path (194,53)(188,47)(196,43)
\path (78,81)(78.92,80.955)(79.88,80.92)(80.88,80.895)(81.92,80.88)
(83,80.875)(84.12,80.88)(85.28,80.895)(86.48,80.92)(87.72,80.955)
(89,81)(90.41,80.83)(91.84,80.72)(93.29,80.67)(94.76,80.68)
(96.25,80.75)(97.76,80.88)(99.29,81.07)(100.84,81.32)(102.41,81.63)
(104,82)(105.835,82.52)(107.64,83.08)(109.415,83.68)(111.16,84.32)
(112.875,85)(114.56,85.72)(116.215,86.48)(117.84,87.28)(119.435,88.12)
(121,89)(122.625,90.055)(124.2,91.12)(125.725,92.195)(127.2,93.28)
(128.625,94.375)(130,95.48)(131.325,96.595)(132.6,97.72)(133.825,98.855)
(135,100)(135.675,101.38)(136.4,102.72)(137.175,104.02)(138,105.28)
(138.875,106.5)(139.8,107.68)(140.775,108.82)(141.8,109.92)(142.875,110.98)
(144,112)(145.535,113.07)(147.04,114.08)(148.515,115.03)(149.96,115.92)
(151.375,116.75)(152.76,117.52)(154.115,118.23)(155.44,118.88)(156.735,119.47)
(158,120)(159.37,120.29)(160.68,120.56)(161.93,120.81)(163.12,121.04)
(164.25,121.25)(165.32,121.44)(166.33,121.61)(167.28,121.76)(168.17,121.89)
(169,122)(169.365,122)(169.76,122)(170.185,122)(170.64,122)
(171.125,122)(171.64,122)(172.185,122)(172.76,122)(173.365,122)
(174,122)(174.845,122.18)(175.68,122.32)(176.505,122.42)(177.32,122.48)
(178.125,122.5)(178.92,122.48)(179.705,122.42)(180.48,122.32)(181.245,122.18)
(182,122)(182.52,121.78)(183.08,121.52)(183.68,121.22)(184.32,120.88)
(185,120.5)(185.72,120.08)(186.48,119.62)(187.28,119.12)(188.12,118.58)
(189,118)(190.19,117.29)(191.36,116.56)(192.51,115.81)(193.64,115.04)
(194.75,114.25)(195.84,113.44)(196.91,112.61)(197.96,111.76)(198.99,110.89)
(200,110)(200.855,109)(201.72,108)(202.595,107)(203.48,106)
(204.375,105)(205.28,104)(206.195,103)(207.12,102)(208.055,101)
(209,100)(209.82,98.91)(210.68,97.84)(211.58,96.79)(212.52,95.76)
(213.5,94.75)(214.52,93.76)(215.58,92.79)(216.68,91.84)(217.82,90.91)
(219,90)(220.355,88.93)(221.72,87.92)(223.095,86.97)(224.48,86.08)
(225.875,85.25)(227.28,84.48)(228.695,83.77)(230.12,83.12)(231.555,82.53)
(233,82)(234.365,81.665)(235.76,81.36)(237.185,81.085)(238.64,80.84)
(240.125,80.625)(241.64,80.44)(243.185,80.285)(244.76,80.16)(246.365,80.065)
(248,80)(249.665,79.965)(251.36,79.96)(253.085,79.985)(254.84,80.04)
(256.625,80.125)(258.44,80.24)(260.285,80.385)(262.16,80.56)(264.065,80.765)
(266,81)
\path (171,124)(177,119)
\path (171,119)(177,124)
\path (263,83)(269,78)
\path (263,78)(269,83)
\path (78,81)(67,83)(57,88)
\put(71.676,77.147){\arc{27.739}{3.834}{6.002}}
\put(172.557,134.812){\arc{146.115}{.75}{2.751}}
\path (198,106)(206,103)(206,111)
\put(200,114){{\rm\rm {$L_C$}}}
\put(169.5,150.744){\arc{223.49}{.686}{2.456}}
\dottedline{3}(237,15)(238.1,18.055)(239.2,21.12)(240.3,24.195)(241.4,27.28)
(242.5,30.375)(243.6,33.48)(244.7,36.595)(245.8,39.72)(246.9,42.855)
(248,46)(249.1,49.155)(250.2,52.32)(251.3,55.495)(252.4,58.68)
(253.5,61.875)(254.6,65.08)(255.7,68.295)(256.8,71.52)(257.9,74.755)
(259,78)
\path (254.027,71.734)(259,78)(259.219,70.003)
\put(224,-1){{\rm\rm {$L_D$}}}
\put(472,101){\circle{139.714}}
\dottedline{3}(403,99)(541,99)
\path (460,99)(427,99)
\path (427,99)(427.73,98.21)(428.52,97.44)(429.37,96.69)(430.28,95.96)
(431.25,95.25)(432.28,94.56)(433.37,93.89)(434.52,93.24)(435.73,92.61)
(437,92)(438.195,91.23)(439.48,90.52)(440.855,89.87)(442.32,89.28)
(443.875,88.75)(445.52,88.28)(447.255,87.87)(449.08,87.52)(450.995,87.23)
(453,87)(455.905,86.785)(458.72,86.64)(461.445,86.565)(464.08,86.56)
(466.625,86.625)(469.08,86.76)(471.445,86.965)(473.72,87.24)(475.905,87.585)
(478,88)(479.825,88.71)(481.6,89.44)(483.325,90.19)(485,90.96)
(486.625,91.75)(488.2,92.56)(489.725,93.39)(491.2,94.24)(492.625,95.11)
(494,96)(494.875,97.135)(495.8,98.24)(496.775,99.315)(497.8,100.36)
(498.875,101.375)(500,102.36)(501.175,103.315)(502.4,104.24)(503.675,105.135)
(505,106)(506.375,106.835)(507.8,107.64)(509.275,108.415)(510.8,109.16)
(512.375,109.875)(514,110.56)(515.675,111.215)(517.4,111.84)(519.175,112.435)
(521,113)
\put(452,43){{\rm\rm {$y<0$}}}
\put(445,143){{\rm\rm {$y>0$}}}
\put(533,53){{\rm\rm {$\lambda=0$}}}
\put(417,106){{\rm\rm {$A$}}}
\put(494,114){{\rm\rm {$L_C$}}}
\put(444,108){{\rm\rm {$L_D$}}}
\path (447,103)(439,99)(447,96)
\path (465,91)(457,87)(465,83)
\path (458,102)(464,97)
\path (458,97)(465,102)
\path (520,116)(526,111)
\path (520,111)(527,116)
\path (102,105)(108,109)
\path (102,109)(106,105)
\thicklines
\path (422,102)(426,100)(422,97)
\thinlines
\end{picture}

\end{center}
\caption{}\label{cutf3}
\end{figure}

The cut-locus has two branches $L_C$ and $L_D$ corresponding
respectively to oscillating and rotating trajectories. They
ramify on the abnormal direction $A$ which is not contained in
the cut-locus. The extremities of the branches $L_C$ and $L_D$
are conjugate points corresponding respectively to $\theta
(0)=\pi$ and $\theta (0)=0$. The branch $L_C$ has only one
intersection with $y=0$ which corresponds approximatively to
$\theta (0)=\f{\pi}{2}$. \it{The cut-locus is not subanalytic at
$A$ but belongs to the log-exp category}, and thus
from \cite{DriesMiller}~:

\begin{lem}
The cut-locus is $C^\infty$-stratifiable.
\end{lem}

To generalize this analysis we must observe the following. The
respective positions of the branches $C_1, D_2$ are given in
Section \ref{section4.2.13}.
Here the curve $D_2$ is above and hence the
rotating trajectories are optimal up to the second intersection.
Also the integrability of the geodesic flow is not crucial and in
general the branch $L_D$ is not contained in the plane $y=0$. We
make the following conjecture.

\begin{conj}
\begin{enumerate}
\item
In the Martinet case the cut-locus is $C^1$-stratifiable.
\item
In the generic case the cut-locus has two branches in the
northern hemisphere ramifying at the end-point of the abnormal
direction.
\end{enumerate}
\end{conj}


\section{Some extensions of Martinet SR geometry and microlocal
analysis of the singularity of the SR sphere in the abnormal
direction}
\subsection{Non properness and Tangency Theorem}
\label{sectiontangent}
This analysis is based on the sub-Riemannian Martinet case, where
it was shown in the previous Section
that the exponential mapping is not
proper and that in the generic case the sphere is
tangent to the abnormal direction.
This fact is actually general and we have the following results
(see \cite{trelatJDCS}).\\

Consider a smooth
sub-Riemannian structure $(M,\Delta,g)$ where $M$ is a
Riemannian $n$-dimensional manifold, $n\geq 3$, $\Delta$ is
a rank $m$ distribution on $M$, and $g$ is a metric on $\Delta$.
Let $q_0\in M$~; our point of view is local and we can assume
that $M=\R^n$ and $q_0=0$. Suppose there exists a strict (in the
sense of definition \ref{defistrict})
abnormal trajectory
$\gamma$ passing through $0$. Up to
reparametrizing we can assume that
$\Delta=\rm{Span }\{F_1,\ldots,F_m\}$ where the system of $F_i$'s
is $g$-orthonormal. Then the sub-Riemannian problem is equivalent
to the \it{time-optimal problem} for the system~:
\begin{equation}\label{srsystem}
\dot{q}=\sum_{i=1}^mu_iF_i(q),\ q(0)=0
\end{equation}
where the controls satisfy the \it{constraint}
$\sum_{i=1}^mu_i^2\leq 1$.
Suppose further that $\gamma$ is associated to an unique strictly
abnormal control. Then~:

\begin{thm} \label{nonproper}
The exponential mapping is not proper near $\gamma$.
\end{thm}

\begin{proof}
Set $A=\gamma(r), r>0$.
Let $(A_n)$ be a sequence of end-points of minimizing
normal geodesics $q_n$ converging to $A$.
To each geodesic $q_n$ is associated a control $u_n$ and an
adjoint vector $(p_n,p^0_n)$.
As $q_n$ is normal we may suppose that
$p^0_n=-\inv{2}$. Let $\psi_n$ the end-point of the adjoint
vector $p_n$. Then if $E$ denotes the end-point mapping and $C$
is the cost (here the cost is quadratic in the control),
we have the following \it{Lagrange multiplier
equality}~:
$$\psi_n.dE(u_n)=\inv{2}dC(u_n)$$
If the sequence $\psi_n$ were bounded then up to a subsequence it
would converge to some $\psi\in\R^n$. Now since the $u_n$ are
minimizing the sequence $(u_n)$ is bounded in $L^2$, hence up to
a subsequence it converges weakly to some $u\in L^2$. Using the
regularity properties of the end-point mapping (see
\cite{trelatJDCS}), we can pass through the limit in the previous
equality and we get~:
$$\psi.dE(u)=\inv{2}dC(u)$$
and on the other part~: $A=E(u)$. It is not difficult to see that
$u$ has to be minimizing, and then we get a contradiction with
the fact that $\gamma$ is strict.
\end{proof}

\begin{rem}
Conversely if the exponential mapping
is not proper then actually there exists an abnormal
minimizer. This shows the interaction between abnormal and normal
minimizers. In a sense normal extremals recognize abnormal
extremals. This phenomenon of non-properness is characteristic
for abnormality.
\end{rem}

This non-properness is actually responsible for a phenomenon of
\it{tangency} described in the following Theorem (see
\cite{trelatCOCV} for a more general statement):

\begin{thm}\label{thmtangent}
Consider the SR system (\ref{srsystem}) and suppose there exists
a minimizing geodesic $\gamma$ associated to an unique strictly
abnormal control $u$. Let $A\in S(0,r)$ be the end-point of
$\gamma$.
Assume $(\sigma(\tau))_{0<\tau\leq 1}$ is a $C^1$ curve on
$S(0,r)$ such that $\displaystyle
{\lim_{\tau\rightarrow 0}\sigma(\tau)=A}$.
Then~:
$\displaystyle
{\lim_{\tau\rightarrow 0}\sigma'(\tau) \in \rm{Im
}dE(u) }$.
\end{thm}

\begin{figure}[h]
\epsfxsize=9cm
\epsfysize=7cm
\centerline{\epsffile{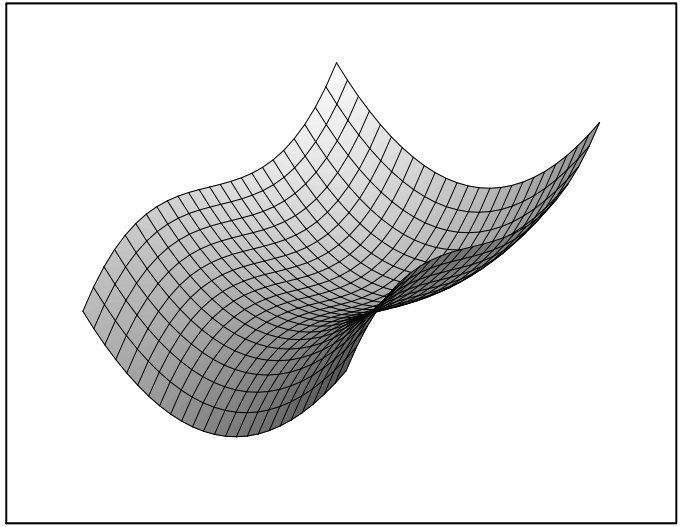}}
\caption{Tangency phenomenon}\label{figtang}
\end{figure}

In particular if $S(0,r)$ is $C^1$-stratifiable near $A$ then the
strata of $S(0,r)$ are tangent at $A$ to the hyperplane $\rm{Im
}dE(u)$ (see Fig. \ref{figtang}). Moreover if $B$ is a
$C^1$-branch of the cut-locus ramifying at $A$ then
$B$ is tangent at $A$ to this hyperplane.


\begin{proof}
For each $\tau$ the point $\sigma(\tau)$ is the
end-point of a minimizing geodesic, and we denote by $(p_\tau,
p^0_\tau)$
(resp. $u_\tau$) an associated adjoint vector (resp. an
associated control). Let $(\psi_\tau,\psi^0_\tau)$ be the
end-point of this adjoint vector. We may suppose that it is
unitary in $\R^n\times\R$. We have~:
$$\psi_\tau.dE(u_\tau)=-\psi^0_\tau dC(u_\tau)$$
Using the same reasoning as in
the Proof of Theorem \ref{nonproper} we get that
$||\psi^0_\tau||\rightarrow 0$ as $\tau\rightarrow 0$.
To conclude it suffices to show that
$\psi_\tau$ is normal to the curve $\sigma(\tau)$. Indeed the
previous equality implies~:
$$\psi_\tau.dE(u_\tau) .\f{du_\tau}{d\tau}=-\psi^0_\tau
dC(u_\tau) .\f{du_\tau}{d\tau}$$
But $C(u_\tau)$ is constant (equal to $r$) and
$E(u_\tau)=\sigma(u_\tau)$, and thus~:
$$\psi_\tau.\f{d\sigma(\tau)}{d\tau}=0$$
which ends the proof.
\end{proof}


\subsection{The tangential case}
\subsubsection{Preliminaries}
In this Subsection we shall make a brief analysis of the
so-called tangential case. According to Section 2.2.2 the
distribution $D=\rm{Ker } \omega$ can be reduced \cite{Zh}
to one of the normal forms~:

\begin{itemize}
\item \underline{elliptic case :} 
$\omega_e=dy-(\varepsilon xy+\f{x^3}{3}+xz^2+mx^3z^2)dz$
\item \underline{hyperbolic case :}
$\omega_h=dy-(\varepsilon xy+x^2z+mx^3z^2)dz$
\end{itemize}

where $\varepsilon=\pm 1$. The parameter $\varepsilon$ is a
deformation
parameter whose introduction will be justified later.

A general metric $g$ is then defined by :
$a(q)dx^2+2b(q)dxdz+c(q)dz^2$ where $a,b,c$ can be taken as
constant in the nilpotent approximation of order $-1$. Our study
is far to be complete and we shall describe briefly the case
$g=dx^2+dz^2$.

The general case of order $-1$ depends on a parameter $\lambda$
and corresponds to a 6-dimensional nilpotent Lie algebra. It
contains both elliptic and hyperbolic cases. It is the Lie
algebra generated by $F_1$, $F_2$ with the following Lie brackets
relations~:
$$F_3=[F_1,F_2]\ , \ \ F_4=[F_3,F_1]$$
$$F_5=[F_3,F_2]\ , \ \ F_6=[F_4,F_1]$$
$$[F_5,F_2]=\lambda F_6 \ , $$
\noindent and all other Lie brackets are $0$.

Introducing $P_i=<p,F_i>$ the geodesic equations are given by~:
\begin{equation*}
\begin{split}
\dot{P_1}&=P_3P_2\ , \ \ \dot{P_2}=-P_3P_1\ , \ \ \dot{P_3}
=P_4P_1+P_5P_2\ , \\
\dot{P_4}&=P_6P_1\ , \ \ \dot{P_5}=\lambda P_6P_2
\end{split}
\end{equation*}

\no and $P_6$ is a \it{Casimir first integral}. The value
$\lambda=0$ represents the bifurcation between the two cases.


\subsubsection{Abnormal geodesics}
\no{\bf Elliptic case}
The abnormal geodesics are contained in the Martinet surface~:
$$\varepsilon y+x^2+z^2+3mx^2z^2=0$$
and are solutions of the equations~:
\begin{equation*}
\begin{split}
\dot{x}&=(2z+6mx^2z)-\varepsilon (\f{2x^3}{3}+2mx^3z^2) \\
\dot{z}&=-(2x+6mx^2z) 
\end{split}
\end{equation*}

\begin{figure}[h]
\epsfxsize=7cm
\epsfysize=7cm
\centerline{\epsffile{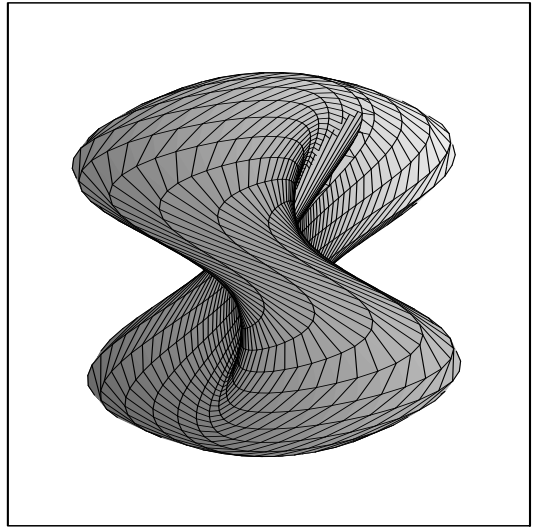}}
\caption{Elliptic case}\label{elliptic}
\end{figure}

From \cite{Zh}, the singularity $x=z=0$ is a weak focus and a spiral
passing through $0$ is with infinite length. Since any minimizer
is smooth \it{no piece of abnormal geodesic is a minimizer} when
computing the distance to $0$. Using the general result of
\cite{Ag},
the sphere of small radius is the image by the exponential
mapping of a compact set and it is subanalytic. This is also
clearly shown by numerical simulations and the sphere is
represented on Fig. \ref{elliptic}.


By taking $\varepsilon =0$, the Martinet surface becomes~:
$x^2+z^2+3mx^2z^2=0$ and reduces near $0$ to~: $x=z=0$. Hence the
spiral disappears. Since the weight of $x,z$ is one and the
weight of $y$ is four, it corresponds to the nilpotent
approximation of order $-1$ where $m$ is $0$.
\\

\no{\bf Hyperbolic case}
The Martinet surface is given by the equation~:
$$\varepsilon y+2xz+3mx^2z^2=0$$
\noindent and the abnormal geodesics are solutions of~:
\begin{equation*}
\begin{split}
\dot{x}&=2x-x^2z(\varepsilon -6m)-2mx^3z^2 \\
\dot{z}&=-(2z+6mx^2z^2)
\end{split}
\end{equation*}

\begin{figure}[h]
\epsfxsize=7cm
\epsfysize=7cm
\centerline{\epsffile{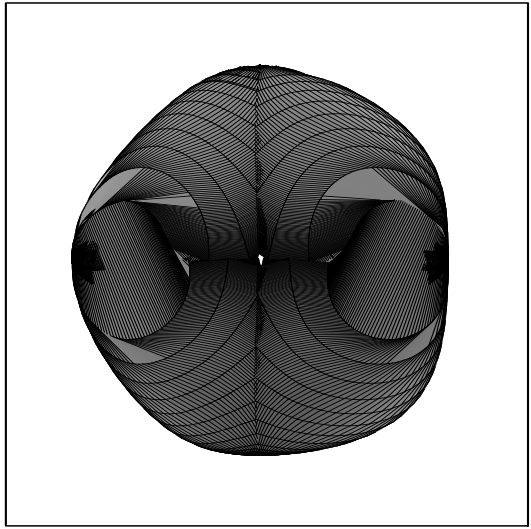}}
\caption{Hyperbolic case}\label{trou}
\end{figure}

\no and the singularity at $x=z=0$ is a saddle point. \it{The two
lines $x=0$ and $z=0$ are optimal for the metric $dx^2+dz^2$}.
Hence they play a role when computing the distance to $0$.
Numerical simulations show that \it{the sphere is not the image
of a compact set}. This can be seen on Fig. \ref{trou} because the
sphere cannot be numerically represented in the abnormal
direction (there is a hole).
It is similar to the situation encountered in
the Martinet case. The sphere is pinched in both abnormal
directions.

The nilpotent approximation of order $-1$ is obtained by taking
$\varepsilon =0$ and $m=0$. The Martinet surface becomes~: $xz=0$
and the two lines $x=0$ and $z=0$ remain abnormal geodesics.


\subsubsection{Normal geodesics}
\no{\bf Elliptic case}
We take the frame~:
$$F_1=\f{\partial}{\partial x}\ ,\ \ F_2=\f{\partial}{\partial
z}+(\varepsilon xy+\f{x^3}{3}+xz^2+mx^3z^2)
\f{\partial}{\partial y}\ ,\ \
F_3=\f{\partial}{\partial y}$$
\no and we introduce~: $P_i=<p,F_i(q)>$. The geodesics equations
are~:
\begin{equation*}
\begin{split}
\dot{x}&=P_1 \\
\dot{y}&=P_2(\varepsilon xy+\f{x^3}{3}+xz^2+mx^3z^2) \\
\dot{z}&=P_2 \\
\dot{P_1}&=-(\varepsilon y+x^2+z^2+3mx^2z^2)P_2P_3 \\
\dot{P_2}&=(\varepsilon y+x^2+z^2+3mx^2z^2)P_1P_3 \\
\dot{P_3}&=-\varepsilon xP_3 
\end{split}
\end{equation*}

\no and they can be truncated at order $-1$ by making
$\varepsilon =m=0$. In this case $P_3$ is a first integral and we
can set~: $P_3=\lambda$. Moreover if we introduce~: $P_1=\sin
\theta$, $P_2=\cos \theta$, the equations become~:
$$\dot{x}=\sin \theta\ , \ \ \dot{z}=\cos \theta\ ,\ \
\dot{y}=\cos \theta \ (\f{x^3}{3}+xz^2)$$
$$\dot{\theta}=-(x^2+z^2)P_3\ , \ \ P_3=\lambda $$

\no They can be projected onto the space $(x,y,\theta )$ and the
foliation is defined by~:
$$\dot{x}=\sin \theta\ , \ \ \dot{z}=\cos \theta\ ,\ \
\dot{\theta}=-(x^2+z^2)\lambda $$

It is not Liouville-integrable but the equations can be
integrated by quadratures, see \cite{Pel}. Using polar coordinates~:
$x=r\cos \psi, z=r\sin \psi$, it becomes~:
$$\dot{x}=\f{\cos (\theta +\psi)}{r} \ , \ \ \dot{r}=\sin
(\theta+\psi)\ ,\ \ \dot{\theta}=-\lambda r^2$$
The important property is the following~:

\begin{lem}
The sign of $\dot{\theta}$ is constant and reparametrizing the
equation can be rewritten~: $\dot{\theta}=1$ as in the contact
case.
\end{lem}

\paragraph{Numerical simulations \\ }
The geodesics equations can be integrated numerically. The
projections in the plane $(x,z)$ of the geodesics starting from
$0$ are flowers with three petals (they are circles in the
Heisenberg case), see Fig. \ref{f19}.


\unitlength=.25mm
\makeatletter
\def\shade{\@ifnextchar[{\shade@special}{\@killglue\special{sh}\ignorespaces}}
\def\shade@special[#1]{\@killglue\special{sh #1}\ignorespaces}
\makeatother

\begin{figure}[h]
\begin{center}

\begin{picture}(378,120)(142,-5)
\thinlines
\typeout{\space\space\space eepic-ture exported by 'qfig'.}
\font\FonttenBI=cmbxti10\relax
\font\FonttwlBI=cmbxti10 scaled \magstep1\relax
\path (142,47)(275,47)
\path (210,0)(210,106)
\path (207.264,98.482)(210,106)(212.736,98.482)
\path (267.483,49.736)(275,47)(267.483,44.264)
\put(274,30){{\rm\rm {$x$}}}
\put(215,98){{\rm\rm {$z$}}}
\path (348,47)(481,47)
\path (416,0)(416,106)
\put(480,30){{\rm\rm {$x$}}}
\put(421,98){{\rm\rm {$z$}}}
\path (473.483,49.736)(481,47)(473.483,44.264)
\path (413.264,98.482)(416,106)(418.736,98.482)
\put(193,47){\circle{34}}
\path (416,47)(416.955,47.735)(417.92,48.44)(418.895,49.115)(419.88,49.76)
(420.875,50.375)(421.88,50.96)(422.895,51.515)(423.92,52.04)(424.955,52.535)
(426,53)(427.325,53.57)(428.6,54.08)(429.825,54.53)(431,54.92)
(432.125,55.25)(433.2,55.52)(434.225,55.73)(435.2,55.88)(436.125,55.97)
(437,56)(437.78,55.835)(438.52,55.64)(439.22,55.415)(439.88,55.16)
(440.5,54.875)(441.08,54.56)(441.62,54.215)(442.12,53.84)(442.58,53.435)
(443,53)(443.38,52.355)(443.72,51.72)(444.02,51.095)(444.28,50.48)
(444.5,49.875)(444.68,49.28)(444.82,48.695)(444.92,48.12)(444.98,47.555)
(445,47)(445.07,46.455)(445.08,45.92)(445.03,45.395)(444.92,44.88)
(444.75,44.375)(444.52,43.88)(444.23,43.395)(443.88,42.92)(443.47,42.455)
(443,42)(442.2,41.375)(441.4,40.8)(440.6,40.275)(439.8,39.8)
(439,39.375)(438.2,39)(437.4,38.675)(436.6,38.4)(435.8,38.175)
(435,38)(434.335,37.785)(433.64,37.64)(432.915,37.565)(432.16,37.56)
(431.375,37.625)(430.56,37.76)(429.715,37.965)(428.84,38.24)(427.935,38.585)
(427,39)(425.945,39.71)(424.88,40.44)(423.805,41.19)(422.72,41.96)
(421.625,42.75)(420.52,43.56)(419.405,44.39)(418.28,45.24)(417.145,46.11)
(416,47)(414.44,48.045)(412.96,49.08)(411.56,50.105)(410.24,51.12)
(409,52.125)(407.84,53.12)(406.76,54.105)(405.76,55.08)(404.84,56.045)
(404,57)(403.285,57.99)(402.64,58.96)(402.065,59.91)(401.56,60.84)
(401.125,61.75)(400.76,62.64)(400.465,63.51)(400.24,64.36)(400.085,65.19)
(400,66)(400.165,66.97)(400.36,67.88)(400.585,68.73)(400.84,69.52)
(401.125,70.25)(401.44,70.92)(401.785,71.53)(402.16,72.08)(402.565,72.57)
(403,73)(403.6,73.505)(404.2,73.92)(404.8,74.245)(405.4,74.48)
(406,74.625)(406.6,74.68)(407.2,74.645)(407.8,74.52)(408.4,74.305)
(409,74)(409.825,73.245)(410.6,72.48)(411.325,71.705)(412,70.92)
(412.625,70.125)(413.2,69.32)(413.725,68.505)(414.2,67.68)(414.625,66.845)
(415,66)(415.145,65.1)(415.28,64.2)(415.405,63.3)(415.52,62.4)
(415.625,61.5)(415.72,60.6)(415.805,59.7)(415.88,58.8)(415.945,57.9)
(416,57)(416.045,56.37)(416.08,55.68)(416.105,54.93)(416.12,54.12)
(416.125,53.25)(416.12,52.32)(416.105,51.33)(416.08,50.28)(416.045,49.17)
(416,48)(415.99,46.05)(415.96,44.2)(415.91,42.45)(415.84,40.8)
(415.75,39.25)(415.64,37.8)(415.51,36.45)(415.36,35.2)(415.19,34.05)
(415,33)(414.97,32.365)(414.88,31.76)(414.73,31.185)(414.52,30.64)
(414.25,30.125)(413.92,29.64)(413.53,29.185)(413.08,28.76)(412.57,28.365)
(412,28)(410.92,27.575)(409.88,27.2)(408.88,26.875)(407.92,26.6)
(407,26.375)(406.12,26.2)(405.28,26.075)(404.48,26)(403.72,25.975)
(403,26)(402.23,26.075)(401.52,26.2)(400.87,26.375)(400.28,26.6)
(399.75,26.875)(399.28,27.2)(398.87,27.575)(398.52,28)(398.23,28.475)
(398,29)(397.74,29.89)(397.56,30.76)(397.46,31.61)(397.44,32.44)
(397.5,33.25)(397.64,34.04)(397.86,34.81)(398.16,35.56)(398.54,36.29)
(399,37)(399.945,37.69)(400.88,38.36)(401.805,39.01)(402.72,39.64)
(403.625,40.25)(404.52,40.84)(405.405,41.41)(406.28,41.96)(407.145,42.49)
(408,43)(408.845,43.49)(409.68,43.96)(410.505,44.41)(411.32,44.84)
(412.125,45.25)(412.92,45.64)(413.705,46.01)(414.48,46.36)(415.245,46.69)
(416,47)
\end{picture}
\end{center}
\caption{}\label{f19}
\end{figure}


\no{\bf Hyperbolic case}
We take the frame~:
$$F_1=\f{\partial}{\partial x}\ ,\ \ F_2=\f{\partial}{\partial
z}+(\varepsilon xy+x^2z+mx^3z^2)\f{\partial}{\partial y}\ ,\ \
F_3=\f{\partial}{\partial y}$$
\no and we introduce~: $P_i=<p,F_i(q)>$. The geodesics equations
are~:
\begin{equation*}
\begin{split}
\dot{x}&=P_1,\
\dot{y}=P_2(\varepsilon xy+x^2z+mx^3z^2),\
\dot{z}=P_2 \\
\dot{P_1}&=-(\varepsilon y+2xz+3mx^2z^2)P_2P_3 \\
\dot{P_2}&=(\varepsilon y+2xz+3mx^2z^2)P_1P_3 \\
\dot{P_3}&=-\varepsilon xP_2P_3 
\end{split}
\end{equation*}

\no and they can be truncated at order $-1$ by making
$\varepsilon =m=0$. In this case $P_3$ is a first integral and we
can set~: $P_3=\lambda$. Moreover if we introduce~: $P_1=\sin
\theta$, $P_2=\cos \theta$, the equations become~:
$$\dot{x}=\sin \theta\ , \ \ \dot{z}=\cos \theta\ ,\ \
\dot{y}=\cos \theta \ x^2z$$
$$\dot{\theta}=-2xzP_3\ , \ \ P_3=\lambda $$

\no They can be projected onto the space $(x,z,\theta )$.

\setlength{\unitlength}{0.5mm}
\begin{figure}[h]
\begin{center} 

\begin{picture}(180,100)
\thinlines
\drawvector{10.0}{48.0}{74.0}{1}{0}
\drawvector{96.0}{48.0}{74.0}{1}{0}
\drawvector{46.0}{16.0}{68.0}{0}{1}
\drawvector{100.0}{32.0}{52.0}{0}{1}
\path(12.0,48.0)(12.0,48.0)(12.27,47.95)(12.55,47.9)(12.82,47.86)(13.1,47.81)(13.38,47.77)(13.64,47.7)(13.92,47.65)(14.18,47.59)
\path(14.18,47.59)(14.44,47.54)(14.71,47.47)(14.97,47.4)(15.23,47.34)(15.5,47.27)(15.76,47.2)(16.01,47.13)(16.27,47.04)(16.52,46.97)
\path(16.52,46.97)(16.78,46.88)(17.03,46.79)(17.28,46.7)(17.52,46.63)(17.77,46.52)(18.01,46.43)(18.25,46.34)(18.5,46.25)(18.72,46.13)
\path(18.72,46.13)(18.96,46.04)(19.2,45.93)(19.44,45.83)(19.67,45.72)(19.9,45.59)(20.13,45.49)(20.36,45.36)(20.59,45.25)(20.8,45.13)
\path(20.8,45.13)(21.04,45.0)(21.26,44.86)(21.47,44.74)(21.7,44.61)(21.9,44.47)(22.12,44.34)(22.34,44.2)(22.55,44.06)(22.77,43.9)
\path(22.77,43.9)(22.96,43.77)(23.18,43.61)(23.38,43.45)(23.59,43.31)(23.79,43.15)(23.98,43.0)(24.19,42.83)(24.38,42.66)(24.59,42.5)
\path(24.59,42.5)(24.78,42.34)(24.96,42.16)(25.17,41.99)(25.36,41.81)(25.54,41.63)(25.72,41.45)(25.9,41.27)(26.1,41.09)(26.28,40.9)
\path(26.28,40.9)(26.45,40.7)(26.63,40.52)(26.8,40.33)(26.98,40.13)(27.15,39.93)(27.34,39.72)(27.51,39.52)(27.67,39.31)(27.84,39.11)
\path(27.84,39.11)(28.01,38.88)(28.17,38.68)(28.32,38.45)(28.48,38.25)(28.64,38.02)(28.8,37.79)(28.96,37.56)(29.12,37.34)(29.27,37.11)
\path(29.27,37.11)(29.43,36.88)(29.57,36.65)(29.71,36.4)(29.87,36.16)(30.01,35.93)(30.15,35.68)(30.29,35.43)(30.44,35.18)(30.57,34.93)
\path(30.57,34.93)(30.7,34.68)(30.85,34.41)(30.97,34.15)(31.12,33.9)(31.25,33.63)(31.37,33.36)(31.5,33.09)(31.62,32.81)(31.75,32.54)
\path(31.75,32.54)(31.87,32.27)(31.98,32.0)(32.0,32.0)
\path(46.0,48.0)(46.0,48.0)(45.63,47.95)(45.27,47.9)(44.93,47.86)(44.56,47.81)(44.22,47.77)(43.88,47.7)(43.52,47.65)(43.18,47.59)
\path(43.18,47.59)(42.84,47.54)(42.52,47.47)(42.18,47.4)(41.84,47.34)(41.52,47.27)(41.18,47.2)(40.86,47.13)(40.54,47.04)(40.22,46.97)
\path(40.22,46.97)(39.9,46.88)(39.59,46.79)(39.27,46.7)(38.95,46.63)(38.65,46.52)(38.34,46.43)(38.04,46.34)(37.75,46.25)(37.45,46.13)
\path(37.45,46.13)(37.15,46.04)(36.86,45.93)(36.56,45.83)(36.27,45.72)(35.99,45.59)(35.7,45.49)(35.41,45.36)(35.13,45.25)(34.86,45.13)
\path(34.86,45.13)(34.59,45.0)(34.31,44.86)(34.04,44.74)(33.77,44.61)(33.52,44.47)(33.25,44.34)(32.99,44.2)(32.72,44.06)(32.47,43.9)
\path(32.47,43.9)(32.22,43.77)(31.96,43.61)(31.72,43.45)(31.47,43.31)(31.23,43.15)(31.0,43.0)(30.76,42.83)(30.52,42.66)(30.29,42.5)
\path(30.29,42.5)(30.04,42.34)(29.82,42.16)(29.6,41.99)(29.37,41.81)(29.14,41.63)(28.93,41.45)(28.71,41.27)(28.5,41.09)(28.29,40.9)
\path(28.29,40.9)(28.07,40.7)(27.87,40.52)(27.67,40.33)(27.45,40.13)(27.26,39.93)(27.05,39.72)(26.87,39.52)(26.68,39.31)(26.47,39.11)
\path(26.47,39.11)(26.29,38.88)(26.11,38.68)(25.93,38.45)(25.75,38.25)(25.56,38.02)(25.38,37.79)(25.21,37.56)(25.04,37.34)(24.87,37.11)
\path(24.87,37.11)(24.7,36.88)(24.54,36.65)(24.37,36.4)(24.21,36.16)(24.06,35.93)(23.9,35.68)(23.76,35.43)(23.61,35.18)(23.45,34.93)
\path(23.45,34.93)(23.31,34.68)(23.17,34.41)(23.03,34.15)(22.88,33.9)(22.76,33.63)(22.62,33.36)(22.48,33.09)(22.37,32.81)(22.23,32.54)
\path(22.23,32.54)(22.12,32.27)(22.0,32.0)(22.0,32.0)
\drawarc{27.0}{30.0}{10.77}{5.9}{3.52}
\path(46.0,48.0)(46.0,48.0)(46.38,48.04)(46.79,48.08)(47.18,48.13)(47.56,48.16)(47.95,48.22)(48.34,48.27)(48.72,48.33)(49.09,48.38)
\path(49.09,48.38)(49.47,48.45)(49.84,48.52)(50.2,48.58)(50.56,48.65)(50.91,48.72)(51.27,48.79)(51.63,48.86)(51.99,48.93)(52.33,49.02)
\path(52.33,49.02)(52.68,49.09)(53.02,49.18)(53.36,49.27)(53.68,49.36)(54.02,49.45)(54.34,49.54)(54.66,49.65)(55.0,49.75)(55.31,49.84)
\path(55.31,49.84)(55.63,49.95)(55.93,50.06)(56.25,50.15)(56.54,50.27)(56.86,50.38)(57.15,50.5)(57.45,50.61)(57.75,50.74)(58.04,50.86)
\path(58.04,50.86)(58.31,50.99)(58.59,51.11)(58.88,51.25)(59.15,51.38)(59.43,51.5)(59.7,51.65)(59.97,51.79)(60.24,51.93)(60.5,52.08)
\path(60.5,52.08)(60.75,52.22)(61.0,52.36)(61.25,52.52)(61.5,52.68)(61.75,52.84)(61.99,52.99)(62.22,53.15)(62.47,53.31)(62.7,53.49)
\path(62.7,53.49)(62.93,53.65)(63.15,53.81)(63.38,54.0)(63.59,54.16)(63.81,54.34)(64.02,54.52)(64.23,54.7)(64.44,54.9)(64.63,55.09)
\path(64.63,55.09)(64.83,55.27)(65.04,55.47)(65.23,55.65)(65.43,55.86)(65.61,56.06)(65.8,56.25)(65.98,56.47)(66.15,56.66)(66.33,56.88)
\path(66.33,56.88)(66.5,57.09)(66.66,57.31)(66.83,57.52)(66.98,57.74)(67.15,57.97)(67.3,58.18)(67.45,58.41)(67.61,58.63)(67.75,58.86)
\path(67.75,58.86)(67.9,59.11)(68.04,59.34)(68.16,59.58)(68.3,59.81)(68.43,60.06)(68.55,60.31)(68.68,60.56)(68.8,60.81)(68.91,61.06)
\path(68.91,61.06)(69.02,61.31)(69.15,61.56)(69.25,61.83)(69.36,62.09)(69.45,62.36)(69.55,62.61)(69.65,62.88)(69.73,63.16)(69.83,63.43)
\path(69.83,63.43)(69.91,63.72)(70.0,63.99)(70.0,64.0)
\path(80.0,48.0)(80.0,48.0)(79.59,48.11)(79.19,48.24)(78.8,48.36)(78.43,48.47)(78.05,48.61)(77.66,48.72)(77.29,48.84)(76.91,48.97)
\path(76.91,48.97)(76.55,49.11)(76.19,49.22)(75.83,49.36)(75.48,49.49)(75.12,49.61)(74.79,49.75)(74.44,49.88)(74.11,50.02)(73.76,50.15)
\path(73.76,50.15)(73.44,50.27)(73.12,50.41)(72.8,50.56)(72.48,50.68)(72.16,50.83)(71.84,50.97)(71.55,51.11)(71.25,51.25)(70.94,51.38)
\path(70.94,51.38)(70.65,51.52)(70.36,51.66)(70.08,51.81)(69.8,51.95)(69.51,52.09)(69.23,52.24)(68.97,52.38)(68.7,52.54)(68.44,52.68)
\path(68.44,52.68)(68.19,52.83)(67.93,52.97)(67.68,53.13)(67.44,53.27)(67.19,53.43)(66.95,53.59)(66.72,53.74)(66.48,53.88)(66.26,54.04)
\path(66.26,54.04)(66.05,54.2)(65.83,54.36)(65.61,54.52)(65.4,54.68)(65.19,54.84)(65.0,54.99)(64.8,55.15)(64.59,55.31)(64.41,55.47)
\path(64.41,55.47)(64.23,55.63)(64.05,55.79)(63.86,55.97)(63.68,56.13)(63.52,56.29)(63.36,56.47)(63.2,56.63)(63.04,56.79)(62.88,56.97)
\path(62.88,56.97)(62.72,57.13)(62.59,57.31)(62.45,57.47)(62.31,57.65)(62.16,57.83)(62.04,58.0)(61.91,58.18)(61.79,58.34)(61.68,58.52)
\path(61.68,58.52)(61.56,58.7)(61.45,58.88)(61.34,59.06)(61.25,59.24)(61.15,59.43)(61.04,59.61)(60.95,59.79)(60.88,59.97)(60.79,60.15)
\path(60.79,60.15)(60.72,60.34)(60.63,60.52)(60.56,60.7)(60.5,60.9)(60.45,61.08)(60.38,61.27)(60.33,61.45)(60.27,61.65)(60.24,61.84)
\path(60.24,61.84)(60.2,62.02)(60.15,62.22)(60.11,62.41)(60.09,62.61)(60.06,62.81)(60.04,63.0)(60.02,63.2)(60.0,63.4)(60.0,63.59)
\path(60.0,63.59)(60.0,63.79)(60.0,63.99)(60.0,64.0)
\drawarc{65.0}{64.0}{10.0}{3.14}{6.28}
\drawcenteredtext{50.0}{80.0}{$z$}
\drawcenteredtext{80.0}{44.0}{$x$}
\drawcenteredtext{96.0}{80.0}{$z$}
\drawcenteredtext{166.0}{44.0}{$x$}
\path(100.0,80.0)(100.0,80.0)(100.34,79.44)(100.7,78.88)(101.06,78.34)(101.43,77.81)(101.79,77.29)(102.13,76.77)(102.5,76.26)(102.84,75.76)
\path(102.84,75.76)(103.19,75.27)(103.55,74.79)(103.91,74.31)(104.26,73.84)(104.61,73.38)(104.95,72.94)(105.3,72.5)(105.65,72.05)(106.0,71.62)
\path(106.0,71.62)(106.34,71.2)(106.69,70.8)(107.04,70.4)(107.37,70.0)(107.72,69.61)(108.05,69.23)(108.4,68.86)(108.75,68.5)(109.08,68.13)
\path(109.08,68.13)(109.41,67.79)(109.76,67.44)(110.09,67.12)(110.44,66.79)(110.76,66.48)(111.11,66.16)(111.44,65.87)(111.76,65.58)(112.09,65.3)
\path(112.09,65.3)(112.44,65.01)(112.76,64.75)(113.09,64.48)(113.43,64.23)(113.75,64.0)(114.08,63.75)(114.41,63.52)(114.73,63.31)(115.05,63.09)
\path(115.05,63.09)(115.37,62.9)(115.7,62.7)(116.02,62.5)(116.34,62.33)(116.66,62.15)(117.0,62.0)(117.3,61.84)(117.62,61.68)(117.94,61.54)
\path(117.94,61.54)(118.26,61.41)(118.58,61.29)(118.9,61.18)(119.22,61.06)(119.52,60.97)(119.83,60.88)(120.15,60.79)(120.47,60.72)(120.77,60.65)
\path(120.77,60.65)(121.08,60.59)(121.4,60.54)(121.69,60.5)(122.01,60.45)(122.31,60.43)(122.62,60.4)(122.93,60.4)(123.23,60.38)(123.54,60.4)
\path(123.54,60.4)(123.83,60.4)(124.13,60.43)(124.44,60.45)(124.73,60.5)(125.04,60.54)(125.33,60.59)(125.63,60.65)(125.94,60.72)(126.23,60.79)
\path(126.23,60.79)(126.52,60.88)(126.83,60.97)(127.12,61.06)(127.41,61.18)(127.69,61.29)(128.0,61.41)(128.28,61.54)(128.58,61.68)(128.86,61.84)
\path(128.86,61.84)(129.14,61.99)(129.44,62.15)(129.72,62.33)(130.02,62.5)(130.3,62.7)(130.58,62.88)(130.86,63.09)(131.14,63.31)(131.42,63.52)
\path(131.16,63.36)(131.16,63.36)(131.21,63.45)(131.25,63.52)(131.3,63.61)(131.33,63.7)(131.38,63.79)(131.41,63.86)(131.46,63.95)(131.5,64.02)
\path(131.5,64.02)(131.52,64.11)(131.57,64.19)(131.61,64.27)(131.63,64.36)(131.67,64.43)(131.71,64.51)(131.74,64.58)(131.77,64.66)(131.8,64.75)
\path(131.8,64.75)(131.83,64.83)(131.86,64.9)(131.88,64.98)(131.91,65.05)(131.94,65.12)(131.97,65.2)(132.0,65.27)(132.02,65.36)(132.03,65.43)
\path(123.06,68.52)(123.05,68.51)(123.19,68.66)(123.33,68.8)(123.47,68.93)(123.59,69.05)(123.73,69.18)(123.86,69.3)(123.98,69.41)(124.12,69.54)
\path(124.12,69.54)(124.25,69.65)(124.37,69.76)(124.51,69.87)(124.62,69.97)(124.76,70.06)(124.87,70.16)(125.01,70.26)(125.12,70.34)(125.26,70.44)
\path(125.26,70.44)(125.37,70.51)(125.5,70.61)(125.62,70.68)(125.73,70.76)(125.86,70.83)(125.98,70.9)(126.09,70.95)(126.22,71.02)(126.33,71.08)
\path(126.33,71.08)(126.44,71.13)(126.55,71.19)(126.68,71.23)(126.79,71.29)(126.9,71.33)(127.01,71.37)(127.12,71.41)(127.23,71.44)(127.34,71.48)
\path(127.34,71.48)(127.45,71.5)(127.56,71.52)(127.68,71.55)(127.77,71.56)(127.88,71.58)(127.98,71.59)(128.1,71.61)(128.19,71.61)(128.3,71.61)
\path(128.3,71.61)(128.41,71.61)(128.5,71.61)(128.61,71.59)(128.72,71.58)(128.82,71.58)(128.91,71.55)(129.0,71.54)(129.11,71.51)(129.21,71.48)
\path(129.21,71.48)(129.3,71.47)(129.39,71.43)(129.5,71.4)(129.58,71.36)(129.67,71.31)(129.77,71.26)(129.86,71.23)(129.96,71.18)(130.05,71.12)
\path(130.05,71.12)(130.13,71.05)(130.22,71.0)(130.32,70.94)(130.41,70.87)(130.49,70.8)(130.58,70.73)(130.66,70.65)(130.75,70.58)(130.83,70.48)
\path(130.83,70.48)(130.91,70.41)(131.0,70.31)(131.08,70.23)(131.16,70.12)(131.25,70.02)(131.33,69.93)(131.39,69.83)(131.47,69.72)(131.55,69.61)
\path(132.03,65.43)(132.05,65.5)(132.08,65.58)(132.11,65.65)(132.13,65.72)(132.14,65.79)(132.16,65.86)(132.17,65.93)(132.19,66.0)(132.21,66.06)
\path(132.21,66.06)(132.22,66.13)(132.24,66.2)(132.25,66.27)(132.27,66.34)(132.27,66.41)(132.28,66.48)(132.3,66.55)(132.3,66.62)(132.3,66.68)
\path(123.0,68.59)(122.69,68.19)(122.41,67.8)(122.15,67.4)(121.9,67.0)(121.66,66.61)(121.44,66.23)(121.23,65.83)(121.05,65.47)(120.87,65.08)
\path(120.87,65.08)(120.73,64.72)(120.59,64.36)(120.47,64.0)(120.37,63.63)(120.27,63.27)(120.2,62.93)(120.16,62.59)(120.12,62.25)(120.09,61.9)
\path(120.09,61.9)(120.08,61.56)(120.09,61.24)(120.12,60.91)(120.18,60.59)(120.23,60.27)(120.3,59.95)(120.4,59.65)(120.51,59.34)(120.63,59.04)
\path(120.63,59.04)(120.77,58.75)(120.94,58.45)(121.12,58.15)(121.3,57.88)(121.51,57.59)(121.73,57.31)(121.98,57.04)(122.23,56.77)(122.51,56.5)
\path(122.51,56.5)(122.8,56.25)(123.11,55.99)(123.43,55.72)(123.76,55.47)(124.12,55.24)(124.48,55.0)(124.87,54.75)(125.29,54.52)(125.7,54.29)
\path(125.7,54.29)(126.15,54.06)(126.59,53.84)(127.06,53.63)(127.55,53.4)(128.05,53.2)(128.58,53.0)(129.11,52.79)(129.66,52.59)(130.24,52.4)
\path(130.24,52.4)(130.82,52.2)(131.41,52.02)(132.03,51.84)(132.66,51.66)(133.32,51.49)(133.99,51.31)(134.67,51.15)(135.36,51.0)(136.08,50.84)
\path(136.08,50.84)(136.82,50.68)(137.57,50.52)(138.33,50.38)(139.11,50.25)(139.91,50.11)(140.72,49.97)(141.55,49.84)(142.41,49.72)(143.27,49.59)
\path(143.27,49.59)(144.16,49.47)(145.05,49.36)(145.96,49.25)(146.88,49.15)(147.83,49.04)(148.8,48.95)(149.77,48.84)(150.77,48.75)(151.78,48.68)
\path(151.78,48.68)(152.8,48.59)(153.86,48.52)(154.91,48.43)(156.0,48.36)(157.08,48.31)(158.19,48.25)(159.33,48.18)(160.47,48.13)(161.63,48.08)
\path(161.63,48.08)(162.8,48.04)(163.99,48.0)(164.0,48.0)
\path(132.3,66.68)(132.32,66.75)(132.33,66.8)(132.33,66.87)(132.33,66.94)(132.33,67.0)(132.33,67.06)(132.33,67.12)(132.35,67.19)(132.35,67.25)
\path(132.35,67.25)(132.33,67.31)(132.33,67.37)(132.33,67.44)(132.33,67.5)(132.33,67.55)(132.33,67.62)(132.32,67.66)(132.32,67.73)(132.3,67.79)
\path(132.3,67.79)(132.3,67.84)(132.28,67.9)(132.27,67.95)(132.27,68.01)(132.25,68.06)(132.25,68.12)(132.22,68.18)(132.22,68.23)(132.19,68.29)
\path(132.19,68.29)(132.19,68.33)(132.16,68.4)(132.14,68.44)(132.13,68.5)(132.11,68.55)(132.08,68.59)(132.07,68.65)(132.05,68.69)(132.02,68.75)
\path(132.02,68.75)(132.0,68.8)(131.97,68.84)(131.96,68.9)(131.92,68.94)(131.89,69.0)(131.86,69.04)(131.85,69.08)(131.82,69.13)(131.78,69.18)
\path(131.78,69.18)(131.75,69.23)(131.72,69.26)(131.69,69.31)(131.66,69.36)(131.61,69.4)(131.58,69.44)(131.55,69.48)(131.5,69.52)(131.47,69.56)
\path(131.47,69.56)(131.42,69.61)(131.38,69.66)(131.39,69.66)
\end{picture}

\end{center}
\caption{}\label{hyp}
\end{figure}

\paragraph{Numerical simulations \\ }
The projected equations~: $\dot{x}=\sin \theta, \dot{z}=\cos
\theta, \dot{\theta}=-2xz\lambda$ can be integrated numerically
and the solutions compared with the pendulum~: $\dot{x}=\sin
\theta, \dot{z}=\cos \theta, \dot{\theta}=-\lambda x$, see
Fig. \ref{hyp}.
The behaviour is \it{quite chaotic} and $\theta$ exhibits
oscillating and dissipative phenomena. Also it shows in the plane
$(x,z)$ a \it{coupling effect between the two abnormal
directions}.

\begin{figure}[h]
\epsfxsize=7cm
\epsfysize=7cm
\centerline{\epsffile{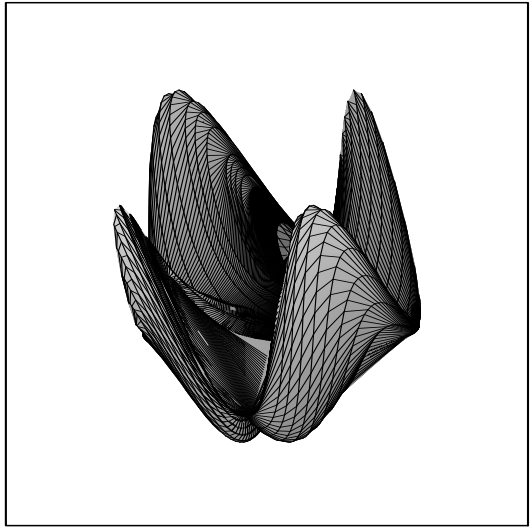}}
\caption{Tangential hyperbolic sphere}\label{figtanghyp}
\end{figure}

The SR sphere in the tangential hyperbolic sphere is represented
on Fig. \ref{figtanghyp}.


\subsubsection{Conclusion}
The elliptic situation is similar to the contact case. The
analysis of the hyperbolic case is intricate. A tool to study
the sphere is to introduce as in the Martinet case a return
mapping by taking the intersections of the geodesics with one of
the planes~: $x=0$ or $z=0$. The non properness of
this application can be checked numerically.


\subsection{The Engel case and left-invariant SR geometry
on nilpotent Lie groups}
If $q=(x,y,z,w)$, we consider the system in $\R^4$~:
$$F_1=\f{\partial}{\partial x}+y\f{\partial}{\partial
z}+\f{y^2}{2}\f{\partial}{\partial w},\esp
F_2=\f{\partial}{\partial y}$$
We have the following relations~:
$F_3=[F_1,F_2]=\f{\partial}{\partial z}+y\f{\partial}{\partial
w}, \esp
F_4=[[F_1,F_2],F_2]=\f{\partial}{\partial w}$
and $[[F_1,F_2],F_1]=0$. Moreover all Lie brackets with length
greater than $4$ are equal to zero. Set~:
$$L_1=\l(\begin{array}{cccc}
0&0&0&0\\
0&0&1&0\\
0&0&0&1\\
0&0&0&0
\end{array}\r),\quad
L_2=\l(\begin{array}{cccc}
0&1&0&0\\
0&0&0&0\\
0&0&0&0\\
0&0&0&0
\end{array}\r)$$
and define the following representation~:
$\rho(F_1)=L_1, \rho(F_2)=L_2$
which allows to identify the previous system in $\R^4$ to the
\it{left-invariant system}
$\dot{R}=(u_1L_1+u_2L_2)R$
on the Engel group $G_e$, here represented by the nilpotent
matrices~:
$$\l(\begin{array}{cccc}
1&q_2&q_3&q_4\\
0&1&q_1&\f{q_1^2}{2}\\
0&0&1&q_1\\
0&0&0&1
\end{array}\r)$$
The weight of $x,y$ is one, the weight of $z$ is two, and the
weight of $w$ is three. For any sub-Riemannian metric on $G_e$,
the approximation of order $-1$ is the flat metric $g=dx^2+dy^2$.
Any sub-Riemannian Martinet metric can be written
$g=adx^2+cdy^2$ and can be lifted on $G_e$.
\\

\subsubsection{Paramatrization of geodesics in the flat case}
Non trivial abnormal extremals are solutions of~:
\begin{equation*}
P_1=P_2=\{P_1,P_2\}=0 ,\esp
u_1\{\{P_1,P_2\},P_1\}+u_2\{\{P_1,P_2\},P_2\}=0
\end{equation*}
Set $p=(p_x,p_y,p_z,p_w)$. We get~:
$$p_x+p_zy+p_w\f{y^2}{2}=p_y=p_z+yp_w=p_wu_2=0$$
This implies $p_w\neq 0$ and thus $u_2=0$. The abnormal flow is
given by~:
$$\dot{x}=u_1,\esp \dot{y}=0,\esp \dot{z}=u_1y,\esp
\dot{w}=u_1\f{y^2}{2}$$
where $|u_1|=1$ if the parameter is the length.

To compute normal extremals, we set $P_i=<p,F_i(q)>,
i=1,2,3,4$ and $H_n=\inv{2}(P_1^2+P_2^2)$~; we get~:
\begin{equation*}
\dot{P}_1=P_2P_3,\esp
\dot{P}_2=-P_1P_3,\esp
\dot{P}_3=P_2P_4,\esp
\dot{P}_4=0
\end{equation*}
Parametrizing by the length $H_n=\inv{2}$, we may set~:
$P_1=\cos\theta, P_2=\sin\theta$, and if $\theta\neq k\pi$ we
get~:
$\dot{\theta}=-P_3,\esp \ddot{\theta}=-P_2P_4$.
Denote $P_4=\lambda$, then this is equivalent to the pendulum
equation~:
$$\ddot{\theta}+\lambda\sin\theta=0$$

Let $L$ denote the abnormal line starting from $0$~: $t\mapsto
(\pm
t,0,0,0)$. \it{It is not strict} and projects onto
$\theta=k\pi$.

In order to obtain an uniform representation of normal geodesics,
we shall use the \it{Weierstrass elliptic function}
${\cal P}$. Indeed the system admits three integrals~:
$P_1^2+P_2^2=1$, and two Casimir functions~:
$-2P_1P_4+P_3^2=C$ et $P_4=\lambda$.
Using $\dot{P}_1=P_2P_3$ we get~:
$\ddot{P}_1=-CP_1-3\lambda P_1^2+\lambda$,
which is equivalent, with $\dot{P}_1\neq 0$ and $\lambda\neq 0$,
to the equation~:
$\dot{P_1}^2=-2\lambda(P_1^3+\f{C}{2\lambda}P_1^2-P_1+
D)$.
Let ${\cal P}(u)$ denote the Weierstrass elliptic function
(cf \cite{L}) solution of~:
$${\cal P}'(u)=-2\sqrt{({\cal P}(u)-e_1)({\cal P}(u)-e_2)({\cal
P}(u)-e_3)}$$
where the complex numbers $e_i$ satisfy $e_1+e_2+e_3=0$. Set
$g_2=-4(e_2e_3+e_3e_1+e_1e_2)$ and $g_3=4e_1e_2e_3$~;
then~:
${\cal P}'(u)=4{\cal P}^3(u)-g_2{\cal P}(u)-g_3$.
The function ${\cal P}(u)$ can be expanded at $0$ in the
following way~:
$${\cal P}(u)=\inv{u^2}+\inv{20}q_2u^2+\inv{28}q_3u^4+\cdots$$
Hence the solution $P_1$ can be written~: $a{\cal
P}(u)+b$. Then we can compute $P_2$ and $P_3$ using
the integrals, and $x,y,z,w$ can be computed by
quadratures. We find again oscillating and rotating solutions of
the pendulum using Jacobi elliptic functions given by the
formulas~:
$$\rm{cn }u=\l(\f{{\cal P}(u)-e_1}{{\cal P}(u)-e_2}\r)^\inv{2},
\esp \rm{dn }u=
\l(\f{{\cal P}(u)-e_2}{{\cal P}(u)-e_3}\r)^\inv{2}$$


\subsubsection{Heisenberg and Martinet flat cases deduced from the
Engel case. Blowing-up in lines}
Note that the two vector fields $\f{\partial}{\partial z}$ and
$\f{\partial}{\partial w}$ commute with $F_1$ and $F_2$. The
Engel case contains the flat contact case and the flat Martinet
case which are given by the following operations~:
\begin{itemize}
\item Setting $p_z=0$, we obtain the geodesics of the
Heisenberg case.
\item Setting $p_w=0$, we obtain the geodesics of the
Martinet flat case.
\end{itemize}

The interpretation is the following.

\begin{lem}
We obtain the Martinet flat case (resp. Heisenberg) by minimizing
the SR distance to the line $(Oz)$
(resp. $(Ow)$).
\end{lem}

Indeed the condition $p_z=0$ (resp. $p_w=0$) corresponds to the
\it{transversality condition}. It may be observed that, since the
SR distance to a line is more regular than to a point, the
SR distance in the Engel case has at least all singularities of
the Heisenberg and Martinet flat cases.

Another way to get the Martinet flat case is to use the following
general fact from \cite{Be}~:

\begin{lem}
The Martinet flat case is isometric to $(G_{e/H}, dx^2+dy^2)$
where $H$ is the following sub-group of
$G_e$~: $\{\rm{exp }t[F_1,F_2]\ /\ t\in\R\}$.
\end{lem}

The Engel case can be imbedded in any dimension, for more details
about left-invariant SR-geometry on nilpotent Lie groups see
\cite{Sach}.

In the Martinet flat case, using the uniform parametrization of the
geodesics by elliptic functions, the sphere is evaluated in any
direction.

Next we give a description of the SR sphere in an abnormal direction
when the flag associated to the distribution $D$ satisfies $D^3\neq D^2$,
see Section 3.4.


\subsection{Microlocal analysis of the singularity of the
SR sphere in the abnormal direction}
The aim of this Section is to stratify the singularity of the
SR sphere in the abnormal direction. We use symplectic
geometry. This leads to a stratification of the solutions of
the Hamilton-Jacobi-Bellman equation viewed in the cotangent
bundle.

\subsubsection{Lagrangian manifolds and generating mapping}
\begin{defi}
Let $(M,\omega)$ be a smooth symplectic manifold and $L\subset
M$ be a smooth regular submanifold. We say that $L$ is
isotropic if the restriction of $\omega$ to $TL$ is equal to
zero, and if $\rm{dim }L=\inv{2}\rm{dim }M$ then $L$ is called
Lagrangian.
\end{defi}

The following result is crucial, see \cite{Mi}.

\begin{prop}
Let $(M,\omega)$ be a $2n$-dimensional manifold and $L\subset M$
be a Lagrangian submanifold. Then there exists Darboux
local coordinates $(q,p)$ and a smooth function
$S(q_I,p_{\bar{I}})$
where $I=\{1,\ldots, m\}$, $\bar{I}=\{m+1,\ldots,n\}$ is a
partition of $\{1,\ldots,n\}$ such that $L$ is given locally
by the equations~:
$$p_I=\f{\partial S}{\partial q_I},\
q_{\bar{I}}=-\f{\partial S}{\partial p_{\bar{I}}}$$
\end{prop}

\begin{defi}
The mapping $S$ which represents locally $L$ is called the
\it{generating mapping} of $L$.
\end{defi}

\begin{defi}
Let $L$ be a Lagrangian manifold and $\Pi$ the standard projection
$(q,p)\mapsto q$ from $TM$ onto $M$. The caustic is the projection
on $M$ of the singularities of $(L,\Pi)$.
\end{defi}


\subsubsection{Lagrangian manifolds and SR normal case}
Consider the SR problem~:
$$\dot{q}=u_1F_1(q)+u_2F_2(q),\ q\in M$$
where the length of $q$ is~:
$$L(q)=\int_0^T (u_1^2(t)+u_2^2(t))^\inv{2}dt$$
We use the notations of Section 3.2. Let $P_i=<p,F_i(q)>,
i=1,2$, the Hamiltonian associated to normal geodesics is given
by $H_n=\inv{2}(P_1^2+P_2^2)$.

Let $t\mapsto \gamma(t), t\in[0,T]$ be a reference one-to-one
normal geodesic. We assume the following~:

\begin{quote}
{\bf Hypothesis~:} We assume that the reference geodesic is
\it{strict}, i.e. there exists an unique lifting $[\tilde{\gamma}]$
of $\gamma$ in the projective bundle $P(T^*M)$.
\end{quote}

\no{\bf Notations}
\begin{itemize}
\item $\rm{exp}_{\gamma(0)}$ is the exponential mapping. If the
geodesics are parametrized by arc-length $H_n=\inv{2}$, it is
defined by $t\mapsto \Pi(\tilde{q}(t))$ where $t\mapsto\tilde{q}(t)$
is a solution of $\vec{H}_n$ starting from $\gamma(0)$ at time
$t=0$.
\item $L_t=\rm{exp }t\vec{H}_n(T^*_{\gamma(0)}M)$, where
$\rm{exp }t\vec{H}_n$ is the local one-parameter group associated
to $\vec{H}_n$.
\end{itemize}

The length of a curve does not depend on the parametrization and
the optimal control problem is \it{parametric}. This induces a
symmetry which has to be taken into account when writing
Hamilton-Jacobi equation in the normal case. Indeed we have~:

\begin{lem}\label{lem5.5}
The solutions of $\vec{H}_n$ satisfy the relation~:
$$q(t,q_1,\lambda p_1)=q(\lambda t,q_1,p_1),\
p(t,q_1,\lambda p_1)=\lambda p(\lambda t,q_1,p_1)$$
\end{lem}

The following results are standard~:

\begin{prop}
\begin{enumerate}
\item $L_0=T^*_{\gamma(0)}M$ is a linear Lagrangian manifold, and
for each $t>0$, $L_t$ is a Lagrangian manifold.
\item The time $t_c$ is conjugate along $\gamma$ if and only if
the projection $\Pi~: L_{t_c}\rightarrow M$ is singular at
$\tilde{\gamma}(t_c)$.
\item Assume that geodesics are parametrized by arc-length $t$,
and let
$$\displaystyle{W=\bigcup_{0<t\leq T}\rm{exp }t\vec{H}_n
(T^*_{\gamma(0)}M\cap (H_n=\inv{2})) }$$
where $T<t_{1c}$ (first
conjugate time along $\gamma$. Then ${\cal E}=\Pi(W)$ is a central
field along $\gamma$.
\end{enumerate}
\end{prop}

\begin{rem}
\begin{itemize}
\item The caustic of $L_t$ is the set of conjugate points which can
be analyzed using Lagrangian singularities.
\item We represent locally $W$ by an Hamilton-Jacobi or wave
function defined as follows. We integrate the normal flow starting
from $\gamma(0)$ and parametrized by arc-length~: $P_1^2+P_2^2=1$.
By setting $P_1(0)=\cos\theta$, this gives us the family of
geodesics~:
$$E~: p=(\theta,\lambda_1,\ldots,\lambda_{n-2},t)\in {\mathbb S}^1
\times \R^{n-1}\mapsto M$$
and $p$ is eliminated by solving the equation $E(p)=q$ near
$\gamma(t)$ using the Implicit Function Theorem. Beyond the
computations need the Preparation Theorem and \it{Legendrian
singularity theory}.
\end{itemize}
\end{rem}

Next we describe the tangent space to the Lagrangian manifold.

\begin{defi}
We denote by $(V_n)$ the variational equation~:
\begin{equation} \label{eq5.1}
\delta\dot{\tilde{q}}=\f{\partial\vec{H}_n}{\partial\tilde{\gamma}}
(\tilde{q}(t))\delta q
\end{equation}
along the reference trajectory $t\mapsto \tilde{\gamma}(t)$.
This Hamiltonian linear equation is called \it{Jacobi equation}.
A \it{Jacobi field} $J(t)=(\delta q(t),\delta p(t))$ is a nontrivial
solution of (\ref{eq5.1}). It is called \it{vertical} if $\delta
q(0)=0$.
\end{defi}

\begin{prop}
\begin{enumerate}
\item Let $L_t=\rm{exp }t\vec{H}_n(T^*_{\gamma(0)}M)$. Then the
space of vertical Jacobi fields is the tangent space to $L_t$ for
$t>0$.
\item Assume we are in the analytic category. Let $J(.)$ be a
vertical Jacobi field and let $\varepsilon\mapsto\alpha
(\varepsilon)$ be an analytic curve such that $\alpha(0) =J(0)$.
If $Y$ is an analytic vector field on $T^*M$ such that
$Y(\tilde{\gamma}(0))=\alpha(0)$, then $t\mapsto J(t)$ is given
for $t$ small by the Baker-Campbell-Hausdorff formula~:
$$J(t)=\sum_{n\geq 0} \f{t^n}{n!}ad^n\vec{H}_n(Y)(\tilde{\gamma}(t))$$
\end{enumerate}
\end{prop}

A consequence of Lemma \ref{lem5.5} is~:

\begin{lem}
Let $\tilde{\gamma}(0)=(q_0,p_0)$ and consider the curve
$\alpha(\varepsilon)=(q_0,p_0+\varepsilon p_0)$. Then it is a
vertical curve, and if $J_1$ is the associated Jacobi field then
$\Pi(J_1(t))=t\dot{\gamma}(t)$.
\end{lem}


\subsubsection{Isotropic manifolds and SR abnormal case}
Consider the system $ \dot{q}=u_1F_1(q)+u_2F_2(q)$. According to
Section 3, the abnormal geodesics are solutions of the
equations~:
$$\f{dq}{dt}=\f{\partial H_a}{\partial p},\
\f{dp}{dt}=-\f{\partial H_a}{\partial q}$$
where $H_a=u_1P_1+u_2P_2$. They are contained in~:
$$P_1=P_2=\{P_1,P_2\}=0$$
and the abnormal controls are computed using~:
$$u_1\{\{P_1,P_2\},P_1\}+u_2\{\{P_1,P_2\},P_2\}=0$$

\no{\bf Assumptions}
Let $t\mapsto \gamma(t),t\in[-T,T]$ be a one-to-one abnormal
reference geodesic. One may assume that it corresponds to the
control $u_2=0$. We suppose that the following conditions are
satisfied along $\gamma$ for the couple $(F_1,F_2)$, see Section
2.5.
\begin{itemize}
\item[($H_1$)] The first order Pontryagin's cone $K(t)=\rm{Span
}\{ad^kF_1.{F_2}_{|\gamma}\ /\ k\in\N\}$ has codimension one and is
generated by $\{F_2,\ldots,ad^{n-2}F_1.F_2 \}_{|\gamma(t)}$.
\item[($H_2$)] If $n\geq 3$, for each $t$,
$F_1(\gamma(t))\notin\rm{Span }\{ad^kF_1.{F_2}_{|\gamma}\ /\
k=0\ldots n-3\}$.
\item[($H_3$)] $\{P_2,\{P_1,P_2\}\}\neq 0$ along $\gamma$.
\end{itemize}

\no{\bf Notations}
\begin{itemize}
\item Under the previous assumptions $\gamma$ admits an unique
lifting $[\tilde{\gamma}]=(\gamma,p_\gamma)$ in $P(T^*M)$. One
may identify locally $M$ to a neighborhood $U$ of $\gamma(0)=0$
in $\R^n$. Let $V$ be a neighborhood of $p_\gamma$ in
$P(T_0^*U)$. We can choose $V$ small enough such that all
abnormal geodesics starting from $\{0\}\times V$ satisfy the
assumptions $(H_1-H_3)$. We denote by $\Sigma_r$ the sector of
$U$ covered by abnormal geodesics with length $\leq r$ and
starting from $\{0\}\times V$. This defines a mapping denoted
Exp. The construction is represented on Fig. \ref{lafig1}.

\setlength{\unitlength}{0.3mm}
\begin{figure}[h]
\begin{center}

\begin{picture}(180,100)
\thinlines
\drawarc{-18.0}{48.0}{145.41}{5.73}{0.57}
\drawdot{-18.0}{48.0}
\drawthickdot{-18.0}{48.0}
\drawpath{-18.0}{48.0}{44.0}{86.0}
\drawpath{-18.0}{48.0}{52.0}{68.0}
\drawpath{-18.0}{48.0}{52.0}{30.0}
\drawpath{-18.0}{48.0}{44.0}{10.0}
\drawpath{-18.06}{48.08}{54.79}{48.08}
\drawdot{54.79}{48.08}
\drawthickdot{54.79}{48.08}
\drawlefttext{61.75}{48.08}{$p_\gamma$}
\drawcenteredtext{20.36}{-6.61}{$P(T^*_0U)$}
\path(77.05,60.43)(77.05,60.43)(77.51,60.72)(78.01,60.99)(78.5,61.25)(78.98,61.52)(79.48,61.77)(79.94,62.02)(80.44,62.27)(80.93,62.5)
\path(80.93,62.5)(81.41,62.75)(81.91,62.97)(82.37,63.2)(82.87,63.4)(83.37,63.61)(83.86,63.83)(84.33,64.01)(84.83,64.22)(85.3,64.41)
\path(85.3,64.41)(85.8,64.58)(86.29,64.76)(86.76,64.94)(87.26,65.08)(87.76,65.26)(88.23,65.41)(88.73,65.55)(89.22,65.69)(89.69,65.83)
\path(89.69,65.83)(90.19,65.94)(90.68,66.08)(91.16,66.19)(91.66,66.3)(92.15,66.41)(92.62,66.51)(93.12,66.62)(93.62,66.69)(94.11,66.79)
\path(94.11,66.79)(94.58,66.87)(95.08,66.94)(95.55,67.01)(96.05,67.05)(96.55,67.12)(97.04,67.16)(97.51,67.19)(98.01,67.25)(98.51,67.26)
\path(98.51,67.26)(99.0,67.3)(99.48,67.3)(99.98,67.33)(100.47,67.33)(100.94,67.33)(101.44,67.33)(101.94,67.33)(102.43,67.3)(102.91,67.3)
\path(102.91,67.3)(103.41,67.26)(103.91,67.23)(104.4,67.19)(104.87,67.15)(105.37,67.08)(105.87,67.05)(106.36,66.98)(106.83,66.91)(107.33,66.83)
\path(107.33,66.83)(107.83,66.76)(108.33,66.68)(108.8,66.58)(109.3,66.48)(109.8,66.37)(110.3,66.26)(110.76,66.16)(111.26,66.05)(111.76,65.91)
\path(111.76,65.91)(112.26,65.79)(112.75,65.66)(113.25,65.51)(113.73,65.36)(114.23,65.19)(114.72,65.05)(115.22,64.87)(115.72,64.69)(116.19,64.51)
\path(116.19,64.51)(116.69,64.33)(117.19,64.16)(117.69,63.97)(118.16,63.77)(118.66,63.56)(119.16,63.34)(119.66,63.13)(120.16,62.9)(120.62,62.65)
\path(120.62,62.65)(121.12,62.43)(121.62,62.2)(122.12,61.95)(122.62,61.7)(123.12,61.43)(123.61,61.15)(124.11,60.9)(124.61,60.63)(125.08,60.34)
\path(125.08,60.34)(125.58,60.06)(126.08,59.77)(126.08,59.77)
\drawvector{124.61}{60.63}{1.49}{2}{-1}
\drawvector{184.0}{44.0}{30.0}{-1}{0}
\drawvector{184.0}{44.0}{28.0}{1}{0}
\drawvector{184.0}{44.0}{28.0}{-3}{2}
\drawvector{184.0}{44.0}{26.0}{3}{2}
\drawvector{184.0}{44.0}{26.0}{3}{-2}
\drawvector{184.0}{44.0}{24.0}{-3}{-2}
\drawpath{156.0}{62.0}{146.0}{68.0}
\drawpath{154.0}{44.0}{144.0}{44.0}
\drawpath{160.0}{28.0}{148.0}{20.0}
\drawpath{210.0}{26.0}{222.0}{18.0}
\drawpath{212.0}{44.0}{222.0}{44.0}
\drawpath{210.0}{62.0}{224.0}{72.0}
\drawcenteredtext{184.0}{4.0}{$\Sigma_r$}
\drawcenteredtext{100.0}{74.0}{Exp}
\end{picture}

\end{center}
\caption{}\label{lafig1}
\end{figure}

\item On $\Omega=T^*M\backslash(\{\{P_1,P_2\},P_2\}=0)$ let
$\hat{H}_a$ be the Hamiltonian $\hat{H}_a=P_1+\hat{u}P_2$ where
$\hat{u}=-\f{\{\{P_1,P_2\},P_1\}} { \{\{P_1,P_2\},P_2\}  }$, and
let $\rm{exp }t\vec{\hat{H}}_a$ be the one-parameter local group.
We denote by $L_t^a$ (resp. $I_t^a$) the image of
$T^*_{\gamma(0)}M$ (resp. $T^*_{\gamma(0)}M \cap
(P_1=P_2=\{P_1,P_2\}=0)$).
\end{itemize}

\begin{lem}
On $\Omega$, $L_t^a$ is a Lagrangian submanifold and $I_t^a$ is
an isotropic submanifold.
\end{lem}


\subsubsection{The smooth abnormal sector of the SR sphere}
\begin{lem}
Consider the SR problem where $\gamma$ is an abnormal reference
trajectory and assume $(H_1-H_3)$.
It can be identified to a trajectory of $F_1$ where
the system $(F_1,F_2)$ is orthonormal. Then the abnormal geodesic
is strict, and there exists $r>0$ such that if the length of
$\gamma$ is less than $r$ then $\gamma$ is a global minimizer.
\end{lem}

\begin{proof}
Under assumption $(H_1)$ the first order Pontryagin's cone along
$\gamma$ has codimension one, and from $(H_3)$ $\gamma$ is not a
normal geodesic. The optimality assertion follows from
\cite{AS2}, see also \cite{LS}.
\end{proof}

Hence the end-point of $\gamma$ belongs to the sphere. Moreover
$r$ can be estimated and the estimate is uniform for each
abnormal geodesic $C^1$-close to $\gamma$. Therefore we have~:

\begin{prop}
For $r$ small enough $\Sigma_r$ is a sector of the SR ball
homeomorphic to $C\cup -C$, where $C$ is a positive cone of
dimension $n-3$ if $n\geq 4$ and $1$ if $n=3$. Its intersection
with the sphere consists of two ($C^\infty$ or $C^\omega$)
surfaces of dimension $n-4$ if $n\geq 4$ and reduced to two
points if $n=3$.
\end{prop}


\subsubsection{Gluing both normal and abnormal parts}
The tangent space to the sphere near the abnormal directions is
described by the results of Section \ref{sectiontangent}, namely
Theorem \ref{thmtangent}.

Let $A$ be the end-point of the abnormal trajectory and let
$K(r)$ be the first order Pontryagin's cone evaluated at $A$, $r$
small enough. Let $\varepsilon\mapsto\alpha(\varepsilon)$ be a
$C^1$ curve on the sphere $S(0,r)$, $\alpha(0)=A$,
$\varepsilon\geq 0$. Assume the following~:
\begin{enumerate}
\item $\alpha(\varepsilon)\subset S(0,r)\backslash\Sigma_r$ for
$\varepsilon\neq 0$.
\item $\alpha(\varepsilon)\cap L=\emptyset$, where $L$ is the
cut-locus for geodesics starting from $0$.
\end{enumerate}
Then the tangent space to the sphere evaluated at
$\alpha(\varepsilon)$ tends to $K(r)$ when
$\varepsilon\rightarrow 0$ (see Fig. \ref{mafig}).

\setlength{\unitlength}{0.4mm}
\begin{figure}[h]
\begin{center}

\begin{picture}(180,100)
\thicklines
\path(82.0,74.0)(82.0,74.0)(81.62,73.23)(81.3,72.48)(80.97,71.72)(80.65,70.97)(80.33,70.22)(80.04,69.44)(79.76,68.69)(79.48,67.97)
\path(79.48,67.97)(79.22,67.22)(78.97,66.47)(78.73,65.73)(78.51,64.98)(78.3,64.25)(78.08,63.5)(77.9,62.77)(77.72,62.04)(77.55,61.31)
\path(77.55,61.31)(77.37,60.56)(77.25,59.84)(77.12,59.11)(76.98,58.38)(76.87,57.65)(76.76,56.93)(76.69,56.22)(76.62,55.5)(76.55,54.77)
\path(76.55,54.77)(76.5,54.06)(76.44,53.34)(76.43,52.63)(76.41,51.9)(76.41,51.2)(76.41,50.49)(76.43,49.79)(76.44,49.08)(76.5,48.38)
\path(76.5,48.38)(76.55,47.65)(76.62,46.97)(76.69,46.27)(76.76,45.56)(76.87,44.88)(76.98,44.18)(77.11,43.49)(77.23,42.79)(77.37,42.09)
\path(77.37,42.09)(77.54,41.4)(77.69,40.72)(77.87,40.04)(78.08,39.36)(78.26,38.68)(78.48,38.0)(78.72,37.31)(78.94,36.63)(79.19,35.95)
\path(79.19,35.95)(79.47,35.29)(79.73,34.61)(80.01,33.93)(80.3,33.27)(80.62,32.61)(80.94,31.94)(81.26,31.28)(81.62,30.61)(81.97,29.95)
\path(81.97,29.95)(82.33,29.29)(82.69,28.62)(83.08,27.95)(83.5,27.29)(83.91,26.67)(84.33,26.01)(84.76,25.36)(85.19,24.7)(85.66,24.04)
\path(85.66,24.04)(86.12,23.42)(86.62,22.78)(87.12,22.12)(87.62,21.5)(88.12,20.86)(88.66,20.2)(89.19,19.55)(89.75,18.95)(90.3,18.29)
\path(90.3,18.29)(90.87,17.68)(91.47,17.04)(92.05,16.43)(92.68,15.8)(93.3,15.18)(93.93,14.55)(94.58,13.93)(95.23,13.31)(95.9,12.68)
\path(95.9,12.68)(96.55,12.06)(97.26,11.46)(97.97,10.85)(98.68,10.22)(99.4,9.6)(100.12,9.02)(100.87,8.39)(101.65,7.8)(102.41,7.19)
\path(102.41,7.19)(103.19,6.59)(103.98,6.0)(104.0,6.0)
\thinlines
\path(78.0,64.0)(78.0,64.0)(77.23,63.99)(76.48,63.97)(75.73,63.95)(74.97,63.93)(74.22,63.9)(73.48,63.86)(72.73,63.83)(71.98,63.77)
\path(71.98,63.77)(71.25,63.72)(70.51,63.65)(69.76,63.58)(69.05,63.5)(68.3,63.4)(67.58,63.33)(66.87,63.22)(66.12,63.11)(65.41,63.0)
\path(65.41,63.0)(64.69,62.88)(63.99,62.77)(63.27,62.63)(62.56,62.5)(61.86,62.34)(61.15,62.2)(60.45,62.04)(59.75,61.86)(59.04,61.7)
\path(59.04,61.7)(58.34,61.52)(57.65,61.33)(56.95,61.13)(56.27,60.93)(55.59,60.72)(54.9,60.5)(54.22,60.29)(53.54,60.06)(52.86,59.83)
\path(52.86,59.83)(52.18,59.59)(51.52,59.34)(50.84,59.09)(50.18,58.81)(49.52,58.56)(48.84,58.27)(48.18,58.0)(47.52,57.7)(46.88,57.4)
\path(46.88,57.4)(46.22,57.11)(45.56,56.79)(44.93,56.47)(44.27,56.15)(43.63,55.83)(43.0,55.5)(42.36,55.15)(41.72,54.79)(41.09,54.43)
\path(41.09,54.43)(40.45,54.08)(39.83,53.7)(39.2,53.33)(38.56,52.95)(37.95,52.56)(37.33,52.15)(36.72,51.75)(36.09,51.34)(35.49,50.93)
\path(35.49,50.93)(34.88,50.5)(34.27,50.06)(33.65,49.63)(33.06,49.18)(32.45,48.72)(31.86,48.27)(31.27,47.81)(30.68,47.34)(30.05,46.86)
\path(30.05,46.86)(29.5,46.36)(28.88,45.88)(28.3,45.38)(27.75,44.86)(27.17,44.36)(26.59,43.84)(26.02,43.31)(25.44,42.77)(24.87,42.24)
\path(24.87,42.24)(24.29,41.68)(23.71,41.13)(23.18,40.56)(22.62,40.0)(22.04,39.43)(21.51,38.84)(20.95,38.25)(20.38,37.65)(19.86,37.06)
\path(19.86,37.06)(19.29,36.45)(18.77,35.84)(18.2,35.22)(17.69,34.59)(17.13,33.95)(16.62,33.31)(16.09,32.65)(15.56,32.0)(15.02,31.34)
\path(15.02,31.34)(14.52,30.67)(14.0,30.0)(14.0,30.0)
\path(78.0,64.0)(78.0,64.0)(78.98,64.04)(79.98,64.08)(80.97,64.12)(81.94,64.19)(82.91,64.26)(83.87,64.33)(84.83,64.4)(85.8,64.48)
\path(85.8,64.48)(86.75,64.55)(87.69,64.65)(88.62,64.75)(89.55,64.83)(90.48,64.94)(91.41,65.05)(92.3,65.18)(93.23,65.3)(94.12,65.43)
\path(94.12,65.43)(95.01,65.55)(95.91,65.69)(96.8,65.83)(97.66,65.98)(98.54,66.12)(99.41,66.29)(100.26,66.44)(101.12,66.62)(101.97,66.79)
\path(101.97,66.79)(102.8,66.97)(103.62,67.15)(104.47,67.33)(105.3,67.54)(106.11,67.73)(106.91,67.94)(107.73,68.15)(108.51,68.36)(109.3,68.58)
\path(109.3,68.58)(110.11,68.8)(110.87,69.01)(111.66,69.26)(112.43,69.51)(113.19,69.76)(113.94,70.01)(114.69,70.26)(115.44,70.51)(116.19,70.79)
\path(116.19,70.79)(116.91,71.05)(117.65,71.33)(118.37,71.62)(119.08,71.91)(119.79,72.19)(120.48,72.48)(121.19,72.8)(121.87,73.11)(122.55,73.41)
\path(122.55,73.41)(123.25,73.73)(123.91,74.05)(124.58,74.37)(125.25,74.72)(125.9,75.05)(126.55,75.41)(127.19,75.75)(127.83,76.11)(128.46,76.47)
\path(128.46,76.47)(129.08,76.83)(129.71,77.19)(130.32,77.58)(130.91,77.94)(131.52,78.33)(132.11,78.73)(132.71,79.12)(133.27,79.51)(133.86,79.94)
\path(133.86,79.94)(134.44,80.33)(135.0,80.76)(135.57,81.19)(136.11,81.62)(136.66,82.05)(137.21,82.48)(137.74,82.93)(138.27,83.37)(138.77,83.83)
\path(138.77,83.83)(139.3,84.29)(139.82,84.76)(140.33,85.23)(140.83,85.69)(141.32,86.18)(141.8,86.66)(142.27,87.15)(142.75,87.65)(143.22,88.15)
\path(143.22,88.15)(143.69,88.65)(144.13,89.16)(144.6,89.68)(145.05,90.19)(145.49,90.73)(145.91,91.26)(146.35,91.8)(146.77,92.33)(147.16,92.87)
\path(147.16,92.87)(147.58,93.44)(147.99,93.98)(148.0,94.0)
\path(76.0,-8.0)(76.0,-8.0)(75.04,-8.27)(74.08,-8.52)(73.15,-8.77)(72.19,-9.01)(71.26,-9.22)(70.36,-9.43)(69.44,-9.63)(68.51,-9.81)
\path(68.51,-9.81)(67.62,-9.97)(66.73,-10.13)(65.83,-10.27)(64.94,-10.39)(64.08,-10.52)(63.22,-10.6)(62.36,-10.71)(61.5,-10.77)(60.65,-10.85)
\path(60.65,-10.85)(59.81,-10.89)(58.97,-10.93)(58.15,-10.96)(57.33,-10.96)(56.52,-10.96)(55.7,-10.93)(54.9,-10.89)(54.11,-10.85)(53.33,-10.81)
\path(53.33,-10.81)(52.54,-10.72)(51.77,-10.64)(51.0,-10.56)(50.25,-10.46)(49.5,-10.31)(48.75,-10.18)(48.02,-10.05)(47.29,-9.89)(46.56,-9.71)
\path(46.56,-9.71)(45.84,-9.52)(45.13,-9.31)(44.4,-9.1)(43.72,-8.88)(43.04,-8.64)(42.34,-8.38)(41.65,-8.1)(41.0,-7.82)(40.34,-7.53)
\path(40.34,-7.53)(39.68,-7.23)(39.02,-6.9)(38.38,-6.57)(37.75,-6.23)(37.11,-5.86)(36.5,-5.5)(35.88,-5.11)(35.27,-4.71)(34.65,-4.3)
\path(34.65,-4.3)(34.06,-3.84)(33.47,-3.43)(32.9,-2.98)(32.31,-2.5)(31.75,-2.0)(31.19,-1.51)(30.62,-1.01)(30.09,-0.5)(29.54,0.0)
\path(29.54,0.0)(29.01,0.55)(28.45,1.11)(27.95,1.65)(27.45,2.25)(26.94,2.83)(26.44,3.47)(25.94,4.09)(25.45,4.73)(24.95,5.38)
\path(24.95,5.38)(24.5,6.05)(24.03,6.73)(23.54,7.42)(23.12,8.1)(22.67,8.84)(22.2,9.56)(21.79,10.31)(21.37,11.06)(20.95,11.81)
\path(20.95,11.81)(20.54,12.6)(20.12,13.39)(19.71,14.22)(19.35,15.02)(18.95,15.88)(18.55,16.7)(18.2,17.59)(17.84,18.45)(17.46,19.35)
\path(17.46,19.35)(17.12,20.25)(16.79,21.17)(16.45,22.1)(16.12,23.04)(15.8,23.96)(15.47,24.95)(15.17,25.94)(14.85,26.93)(14.56,27.94)
\path(14.56,27.94)(14.27,28.95)(14.0,29.96)(14.0,30.0)
\path(76.0,-8.0)(76.0,-8.0)(76.26,-7.48)(76.55,-6.96)(76.83,-6.46)(77.12,-5.96)(77.43,-5.46)(77.73,-4.96)(78.01,-4.48)(78.33,-4.0)
\path(78.33,-4.0)(78.62,-3.5)(78.94,-3.03)(79.26,-2.56)(79.58,-2.13)(79.91,-1.64)(80.23,-1.22)(80.55,-0.75)(80.87,-0.34)(81.22,0.07)
\path(81.22,0.07)(81.55,0.5)(81.87,0.93)(82.23,1.36)(82.58,1.75)(82.93,2.18)(83.26,2.55)(83.62,2.98)(84.0,3.34)(84.36,3.75)
\path(84.36,3.75)(84.72,4.13)(85.08,4.51)(85.44,4.88)(85.83,5.25)(86.19,5.61)(86.58,5.96)(86.98,6.32)(87.36,6.67)(87.76,7.01)
\path(87.76,7.01)(88.15,7.34)(88.55,7.67)(88.94,8.0)(89.33,8.31)(89.75,8.63)(90.16,8.93)(90.58,9.25)(90.98,9.55)(91.41,9.84)
\path(91.41,9.84)(91.83,10.13)(92.26,10.39)(92.69,10.68)(93.12,10.96)(93.55,11.22)(93.98,11.47)(94.44,11.75)(94.87,12.0)(95.33,12.25)
\path(95.33,12.25)(95.76,12.47)(96.23,12.72)(96.69,12.96)(97.15,13.18)(97.62,13.39)(98.08,13.6)(98.55,13.81)(99.01,14.02)(99.51,14.22)
\path(99.51,14.22)(99.98,14.43)(100.47,14.63)(100.94,14.81)(101.44,14.97)(101.94,15.14)(102.43,15.31)(102.93,15.5)(103.43,15.64)(103.94,15.81)
\path(103.94,15.81)(104.44,15.96)(104.94,16.1)(105.48,16.21)(105.98,16.37)(106.51,16.5)(107.04,16.62)(107.55,16.71)(108.08,16.85)(108.62,16.95)
\path(108.62,16.95)(109.16,17.04)(109.69,17.12)(110.26,17.21)(110.8,17.3)(111.33,17.38)(111.91,17.46)(112.47,17.54)(113.01,17.62)(113.58,17.68)
\path(114.0,18.0)(114.0,18.0)(114.3,18.04)(114.62,18.05)(114.94,18.12)(115.26,18.13)(115.58,18.2)(115.9,18.21)(116.22,18.28)(116.51,18.3)
\path(116.51,18.3)(116.83,18.37)(117.15,18.42)(117.47,18.45)(117.76,18.5)(118.08,18.54)(118.4,18.59)(118.69,18.62)(119.01,18.69)(119.3,18.7)
\path(119.3,18.7)(119.62,18.78)(119.93,18.8)(120.23,18.87)(120.54,18.92)(120.83,18.95)(121.12,19.02)(121.44,19.04)(121.75,19.12)(122.04,19.17)
\path(122.04,19.17)(122.33,19.2)(122.62,19.27)(122.94,19.29)(123.23,19.37)(123.51,19.43)(123.83,19.45)(124.12,19.53)(124.41,19.59)(124.69,19.62)
\path(124.69,19.62)(125.0,19.69)(125.29,19.75)(125.58,19.79)(125.87,19.86)(126.16,19.92)(126.44,19.95)(126.73,20.03)(127.01,20.05)(127.3,20.12)
\path(127.3,20.12)(127.58,20.2)(127.87,20.26)(128.13,20.29)(128.41,20.37)(128.71,20.44)(129.0,20.5)(129.27,20.54)(129.55,20.62)(129.83,20.68)
\path(129.83,20.68)(130.11,20.71)(130.38,20.79)(130.66,20.86)(130.94,20.92)(131.21,20.96)(131.47,21.04)(131.75,21.11)(132.02,21.18)(132.3,21.21)
\path(141.52,23.79)(141.52,23.79)(141.69,23.87)(141.86,23.95)(142.02,24.04)(142.21,24.12)(142.38,24.21)(142.55,24.29)(142.72,24.37)(142.88,24.46)
\path(142.88,24.46)(143.05,24.54)(143.21,24.62)(143.38,24.7)(143.52,24.79)(143.69,24.87)(143.86,24.95)(144.02,25.04)(144.16,25.12)(144.33,25.2)
\path(153.41,31.04)(154.07,31.46)(154.71,31.92)(155.35,32.36)(156.0,32.81)(156.66,33.25)(157.3,33.72)(157.97,34.18)(158.63,34.65)(159.3,35.13)
\path(159.3,35.13)(159.97,35.61)(160.63,36.11)(161.32,36.61)(161.99,37.11)(162.66,37.63)(163.36,38.15)(164.05,38.68)(164.74,39.2)(165.41,39.75)
\path(165.41,39.75)(166.13,40.29)(166.83,40.86)(167.52,41.4)(168.24,41.99)(168.94,42.56)(169.66,43.13)(170.38,43.72)(171.11,44.31)(171.83,44.93)
\path(171.83,44.93)(172.55,45.52)(173.27,46.15)(174.02,46.77)(174.75,47.4)(175.49,48.04)(176.22,48.68)(176.97,49.33)(177.72,49.99)(178.47,50.65)
\path(178.47,50.65)(179.24,51.31)(179.99,51.99)(180.0,52.0)
\path(180.0,52.0)(180.0,52.0)(179.38,52.27)(178.8,52.56)(178.22,52.84)(177.63,53.13)(177.07,53.43)(176.5,53.72)(175.91,54.02)(175.36,54.31)
\path(175.36,54.31)(174.82,54.63)(174.27,54.93)(173.72,55.24)(173.19,55.56)(172.66,55.86)(172.13,56.18)(171.63,56.5)(171.11,56.83)(170.6,57.15)
\path(170.6,57.15)(170.1,57.49)(169.61,57.81)(169.11,58.15)(168.63,58.49)(168.13,58.83)(167.66,59.18)(167.21,59.52)(166.75,59.86)(166.27,60.22)
\path(166.27,60.22)(165.83,60.58)(165.38,60.93)(164.94,61.29)(164.52,61.65)(164.08,62.02)(163.66,62.38)(163.24,62.75)(162.83,63.13)(162.41,63.5)
\path(162.41,63.5)(162.02,63.88)(161.63,64.26)(161.24,64.66)(160.85,65.04)(160.47,65.43)(160.1,65.83)(159.72,66.22)(159.36,66.62)(159.02,67.01)
\path(159.02,67.01)(158.66,67.43)(158.32,67.83)(157.97,68.25)(157.63,68.66)(157.32,69.08)(157.0,69.48)(156.66,69.91)(156.36,70.33)(156.05,70.76)
\path(156.05,70.76)(155.75,71.19)(155.47,71.62)(155.16,72.05)(154.88,72.5)(154.61,72.94)(154.33,73.37)(154.08,73.83)(153.8,74.26)(153.55,74.73)
\path(153.55,74.73)(153.3,75.19)(153.05,75.65)(152.83,76.11)(152.58,76.55)(152.36,77.04)(152.13,77.51)(151.91,77.98)(151.72,78.44)(151.5,78.93)
\path(151.5,78.93)(151.3,79.41)(151.11,79.9)(150.91,80.37)(150.75,80.87)(150.57,81.36)(150.38,81.86)(150.22,82.33)(150.07,82.83)(149.91,83.33)
\path(149.91,83.33)(149.77,83.86)(149.61,84.37)(149.47,84.87)(149.35,85.37)(149.22,85.91)(149.1,86.43)(148.99,86.94)(148.88,87.48)(148.77,88.0)
\path(148.77,88.0)(148.66,88.51)(148.58,89.05)(148.49,89.58)(148.41,90.12)(148.33,90.69)(148.27,91.23)(148.19,91.76)(148.13,92.33)(148.08,92.87)
\path(148.08,92.87)(148.02,93.44)(148.0,93.98)(148.0,94.0)
\path(134.0,80.0)(134.0,80.0)(134.02,79.48)(134.08,78.94)(134.13,78.44)(134.19,77.93)(134.25,77.41)(134.3,76.9)(134.38,76.37)(134.46,75.87)
\path(134.46,75.87)(134.52,75.37)(134.61,74.87)(134.69,74.37)(134.77,73.87)(134.88,73.37)(134.99,72.87)(135.08,72.37)(135.19,71.87)(135.3,71.37)
\path(135.3,71.37)(135.41,70.87)(135.55,70.4)(135.66,69.91)(135.8,69.43)(135.94,68.94)(136.08,68.44)(136.22,67.98)(136.36,67.5)(136.52,67.01)
\path(136.52,67.01)(136.66,66.54)(136.83,66.05)(137.0,65.58)(137.16,65.11)(137.35,64.62)(137.52,64.16)(137.71,63.7)(137.88,63.24)(138.08,62.77)
\path(138.08,62.77)(138.27,62.31)(138.49,61.84)(138.69,61.38)(138.88,60.93)(139.11,60.47)(139.33,60.02)(139.55,59.56)(139.77,59.11)(140.0,58.65)
\path(140.0,58.65)(140.25,58.22)(140.49,57.77)(140.72,57.31)(140.97,56.88)(141.24,56.43)(141.5,56.0)(141.75,55.56)(142.02,55.11)(142.27,54.68)
\path(142.27,54.68)(142.57,54.25)(142.85,53.81)(143.13,53.38)(143.41,52.95)(143.72,52.52)(144.0,52.09)(144.3,51.68)(144.61,51.25)(144.91,50.83)
\path(144.91,50.83)(145.25,50.4)(145.57,49.99)(145.88,49.58)(146.22,49.15)(146.55,48.75)(146.88,48.33)(147.22,47.9)(147.57,47.52)(147.91,47.11)
\path(147.91,47.11)(148.27,46.7)(148.63,46.29)(149.0,45.9)(149.36,45.5)(149.74,45.09)(150.11,44.7)(150.5,44.29)(150.88,43.9)(151.27,43.52)
\path(151.27,43.52)(151.66,43.11)(152.07,42.72)(152.47,42.34)(152.88,41.95)(153.27,41.58)(153.71,41.18)(154.13,40.81)(154.55,40.43)(154.97,40.04)
\path(154.97,40.04)(155.41,39.68)(155.85,39.29)(156.3,38.93)(156.74,38.54)(157.19,38.18)(157.63,37.81)(158.11,37.45)(158.57,37.08)(159.02,36.72)
\path(159.02,36.72)(159.52,36.36)(159.99,36.0)(160.0,36.0)
\path(114.0,70.0)(114.0,70.0)(114.04,69.36)(114.08,68.72)(114.12,68.08)(114.18,67.44)(114.23,66.83)(114.3,66.19)(114.36,65.58)(114.43,64.98)
\path(114.43,64.98)(114.5,64.36)(114.55,63.75)(114.65,63.15)(114.73,62.54)(114.8,61.95)(114.91,61.34)(115.0,60.75)(115.08,60.15)(115.19,59.58)
\path(115.19,59.58)(115.3,58.99)(115.4,58.4)(115.51,57.83)(115.62,57.25)(115.75,56.68)(115.87,56.11)(115.98,55.56)(116.12,55.0)(116.25,54.43)
\path(116.25,54.43)(116.37,53.88)(116.51,53.33)(116.66,52.77)(116.8,52.24)(116.94,51.68)(117.12,51.15)(117.26,50.61)(117.44,50.08)(117.58,49.56)
\path(117.58,49.56)(117.76,49.02)(117.94,48.5)(118.11,47.99)(118.29,47.47)(118.48,46.95)(118.66,46.43)(118.83,45.93)(119.04,45.43)(119.23,44.93)
\path(119.23,44.93)(119.44,44.43)(119.62,43.93)(119.83,43.45)(120.05,42.95)(120.26,42.47)(120.48,42.0)(120.72,41.52)(120.94,41.04)(121.16,40.56)
\path(121.16,40.56)(121.4,40.09)(121.62,39.63)(121.87,39.15)(122.12,38.7)(122.37,38.25)(122.62,37.79)(122.87,37.36)(123.12,36.9)(123.37,36.47)
\path(123.37,36.47)(123.66,36.02)(123.93,35.59)(124.19,35.15)(124.48,34.72)(124.76,34.29)(125.04,33.86)(125.3,33.45)(125.61,33.04)(125.91,32.61)
\path(125.91,32.61)(126.19,32.2)(126.51,31.79)(126.8,31.37)(127.12,31.0)(127.43,30.6)(127.75,30.2)(128.07,29.79)(128.38,29.42)(128.71,29.04)
\path(128.71,29.04)(129.02,28.62)(129.38,28.27)(129.72,27.87)(130.05,27.52)(130.38,27.13)(130.75,26.79)(131.1,26.43)(131.44,26.04)(131.8,25.7)
\path(131.8,25.7)(132.16,25.36)(132.52,25.0)(132.91,24.63)(133.27,24.29)(133.66,23.95)(134.02,23.62)(134.41,23.29)(134.8,22.95)(135.19,22.62)
\path(132.3,21.21)(132.57,21.29)(132.83,21.37)(133.11,21.44)(133.36,21.51)(133.63,21.54)(133.91,21.62)(134.16,21.7)(134.41,21.77)(134.69,21.84)
\path(96.0,66.0)(96.0,66.0)(95.87,65.48)(95.76,64.94)(95.66,64.44)(95.55,63.9)(95.44,63.4)(95.37,62.88)(95.26,62.36)(95.19,61.86)
\path(95.19,61.86)(95.12,61.34)(95.05,60.84)(94.98,60.31)(94.93,59.81)(94.87,59.29)(94.8,58.79)(94.76,58.29)(94.73,57.77)(94.69,57.27)
\path(94.69,57.27)(94.68,56.75)(94.65,56.25)(94.62,55.75)(94.62,55.25)(94.61,54.75)(94.61,54.25)(94.61,53.75)(94.62,53.25)(94.62,52.75)
\path(94.62,52.75)(94.65,52.25)(94.66,51.75)(94.69,51.25)(94.73,50.75)(94.76,50.25)(94.8,49.75)(94.87,49.27)(94.91,48.77)(94.98,48.29)
\path(94.98,48.29)(95.04,47.79)(95.11,47.29)(95.19,46.81)(95.26,46.31)(95.36,45.84)(95.44,45.34)(95.54,44.86)(95.62,44.36)(95.75,43.88)
\path(95.75,43.88)(95.86,43.4)(95.98,42.9)(96.08,42.43)(96.23,41.95)(96.36,41.47)(96.48,41.0)(96.62,40.52)(96.79,40.04)(96.94,39.56)
\path(96.94,39.56)(97.08,39.08)(97.26,38.61)(97.43,38.13)(97.58,37.65)(97.76,37.18)(97.97,36.7)(98.15,36.24)(98.33,35.75)(98.55,35.29)
\path(98.55,35.29)(98.75,34.81)(98.94,34.34)(99.18,33.88)(99.4,33.4)(99.62,32.95)(99.86,32.47)(100.08,32.02)(100.33,31.54)(100.58,31.09)
\path(100.58,31.09)(100.83,30.62)(101.08,30.17)(101.33,29.7)(101.62,29.25)(101.87,28.79)(102.16,28.3)(102.44,27.87)(102.73,27.38)(103.01,26.95)
\path(103.01,26.95)(103.33,26.5)(103.62,26.04)(103.94,25.59)(104.26,25.12)(104.58,24.69)(104.9,24.2)(105.23,23.78)(105.55,23.3)(105.91,22.87)
\path(105.91,22.87)(106.25,22.44)(106.61,21.96)(106.94,21.54)(107.3,21.09)(107.69,20.62)(108.05,20.2)(108.44,19.76)(108.8,19.29)(109.19,18.87)
\path(109.19,18.87)(109.58,18.44)(109.98,18.0)(110.0,18.0)
\path(60.0,62.0)(60.0,62.0)(59.79,61.38)(59.61,60.79)(59.43,60.18)(59.25,59.59)(59.09,58.99)(58.93,58.38)(58.79,57.79)(58.65,57.18)
\path(58.65,57.18)(58.52,56.58)(58.4,55.97)(58.27,55.36)(58.15,54.77)(58.06,54.15)(57.97,53.56)(57.9,52.95)(57.81,52.34)(57.75,51.74)
\path(57.75,51.74)(57.68,51.13)(57.63,50.52)(57.59,49.9)(57.56,49.31)(57.52,48.7)(57.5,48.09)(57.5,47.47)(57.5,46.86)(57.5,46.25)
\path(57.5,46.25)(57.5,45.65)(57.52,45.04)(57.56,44.43)(57.59,43.81)(57.63,43.2)(57.68,42.59)(57.75,41.97)(57.81,41.36)(57.88,40.75)
\path(57.88,40.75)(57.97,40.13)(58.06,39.52)(58.15,38.9)(58.27,38.29)(58.38,37.68)(58.52,37.06)(58.65,36.43)(58.79,35.83)(58.93,35.2)
\path(58.93,35.2)(59.09,34.59)(59.25,33.97)(59.43,33.34)(59.61,32.72)(59.79,32.11)(59.99,31.5)(60.2,30.87)(60.4,30.25)(60.63,29.62)
\path(60.63,29.62)(60.86,29.01)(61.09,28.37)(61.34,27.77)(61.59,27.12)(61.84,26.52)(62.11,25.87)(62.38,25.28)(62.68,24.62)(62.97,24.03)
\path(62.97,24.03)(63.27,23.37)(63.58,22.78)(63.88,22.12)(64.22,21.52)(64.55,20.87)(64.87,20.27)(65.23,19.62)(65.58,19.02)(65.94,18.37)
\path(65.94,18.37)(66.33,17.76)(66.69,17.12)(67.08,16.5)(67.48,15.85)(67.9,15.22)(68.3,14.6)(68.73,13.97)(69.16,13.35)(69.58,12.72)
\path(69.58,12.72)(70.04,12.06)(70.48,11.43)(70.94,10.81)(71.41,10.18)(71.87,9.55)(72.37,8.92)(72.87,8.27)(73.37,7.65)(73.87,7.01)
\path(73.87,7.01)(74.37,6.38)(74.91,5.73)(75.44,5.09)(75.98,4.46)(76.54,3.8)(77.08,3.19)(77.66,2.52)(78.23,1.88)(78.8,1.25)
\path(78.8,1.25)(79.4,0.62)(79.98,0.0)(80.0,0.0)
\path(34.0,50.0)(34.0,50.0)(33.95,49.31)(33.93,48.63)(33.9,47.95)(33.88,47.29)(33.88,46.61)(33.88,45.93)(33.88,45.27)(33.88,44.61)
\path(33.88,44.61)(33.9,43.93)(33.93,43.27)(33.97,42.61)(34.0,41.95)(34.04,41.29)(34.09,40.63)(34.15,39.97)(34.22,39.31)(34.29,38.65)
\path(34.29,38.65)(34.38,38.0)(34.45,37.36)(34.54,36.72)(34.65,36.06)(34.75,35.4)(34.86,34.77)(34.99,34.13)(35.11,33.5)(35.25,32.86)
\path(35.25,32.86)(35.38,32.22)(35.54,31.55)(35.68,30.95)(35.86,30.29)(36.02,29.68)(36.2,29.04)(36.38,28.43)(36.56,27.79)(36.75,27.18)
\path(36.75,27.18)(36.95,26.54)(37.15,25.93)(37.38,25.29)(37.61,24.69)(37.83,24.05)(38.06,23.45)(38.31,22.85)(38.56,22.2)(38.81,21.62)
\path(38.81,21.62)(39.08,21.02)(39.34,20.38)(39.63,19.79)(39.9,19.2)(40.2,18.6)(40.49,18.0)(40.79,17.37)(41.11,16.79)(41.43,16.2)
\path(41.43,16.2)(41.75,15.6)(42.08,15.02)(42.4,14.42)(42.75,13.81)(43.11,13.25)(43.47,12.64)(43.83,12.06)(44.2,11.47)(44.58,10.89)
\path(44.58,10.89)(44.97,10.31)(45.36,9.75)(45.75,9.18)(46.15,8.6)(46.58,8.02)(47.0,7.44)(47.4,6.88)(47.84,6.32)(48.29,5.75)
\path(48.29,5.75)(48.74,5.17)(49.18,4.61)(49.65,4.05)(50.11,3.5)(50.59,2.94)(51.06,2.38)(51.56,1.82)(52.04,1.25)(52.54,0.72)
\path(52.54,0.72)(53.06,0.15)(53.58,-0.37)(54.09,-0.92)(54.63,-1.47)(55.15,-2.0)(55.7,-2.53)(56.25,-3.07)(56.79,-3.64)(57.36,-4.17)
\path(57.36,-4.17)(57.93,-4.71)(58.5,-5.25)(59.09,-5.78)(59.68,-6.32)(60.27,-6.84)(60.88,-7.36)(61.49,-7.9)(62.11,-8.43)(62.72,-8.93)
\path(62.72,-8.93)(63.36,-9.47)(63.99,-9.97)(64.0,-10.0)
\path(18.0,18.0)(18.0,18.0)(18.62,18.59)(19.27,19.18)(19.88,19.77)(20.54,20.35)(21.18,20.92)(21.79,21.46)(22.45,22.04)(23.05,22.6)
\path(23.05,22.6)(23.7,23.12)(24.3,23.69)(24.95,24.2)(25.59,24.76)(26.2,25.29)(26.84,25.79)(27.45,26.29)(28.05,26.8)(28.7,27.3)
\path(28.7,27.3)(29.29,27.79)(29.94,28.29)(30.54,28.79)(31.17,29.27)(31.78,29.71)(32.4,30.2)(33.0,30.67)(33.61,31.12)(34.22,31.54)
\path(34.22,31.54)(34.84,32.0)(35.43,32.43)(36.04,32.86)(36.65,33.29)(37.25,33.7)(37.86,34.11)(38.45,34.52)(39.06,34.93)(39.65,35.31)
\path(39.65,35.31)(40.25,35.7)(40.84,36.09)(41.45,36.45)(42.04,36.83)(42.63,37.18)(43.22,37.54)(43.81,37.9)(44.4,38.25)(44.99,38.59)
\path(44.99,38.59)(45.58,38.9)(46.15,39.25)(46.75,39.56)(47.33,39.88)(47.9,40.18)(48.49,40.49)(49.06,40.79)(49.65,41.08)(50.22,41.36)
\path(50.22,41.36)(50.81,41.65)(51.38,41.9)(51.95,42.18)(52.52,42.45)(53.09,42.7)(53.65,42.95)(54.22,43.18)(54.79,43.43)(55.36,43.65)
\path(55.36,43.65)(55.93,43.88)(56.5,44.11)(57.06,44.31)(57.61,44.52)(58.18,44.72)(58.74,44.9)(59.29,45.11)(59.84,45.29)(60.4,45.47)
\path(60.4,45.47)(60.95,45.63)(61.52,45.81)(62.06,45.97)(62.61,46.11)(63.15,46.27)(63.72,46.4)(64.26,46.54)(64.8,46.65)(65.33,46.79)
\path(65.33,46.79)(65.9,46.9)(66.44,47.02)(66.98,47.13)(67.51,47.22)(68.05,47.31)(68.58,47.4)(69.12,47.49)(69.66,47.56)(70.19,47.63)
\path(70.19,47.63)(70.73,47.68)(71.26,47.75)(71.8,47.79)(72.33,47.84)(72.83,47.88)(73.37,47.9)(73.91,47.95)(74.43,47.97)(74.94,47.97)
\path(74.94,47.97)(75.47,47.99)(75.98,48.0)(76.0,48.0)
\path(76.0,48.0)(76.0,48.0)(77.11,48.08)(78.22,48.15)(79.3,48.25)(80.41,48.34)(81.51,48.43)(82.58,48.54)(83.66,48.63)(84.73,48.75)
\path(84.73,48.75)(85.8,48.86)(86.86,48.97)(87.9,49.09)(88.94,49.2)(89.98,49.34)(91.01,49.47)(92.01,49.59)(93.04,49.74)(94.05,49.88)
\path(94.05,49.88)(95.05,50.02)(96.05,50.15)(97.04,50.31)(98.01,50.47)(98.98,50.63)(99.94,50.79)(100.91,50.95)(101.87,51.11)(102.8,51.29)
\path(102.8,51.29)(103.76,51.47)(104.69,51.65)(105.62,51.83)(106.54,52.02)(107.44,52.2)(108.33,52.4)(109.25,52.59)(110.12,52.79)(111.01,53.0)
\path(111.01,53.0)(111.91,53.2)(112.76,53.4)(113.65,53.63)(114.5,53.84)(115.33,54.06)(116.19,54.29)(117.04,54.52)(117.87,54.75)(118.69,55.0)
\path(118.69,55.0)(119.51,55.24)(120.3,55.47)(121.12,55.72)(121.91,55.97)(122.69,56.24)(123.48,56.49)(124.26,56.75)(125.04,57.02)(125.8,57.29)
\path(125.8,57.29)(126.55,57.56)(127.3,57.84)(128.05,58.11)(128.77,58.4)(129.52,58.68)(130.24,58.97)(130.94,59.27)(131.66,59.56)(132.36,59.86)
\path(132.36,59.86)(133.05,60.18)(133.75,60.49)(134.41,60.79)(135.1,61.11)(135.77,61.43)(136.41,61.75)(137.08,62.08)(137.72,62.4)(138.38,62.75)
\path(138.38,62.75)(139.0,63.09)(139.63,63.43)(140.25,63.77)(140.86,64.12)(141.47,64.47)(142.08,64.83)(142.66,65.19)(143.25,65.55)(143.83,65.91)
\path(143.83,65.91)(144.41,66.26)(144.97,66.66)(145.52,67.04)(146.08,67.41)(146.63,67.8)(147.16,68.19)(147.69,68.58)(148.22,68.97)(148.74,69.37)
\path(148.74,69.37)(149.25,69.76)(149.75,70.18)(150.25,70.58)(150.75,71.0)(151.22,71.41)(151.71,71.83)(152.16,72.26)(152.63,72.69)(153.1,73.12)
\path(153.1,73.12)(153.55,73.55)(153.99,73.98)(154.0,74.0)
\path(80.0,34.0)(80.0,34.0)(79.3,33.95)(78.62,33.9)(77.97,33.84)(77.3,33.79)(76.62,33.72)(75.97,33.65)(75.3,33.58)(74.66,33.5)
\path(74.66,33.5)(74.0,33.4)(73.36,33.31)(72.69,33.22)(72.05,33.11)(71.43,33.0)(70.79,32.88)(70.16,32.77)(69.51,32.63)(68.9,32.5)
\path(68.9,32.5)(68.26,32.36)(67.65,32.22)(67.04,32.06)(66.41,31.92)(65.8,31.76)(65.19,31.59)(64.58,31.42)(64.0,31.25)(63.4,31.04)
\path(63.4,31.04)(62.79,30.87)(62.2,30.68)(61.61,30.45)(61.04,30.28)(60.45,30.04)(59.86,29.85)(59.29,29.62)(58.72,29.37)(58.15,29.17)
\path(58.15,29.17)(57.59,28.93)(57.02,28.68)(56.47,28.43)(55.9,28.18)(55.36,27.92)(54.79,27.62)(54.25,27.37)(53.7,27.1)(53.15,26.79)
\path(53.15,26.79)(52.63,26.53)(52.09,26.2)(51.56,25.93)(51.04,25.62)(50.52,25.29)(50.0,25.0)(49.47,24.67)(48.95,24.34)(48.45,24.01)
\path(48.45,24.01)(47.93,23.67)(47.43,23.3)(46.93,22.95)(46.43,22.62)(45.93,22.26)(45.43,21.87)(44.95,21.52)(44.47,21.12)(43.99,20.75)
\path(43.99,20.75)(43.5,20.36)(43.02,19.95)(42.56,19.54)(42.08,19.13)(41.61,18.75)(41.15,18.3)(40.68,17.87)(40.24,17.45)(39.77,17.04)
\path(39.77,17.04)(39.33,16.6)(38.88,16.12)(38.43,15.68)(38.0,15.25)(37.56,14.77)(37.11,14.31)(36.68,13.84)(36.25,13.35)(35.84,12.88)
\path(35.84,12.88)(35.4,12.38)(34.99,11.89)(34.58,11.39)(34.15,10.88)(33.75,10.35)(33.34,9.85)(32.95,9.31)(32.54,8.77)(32.15,8.26)
\path(32.15,8.26)(31.76,7.71)(31.36,7.17)(30.95,6.61)(30.59,6.05)(30.2,5.48)(29.84,4.92)(29.45,4.34)(29.09,3.75)(28.7,3.18)
\path(28.7,3.18)(28.36,2.56)(28.0,2.0)(28.0,2.0)
\path(80.0,34.0)(80.0,34.0)(81.15,34.08)(82.3,34.15)(83.44,34.25)(84.58,34.34)(85.69,34.45)(86.83,34.54)(87.94,34.65)(89.05,34.77)
\path(89.05,34.77)(90.16,34.88)(91.25,35.02)(92.33,35.13)(93.43,35.27)(94.5,35.4)(95.55,35.54)(96.62,35.68)(97.68,35.84)(98.73,35.99)
\path(98.73,35.99)(99.76,36.15)(100.8,36.31)(101.83,36.47)(102.86,36.65)(103.87,36.81)(104.87,37.0)(105.87,37.18)(106.87,37.36)(107.86,37.56)
\path(107.86,37.56)(108.83,37.75)(109.8,37.95)(110.76,38.15)(111.73,38.38)(112.69,38.59)(113.62,38.81)(114.55,39.02)(115.5,39.25)(116.43,39.49)
\path(116.43,39.49)(117.33,39.72)(118.26,39.97)(119.16,40.2)(120.05,40.45)(120.94,40.7)(121.83,40.97)(122.72,41.24)(123.58,41.5)(124.44,41.77)
\path(124.44,41.77)(125.3,42.04)(126.16,42.33)(127.0,42.61)(127.83,42.9)(128.66,43.2)(129.49,43.49)(130.3,43.79)(131.11,44.09)(131.91,44.4)
\path(131.91,44.4)(132.72,44.72)(133.5,45.04)(134.27,45.36)(135.07,45.7)(135.83,46.04)(136.6,46.36)(137.35,46.7)(138.1,47.06)(138.85,47.4)
\path(138.85,47.4)(139.58,47.77)(140.3,48.13)(141.02,48.49)(141.74,48.86)(142.44,49.22)(143.13,49.61)(143.85,49.99)(144.52,50.36)(145.22,50.77)
\path(145.22,50.77)(145.88,51.15)(146.55,51.56)(147.22,51.95)(147.86,52.36)(148.52,52.77)(149.16,53.2)(149.77,53.61)(150.41,54.04)(151.02,54.47)
\path(151.02,54.47)(151.63,54.9)(152.25,55.34)(152.85,55.79)(153.44,56.24)(154.02,56.68)(154.61,57.15)(155.16,57.61)(155.75,58.06)(156.3,58.54)
\path(156.3,58.54)(156.85,59.0)(157.38,59.49)(157.94,59.97)(158.47,60.45)(158.99,60.95)(159.5,61.45)(160.02,61.95)(160.52,62.45)(161.02,62.95)
\path(161.02,62.95)(161.5,63.47)(161.99,63.99)(162.0,64.0)
\path(48.0,-10.0)(48.0,-10.0)(48.31,-9.52)(48.63,-9.02)(48.95,-8.56)(49.29,-8.1)(49.61,-7.63)(49.93,-7.17)(50.27,-6.71)(50.61,-6.26)
\path(50.61,-6.26)(50.93,-5.82)(51.27,-5.38)(51.61,-4.92)(51.95,-4.48)(52.29,-4.05)(52.63,-3.63)(52.97,-3.2)(53.31,-2.75)(53.65,-2.33)
\path(53.65,-2.33)(54.0,-1.94)(54.36,-1.5)(54.72,-1.11)(55.06,-0.68)(55.4,-0.31)(55.77,0.07)(56.13,0.46)(56.5,0.87)(56.86,1.25)
\path(56.86,1.25)(57.22,1.62)(57.58,2.0)(57.95,2.4)(58.31,2.75)(58.68,3.15)(59.04,3.5)(59.43,3.84)(59.79,4.23)(60.15,4.59)
\path(60.15,4.59)(60.54,4.94)(60.93,5.28)(61.31,5.63)(61.68,5.98)(62.06,6.3)(62.45,6.65)(62.84,6.98)(63.22,7.3)(63.61,7.63)
\path(63.61,7.63)(64.01,7.94)(64.41,8.27)(64.8,8.56)(65.19,8.89)(65.58,9.18)(65.98,9.47)(66.4,9.77)(66.8,10.09)(67.19,10.38)
\path(67.19,10.38)(67.61,10.67)(68.01,10.93)(68.41,11.22)(68.83,11.51)(69.25,11.77)(69.66,12.05)(70.05,12.31)(70.48,12.56)(70.91,12.84)
\path(70.91,12.84)(71.33,13.09)(71.75,13.34)(72.16,13.59)(72.58,13.81)(73.01,14.06)(73.44,14.31)(73.87,14.55)(74.3,14.77)(74.75,15.0)
\path(74.75,15.0)(75.18,15.22)(75.62,15.43)(76.05,15.64)(76.48,15.85)(76.94,16.05)(77.37,16.28)(77.8,16.45)(78.26,16.68)(78.69,16.87)
\path(78.69,16.87)(79.16,17.04)(79.61,17.25)(80.05,17.43)(80.51,17.61)(80.97,17.79)(81.43,17.95)(81.87,18.12)(82.33,18.29)(82.8,18.45)
\path(82.8,18.45)(83.26,18.61)(83.73,18.77)(84.19,18.92)(84.66,19.04)(85.12,19.2)(85.61,19.35)(86.08,19.46)(86.55,19.62)(87.04,19.75)
\path(87.04,19.75)(87.51,19.87)(87.98,19.96)(88.0,20.0)
\path(88.0,20.0)(88.0,20.0)(89.15,20.2)(90.3,20.37)(91.44,20.61)(92.58,20.79)(93.69,21.04)(94.83,21.26)(95.94,21.45)(97.05,21.7)
\path(97.05,21.7)(98.16,21.94)(99.25,22.18)(100.33,22.38)(101.43,22.62)(102.5,22.87)(103.55,23.12)(104.62,23.37)(105.68,23.63)(106.73,23.92)
\path(106.73,23.92)(107.76,24.18)(108.8,24.44)(109.83,24.7)(110.86,24.96)(111.87,25.27)(112.87,25.54)(113.87,25.8)(114.87,26.12)(115.86,26.38)
\path(115.86,26.38)(116.83,26.7)(117.8,27.01)(118.76,27.29)(119.73,27.61)(120.69,27.92)(121.62,28.21)(122.55,28.54)(123.5,28.87)(124.43,29.2)
\path(124.43,29.2)(125.33,29.53)(126.26,29.86)(127.16,30.19)(128.05,30.53)(128.94,30.87)(129.83,31.2)(130.72,31.54)(131.58,31.92)(132.44,32.27)
\path(132.44,32.27)(133.3,32.63)(134.16,33.0)(135.0,33.36)(135.83,33.74)(136.66,34.11)(137.49,34.49)(138.3,34.88)(139.11,35.25)(139.91,35.65)
\path(139.91,35.65)(140.72,36.04)(141.5,36.43)(142.27,36.84)(143.07,37.24)(143.83,37.65)(144.6,38.06)(145.35,38.47)(146.1,38.88)(146.85,39.31)
\path(146.85,39.31)(147.58,39.74)(148.3,40.15)(149.02,40.59)(149.74,41.04)(150.44,41.47)(151.13,41.9)(151.85,42.36)(152.52,42.81)(153.22,43.27)
\path(153.22,43.27)(153.88,43.72)(154.55,44.18)(155.22,44.65)(155.86,45.11)(156.52,45.59)(157.16,46.06)(157.77,46.54)(158.41,47.02)(159.02,47.5)
\path(159.02,47.5)(159.63,48.0)(160.25,48.5)(160.85,49.0)(161.44,49.5)(162.02,50.0)(162.61,50.5)(163.16,51.02)(163.75,51.52)(164.3,52.04)
\path(164.3,52.04)(164.85,52.56)(165.38,53.09)(165.94,53.63)(166.47,54.15)(166.99,54.7)(167.5,55.24)(168.02,55.77)(168.52,56.33)(169.02,56.88)
\path(169.02,56.88)(169.5,57.43)(169.99,57.99)(170.0,58.0)
\path(134.69,21.84)(134.96,21.88)(135.22,21.95)(135.47,22.04)(135.75,22.12)(136.0,22.19)(136.25,22.26)(136.52,22.3)(136.77,22.37)(137.02,22.45)
\path(137.02,22.45)(137.27,22.54)(137.55,22.62)(137.8,22.69)(138.05,22.77)(138.3,22.84)(138.55,22.88)(138.8,22.96)(139.05,23.04)(139.3,23.12)
\path(139.3,23.12)(139.55,23.2)(139.8,23.29)(140.05,23.37)(140.3,23.44)(140.52,23.52)(140.77,23.6)(141.02,23.68)(141.27,23.76)(141.5,23.84)
\path(144.33,25.2)(144.49,25.29)(144.63,25.37)(144.77,25.45)(144.94,25.53)(145.08,25.61)(145.24,25.69)(145.38,25.77)(145.52,25.84)(145.66,25.92)
\path(145.66,25.92)(145.83,26.0)(145.97,26.04)(146.11,26.12)(146.25,26.2)(146.38,26.29)(146.52,26.37)(146.66,26.45)(146.8,26.53)(146.94,26.6)
\path(146.94,26.6)(147.08,26.68)(147.21,26.75)(147.33,26.79)(147.47,26.87)(147.6,26.95)(147.72,27.04)(147.86,27.11)(147.97,27.19)(148.11,27.26)
\path(148.11,27.26)(148.22,27.3)(148.35,27.37)(148.47,27.45)(148.58,27.54)(148.71,27.61)(148.83,27.68)(148.94,27.75)(149.05,27.79)(149.16,27.87)
\path(149.16,27.87)(149.27,27.95)(149.38,28.02)(149.5,28.09)(149.61,28.13)(149.72,28.2)(149.83,28.29)(149.94,28.36)(150.05,28.43)(150.13,28.46)
\path(150.13,28.46)(150.25,28.54)(150.35,28.62)(150.44,28.69)(150.55,28.75)(150.63,28.79)(150.75,28.87)(150.83,28.95)(150.94,29.01)(151.02,29.04)
\path(151.02,29.04)(151.11,29.12)(151.21,29.2)(151.3,29.26)(151.38,29.29)(151.47,29.37)(151.55,29.45)(151.63,29.51)(151.72,29.54)(151.8,29.62)
\path(151.8,29.62)(151.88,29.69)(151.97,29.75)(152.05,29.79)(152.13,29.87)(152.21,29.93)(152.27,29.96)(152.36,30.04)(152.41,30.1)(152.5,30.13)
\path(152.5,30.13)(152.57,30.2)(152.63,30.27)(152.71,30.3)(152.77,30.37)(152.83,30.44)(152.91,30.5)(152.97,30.54)(153.02,30.61)(153.1,30.67)
\path(153.1,30.67)(153.16,30.7)(153.22,30.77)(153.22,30.78)
\drawthickdot{79.87}{33.93}
\drawhollowdot{79.65}{34.15}
\drawlefttext{82.0}{40.0}{$A$}
\drawdashline{146.0}{-2.0}{-24.0}{-2.0}
\drawdashline{-24.0}{-2.0}{36.0}{72.0}
\drawdashline{146.0}{-2.0}{212.0}{-2.0}
\drawdotline{80.0}{34.0}{194.0}{34.0}
\drawthickdot{194.0}{34.0}
\drawhollowdot{194.0}{34.0}
\drawlefttext{196.0}{40.0}{$0$}
\drawrighttext{176.0}{28.0}{$\gamma$}
\drawcenteredtext{84.0}{86.0}{abnormal isotropic}
\drawcenteredtext{156.0}{76.0}{normal Lagrangian}
\drawlefttext{182.0}{4.0}{$K(r)$}
\drawcenteredtext{84.0}{80.0}{submanifold}
\drawcenteredtext{156.0}{70.0}{manifold}
\end{picture}

\end{center}
\caption{}\label{mafig}
\end{figure}


\subsubsection{Lagrangian splitting and the Martinet sector}
\begin{defi}
We call \it{Martinet sector} of the Martinet sphere the trace of
the ball $B(0,r)$ with the Martinet plane identified to $y=0$.
\end{defi}

A precise description is obtained if we use the pendulum
representation of Section 4, where the metric is truncated to
order $0$~: $g=(1+\alpha y)^2dx^2+(1+\beta x+\gamma y)^2dy^2$.
The abnormal geodesic is strict if and only if $\alpha\neq 0$.
The pendulum equation is~:
$$\theta''+\sin\theta+\varepsilon\beta\cos\theta\
\theta'+\varepsilon^2\alpha\sin\theta
(\alpha\cos\theta-\beta\sin\theta) =0$$
where $\varepsilon=\inv{\rala}$ is a \it{parameter}. Cutting
by $y=0$ induces a \it{one-parameter} section~:
$$S~:\ \theta'=\varepsilon(\alpha\cos\theta+\beta\sin\theta)$$
The trace of the sphere with the Martinet plane near the
end-point $A=(-r,0,0)$ of the abnormal direction is described in
Section 4. We take the first and second intersections of the
pendulum trajectories with $S$, see Fig. \ref{lafig2}.

\setlength{\unitlength}{0.5mm}
\begin{figure}[h]
\begin{center}

\begin{picture}(220,120)
\thinlines
\drawvector{-20.0}{60.0}{142.0}{1}{0}
\drawvector{50.0}{16.0}{90.0}{0}{1}
\drawdashline{0.0}{22.0}{0.0}{102.0}
\drawdashline{100.0}{22.0}{100.0}{102.0}
\drawlefttext{50.9}{56.69}{$0$}
\drawcenteredtext{124.47}{63.44}{$\theta$}
\drawlefttext{51.55}{103.22}{$\theta'$}
\drawlefttext{103.08}{56.69}{$+\pi$}
\drawdotline{0.0}{60.0}{20.0}{30.0}
\thicklines
\drawarc{50.0}{100.0}{128.05}{0.67}{2.46}
\drawarc{50.0}{20.0}{128.05}{3.75}{5.59}
\thinlines
\drawdotline{-20.0}{30.0}{0.0}{60.0}
\drawdotline{0.0}{60.0}{20.0}{90.0}
\drawdotline{-20.0}{90.0}{0.0}{60.0}
\drawdotline{100.0}{60.0}{120.0}{90.0}
\drawdotline{80.0}{90.0}{100.0}{60.0}
\drawdotline{100.0}{60.0}{80.0}{30.0}
\drawdotline{120.0}{30.0}{100.0}{60.0}
\drawrighttext{-3.13}{56.48}{$-\pi$}
\thicklines
\drawarc{139.77}{87.16}{96.33}{1.9}{2.5}
\drawarc{138.0}{34.29}{91.69}{3.73}{4.4}
\thinlines
\path(98.09,51.0)(98.08,51.0)(98.01,51.0)(97.94,50.99)(97.86,50.99)(97.79,50.97)(97.7,50.97)(97.62,50.95)(97.55,50.95)(97.48,50.93)
\drawarc{50.0}{22.0}{87.72}{3.95}{5.46}
\drawarc{104.0}{74.0}{62.47}{0.86}{2.44}
\drawarc{-5.63}{76.76}{68.58}{0.71}{2.13}
\path(97.48,50.93)(97.4,50.93)(97.31,50.9)(97.23,50.9)(97.16,50.88)(97.08,50.86)(97.01,50.84)(96.93,50.84)(96.84,50.81)(96.76,50.79)
\path(96.76,50.79)(96.69,50.77)(96.61,50.75)(96.52,50.74)(96.44,50.72)(96.37,50.68)(96.29,50.66)(96.2,50.65)(96.12,50.61)(96.05,50.59)
\path(96.05,50.59)(95.97,50.56)(95.88,50.54)(95.8,50.5)(95.73,50.47)(95.65,50.45)(95.56,50.41)(95.48,50.38)(95.41,50.34)(95.33,50.31)
\path(95.33,50.31)(95.25,50.27)(95.16,50.25)(95.08,50.2)(95.0,50.16)(94.91,50.13)(94.83,50.09)(94.75,50.06)(94.66,50.0)(94.58,49.97)
\path(94.58,49.97)(94.51,49.93)(94.41,49.88)(94.33,49.84)(94.26,49.79)(94.18,49.75)(94.08,49.7)(94.01,49.65)(93.93,49.61)(93.83,49.56)
\path(93.83,49.56)(93.76,49.5)(93.66,49.45)(93.58,49.4)(93.51,49.34)(93.41,49.29)(93.33,49.25)(93.25,49.18)(93.16,49.13)(93.08,49.08)
\path(93.08,49.08)(93.0,49.02)(92.91,48.95)(92.83,48.9)(92.73,48.84)(92.66,48.77)(92.56,48.72)(92.48,48.65)(92.4,48.59)(92.3,48.52)
\path(92.3,48.52)(92.23,48.45)(92.13,48.4)(92.05,48.33)(91.97,48.25)(91.87,48.18)(91.8,48.11)(91.7,48.04)(91.62,47.97)(91.54,47.9)
\path(91.54,47.9)(91.44,47.83)(91.36,47.75)(91.26,47.68)(91.19,47.61)(91.09,47.52)(91.01,47.45)(90.91,47.38)(90.83,47.29)(90.73,47.22)
\path(90.73,47.22)(90.66,47.13)(90.56,47.06)(90.48,46.97)(90.38,46.88)(90.3,46.81)(90.2,46.72)(90.12,46.63)(90.02,46.54)(89.94,46.45)
\path(89.94,46.45)(89.84,46.36)(89.76,46.29)(89.76,46.29)
\path(89.76,46.29)(89.76,46.29)(89.43,45.99)(89.09,45.7)(88.76,45.41)(88.41,45.13)(88.08,44.86)(87.73,44.58)(87.4,44.29)(87.05,44.02)
\path(87.05,44.02)(86.72,43.75)(86.37,43.49)(86.01,43.22)(85.68,42.95)(85.33,42.7)(84.98,42.45)(84.62,42.18)(84.26,41.93)(83.91,41.68)
\path(83.91,41.68)(83.55,41.43)(83.19,41.2)(82.83,40.95)(82.48,40.7)(82.12,40.47)(81.76,40.24)(81.4,40.0)(81.02,39.77)(80.66,39.54)
\path(80.66,39.54)(80.3,39.33)(79.93,39.11)(79.55,38.88)(79.19,38.66)(78.81,38.45)(78.44,38.25)(78.06,38.04)(77.69,37.83)(77.3,37.63)
\path(98.01,50.93)(98.19,50.98)(98.36,51.04)(98.51,51.08)(98.69,51.12)(98.86,51.18)(99.05,51.19)(99.22,51.23)(99.36,51.25)(99.55,51.26)
\path(99.55,51.26)(99.73,51.29)(99.91,51.29)(100.08,51.29)(100.25,51.29)(100.44,51.26)(100.61,51.25)(100.76,51.23)(100.94,51.19)(101.11,51.15)
\path(101.11,51.15)(101.3,51.12)(101.48,51.08)(101.66,51.04)(101.83,50.97)(102.01,50.9)(102.19,50.83)(102.36,50.76)(102.55,50.69)(102.73,50.61)
\path(102.73,50.61)(102.91,50.51)(103.11,50.43)(103.26,50.3)(103.47,50.22)(103.61,50.08)(103.83,49.98)(104.0,49.86)(104.19,49.72)(104.38,49.58)
\path(104.38,49.58)(104.55,49.44)(104.75,49.3)(104.91,49.15)(105.11,49.0)(105.3,48.83)(105.47,48.65)(105.66,48.48)(105.86,48.3)(106.02,48.12)
\path(106.02,48.12)(106.22,47.93)(106.41,47.73)(106.6,47.51)(106.77,47.3)(106.97,47.11)(107.16,46.87)(107.36,46.65)(107.55,46.43)(107.74,46.19)
\path(107.74,46.19)(107.91,45.94)(108.11,45.69)(108.3,45.44)(108.5,45.19)(108.69,44.93)(108.88,44.65)(109.07,44.37)(109.27,44.08)(109.46,43.79)
\path(109.46,43.79)(109.63,43.51)(109.85,43.19)(109.85,43.19)
\path(77.3,37.63)(76.93,37.43)(76.55,37.22)(76.16,37.02)(75.79,36.84)(75.4,36.65)(75.01,36.45)(74.62,36.27)(74.23,36.09)(73.84,35.9)
\path(15.79,49.97)(15.79,49.97)(15.15,50.36)(14.55,50.76)(14.0,51.18)(13.44,51.58)(12.94,52.0)(12.47,52.4)(12.04,52.8)(11.63,53.22)
\path(11.63,53.22)(11.24,53.62)(10.88,54.04)(10.55,54.44)(10.27,54.86)(10.0,55.26)(9.75,55.68)(9.54,56.08)(9.36,56.5)(9.19,56.9)
\path(9.19,56.9)(9.05,57.3)(8.94,57.72)(8.86,58.12)(8.8,58.54)(8.77,58.93)(8.77,59.33)(8.79,59.73)(8.83,60.12)(8.88,60.54)
\path(8.88,60.54)(8.97,60.93)(9.08,61.33)(9.19,61.72)(9.36,62.11)(9.52,62.5)(9.72,62.87)(9.93,63.25)(10.15,63.65)(10.4,64.01)
\path(10.4,64.01)(10.68,64.4)(10.97,64.76)(11.27,65.12)(11.61,65.5)(11.94,65.86)(12.3,66.22)(12.68,66.58)(13.08,66.94)(13.49,67.29)
\path(13.49,67.29)(13.9,67.62)(14.36,67.98)(14.8,68.33)(15.29,68.65)(15.77,68.98)(16.27,69.3)(16.77,69.62)(17.29,69.94)(17.83,70.26)
\path(17.83,70.26)(18.38,70.58)(18.94,70.87)(19.52,71.18)(20.08,71.48)(20.69,71.76)(21.29,72.05)(21.9,72.33)(22.52,72.61)(23.15,72.87)
\path(23.15,72.87)(23.8,73.15)(24.43,73.4)(25.11,73.65)(25.75,73.9)(26.43,74.15)(27.11,74.37)(27.79,74.61)(28.5,74.83)(29.18,75.04)
\path(29.18,75.04)(29.86,75.25)(30.58,75.44)(31.29,75.65)(32.0,75.83)(32.72,76.01)(33.43,76.19)(34.18,76.37)(34.9,76.51)(35.61,76.69)
\path(35.61,76.69)(36.36,76.83)(37.08,76.97)(37.79,77.08)(38.54,77.22)(39.25,77.33)(40.01,77.44)(40.75,77.55)(41.48,77.62)(42.19,77.72)
\path(42.19,77.72)(42.93,77.8)(43.65,77.86)(44.37,77.91)(45.08,77.97)(45.8,78.01)(46.51,78.04)(47.22,78.05)(47.91,78.05)(48.62,78.08)
\path(48.62,78.08)(49.3,78.05)(49.98,78.05)(50.0,78.05)
\drawvector{47.91}{78.05}{2.0}{1}{0}
\path(84.41,49.75)(84.41,49.75)(85.08,50.15)(85.73,50.58)(86.33,51.01)(86.91,51.43)(87.44,51.86)(87.97,52.29)(88.44,52.72)(88.87,53.12)
\path(88.87,53.12)(89.3,53.55)(89.69,53.97)(90.01,54.4)(90.36,54.83)(90.66,55.25)(90.91,55.65)(91.12,56.08)(91.36,56.5)(91.51,56.9)
\path(91.51,56.9)(91.66,57.33)(91.8,57.73)(91.87,58.15)(91.94,58.55)(92.0,58.94)(92.01,59.36)(92.0,59.76)(91.97,60.15)(91.91,60.55)
\path(91.91,60.55)(91.83,60.97)(91.72,61.36)(91.61,61.75)(91.44,62.15)(91.26,62.54)(91.08,62.9)(90.86,63.29)(90.62,63.68)(90.37,64.05)
\path(90.37,64.05)(90.08,64.43)(89.8,64.79)(89.47,65.15)(89.12,65.51)(88.76,65.87)(88.41,66.25)(88.01,66.58)(87.62,66.94)(87.19,67.29)
\path(87.19,67.29)(86.76,67.62)(86.3,67.97)(85.83,68.3)(85.33,68.65)(84.83,68.97)(84.33,69.29)(83.8,69.61)(83.26,69.93)(82.69,70.23)
\path(82.69,70.23)(82.12,70.54)(81.55,70.83)(80.98,71.12)(80.37,71.43)(79.76,71.69)(79.12,71.98)(78.51,72.26)(77.86,72.54)(77.19,72.79)
\path(77.19,72.79)(76.55,73.05)(75.87,73.3)(75.19,73.55)(74.51,73.8)(73.83,74.04)(73.12,74.26)(72.44,74.5)(71.73,74.72)(71.01,74.93)
\path(71.01,74.93)(70.3,75.12)(69.58,75.33)(68.83,75.51)(68.12,75.72)(67.37,75.9)(66.66,76.08)(65.91,76.23)(65.18,76.4)(64.44,76.55)
\path(64.44,76.55)(63.69,76.69)(62.94,76.83)(62.19,76.98)(61.47,77.08)(60.73,77.22)(59.98,77.33)(59.25,77.44)(58.51,77.51)(57.76,77.62)
\path(57.76,77.62)(57.05,77.69)(56.3,77.76)(55.58,77.83)(54.87,77.87)(54.16,77.94)(53.44,77.98)(52.75,78.01)(52.05,78.04)(51.36,78.05)
\path(51.36,78.05)(50.66,78.05)(50.0,78.05)(50.0,78.05)
\drawvector{100.61}{51.3}{1.11}{-1}{0}
\thicklines
\drawvector{48.0}{84.0}{2.0}{1}{0}
\drawvector{52.0}{36.0}{2.0}{-1}{0}
\thinlines
\drawcenteredtext{134.0}{52.0}{section $S$}
\path(73.84,35.9)(73.45,35.72)(73.06,35.54)(72.66,35.38)(72.27,35.2)(71.87,35.04)(71.48,34.88)(71.08,34.7)(70.69,34.54)(70.27,34.4)
\path(70.27,34.4)(69.87,34.24)(69.48,34.09)(69.06,33.93)(68.66,33.79)(68.26,33.65)(67.84,33.5)(67.44,33.36)(67.02,33.22)(66.61,33.09)
\drawvector{170.0}{60.0}{76.0}{1}{0}
\drawvector{226.0}{28.0}{66.0}{0}{1}
\drawthickdot{194.0}{60.0}
\drawhollowdot{194.0}{60.0}
\path(194.0,60.0)(194.0,60.0)(194.47,60.0)(194.94,60.0)(195.41,60.0)(195.88,60.01)(196.35,60.04)(196.8,60.05)(197.27,60.08)(197.72,60.11)
\path(197.72,60.11)(198.16,60.16)(198.61,60.19)(199.05,60.23)(199.5,60.26)(199.91,60.33)(200.36,60.37)(200.77,60.44)(201.21,60.5)(201.63,60.55)
\path(201.63,60.55)(202.05,60.62)(202.47,60.72)(202.88,60.79)(203.27,60.87)(203.66,60.94)(204.08,61.05)(204.47,61.15)(204.86,61.25)(205.25,61.33)
\path(205.25,61.33)(205.63,61.44)(206.02,61.55)(206.38,61.68)(206.77,61.79)(207.13,61.91)(207.5,62.04)(207.86,62.16)(208.22,62.3)(208.58,62.44)
\path(208.58,62.44)(208.94,62.58)(209.27,62.72)(209.63,62.87)(209.97,63.04)(210.3,63.19)(210.63,63.36)(210.97,63.51)(211.3,63.69)(211.63,63.86)
\path(211.63,63.86)(211.94,64.04)(212.27,64.22)(212.58,64.41)(212.88,64.58)(213.19,64.8)(213.5,64.98)(213.77,65.19)(214.08,65.4)(214.38,65.61)
\path(214.38,65.61)(214.66,65.83)(214.94,66.04)(215.22,66.26)(215.5,66.48)(215.77,66.72)(216.02,66.94)(216.32,67.19)(216.57,67.44)(216.83,67.68)
\path(216.83,67.68)(217.08,67.93)(217.33,68.19)(217.58,68.44)(217.82,68.69)(218.07,68.97)(218.3,69.23)(218.52,69.51)(218.75,69.79)(219.0,70.08)
\path(219.0,70.08)(219.22,70.36)(219.44,70.65)(219.66,70.94)(219.86,71.23)(220.07,71.55)(220.27,71.83)(220.47,72.16)(220.66,72.47)(220.86,72.79)
\path(220.86,72.79)(221.07,73.11)(221.25,73.44)(221.41,73.76)(221.6,74.11)(221.77,74.44)(221.94,74.79)(222.13,75.12)(222.27,75.47)(222.44,75.83)
\path(222.44,75.83)(222.6,76.19)(222.75,76.55)(222.91,76.91)(223.07,77.29)(223.19,77.66)(223.35,78.04)(223.47,78.43)(223.6,78.8)(223.75,79.19)
\path(223.75,79.19)(223.85,79.58)(224.0,79.98)(224.0,80.0)
\path(194.0,60.0)(194.0,60.0)(193.66,59.98)(193.36,59.98)(193.05,59.97)(192.74,59.97)(192.41,59.94)(192.13,59.93)(191.82,59.91)(191.52,59.87)
\path(191.52,59.87)(191.22,59.83)(190.94,59.8)(190.63,59.76)(190.36,59.73)(190.07,59.69)(189.77,59.62)(189.5,59.58)(189.22,59.51)(188.96,59.47)
\path(188.96,59.47)(188.69,59.41)(188.41,59.33)(188.16,59.26)(187.88,59.19)(187.63,59.11)(187.38,59.04)(187.11,58.94)(186.86,58.86)(186.61,58.76)
\path(186.61,58.76)(186.38,58.68)(186.13,58.58)(185.88,58.47)(185.66,58.37)(185.41,58.26)(185.19,58.15)(184.96,58.01)(184.72,57.91)(184.5,57.79)
\path(184.5,57.79)(184.27,57.66)(184.07,57.51)(183.86,57.4)(183.63,57.25)(183.44,57.11)(183.22,56.97)(183.02,56.8)(182.82,56.66)(182.63,56.5)
\path(182.63,56.5)(182.41,56.33)(182.24,56.19)(182.05,56.01)(181.86,55.83)(181.66,55.66)(181.5,55.5)(181.32,55.3)(181.13,55.12)(180.97,54.94)
\path(180.97,54.94)(180.8,54.75)(180.63,54.55)(180.47,54.33)(180.3,54.15)(180.13,53.94)(179.99,53.72)(179.83,53.51)(179.66,53.3)(179.52,53.08)
\path(179.52,53.08)(179.38,52.83)(179.25,52.61)(179.11,52.37)(178.97,52.15)(178.83,51.91)(178.71,51.66)(178.58,51.43)(178.46,51.18)(178.33,50.91)
\path(178.33,50.91)(178.21,50.66)(178.1,50.4)(177.97,50.12)(177.86,49.86)(177.75,49.58)(177.66,49.3)(177.55,49.04)(177.44,48.75)(177.36,48.47)
\path(177.36,48.47)(177.25,48.19)(177.16,47.87)(177.08,47.58)(176.99,47.29)(176.91,46.98)(176.83,46.68)(176.75,46.36)(176.66,46.05)(176.6,45.73)
\path(176.6,45.73)(176.52,45.41)(176.47,45.08)(176.38,44.75)(176.33,44.43)(176.27,44.08)(176.22,43.75)(176.16,43.41)(176.13,43.05)(176.08,42.69)
\path(176.08,42.69)(176.02,42.33)(176.0,42.0)(176.0,42.0)
\drawlefttext{226.0}{76.0}{$C_1$}
\drawlefttext{180.0}{44.0}{$D_1$}
\drawcenteredtext{246.0}{66.0}{$x$}
\drawlefttext{230.0}{92.0}{$z$}
\path(142.0,84.0)(142.0,84.0)(142.38,84.35)(142.8,84.69)(143.19,85.02)(143.6,85.38)(144.0,85.69)(144.38,86.02)(144.8,86.33)(145.21,86.63)
\path(145.21,86.63)(145.61,86.94)(146.02,87.22)(146.41,87.5)(146.82,87.8)(147.22,88.07)(147.63,88.33)(148.02,88.58)(148.44,88.83)(148.85,89.07)
\path(148.85,89.07)(149.25,89.3)(149.66,89.52)(150.08,89.75)(150.47,89.97)(150.88,90.16)(151.3,90.36)(151.71,90.55)(152.11,90.75)(152.52,90.91)
\path(152.52,90.91)(152.94,91.08)(153.35,91.25)(153.75,91.41)(154.16,91.55)(154.58,91.69)(155.0,91.83)(155.41,91.94)(155.83,92.07)(156.24,92.16)
\path(156.24,92.16)(156.63,92.27)(157.07,92.38)(157.47,92.47)(157.88,92.55)(158.3,92.63)(158.72,92.69)(159.13,92.75)(159.55,92.82)(159.97,92.86)
\path(159.97,92.86)(160.38,92.88)(160.82,92.94)(161.24,92.94)(161.66,92.97)(162.08,92.97)(162.5,93.0)(162.91,92.97)(163.33,92.97)(163.75,92.94)
\path(163.75,92.94)(164.16,92.94)(164.6,92.91)(165.02,92.86)(165.44,92.82)(165.86,92.75)(166.27,92.69)(166.71,92.63)(167.13,92.55)(167.55,92.47)
\path(167.55,92.47)(167.99,92.38)(168.41,92.27)(168.83,92.19)(169.27,92.07)(169.69,91.94)(170.11,91.83)(170.55,91.69)(170.97,91.55)(171.38,91.41)
\path(171.38,91.41)(171.83,91.25)(172.25,91.08)(172.69,90.91)(173.11,90.75)(173.55,90.55)(173.97,90.36)(174.41,90.16)(174.83,89.97)(175.27,89.75)
\path(175.27,89.75)(175.71,89.52)(176.13,89.3)(176.57,89.07)(177.0,88.83)(177.44,88.58)(177.86,88.33)(178.3,88.07)(178.74,87.8)(179.16,87.5)
\path(179.16,87.5)(179.61,87.22)(180.05,86.94)(180.49,86.63)(180.91,86.33)(181.36,86.02)(181.8,85.69)(182.24,85.38)(182.66,85.02)(183.11,84.69)
\path(183.11,84.69)(183.55,84.35)(183.99,84.0)(184.0,84.0)
\drawvector{182.66}{85.02}{1.32}{1}{-1}
\drawrighttext{196.0}{66.0}{$A$}
\path(66.61,33.09)(66.19,32.95)(65.77,32.83)(65.37,32.7)(64.94,32.58)(64.52,32.45)(64.11,32.34)(63.68,32.22)(63.25,32.11)(62.84,32.0)
\path(62.84,32.0)(62.4,31.88)(61.97,31.78)(61.54,31.68)(61.11,31.57)(60.68,31.47)(60.25,31.37)(59.81,31.29)(59.38,31.2)(58.95,31.11)
\path(58.95,31.11)(58.5,31.02)(58.06,30.94)(57.63,30.86)(57.18,30.78)(56.75,30.7)(56.29,30.62)(55.86,30.55)(55.4,30.48)(54.95,30.43)
\path(54.95,30.43)(54.5,30.36)(54.06,30.29)(53.61,30.23)(53.15,30.19)(52.7,30.13)(52.25,30.09)(51.79,30.04)(51.33,29.98)(50.86,29.95)
\path(50.86,29.95)(50.4,29.9)(49.95,29.87)(49.95,29.87)
\drawvector{51.33}{29.98}{1.38}{-1}{0}
\path(49.95,29.87)(49.95,29.87)(49.38,29.88)(48.83,29.89)(48.27,29.9)(47.72,29.93)(47.18,29.95)(46.65,29.97)(46.09,30.01)(45.56,30.04)
\path(45.56,30.04)(45.02,30.06)(44.5,30.11)(43.97,30.14)(43.43,30.19)(42.91,30.22)(42.4,30.28)(41.88,30.32)(41.36,30.37)(40.84,30.44)
\path(40.84,30.44)(40.34,30.5)(39.84,30.55)(39.33,30.62)(38.83,30.69)(38.33,30.77)(37.83,30.84)(37.34,30.92)(36.84,31.0)(36.36,31.07)
\path(36.36,31.07)(35.88,31.15)(35.38,31.25)(34.9,31.34)(34.43,31.43)(33.95,31.53)(33.49,31.62)(33.02,31.72)(32.54,31.84)(32.09,31.95)
\path(32.09,31.95)(31.62,32.06)(31.17,32.16)(30.71,32.29)(30.27,32.4)(29.81,32.52)(29.37,32.65)(28.92,32.77)(28.47,32.91)(28.04,33.04)
\path(28.04,33.04)(27.61,33.18)(27.18,33.33)(26.73,33.47)(26.31,33.61)(25.88,33.75)(25.46,33.9)(25.04,34.06)(24.62,34.22)(24.21,34.38)
\path(24.21,34.38)(23.8,34.54)(23.39,34.7)(22.98,34.88)(22.59,35.06)(22.19,35.22)(21.79,35.4)(21.39,35.59)(21.01,35.77)(20.62,35.95)
\path(20.62,35.95)(20.22,36.13)(19.85,36.34)(19.46,36.52)(19.09,36.72)(18.7,36.93)(18.34,37.13)(17.96,37.34)(17.6,37.54)(17.23,37.75)
\path(17.23,37.75)(16.87,37.97)(16.52,38.18)(16.15,38.4)(15.81,38.63)(15.46,38.86)(15.1,39.09)(14.76,39.31)(14.42,39.56)(14.07,39.79)
\path(14.07,39.79)(13.75,40.04)(13.4,40.27)(13.07,40.52)(12.75,40.77)(12.43,41.02)(12.1,41.27)(11.77,41.54)(11.46,41.79)(11.14,42.06)
\path(11.14,42.06)(10.84,42.33)(10.52,42.59)(10.22,42.86)(9.92,43.15)(9.6,43.43)(9.31,43.7)(9.02,43.99)(8.72,44.27)(8.43,44.56)
\path(8.43,44.56)(8.14,44.86)(7.86,45.15)(7.86,45.15)
\drawcenteredtext{164.0}{98.0}{exp}
\path(194.0,60.0)(194.0,60.0)(194.35,60.0)(194.71,60.0)(195.05,60.0)(195.41,60.02)(195.75,60.04)(196.1,60.06)(196.44,60.08)(196.77,60.11)
\path(196.77,60.11)(197.11,60.13)(197.44,60.18)(197.75,60.2)(198.08,60.25)(198.39,60.29)(198.72,60.34)(199.02,60.4)(199.35,60.45)(199.64,60.52)
\path(199.64,60.52)(199.96,60.58)(200.25,60.63)(200.55,60.72)(200.85,60.79)(201.13,60.86)(201.42,60.95)(201.71,61.02)(202.0,61.11)(202.27,61.2)
\path(202.27,61.2)(202.55,61.31)(202.82,61.4)(203.08,61.5)(203.36,61.61)(203.61,61.72)(203.88,61.84)(204.13,61.95)(204.38,62.08)(204.63,62.2)
\path(204.63,62.2)(204.88,62.33)(205.11,62.45)(205.36,62.59)(205.6,62.72)(205.83,62.86)(206.07,63.02)(206.28,63.16)(206.52,63.31)(206.74,63.47)
\path(206.74,63.47)(206.94,63.63)(207.16,63.79)(207.38,63.97)(207.58,64.13)(207.78,64.31)(208.0,64.5)(208.19,64.68)(208.38,64.86)(208.58,65.05)
\path(208.58,65.05)(208.77,65.23)(208.94,65.44)(209.13,65.63)(209.32,65.83)(209.49,66.05)(209.66,66.26)(209.83,66.47)(210.0,66.69)(210.16,66.91)
\path(210.16,66.91)(210.32,67.13)(210.47,67.37)(210.63,67.59)(210.78,67.83)(210.92,68.08)(211.08,68.31)(211.22,68.55)(211.35,68.8)(211.49,69.06)
\path(211.49,69.06)(211.61,69.33)(211.75,69.58)(211.86,69.84)(212.0,70.12)(212.11,70.38)(212.22,70.66)(212.33,70.94)(212.44,71.23)(212.55,71.51)
\path(212.55,71.51)(212.66,71.8)(212.75,72.09)(212.85,72.4)(212.94,72.69)(213.02,73.0)(213.11,73.3)(213.19,73.62)(213.27,73.93)(213.36,74.25)
\path(213.36,74.25)(213.44,74.56)(213.5,74.9)(213.57,75.23)(213.63,75.55)(213.69,75.9)(213.75,76.23)(213.8,76.58)(213.86,76.93)(213.91,77.27)
\path(213.91,77.27)(213.94,77.63)(214.0,77.98)(214.0,78.0)
\drawrighttext{212.0}{78.0}{$D_2$}
\end{picture}

\end{center}
\caption{($a\neq 0$)}\label{lafig2}
\end{figure}

In the conservative case,
the curves $D_1$, $D_2$ correspond to \it{oscillating}
trajectories of the pendulum, and the curve $C_1$ corresponds to
\it{rotating} trajectories.

Only one of the curves $C_1,D_2$ belongs to the sphere (this is
$D_2$ on the figure) and their respective positions depend on
the Gauss curvature of the restriction of the metric $g$ to the
plane $(x,y)$.

Contacts are the following.

\begin{prop}
Let $Z=\f{z}{r^3}$ and $X=\f{x+r}{2r}$. Then~:
\begin{itemize}
\item[$C_1, D_2$]~: $Z=(\inv{6}+\rm{O}(r))X^3+\rm{o}(X^3)$.
\item[$D_1$]~: $Z=-\f{2}{r^2\alpha^2}X^2+\rm{o}(X^2)$.
\end{itemize}
\end{prop}

This allows to describe the Martinet sector in the ball.

\begin{prop}
In the strict case $\alpha\neq 0$ the Martinet sector has the
following properties.
\begin{enumerate}
\item It is the image by the exponential mapping of a non compact
subset of the cylinder~:
$(\theta(0), \lambda)$, $\lambda\rightarrow \infty$.
\item It is homeomorphic to a conic sector centered on the
abnormal line.
\item It is foliated by leaves $D_1,E_1$ in the spheres
$S(0,\varepsilon)$, $\epsilon\leq r$, which glue according to
Fig. \ref{lafig3}.
\end{enumerate}
\end{prop}

\setlength{\unitlength}{0.35mm}
\begin{figure}[h]
\begin{center}

\begin{picture}(180,100)
\thinlines
\drawpath{2.25}{47.86}{177.97}{48.08}
\drawthickdot{2.0}{47.86}
\path(1.77,47.63)(1.76,47.63)(3.94,47.83)(6.07,48.02)(8.21,48.22)(10.31,48.43)(12.43,48.63)(14.56,48.84)(16.63,49.04)(18.75,49.25)
\path(18.75,49.25)(20.8,49.47)(22.88,49.68)(24.95,49.88)(27.03,50.09)(29.09,50.31)(31.12,50.52)(33.15,50.75)(35.18,50.97)(37.2,51.18)
\path(37.2,51.18)(39.22,51.4)(41.22,51.63)(43.22,51.86)(45.22,52.08)(47.2,52.31)(49.15,52.54)(51.13,52.77)(53.09,53.0)(55.04,53.22)
\path(55.04,53.22)(56.97,53.45)(58.9,53.7)(60.84,53.93)(62.75,54.15)(64.66,54.4)(66.55,54.65)(68.44,54.88)(70.33,55.13)(72.22,55.38)
\path(72.22,55.38)(74.08,55.61)(75.94,55.86)(77.8,56.11)(79.66,56.36)(81.48,56.61)(83.3,56.86)(85.15,57.11)(86.94,57.38)(88.76,57.63)
\path(88.76,57.63)(90.55,57.88)(92.36,58.15)(94.12,58.4)(95.91,58.65)(97.68,58.93)(99.44,59.2)(101.19,59.47)(102.93,59.74)(104.66,60.0)
\path(104.66,60.0)(106.4,60.27)(108.12,60.54)(109.83,60.81)(111.51,61.09)(113.23,61.36)(114.91,61.63)(116.58,61.9)(118.26,62.2)(119.93,62.47)
\path(119.93,62.47)(121.58,62.75)(123.23,63.04)(124.87,63.33)(126.51,63.61)(128.13,63.9)(129.75,64.19)(131.36,64.48)(132.96,64.76)(134.55,65.05)
\path(134.55,65.05)(136.13,65.36)(137.71,65.66)(139.27,65.94)(140.83,66.25)(142.38,66.55)(143.94,66.83)(145.47,67.15)(147.0,67.44)(148.52,67.76)
\path(148.52,67.76)(150.02,68.05)(151.55,68.37)(153.02,68.69)(154.52,69.0)(156.02,69.3)(157.49,69.62)(158.96,69.94)(160.41,70.25)(161.86,70.55)
\path(161.86,70.55)(163.3,70.87)(164.74,71.19)(166.16,71.51)(167.58,71.83)(169.0,72.16)(170.38,72.5)(171.77,72.83)(173.16,73.15)(174.55,73.48)
\path(174.55,73.48)(175.91,73.8)(177.27,74.12)(177.3,74.15)
\path(1.77,48.08)(1.76,48.08)(4.73,47.72)(7.67,47.38)(10.56,47.06)(13.43,46.72)(16.29,46.38)(19.12,46.04)(21.93,45.7)(24.7,45.38)
\path(24.7,45.38)(27.45,45.04)(30.2,44.72)(32.9,44.38)(35.59,44.06)(38.25,43.74)(40.88,43.4)(43.5,43.09)(46.08,42.75)(48.65,42.43)
\path(48.65,42.43)(51.18,42.11)(53.7,41.79)(56.18,41.47)(58.65,41.15)(61.09,40.84)(63.52,40.52)(65.91,40.22)(68.26,39.9)(70.62,39.59)
\path(70.62,39.59)(72.94,39.27)(75.25,38.97)(77.51,38.65)(79.76,38.36)(82.0,38.04)(84.19,37.75)(86.37,37.43)(88.51,37.13)(90.66,36.84)
\path(90.66,36.84)(92.76,36.54)(94.83,36.24)(96.9,35.93)(98.93,35.65)(100.94,35.34)(102.91,35.06)(104.87,34.75)(106.8,34.47)(108.73,34.18)
\path(108.73,34.18)(110.62,33.88)(112.48,33.59)(114.3,33.31)(116.12,33.02)(117.93,32.75)(119.69,32.45)(121.44,32.18)(123.16,31.87)(124.83,31.62)
\path(124.83,31.62)(126.51,31.34)(128.16,31.04)(129.77,30.78)(131.38,30.51)(132.96,30.21)(134.5,29.95)(136.02,29.69)(137.52,29.42)(139.02,29.12)
\path(139.02,29.12)(140.47,28.87)(141.88,28.61)(143.3,28.35)(144.69,28.05)(146.02,27.79)(147.38,27.54)(148.69,27.29)(149.97,27.03)(151.22,26.77)
\path(151.22,26.77)(152.47,26.52)(153.66,26.26)(154.86,26.0)(156.02,25.75)(157.16,25.5)(158.27,25.21)(159.38,24.96)(160.44,24.71)(161.49,24.46)
\path(161.49,24.46)(162.5,24.25)(163.5,24.0)(164.47,23.76)(165.41,23.51)(166.35,23.27)(167.25,23.03)(168.11,22.79)(168.97,22.54)(169.8,22.29)
\path(169.8,22.29)(170.6,22.04)(171.38,21.84)(172.13,21.6)(172.86,21.37)(173.57,21.12)(174.25,20.87)(174.91,20.67)(175.52,20.44)(176.13,20.2)
\path(176.13,20.2)(176.72,19.96)(177.27,19.77)(177.3,19.77)
\drawcenteredtext{188.36}{47.4}{abnormal line}
\path(140.22,48.08)(140.22,48.08)(140.61,48.09)(141.0,48.11)(141.41,48.15)(141.8,48.18)(142.19,48.22)(142.57,48.25)(142.96,48.31)(143.33,48.36)
\path(143.33,48.36)(143.72,48.4)(144.1,48.47)(144.47,48.54)(144.85,48.61)(145.22,48.68)(145.58,48.75)(145.94,48.84)(146.3,48.9)(146.66,49.0)
\path(146.66,49.0)(147.02,49.11)(147.38,49.2)(147.74,49.31)(148.08,49.4)(148.44,49.52)(148.77,49.65)(149.13,49.77)(149.47,49.9)(149.8,50.02)
\path(149.8,50.02)(150.13,50.15)(150.47,50.29)(150.8,50.45)(151.13,50.59)(151.47,50.75)(151.77,50.9)(152.11,51.06)(152.41,51.22)(152.75,51.4)
\path(152.75,51.4)(153.05,51.56)(153.36,51.75)(153.66,51.93)(153.99,52.11)(154.27,52.31)(154.6,52.5)(154.88,52.7)(155.19,52.9)(155.49,53.13)
\path(155.49,53.13)(155.77,53.34)(156.07,53.56)(156.36,53.77)(156.63,54.0)(156.91,54.24)(157.19,54.47)(157.47,54.72)(157.75,54.95)(158.02,55.2)
\path(158.02,55.2)(158.3,55.47)(158.57,55.72)(158.83,55.99)(159.1,56.25)(159.36,56.52)(159.61,56.79)(159.86,57.08)(160.13,57.36)(160.38,57.65)
\path(160.38,57.65)(160.63,57.95)(160.86,58.25)(161.11,58.54)(161.36,58.84)(161.6,59.15)(161.83,59.47)(162.07,59.79)(162.3,60.13)(162.52,60.45)
\path(162.52,60.45)(162.75,60.79)(162.97,61.11)(163.19,61.47)(163.41,61.81)(163.63,62.15)(163.85,62.52)(164.05,62.88)(164.27,63.24)(164.47,63.61)
\path(164.47,63.61)(164.69,63.97)(164.88,64.36)(165.08,64.73)(165.27,65.12)(165.47,65.51)(165.66,65.91)(165.86,66.3)(166.05,66.72)(166.24,67.12)
\path(166.24,67.12)(166.41,67.54)(166.6,67.94)(166.77,68.37)(166.94,68.8)(167.13,69.23)(167.3,69.66)(167.47,70.11)(167.63,70.55)(167.8,71.0)
\path(167.8,71.0)(167.96,71.44)(168.11,71.9)(168.11,71.91)
\path(139.77,48.31)(139.77,48.31)(139.52,48.29)(139.27,48.27)(139.02,48.25)(138.8,48.24)(138.57,48.22)(138.33,48.2)(138.1,48.15)(137.86,48.13)
\path(137.86,48.13)(137.63,48.09)(137.41,48.06)(137.19,48.02)(136.97,47.97)(136.75,47.93)(136.52,47.88)(136.3,47.81)(136.08,47.75)(135.86,47.7)
\path(135.86,47.7)(135.66,47.63)(135.44,47.56)(135.24,47.5)(135.02,47.43)(134.82,47.36)(134.61,47.27)(134.41,47.2)(134.21,47.11)(134.0,47.02)
\path(134.0,47.02)(133.8,46.93)(133.61,46.84)(133.41,46.74)(133.22,46.63)(133.02,46.54)(132.83,46.43)(132.63,46.31)(132.47,46.2)(132.27,46.09)
\path(132.27,46.09)(132.1,45.97)(131.91,45.84)(131.72,45.72)(131.55,45.59)(131.38,45.47)(131.19,45.33)(131.02,45.18)(130.86,45.04)(130.69,44.9)
\path(130.69,44.9)(130.52,44.75)(130.35,44.61)(130.16,44.45)(130.02,44.29)(129.86,44.13)(129.69,43.97)(129.52,43.81)(129.38,43.65)(129.22,43.47)
\path(129.22,43.47)(129.07,43.29)(128.91,43.11)(128.77,42.93)(128.61,42.75)(128.47,42.56)(128.33,42.38)(128.16,42.18)(128.02,41.99)(127.9,41.79)
\path(127.9,41.79)(127.76,41.59)(127.62,41.38)(127.48,41.15)(127.36,40.95)(127.22,40.74)(127.08,40.52)(126.94,40.29)(126.83,40.08)(126.69,39.84)
\path(126.69,39.84)(126.58,39.61)(126.47,39.38)(126.33,39.15)(126.23,38.9)(126.11,38.65)(125.98,38.4)(125.87,38.15)(125.76,37.9)(125.66,37.65)
\path(125.66,37.65)(125.55,37.4)(125.44,37.13)(125.33,36.88)(125.23,36.61)(125.12,36.34)(125.01,36.06)(124.93,35.79)(124.83,35.5)(124.73,35.22)
\path(124.73,35.22)(124.62,34.93)(124.55,34.65)(124.44,34.36)(124.37,34.06)(124.26,33.75)(124.19,33.45)(124.11,33.15)(124.01,32.84)(123.94,32.54)
\path(123.94,32.54)(123.87,32.22)(123.8,31.88)(123.8,31.88)
\path(110.3,48.08)(110.3,48.08)(110.55,48.09)(110.8,48.09)(111.05,48.11)(111.33,48.13)(111.58,48.15)(111.83,48.18)(112.08,48.2)(112.33,48.24)
\path(112.33,48.24)(112.58,48.27)(112.83,48.31)(113.08,48.36)(113.33,48.4)(113.58,48.45)(113.8,48.5)(114.05,48.54)(114.3,48.61)(114.55,48.65)
\path(114.55,48.65)(114.79,48.72)(115.01,48.79)(115.26,48.86)(115.48,48.93)(115.73,49.0)(115.94,49.09)(116.19,49.15)(116.41,49.25)(116.65,49.34)
\path(116.65,49.34)(116.87,49.43)(117.11,49.52)(117.33,49.61)(117.55,49.72)(117.76,49.81)(118.0,49.9)(118.22,50.02)(118.44,50.13)(118.65,50.25)
\path(118.65,50.25)(118.87,50.36)(119.08,50.5)(119.3,50.61)(119.51,50.75)(119.72,50.86)(119.93,51.0)(120.12,51.13)(120.33,51.27)(120.55,51.4)
\path(120.55,51.4)(120.76,51.56)(120.94,51.72)(121.16,51.86)(121.37,52.02)(121.55,52.18)(121.76,52.34)(121.94,52.5)(122.16,52.65)(122.33,52.84)
\path(122.33,52.84)(122.55,53.0)(122.73,53.18)(122.93,53.36)(123.12,53.54)(123.3,53.72)(123.5,53.9)(123.68,54.11)(123.87,54.29)(124.05,54.5)
\path(124.05,54.5)(124.23,54.7)(124.41,54.9)(124.58,55.11)(124.76,55.31)(124.94,55.52)(125.12,55.75)(125.3,55.97)(125.47,56.18)(125.65,56.4)
\path(125.65,56.4)(125.8,56.63)(125.98,56.86)(126.16,57.09)(126.3,57.34)(126.48,57.58)(126.65,57.81)(126.8,58.06)(126.98,58.31)(127.12,58.56)
\path(127.12,58.56)(127.3,58.83)(127.44,59.08)(127.62,59.34)(127.76,59.61)(127.93,59.88)(128.08,60.15)(128.22,60.4)(128.38,60.7)(128.52,60.97)
\path(128.52,60.97)(128.66,61.25)(128.83,61.54)(128.97,61.83)(129.11,62.13)(129.25,62.4)(129.38,62.72)(129.52,63.02)(129.66,63.31)(129.82,63.63)
\path(129.82,63.63)(129.96,63.93)(130.1,64.25)(130.1,64.26)
\path(110.51,48.31)(110.51,48.31)(110.33,48.29)(110.16,48.29)(109.98,48.27)(109.8,48.27)(109.65,48.25)(109.48,48.22)(109.3,48.2)(109.12,48.18)
\path(109.12,48.18)(108.97,48.15)(108.8,48.13)(108.62,48.11)(108.48,48.08)(108.3,48.04)(108.15,48.0)(107.98,47.97)(107.83,47.93)(107.66,47.88)
\path(107.66,47.88)(107.51,47.84)(107.33,47.79)(107.19,47.74)(107.04,47.68)(106.87,47.63)(106.73,47.56)(106.58,47.5)(106.43,47.45)(106.26,47.38)
\path(106.26,47.38)(106.12,47.31)(105.97,47.24)(105.83,47.15)(105.68,47.09)(105.51,47.02)(105.37,46.93)(105.23,46.86)(105.08,46.77)(104.94,46.68)
\path(104.94,46.68)(104.8,46.59)(104.66,46.5)(104.51,46.4)(104.37,46.31)(104.25,46.22)(104.12,46.11)(103.98,46.0)(103.83,45.9)(103.69,45.79)
\path(103.69,45.79)(103.58,45.68)(103.44,45.58)(103.3,45.45)(103.18,45.34)(103.05,45.22)(102.93,45.09)(102.8,44.97)(102.66,44.84)(102.54,44.72)
\path(102.54,44.72)(102.41,44.59)(102.3,44.45)(102.16,44.31)(102.05,44.18)(101.93,44.04)(101.8,43.88)(101.69,43.75)(101.55,43.59)(101.44,43.45)
\path(101.44,43.45)(101.33,43.29)(101.22,43.13)(101.11,42.97)(100.98,42.81)(100.87,42.65)(100.76,42.49)(100.66,42.31)(100.55,42.15)(100.44,41.97)
\path(100.44,41.97)(100.33,41.79)(100.23,41.63)(100.12,41.45)(100.01,41.27)(99.91,41.08)(99.8,40.9)(99.69,40.7)(99.58,40.5)(99.5,40.31)
\path(99.5,40.31)(99.4,40.11)(99.3,39.93)(99.19,39.72)(99.11,39.52)(99.01,39.31)(98.91,39.11)(98.83,38.9)(98.73,38.68)(98.62,38.47)
\path(98.62,38.47)(98.55,38.25)(98.44,38.02)(98.37,37.81)(98.29,37.58)(98.19,37.36)(98.11,37.13)(98.01,36.9)(97.94,36.65)(97.86,36.43)
\path(97.86,36.43)(97.76,36.18)(97.69,35.95)(97.69,35.95)
\drawcenteredtext{176.66}{63.59}{$D_1$}
\drawcenteredtext{132.8}{38.4}{$E_1$}
\drawcenteredtext{2.59}{56.15}{$0$}
\end{picture}

\end{center}
\caption{}\label{lafig3}
\end{figure}

\no{\bf The Lagrangian splitting}
In the pendulum representation, the transport of
$T^*_{\gamma(0)}M\cap\{H_n=\inv{2}\}$ by the normal flow has the
following basic interpretation, see Fig. \ref{lafig4}.

\setlength{\unitlength}{0.35mm}
\begin{figure}[h]
\begin{center}

\begin{picture}(180,100)
\thinlines
\path(2.0,48.0)(2.0,48.0)(4.07,48.31)(6.13,48.63)(8.18,48.97)(10.26,49.29)(12.3,49.63)(14.35,49.97)(16.37,50.31)(18.38,50.65)
\path(18.38,50.65)(20.42,51.0)(22.43,51.34)(24.44,51.7)(26.44,52.06)(28.43,52.43)(30.39,52.79)(32.38,53.15)(34.34,53.52)(36.31,53.9)
\path(36.31,53.9)(38.27,54.27)(40.22,54.65)(42.15,55.04)(44.09,55.4)(46.0,55.81)(47.93,56.2)(49.84,56.59)(51.75,57.0)(53.63,57.4)
\path(53.63,57.4)(55.52,57.79)(57.4,58.2)(59.29,58.61)(61.15,59.02)(63.02,59.45)(64.87,59.86)(66.69,60.29)(68.55,60.72)(70.37,61.15)
\path(70.37,61.15)(72.19,61.59)(74.01,62.02)(75.83,62.47)(77.62,62.9)(79.43,63.34)(81.22,63.79)(83.0,64.26)(84.76,64.69)(86.55,65.16)
\path(86.55,65.16)(88.3,65.62)(90.05,66.08)(91.8,66.55)(93.54,67.04)(95.26,67.51)(96.98,67.98)(98.69,68.48)(100.41,68.94)(102.12,69.44)
\path(102.12,69.44)(103.8,69.94)(105.5,70.43)(107.19,70.93)(108.86,71.43)(110.51,71.94)(112.18,72.44)(113.83,72.94)(115.48,73.47)(117.12,73.98)
\path(117.12,73.98)(118.75,74.51)(120.37,75.01)(121.98,75.55)(123.58,76.08)(125.19,76.62)(126.79,77.15)(128.38,77.69)(129.94,78.23)(131.52,78.76)
\path(131.52,78.76)(133.08,79.33)(134.63,79.87)(136.19,80.44)(137.74,80.98)(139.27,81.55)(140.8,82.12)(142.33,82.69)(143.85,83.26)(145.35,83.83)
\path(145.35,83.83)(146.86,84.41)(148.35,84.98)(149.83,85.58)(151.3,86.16)(152.77,86.75)(154.25,87.33)(155.71,87.94)(157.16,88.55)(158.6,89.15)
\path(158.6,89.15)(160.02,89.75)(161.46,90.36)(162.88,90.98)(164.3,91.58)(165.71,92.19)(167.1,92.83)(168.5,93.44)(169.88,94.08)(171.25,94.72)
\path(171.25,94.72)(172.63,95.36)(173.99,95.98)(174.0,96.0)
\path(90.0,66.0)(90.0,66.0)(90.26,65.26)(90.55,64.55)(90.83,63.84)(91.11,63.13)(91.37,62.4)(91.66,61.7)(91.94,60.99)(92.23,60.29)
\path(92.23,60.29)(92.51,59.58)(92.8,58.88)(93.08,58.15)(93.36,57.47)(93.62,56.77)(93.91,56.06)(94.19,55.38)(94.47,54.68)(94.75,53.99)
\path(94.75,53.99)(95.01,53.29)(95.3,52.59)(95.58,51.9)(95.87,51.22)(96.16,50.54)(96.44,49.86)(96.72,49.18)(97.0,48.5)(97.26,47.81)
\path(97.26,47.81)(97.55,47.13)(97.83,46.45)(98.11,45.79)(98.37,45.11)(98.66,44.43)(98.94,43.77)(99.23,43.11)(99.51,42.43)(99.8,41.77)
\path(99.8,41.77)(100.08,41.11)(100.36,40.45)(100.62,39.79)(100.91,39.13)(101.19,38.47)(101.47,37.81)(101.75,37.15)(102.01,36.5)(102.3,35.86)
\path(102.3,35.86)(102.58,35.22)(102.87,34.56)(103.15,33.9)(103.43,33.27)(103.72,32.63)(104.0,32.0)(104.26,31.36)(104.55,30.7)(104.83,30.06)
\path(104.83,30.06)(105.11,29.45)(105.37,28.8)(105.66,28.18)(105.94,27.54)(106.23,26.93)(106.51,26.29)(106.79,25.68)(107.05,25.04)(107.33,24.43)
\path(107.33,24.43)(107.62,23.79)(107.91,23.19)(108.19,22.56)(108.47,21.95)(108.75,21.35)(109.01,20.71)(109.3,20.12)(109.58,19.52)(109.87,18.89)
\path(109.87,18.89)(110.15,18.29)(110.43,17.7)(110.69,17.1)(110.98,16.5)(111.26,15.89)(111.55,15.3)(111.83,14.68)(112.11,14.1)(112.37,13.52)
\path(112.37,13.52)(112.66,12.92)(112.94,12.31)(113.23,11.75)(113.51,11.14)(113.79,10.56)(114.05,9.97)(114.33,9.39)(114.62,8.81)(114.91,8.25)
\path(114.91,8.25)(115.19,7.67)(115.47,7.09)(115.75,6.53)(116.01,5.94)(116.3,5.38)(116.58,4.82)(116.87,4.25)(117.15,3.68)(117.43,3.1)
\path(117.43,3.1)(117.69,2.54)(117.98,2.0)(118.0,2.0)
\path(58.0,-10.0)(58.0,-10.0)(59.59,-9.75)(61.18,-9.51)(62.75,-9.27)(64.3,-9.02)(65.87,-8.77)(67.43,-8.52)(68.97,-8.27)(70.5,-8.02)
\path(70.5,-8.02)(72.01,-7.76)(73.54,-7.51)(75.04,-7.26)(76.51,-7.0)(78.01,-6.73)(79.48,-6.48)(80.94,-6.21)(82.41,-5.94)(83.87,-5.67)
\path(83.87,-5.67)(85.3,-5.42)(86.73,-5.15)(88.16,-4.86)(89.55,-4.59)(90.97,-4.32)(92.36,-4.05)(93.75,-3.75)(95.12,-3.49)(96.48,-3.21)
\path(96.48,-3.21)(97.83,-2.93)(99.19,-2.65)(100.51,-2.34)(101.86,-2.06)(103.16,-1.77)(104.48,-1.5)(105.79,-1.2)(107.08,-0.91)(108.36,-0.62)
\path(108.36,-0.62)(109.62,-0.31)(110.9,0.0)(112.15,0.27)(113.4,0.56)(114.62,0.87)(115.86,1.16)(117.08,1.49)(118.29,1.77)(119.48,2.08)
\path(119.48,2.08)(120.68,2.41)(121.86,2.73)(123.01,3.02)(124.19,3.34)(125.33,3.68)(126.48,3.99)(127.62,4.32)(128.75,4.63)(129.86,4.96)
\path(129.86,4.96)(130.97,5.28)(132.08,5.61)(133.16,5.94)(134.25,6.26)(135.32,6.61)(136.38,6.94)(137.41,7.26)(138.47,7.61)(139.5,7.94)
\path(139.5,7.94)(140.52,8.27)(141.55,8.63)(142.55,8.97)(143.55,9.31)(144.55,9.67)(145.52,10.01)(146.49,10.35)(147.44,10.71)(148.41,11.06)
\path(148.41,11.06)(149.35,11.42)(150.27,11.77)(151.21,12.14)(152.11,12.47)(153.02,12.85)(153.91,13.22)(154.8,13.56)(155.69,13.93)(156.55,14.31)
\path(156.55,14.31)(157.41,14.68)(158.25,15.05)(159.11,15.43)(159.94,15.8)(160.75,16.17)(161.57,16.54)(162.38,16.93)(163.16,17.29)(163.96,17.69)
\path(163.96,17.69)(164.72,18.05)(165.5,18.45)(166.25,18.85)(167.0,19.21)(167.75,19.62)(168.47,20.01)(169.19,20.39)(169.91,20.79)(170.61,21.2)
\path(170.61,21.2)(171.3,21.6)(171.99,21.97)(172.0,22.0)
\path(166.0,88.0)(166.0,88.0)(165.91,87.4)(165.83,86.8)(165.77,86.19)(165.69,85.58)(165.63,85.01)(165.58,84.41)(165.52,83.8)(165.47,83.23)
\path(165.47,83.23)(165.41,82.62)(165.38,82.05)(165.33,81.47)(165.27,80.87)(165.25,80.3)(165.22,79.69)(165.19,79.12)(165.16,78.55)(165.16,77.97)
\path(165.16,77.97)(165.13,77.37)(165.11,76.8)(165.11,76.23)(165.11,75.66)(165.11,75.08)(165.11,74.51)(165.11,73.94)(165.11,73.37)(165.13,72.8)
\path(165.13,72.8)(165.13,72.23)(165.16,71.66)(165.19,71.08)(165.22,70.54)(165.24,69.97)(165.27,69.41)(165.32,68.83)(165.36,68.29)(165.38,67.73)
\path(165.38,67.73)(165.44,67.16)(165.5,66.62)(165.55,66.05)(165.61,65.51)(165.66,64.94)(165.74,64.4)(165.8,63.84)(165.88,63.29)(165.96,62.75)
\path(165.96,62.75)(166.02,62.2)(166.11,61.65)(166.21,61.11)(166.3,60.58)(166.38,60.04)(166.5,59.5)(166.6,58.95)(166.69,58.4)(166.8,57.88)
\path(166.8,57.88)(166.91,57.34)(167.02,56.81)(167.16,56.27)(167.27,55.74)(167.41,55.2)(167.52,54.68)(167.66,54.15)(167.8,53.63)(167.94,53.09)
\path(167.94,53.09)(168.1,52.58)(168.25,52.04)(168.38,51.52)(168.55,51.0)(168.72,50.49)(168.88,49.97)(169.02,49.45)(169.21,48.93)(169.38,48.4)
\path(169.38,48.4)(169.57,47.9)(169.75,47.38)(169.91,46.88)(170.11,46.36)(170.3,45.86)(170.5,45.34)(170.71,44.84)(170.91,44.34)(171.11,43.84)
\path(171.11,43.84)(171.32,43.33)(171.52,42.83)(171.75,42.33)(171.97,41.83)(172.19,41.33)(172.41,40.83)(172.66,40.34)(172.88,39.84)(173.13,39.34)
\path(173.13,39.34)(173.36,38.86)(173.61,38.36)(173.86,37.86)(174.11,37.38)(174.38,36.9)(174.63,36.4)(174.88,35.9)(175.16,35.43)(175.44,34.95)
\path(175.44,34.95)(175.72,34.47)(175.99,34.0)(176.0,34.0)
\path(12.0,42.0)(12.0,42.0)(12.35,41.59)(12.71,41.18)(13.06,40.79)(13.39,40.38)(13.72,39.99)(14.06,39.58)(14.39,39.18)(14.72,38.77)
\path(14.72,38.77)(15.06,38.36)(15.38,37.95)(15.68,37.54)(16.0,37.13)(16.29,36.72)(16.6,36.31)(16.88,35.9)(17.19,35.49)(17.46,35.08)
\path(17.46,35.08)(17.76,34.65)(18.04,34.25)(18.3,33.84)(18.56,33.4)(18.85,33.0)(19.11,32.58)(19.37,32.15)(19.62,31.75)(19.87,31.3)
\path(19.87,31.3)(20.11,30.88)(20.35,30.46)(20.56,30.04)(20.79,29.62)(21.04,29.2)(21.26,28.79)(21.46,28.36)(21.69,27.93)(21.88,27.51)
\path(21.88,27.51)(22.1,27.06)(22.29,26.63)(22.5,26.2)(22.69,25.79)(22.87,25.36)(23.04,24.92)(23.21,24.47)(23.39,24.04)(23.56,23.62)
\path(23.56,23.62)(23.72,23.19)(23.88,22.75)(24.04,22.29)(24.2,21.87)(24.35,21.43)(24.47,21.0)(24.62,20.54)(24.77,20.11)(24.88,19.67)
\path(24.88,19.67)(25.02,19.21)(25.12,18.79)(25.26,18.34)(25.37,17.88)(25.45,17.45)(25.56,17.0)(25.67,16.54)(25.77,16.11)(25.86,15.64)
\path(25.86,15.64)(25.94,15.21)(26.02,14.76)(26.1,14.31)(26.17,13.85)(26.22,13.39)(26.29,12.93)(26.36,12.47)(26.39,12.02)(26.45,11.56)
\path(26.45,11.56)(26.51,11.1)(26.54,10.64)(26.59,10.18)(26.62,9.75)(26.63,9.27)(26.67,8.81)(26.69,8.35)(26.7,7.9)(26.7,7.44)
\path(26.7,7.44)(26.7,6.96)(26.7,6.51)(26.7,6.03)(26.7,5.57)(26.7,5.11)(26.68,4.63)(26.64,4.17)(26.62,3.7)(26.61,3.23)
\path(26.61,3.23)(26.56,2.75)(26.54,2.26)(26.47,1.79)(26.45,1.34)(26.38,0.86)(26.34,0.37)(26.28,-0.07)(26.2,-0.56)(26.13,-1.02)
\path(26.13,-1.02)(26.05,-1.5)(26.0,-1.99)(26.0,-2.0)
\path(20.0,32.0)(20.0,32.0)(21.19,32.27)(22.37,32.52)(23.54,32.79)(24.71,33.02)(25.88,33.25)(27.04,33.47)(28.2,33.68)(29.35,33.88)
\path(29.35,33.88)(30.47,34.08)(31.61,34.25)(32.74,34.4)(33.84,34.58)(34.95,34.72)(36.04,34.86)(37.13,34.97)(38.22,35.09)(39.29,35.18)
\path(39.29,35.18)(40.36,35.29)(41.4,35.36)(42.47,35.43)(43.52,35.49)(44.56,35.54)(45.58,35.58)(46.61,35.59)(47.61,35.61)(48.63,35.61)
\path(48.63,35.61)(49.61,35.61)(50.61,35.59)(51.59,35.56)(52.56,35.54)(53.54,35.49)(54.5,35.43)(55.45,35.34)(56.4,35.27)(57.34,35.18)
\path(57.34,35.18)(58.27,35.08)(59.18,34.95)(60.11,34.84)(61.02,34.7)(61.9,34.56)(62.81,34.4)(63.68,34.22)(64.55,34.04)(65.44,33.86)
\path(65.44,33.86)(66.3,33.65)(67.15,33.45)(68.0,33.22)(68.83,32.99)(69.66,32.75)(70.48,32.5)(71.3,32.22)(72.12,31.95)(72.91,31.67)
\path(72.91,31.67)(73.69,31.37)(74.5,31.04)(75.26,30.72)(76.05,30.39)(76.8,30.05)(77.55,29.7)(78.3,29.36)(79.05,28.96)(79.79,28.6)
\path(79.79,28.6)(80.51,28.2)(81.23,27.79)(81.94,27.37)(82.62,26.95)(83.33,26.51)(84.01,26.05)(84.69,25.61)(85.37,25.12)(86.04,24.63)
\path(86.04,24.63)(86.69,24.14)(87.33,23.64)(87.98,23.12)(88.62,22.62)(89.25,22.06)(89.86,21.54)(90.48,20.96)(91.08,20.39)(91.66,19.84)
\path(91.66,19.84)(92.26,19.25)(92.83,18.63)(93.41,18.03)(93.98,17.39)(94.54,16.78)(95.08,16.12)(95.62,15.47)(96.16,14.81)(96.69,14.14)
\path(96.69,14.14)(97.19,13.46)(97.73,12.76)(98.23,12.05)(98.73,11.31)(99.22,10.6)(99.69,9.85)(100.16,9.1)(100.62,8.35)(101.08,7.57)
\path(101.08,7.57)(101.55,6.78)(101.98,6.0)(102.0,6.0)
\path(123.97,11.06)(123.97,11.06)(123.91,11.92)(123.87,12.77)(123.86,13.63)(123.83,14.47)(123.83,15.31)(123.83,16.14)(123.83,17.0)(123.86,17.8)
\path(123.86,17.8)(123.87,18.63)(123.93,19.45)(123.98,20.28)(124.04,21.06)(124.08,21.87)(124.18,22.68)(124.26,23.45)(124.36,24.25)(124.44,25.03)
\path(124.44,25.03)(124.55,25.79)(124.69,26.56)(124.8,27.34)(124.94,28.1)(125.11,28.85)(125.26,29.6)(125.44,30.35)(125.61,31.06)(125.8,31.8)
\path(125.8,31.8)(125.98,32.54)(126.19,33.25)(126.41,33.97)(126.62,34.68)(126.87,35.4)(127.11,36.09)(127.36,36.79)(127.62,37.47)(127.87,38.15)
\path(127.87,38.15)(128.16,38.84)(128.44,39.52)(128.75,40.18)(129.05,40.84)(129.36,41.5)(129.69,42.15)(130.02,42.81)(130.36,43.45)(130.72,44.09)
\path(130.72,44.09)(131.08,44.72)(131.44,45.34)(131.83,45.97)(132.22,46.58)(132.61,47.18)(133.02,47.79)(133.44,48.4)(133.86,48.99)(134.3,49.58)
\path(134.3,49.58)(134.75,50.15)(135.19,50.74)(135.66,51.31)(136.13,51.88)(136.61,52.45)(137.1,53.0)(137.6,53.56)(138.11,54.09)(138.63,54.63)
\path(138.63,54.63)(139.16,55.18)(139.69,55.7)(140.24,56.24)(140.77,56.75)(141.36,57.27)(141.91,57.77)(142.5,58.27)(143.11,58.77)(143.71,59.27)
\path(143.71,59.27)(144.32,59.75)(144.94,60.25)(145.55,60.72)(146.19,61.2)(146.85,61.65)(147.5,62.13)(148.16,62.58)(148.85,63.02)(149.52,63.47)
\path(149.52,63.47)(150.22,63.9)(150.91,64.36)(151.63,64.79)(152.36,65.19)(153.08,65.62)(153.83,66.04)(154.58,66.44)(155.33,66.86)(156.1,67.26)
\path(156.1,67.26)(156.88,67.65)(157.66,68.04)(158.46,68.41)(159.25,68.79)(160.07,69.16)(160.88,69.51)(161.72,69.87)(162.57,70.25)(163.41,70.58)
\path(163.41,70.58)(164.27,70.94)(165.13,71.26)(165.13,71.29)
\drawpath{48.36}{57.58}{46.56}{56.0}
\drawpath{46.56}{56.0}{48.81}{55.33}
\drawpath{116.3}{74.87}{119.01}{74.43}
\drawpath{119.01}{74.43}{116.76}{72.41}
\drawpath{131.38}{47.02}{133.41}{48.36}
\drawpath{133.41}{48.36}{132.97}{46.11}
\drawpath{100.12}{38.7}{99.9}{41.61}
\drawpath{99.9}{41.61}{101.91}{39.38}
\drawpath{79.19}{30.14}{76.26}{30.14}
\drawpath{76.26}{30.14}{78.3}{27.69}
\drawcenteredtext{179.1}{72.41}{$S_1$}
\drawcenteredtext{7.19}{23.63}{$S_2$}
\drawcenteredtext{173.25}{8.81}{$S$}
\drawvector{108.0}{100.0}{12.0}{-1}{-4}
\drawcenteredtext{108.0}{104.0}{critical geodesic}
\end{picture}

\end{center}
\caption{}\label{lafig4}
\end{figure}

The section splits into two parts $S_1,S_2$ which represent the
splitting of the fiber $T^*_{\gamma(0)}M$ into two Lagrangian
manifolds.


\subsubsection{Microlocal invariants}
The pendulum has two singular points $F=(0,0)$ and $S=(0,\pi)$.
The local analysis is as follows.
\begin{itemize}
\item Near $F$, the linearized system is a \it{focus} whose
eigenvalues are~:
$$\sigma_\pm = -\f{\varepsilon\beta}{2}\pm
i\sqrt{1+\varepsilon^2(\f{\beta^2}{4} -\alpha^2) } $$
and is a perturbation of the linearized pendulum
$\theta''+\theta=0$ of the flat case.
\item Near $S$, the linearized system is a \it{saddle} whose
eigenvalues are~:
$$\eta_\pm=\f{\varepsilon\beta}{2}\pm
\sqrt{1+\varepsilon^2(\f{\beta^2}{4} -\alpha^2) } 
= \pm 1+\f{\varepsilon\beta}{2}+\rm{o}(\varepsilon)$$
where $\varepsilon=\inv{\rala}$. It is a perturbation of the flat
case $\eta_\pm=\pm 1$ which is \it{resonant}.
\end{itemize}

In order to compute the sector we use the \it{spectrum band}
$\eta_\pm$ which is stable by perturbation. When we compute the
sphere we have to compute an \it{averaging}. This is much more
complex.

The linear pendulum appears already in the contact case and the
existence of the focus reflects the existence of a contact sector
in the Martinet sphere.


\subsubsection{The sectors of the Martinet sphere}
In \cite{BC} was described the Martinet sphere in the integrable
case by gluing sectors. We have three kinds of sectors~:
\begin{itemize}
\item A Riemannian sector $R$ located near the equator, image of
$\lambda =0$.
\item A contact sector around $C$, where $C$ is a cut-point.
\item A Martinet sector around $A$, where $A$ is an end-point of
the abnormal line.
\end{itemize}

The sectors are represented on Fig. \ref{lafig5}.

\setlength{\unitlength}{0.35mm}
\begin{figure}[h]
\begin{center}

\begin{picture}(180,100)
\thinlines
\drawcircle{174.0}{48.0}{84.0}{}
\drawcircle{6.0}{48.0}{84.0}{}
\drawcenteredtext{38.0}{88.0}{$\lambda =0$}
\drawpath{-16.0}{48.0}{30.0}{48.0}
\drawthickdot{-16.0}{48.0}
\drawthickdot{30.0}{48.0}
\drawdashline{-16.0}{48.0}{8.0}{62.0}
\drawdashline{8.0}{62.0}{30.0}{48.0}
\drawdashline{30.0}{48.0}{8.0}{36.0}
\drawdashline{8.0}{36.0}{-16.0}{48.0}
\drawvector{52.0}{70.0}{36.0}{-3}{-1}
\drawvector{44.0}{18.0}{36.0}{-4}{3}
\drawlefttext{46.0}{16.0}{cut-locus}
\drawlefttext{56.0}{72.0}{conjugate locus}
\drawcenteredtext{6.0}{-10.0}{generic contact case}
\drawcenteredtext{176.0}{-10.0}{generic Martinet case}
\path(146.0,48.0)(146.0,48.0)(146.6,47.99)(147.19,47.99)(147.77,47.97)(148.38,47.97)(149.0,47.97)(149.6,47.95)(150.19,47.93)(150.77,47.9)
\path(150.77,47.9)(151.38,47.9)(152.0,47.88)(152.6,47.84)(153.19,47.81)(153.77,47.79)(154.38,47.75)(155.0,47.72)(155.6,47.68)(156.19,47.65)
\path(156.19,47.65)(156.77,47.61)(157.38,47.56)(158.0,47.52)(158.6,47.47)(159.19,47.4)(159.77,47.36)(160.38,47.29)(161.0,47.25)(161.6,47.18)
\path(161.6,47.18)(162.19,47.11)(162.77,47.04)(163.38,46.99)(164.0,46.9)(164.6,46.84)(165.19,46.77)(165.77,46.68)(166.38,46.61)(167.0,46.52)
\path(167.0,46.52)(167.6,46.43)(168.19,46.34)(168.77,46.25)(169.38,46.15)(170.0,46.08)(170.6,45.97)(171.19,45.88)(171.77,45.77)(172.38,45.65)
\path(172.38,45.65)(173.0,45.56)(173.58,45.45)(174.19,45.34)(174.77,45.22)(175.38,45.11)(176.0,45.0)(176.58,44.86)(177.19,44.75)(177.77,44.61)
\path(177.77,44.61)(178.38,44.5)(179.0,44.36)(179.58,44.22)(180.19,44.09)(180.77,43.95)(181.38,43.81)(181.99,43.68)(182.58,43.52)(183.19,43.38)
\path(183.19,43.38)(183.77,43.22)(184.38,43.08)(184.99,42.93)(185.58,42.77)(186.19,42.61)(186.77,42.45)(187.38,42.27)(187.99,42.11)(188.58,41.95)
\path(188.58,41.95)(189.19,41.77)(189.77,41.59)(190.38,41.4)(190.99,41.25)(191.58,41.06)(192.19,40.88)(192.77,40.68)(193.38,40.5)(193.99,40.31)
\path(193.99,40.31)(194.58,40.11)(195.19,39.93)(195.77,39.72)(196.38,39.52)(196.99,39.33)(197.58,39.11)(198.19,38.9)(198.77,38.7)(199.38,38.49)
\path(199.38,38.49)(199.99,38.27)(200.58,38.06)(201.19,37.84)(201.77,37.61)(202.38,37.38)(202.99,37.15)(203.58,36.93)(204.19,36.7)(204.77,36.47)
\path(204.77,36.47)(205.38,36.22)(205.99,36.0)(206.0,36.0)
\path(146.0,48.0)(146.0,48.0)(146.35,48.04)(146.71,48.08)(147.07,48.11)(147.41,48.15)(147.77,48.22)(148.11,48.27)(148.47,48.31)(148.8,48.38)
\path(148.8,48.38)(149.13,48.43)(149.5,48.5)(149.83,48.56)(150.16,48.61)(150.5,48.68)(150.83,48.75)(151.16,48.81)(151.5,48.88)(151.83,48.95)
\path(151.83,48.95)(152.13,49.04)(152.47,49.11)(152.8,49.2)(153.11,49.27)(153.41,49.36)(153.75,49.43)(154.05,49.52)(154.36,49.61)(154.66,49.7)
\path(154.66,49.7)(154.99,49.79)(155.27,49.9)(155.58,50.0)(155.88,50.09)(156.19,50.2)(156.49,50.29)(156.77,50.4)(157.08,50.5)(157.36,50.61)
\path(157.36,50.61)(157.66,50.72)(157.94,50.84)(158.22,50.95)(158.5,51.08)(158.77,51.18)(159.07,51.31)(159.35,51.43)(159.63,51.56)(159.88,51.68)
\path(159.88,51.68)(160.16,51.81)(160.44,51.95)(160.71,52.08)(160.97,52.22)(161.22,52.36)(161.5,52.49)(161.75,52.63)(162.0,52.77)(162.27,52.9)
\path(162.27,52.9)(162.52,53.06)(162.77,53.22)(163.02,53.36)(163.27,53.52)(163.5,53.68)(163.75,53.84)(164.0,53.99)(164.22,54.15)(164.47,54.31)
\path(164.47,54.31)(164.71,54.47)(164.94,54.65)(165.16,54.81)(165.38,54.99)(165.63,55.15)(165.85,55.34)(166.07,55.52)(166.27,55.68)(166.5,55.88)
\path(166.5,55.88)(166.72,56.06)(166.94,56.24)(167.16,56.43)(167.36,56.61)(167.58,56.81)(167.77,57.0)(167.99,57.2)(168.19,57.4)(168.38,57.59)
\path(168.38,57.59)(168.58,57.79)(168.77,58.0)(168.99,58.2)(169.16,58.4)(169.36,58.61)(169.55,58.83)(169.75,59.04)(169.91,59.25)(170.11,59.47)
\path(170.11,59.47)(170.27,59.68)(170.47,59.9)(170.63,60.13)(170.83,60.36)(171.0,60.59)(171.16,60.81)(171.33,61.04)(171.5,61.27)(171.66,61.52)
\path(171.66,61.52)(171.83,61.75)(172.0,61.99)(172.0,62.0)
\drawthickdot{146.0}{48.0}
\drawthickdot{172.0}{62.0}
\drawthickdot{206.0}{36.0}
\drawthickdot{182.0}{44.0}
\drawcenteredtext{184.0}{48.0}{$C$}
\drawcenteredtext{146.0}{44.0}{$A$}
\drawcircle{206.0}{76.0}{16.0}{R}
\drawcircle{182.0}{44.0}{20.0}{}
\drawcircle{146.0}{48.0}{20.0}{}
\end{picture}

\end{center}
\caption{}\label{lafig5}
\end{figure}

The microlocal invariants of the SR balls are the following~:
\begin{itemize}
\item Spectrum band corresponding to the focus $F$.
\item Spectrum band corresponding to the saddle $S$.
\item Invariants connected to the family of Riemannian structures
$a(q)dx^2+b(q)dy^2$ induced on the plane $(x,y)$, where $z$ is
taken as a parameter.
\end{itemize}

The conjugate points accumulate in the flat case along the abnormal
direction, see \cite{ABCK}.


\subsubsection{The n-dimensional case}
Our results except the precise asymptotics of Section 4 can be
generalized to the $n$-dimensional case to define a Martinet
sector in the SR ball. Indeed~:
\begin{itemize}
\item $L^2$-compactness of SR minimizers (see \cite{Ag}) allows
to bound the number of oscillations of Lagrangian manifolds. It
appears in our study by taking only the first and second return
mapping to compute the sphere intersected with $y=0$.
\item From \cite{trelatthese}, using a normal form, we can cut
the SR ball by a 2-dimensional plane to identify a Martinet
sector which splits into two curves~: a curve $D_1$ obtained by
using minimizing controls close to the reference abnormal control
in $L^\infty$-topology~; a curve $E_1$ obtained by using controls
close to the reference abnormal control in $L^2$-topology, but
not in $L^\infty$-topology (see Fig. \ref{lafig6}).
\end{itemize}

\setlength{\unitlength}{0.45mm}
\begin{figure}[h]
\begin{center}

\begin{picture}(180,100)
\thinlines
\drawvector{2.0}{48.0}{176.0}{1}{0}
\drawcenteredtext{28.0}{14.0}{$E_1$}
\drawcenteredtext{184.0}{54.0}{abnormal line}
\drawthickdot{90.0}{48.0}
\drawlefttext{88.0}{42.0}{$A$}
\path(90.0,48.0)(90.0,48.0)(90.98,48.04)(91.98,48.09)(92.97,48.15)(93.94,48.22)(94.91,48.29)(95.87,48.36)(96.83,48.45)(97.8,48.56)
\path(97.8,48.56)(98.75,48.65)(99.69,48.77)(100.62,48.88)(101.55,49.02)(102.48,49.15)(103.41,49.29)(104.3,49.45)(105.23,49.61)(106.12,49.77)
\path(106.12,49.77)(107.01,49.95)(107.91,50.13)(108.8,50.31)(109.66,50.5)(110.54,50.7)(111.41,50.93)(112.26,51.13)(113.12,51.36)(113.97,51.59)
\path(113.97,51.59)(114.8,51.84)(115.62,52.09)(116.47,52.34)(117.3,52.61)(118.11,52.88)(118.91,53.15)(119.73,53.45)(120.51,53.75)(121.3,54.04)
\path(121.3,54.04)(122.11,54.36)(122.87,54.68)(123.66,55.0)(124.43,55.33)(125.19,55.65)(125.94,56.02)(126.69,56.38)(127.44,56.74)(128.19,57.11)
\path(128.19,57.11)(128.91,57.49)(129.63,57.88)(130.36,58.27)(131.08,58.65)(131.77,59.08)(132.49,59.49)(133.19,59.9)(133.88,60.34)(134.57,60.79)
\path(134.57,60.79)(135.25,61.24)(135.91,61.68)(136.58,62.15)(137.25,62.61)(137.88,63.09)(138.55,63.58)(139.19,64.05)(139.83,64.55)(140.46,65.08)
\path(140.46,65.08)(141.08,65.58)(141.71,66.12)(142.32,66.65)(142.91,67.19)(143.52,67.73)(144.11,68.29)(144.71,68.83)(145.27,69.41)(145.86,69.98)
\path(145.86,69.98)(146.44,70.55)(147.0,71.16)(147.57,71.76)(148.11,72.37)(148.66,72.98)(149.21,73.61)(149.74,74.23)(150.27,74.87)(150.77,75.51)
\path(150.77,75.51)(151.3,76.16)(151.82,76.83)(152.33,77.48)(152.83,78.16)(153.32,78.83)(153.8,79.54)(154.27,80.23)(154.75,80.94)(155.22,81.65)
\path(155.22,81.65)(155.69,82.37)(156.13,83.08)(156.6,83.83)(157.05,84.58)(157.49,85.33)(157.91,86.08)(158.35,86.86)(158.77,87.62)(159.16,88.41)
\path(159.16,88.41)(159.58,89.19)(159.99,89.98)(160.0,90.0)
\path(90.0,48.0)(90.0,48.0)(89.04,48.02)(88.08,48.06)(87.15,48.08)(86.19,48.09)(85.26,48.11)(84.36,48.11)(83.44,48.11)(82.51,48.09)
\path(82.51,48.09)(81.62,48.08)(80.73,48.06)(79.83,48.02)(78.94,47.99)(78.08,47.93)(77.22,47.88)(76.36,47.83)(75.51,47.75)(74.66,47.68)
\path(74.66,47.68)(73.8,47.61)(72.98,47.52)(72.16,47.43)(71.33,47.34)(70.51,47.22)(69.69,47.11)(68.91,47.0)(68.12,46.86)(67.33,46.74)
\path(67.33,46.74)(66.55,46.59)(65.76,46.45)(65.01,46.29)(64.26,46.13)(63.5,45.97)(62.75,45.79)(62.02,45.61)(61.29,45.4)(60.56,45.22)
\path(60.56,45.22)(59.84,45.02)(59.13,44.81)(58.4,44.61)(57.72,44.38)(57.04,44.15)(56.34,43.9)(55.65,43.68)(55.0,43.43)(54.34,43.15)
\path(54.34,43.15)(53.68,42.9)(53.02,42.63)(52.38,42.36)(51.75,42.08)(51.11,41.79)(50.5,41.5)(49.88,41.18)(49.27,40.88)(48.65,40.56)
\path(48.65,40.56)(48.06,40.24)(47.47,39.9)(46.9,39.56)(46.31,39.22)(45.75,38.88)(45.18,38.52)(44.63,38.15)(44.09,37.77)(43.54,37.4)
\path(43.54,37.4)(43.0,37.02)(42.47,36.63)(41.95,36.22)(41.45,35.81)(40.93,35.4)(40.43,34.99)(39.93,34.56)(39.45,34.13)(38.97,33.7)
\path(38.97,33.7)(38.5,33.25)(38.02,32.79)(37.56,32.34)(37.11,31.87)(36.65,31.37)(36.22,30.92)(35.79,30.43)(35.36,29.94)(34.95,29.44)
\path(34.95,29.44)(34.54,28.93)(34.13,28.37)(33.74,27.87)(33.34,27.36)(32.95,26.79)(32.58,26.29)(32.2,25.7)(31.84,25.19)(31.45,24.62)
\path(31.45,24.62)(31.12,24.04)(30.79,23.45)(30.45,22.87)(30.12,22.29)(29.79,21.7)(29.45,21.11)(29.17,20.5)(28.87,19.87)(28.54,19.26)
\path(28.54,19.26)(28.28,18.62)(28.0,18.0)(28.0,18.0)
\drawrighttext{76.0}{32.0}{$L^2$-sector}
\drawthickdot{4.0}{48.0}
\thicklines
\path(90.0,48.0)(90.0,48.0)(90.98,48.04)(91.98,48.09)(92.97,48.15)(93.94,48.22)(94.91,48.29)(95.87,48.36)(96.83,48.45)(97.8,48.56)
\path(97.8,48.56)(98.75,48.65)(99.69,48.77)(100.62,48.88)(101.55,49.02)(102.48,49.15)(103.41,49.29)(104.3,49.45)(105.23,49.61)(106.12,49.77)
\path(106.12,49.77)(107.01,49.95)(107.91,50.13)(108.8,50.31)(109.66,50.5)(110.54,50.7)(111.41,50.93)(112.26,51.13)(113.12,51.36)(113.97,51.59)
\path(113.97,51.59)(114.8,51.84)(115.62,52.09)(116.47,52.34)(117.3,52.61)(118.11,52.88)(118.91,53.15)(119.73,53.45)(120.51,53.75)(121.3,54.04)
\path(121.3,54.04)(122.11,54.36)(122.87,54.68)(123.66,55.0)(124.43,55.33)(125.19,55.65)(125.94,56.02)(126.69,56.38)(127.44,56.74)(128.19,57.11)
\path(128.19,57.11)(128.91,57.49)(129.63,57.88)(130.36,58.27)(131.08,58.65)(131.77,59.08)(132.49,59.49)(133.19,59.9)(133.88,60.34)(134.57,60.79)
\path(134.57,60.79)(135.25,61.24)(135.91,61.68)(136.58,62.15)(137.25,62.61)(137.88,63.09)(138.55,63.58)(139.19,64.05)(139.83,64.55)(140.46,65.08)
\path(140.46,65.08)(141.08,65.58)(141.71,66.12)(142.32,66.65)(142.91,67.19)(143.52,67.73)(144.11,68.29)(144.71,68.83)(145.27,69.41)(145.86,69.98)
\path(145.86,69.98)(146.44,70.55)(147.0,71.16)(147.57,71.76)(148.11,72.37)(148.66,72.98)(149.21,73.61)(149.74,74.23)(150.27,74.87)(150.77,75.51)
\path(150.77,75.51)(151.3,76.16)(151.82,76.83)(152.33,77.48)(152.83,78.16)(153.32,78.83)(153.8,79.54)(154.27,80.23)(154.75,80.94)(155.22,81.65)
\path(155.22,81.65)(155.69,82.37)(156.13,83.08)(156.6,83.83)(157.05,84.58)(157.49,85.33)(157.91,86.08)(158.35,86.86)(158.77,87.62)(159.16,88.41)
\path(159.16,88.41)(159.58,89.19)(159.99,89.98)(160.0,90.0)
\path(90.0,48.0)(90.0,48.0)(89.04,48.02)(88.08,48.06)(87.15,48.08)(86.19,48.09)(85.26,48.11)(84.36,48.11)(83.44,48.11)(82.51,48.09)
\path(82.51,48.09)(81.62,48.08)(80.73,48.06)(79.83,48.02)(78.94,47.99)(78.08,47.93)(77.22,47.88)(76.36,47.83)(75.51,47.75)(74.66,47.68)
\path(74.66,47.68)(73.8,47.61)(72.98,47.52)(72.16,47.43)(71.33,47.34)(70.51,47.22)(69.69,47.11)(68.91,47.0)(68.12,46.86)(67.33,46.74)
\path(67.33,46.74)(66.55,46.59)(65.76,46.45)(65.01,46.29)(64.26,46.13)(63.5,45.97)(62.75,45.79)(62.02,45.61)(61.29,45.4)(60.56,45.22)
\path(60.56,45.22)(59.84,45.02)(59.13,44.81)(58.4,44.61)(57.72,44.38)(57.04,44.15)(56.34,43.9)(55.65,43.68)(55.0,43.43)(54.34,43.15)
\path(54.34,43.15)(53.68,42.9)(53.02,42.63)(52.38,42.36)(51.75,42.08)(51.11,41.79)(50.5,41.5)(49.88,41.18)(49.27,40.88)(48.65,40.56)
\path(48.65,40.56)(48.06,40.24)(47.47,39.9)(46.9,39.56)(46.31,39.22)(45.75,38.88)(45.18,38.52)(44.63,38.15)(44.09,37.77)(43.54,37.4)
\path(43.54,37.4)(43.0,37.02)(42.47,36.63)(41.95,36.22)(41.45,35.81)(40.93,35.4)(40.43,34.99)(39.93,34.56)(39.45,34.13)(38.97,33.7)
\path(38.97,33.7)(38.5,33.25)(38.02,32.79)(37.56,32.34)(37.11,31.87)(36.65,31.38)(36.22,30.92)(35.79,30.43)(35.36,29.94)(34.95,29.44)
\path(34.95,29.44)(34.54,28.93)(34.13,28.39)(33.74,27.87)(33.34,27.36)(32.95,26.81)(32.58,26.29)(32.2,25.72)(31.84,25.19)(31.47,24.62)
\thinlines
\drawcenteredtext{132.0}{88.0}{$L^\infty$-sector}
\drawcenteredtext{170.0}{88.0}{$D_1$}
\thicklines
\path(31.47,24.62)(31.12,24.04)(30.79,23.46)(30.45,22.88)(30.12,22.29)(29.79,21.7)(29.46,21.11)(29.17,20.5)(28.87,19.87)(28.55,19.26)
\path(28.55,19.26)(28.28,18.62)(28.0,18.0)(28.0,18.0)
\thinlines
\path(4.0,48.0)(4.0,48.0)(5.63,48.04)(7.26,48.08)(8.89,48.13)(10.51,48.15)(12.1,48.22)(13.72,48.27)(15.31,48.33)(16.92,48.38)
\path(16.92,48.38)(18.51,48.45)(20.09,48.52)(21.67,48.58)(23.22,48.65)(24.79,48.72)(26.37,48.79)(27.92,48.86)(29.45,48.93)(31.01,49.02)
\path(31.01,49.02)(32.54,49.09)(34.06,49.18)(35.59,49.27)(37.11,49.36)(38.61,49.45)(40.13,49.54)(41.63,49.65)(43.11,49.75)(44.61,49.84)
\path(44.61,49.84)(46.09,49.95)(47.56,50.06)(49.02,50.15)(50.49,50.27)(51.95,50.38)(53.4,50.5)(54.84,50.61)(56.29,50.74)(57.72,50.86)
\path(57.72,50.86)(59.15,50.99)(60.56,51.11)(61.97,51.25)(63.38,51.38)(64.79,51.5)(66.19,51.65)(67.58,51.79)(68.97,51.93)(70.33,52.08)
\path(70.33,52.08)(71.72,52.22)(73.08,52.36)(74.44,52.52)(75.8,52.68)(77.15,52.84)(78.48,52.99)(79.83,53.15)(81.16,53.31)(82.48,53.49)
\path(82.48,53.49)(83.8,53.65)(85.12,53.81)(86.43,54.0)(87.73,54.15)(89.01,54.34)(90.3,54.52)(91.58,54.7)(92.87,54.9)(94.12,55.09)
\path(94.12,55.09)(95.41,55.27)(96.66,55.47)(97.91,55.65)(99.16,55.86)(100.41,56.06)(101.62,56.25)(102.87,56.47)(104.08,56.65)(105.3,56.88)
\path(105.3,56.88)(106.51,57.09)(107.73,57.31)(108.93,57.52)(110.12,57.74)(111.3,57.97)(112.48,58.18)(113.66,58.4)(114.83,58.63)(115.98,58.86)
\path(115.98,58.86)(117.15,59.11)(118.3,59.34)(119.44,59.58)(120.58,59.81)(121.72,60.06)(122.83,60.31)(123.97,60.56)(125.08,60.81)(126.19,61.06)
\path(126.19,61.06)(127.29,61.31)(128.38,61.56)(129.47,61.83)(130.57,62.09)(131.63,62.36)(132.72,62.61)(133.77,62.88)(134.85,63.15)(135.88,63.43)
\path(135.88,63.43)(136.94,63.72)(137.99,63.99)(138.0,64.0)
\path(4.0,48.0)(4.0,48.0)(5.75,48.24)(7.42,48.47)(9.01,48.72)(10.51,48.97)(11.93,49.22)(13.27,49.49)(14.55,49.74)(15.72,50.0)
\path(15.72,50.0)(16.86,50.25)(17.89,50.52)(18.87,50.77)(19.79,51.04)(20.62,51.29)(21.42,51.56)(22.12,51.81)(22.79,52.06)(23.37,52.33)
\path(23.37,52.33)(23.93,52.58)(24.39,52.83)(24.84,53.08)(25.2,53.33)(25.54,53.56)(25.81,53.81)(26.04,54.04)(26.25,54.27)(26.37,54.5)
\path(26.37,54.5)(26.47,54.72)(26.54,54.93)(26.56,55.15)(26.55,55.36)(26.52,55.56)(26.44,55.75)(26.31,55.93)(26.2,56.11)(26.03,56.27)
\path(26.03,56.27)(25.84,56.45)(25.62,56.59)(25.37,56.75)(25.12,56.88)(24.86,57.02)(24.55,57.13)(24.26,57.25)(23.94,57.34)(23.61,57.43)
\path(23.61,57.43)(23.27,57.52)(22.93,57.59)(22.55,57.63)(22.2,57.68)(21.86,57.72)(21.5,57.75)(21.12,57.75)(20.78,57.75)(20.43,57.72)
\path(20.43,57.72)(20.06,57.7)(19.75,57.65)(19.42,57.59)(19.11,57.52)(18.79,57.43)(18.52,57.33)(18.25,57.2)(18.0,57.08)(17.77,56.9)
\path(17.77,56.9)(17.55,56.75)(17.37,56.56)(17.21,56.36)(17.1,56.15)(17.0,55.9)(16.94,55.65)(16.88,55.38)(16.88,55.09)(16.94,54.77)
\path(16.94,54.77)(17.02,54.45)(17.12,54.09)(17.29,53.72)(17.47,53.34)(17.72,52.93)(18.04,52.49)(18.37,52.04)(18.78,51.56)(19.21,51.06)
\path(19.21,51.06)(19.71,50.54)(20.29,50.0)(20.89,49.43)(21.56,48.84)(22.3,48.24)(23.12,47.59)(23.97,46.93)(24.92,46.25)(25.93,45.54)
\path(25.93,45.54)(27.0,44.81)(28.12,44.04)(29.36,43.25)(30.64,42.45)(32.02,41.61)(33.47,40.74)(35.02,39.84)(36.63,38.93)(38.33,37.97)
\path(38.33,37.97)(40.11,37.0)(41.99,36.0)(42.0,36.0)
\drawpath{90.0}{50.0}{90.0}{46.0}
\end{picture}

\end{center}
\caption{}\label{lafig6}
\end{figure}

\it{The picture explains well the consequence of the existence of
abnormal minimizers in SR geometry. Contrarily to the classical
case we cannot straight the geodesic flow near the abnormal
direction to form a central field.}


\section{Conclusion}
Our analysis explains the role of abnormal geodesics in
SR geometry. It is based on the Martinet case.
Using our gradated normal form of order $0$, the
geodesics foliation is projected onto a one-dimensional foliation
in a plane which corresponds to a \it{one-parameter family of
pendulums}. In this space the abnormal line projects on the
singularities of the foliation. The computation of the sphere in
the abnormal direction is related to the computation of return
mappings evaluated along the separatrices of the pendulum.
We have computed \it{asymptotics}, using techniques similar to
the ones used in the Hilbert's 16th problem. The computations are
complex, because it is a singular perturbation analysis. In these
computations one needs to consider geodesics $C^1$-close to the
abnormal reference trajectory on the one part, and
geodesics which are $C^0$-close, but not $C^1$-close to the
abnormal reference trajectory on the other part. Our asymptotics
are not complete in the latter case and this requires further
studies. Moreover the techniques have to be adapted to analyze
the general case when the geodesics equations are not
projectable. This leads to \it{stability questions about our
asymptotics}.

The projection of the geodesics flow onto a planar foliation,
valid at order $0$ in the Martinet case, is useful to compute
asymptotics but is not crucial from the geometric point of view,
and Martinet geometry 
is representative of SR geometry with abnormal minimizers.
The existence of such minimizers implies \it{hyperbolicity} seen in
the pendulum representation as the behaviors of the geodesics
near the separatrices. The general geometric framework to analyze
SR geometry is Lagrangian manifolds. Here hyperbolicity due to
the existence of abnormal directions is interpretated as a
\it{splitting} of the Lagrangian fiber $T^*_{q_0}M$ when transported
by the normal flow. To construct the sphere in the abnormal direction
we must glue together the projections of several manifolds in the
cotangent space.

The link between the computations of asymptotics and Lagrangian
geometry is the Jacobi fields which allow to compute contacts
for the return mapping underlying our analysis.
In general the evaluation of a return mapping in the analysis
of a differential equation interpretated as a transport problem
is original and source of further studies.

The question of the category of the SR Martinet sphere is still
open. The Martinet sector is homeomorphic to a locally convex
cone and we have given a qualitative description of its
singularities. In the integrable case the SR sphere is log-exp
and this leads to a smooth stratification of the sphere.
In general we
conjecture that the sphere is not log-exp, belongs to some
extended Il'Yashenko's category, and is still $C^1$-stratifiable.
It is an important question connected to \it{Hamilton-Jacobi equation}
which has many applications in physics (optics, quantum theory).
For control theory, SR geometry is part of optimal control. Moreover
our study is related to the stabilization problem.




\end{document}